\documentclass[UTF-8,reqno]{amsart}
\usepackage{setspace}
\usepackage{enumerate}
\usepackage{mhequ}
\usepackage[margin=1.0in]{geometry}
\usepackage{amssymb,url,color, booktabs,nccmath}
\usepackage{graphicx}
\usepackage{tikz}
\usetikzlibrary{arrows.meta,positioning}
\usepackage{mathrsfs}
\usepackage{enumitem}
\usepackage{graphicx}
\usepackage{tikz}
\usetikzlibrary{shapes,snakes}
\usetikzlibrary{calc}
\usetikzlibrary{decorations.shapes}
\usepackage{color}
\usepackage[colorlinks=true]{hyperref}
\hypersetup{
    linkcolor=blue,
    citecolor=red,
    filecolor=magenta,
    urlcolor=cyan
}

\definecolor{blue}{rgb}{0,0,1}
\definecolor{red}{rgb}{1,0,0}
\definecolor{cyan}{rgb}{0,1,1}
\definecolor{magenta}{rgb}{1,0,1}
\definecolor{green}{rgb}{0,0.5,0}
\definecolor{orange}{rgb}{1,0.5,0}
\definecolor{purple}{rgb}{0.5,0,0.5}
\definecolor{darkergreen}{rgb}{0,0.5,0}

\numberwithin{equation}{section}

\newcommand{\be}{\begin{eqnarray}}
\newcommand{\ee}{\end{eqnarray}}
\newcommand{\ce}{\begin{eqnarray*}}
\newcommand{\de}{\end{eqnarray*}}
\newtheorem{theorem}{Theorem}[section]
\newtheorem{lemma}[theorem]{Lemma}
\newtheorem{remark}[theorem]{Remark}
\newtheorem{definition}[theorem]{Definition}
\newtheorem{proposition}[theorem]{Proposition}
\newtheorem{Examples}[theorem]{Example}
\newtheorem{corollary}[theorem]{Corollary}
\newtheorem{assumption}{Assumption}[section]

\newenvironment{nouppercase}{
  
  \renewcommand{\uppercasenonmath}[1]{}}{}

\def\p{\partial}

\def\[{{\Big[}}
\def\]{{\Big]}}
\def\<{{\langle}}
\def\>{{\rangle}}
\def\({{\Big(}}
\def\){{\Big)}}

\def\bx{{\mathbf{x}}}

\def\sgn{\mbox{\rm sgn}}

\def\dif{{\mathord{{\rm d}}}}

\def\min{{\mathord{{\rm min}}}}

\def\no{\nonumber}
\def\={&\!\!=\!\!&}

\def\bB{{\mathbf B}}

\def\cR{{\mathcal R}}

\def\mE{{\mathbb E}}

\def\mN{{\mathbb N}}

\def\mR{{\mathbb R}}

\def\bB{{\mathbf B}}

\def\1{{I}}

\def\sF{{\mathscr F}}

\def\sL{{\mathscr L}}

\def\sS{{\mathscr S}}

\def\geq{\geqslant}
\def\leq{\leqslant}
\def\ge{\geqslant}
\def\le{\leqslant}

\def\iint{\int\!\!\!\int}

\def\p{\partial}

\def\[{{\Big[}}
\def\]{{\Big]}}
\def\<{{\langle}}
\def\>{{\rangle}}
\def\({{\Big(}}
\def\){{\Big)}}

\def\bx{{\mathbf{x}}}

\def\sgn{\mbox{\rm sgn}}

\def\dif{{\mathord{{\rm d}}}}

\def\min{{\mathord{{\rm min}}}}

\def\no{\nonumber}
\def\={&\!\!=\!\!&}
\def\bt{\begin{theorem}}
\def\et{\end{theorem}}
\def\bl{\begin{lemma}}
\def\el{\end{lemma}}
\def\br{\begin{remark}}
\def\er{\end{remark}}
\def\bx{\begin{Examples}}
\def\ex{\end{Examples}}
\def\bd{\begin{definition}}
\def\ed{\end{definition}}
\def\bp{\begin{proposition}}
\def\ep{\end{proposition}}
\def\bc{\begin{corollary}}
\def\ec{\end{corollary}}

\def\geq{\geqslant}
\def\leq{\leqslant}
\def\ge{\geqslant}
\def\le{\leqslant}

\def\iint{\int\!\!\!\int}

\def\<{\langle} \def\>{\rangle}

\allowdisplaybreaks

\tikzset{
        dot/.style={circle,fill=black,inner sep=0pt, outer sep=0.7pt, minimum size=1mm},
        Phi/.style={white!40!red,thick,snake=coil,segment amplitude=0.6pt, segment length=2pt},
         Z/.style={black!40!green,thick,snake=coil,segment amplitude=0.6pt, segment length=2pt},
        C/.style={thick,black!20!blue},
          Cr/.style={thick,black!20!red},
            Cg/.style={thick,black!20!green},
       }

\begin{document}
\setstretch{0.9}

\title[Vlasov-Fokker-Planck-Dean-Kawasaki Equation]{\LARGE Kinetic Theory with Fluctuations: Well-Posedness of The Vlasov--Fokker--Planck--Dean--Kawasaki Equations}

\author[Zimo Hao]{\large Zimo Hao}
\address[Z. Hao]{School of Mathematics and Statistics, Beijing Institute of Technology, Beijing, China}
\email{zimo\_hao@163.com}

\author[Zhengyan Wu]{\large Zhengyan Wu}
\address[Z. Wu]{Department of Mathematics, Technische Universit\"at M\"unchen, Boltzmannstr. 3, 85748 Garching, Germany}
\email{wuzh@cit.tum.de}

\author[Johannes Zimmer]{\large Johannes Zimmer}
\address[J. Zimmer]{Department of Mathematics, Technische Universit\"at M\"unchen, Boltzmannstr. 3, 85748 Garching, Germany}
\email{jz@tum.de}

\begin{abstract}
We study Vlasov--Fokker--Planck--Dean--Kawasaki equations driven by correlated conservative noise. For regular noise coefficients and bounded nonlocal interactions, we establish probabilistically strong existence and uniqueness in a renormalized kinetic framework. For the square-root coefficient, we treat the non-interacting case and construct a probabilistically weak solution. Key challenges
stem from the complexity of the kinetic operator and the irregularity introduced by the conservative noise
with square-root-type coefficients. The proof relies on a novel combination of kinetic semigroup estimates
and the framework of renormalized kinetic solutions.
\end{abstract}


\date{\today}

\begin{nouppercase}
\maketitle
\end{nouppercase}

\setcounter{tocdepth}{1}
\tableofcontents

\section{Introduction}
Kinetic theory offers a fundamental framework for describing the evolution of particle distributions in phase space, serving as a bridge between microscopic dynamics and macroscopic observables. It plays a crucial role in the modeling of dilute gases, plasmas, and interacting particle systems, particularly in far-from-equilibrium regimes. Complementing this, the theory of fluctuating hydrodynamics extends classical continuum models by incorporating stochastic fluctuations, which are indispensable for capturing mesoscopic phenomena. Motivated by these developments, we seek to understand how kinetic equations, when enriched with fluctuation corrections, can give rise to stochastic mesoscopic models that go beyond the scope of deterministic descriptions.

{\color{black}In this work, we analyze a fluctuation-corrected kinetic model, namely the Vlasov--Fokker--Planck--Dean--Kawasaki (VFPDK) equation with correlated noise. We prove strong well-posedness for a regular-coefficient equation with bounded nonlocal interactions and probabilistically weak existence for the square-root equation without interaction. The general formal model takes the form}
\begin{align}\label{SPDE-00}
{\color{black}\partial_tf =\Delta_vf-v\cdot\nabla_xf-\nabla_v\cdot \big(fV\star_x\rho\big)+\nabla_v\cdot (vf)-\nabla_v\cdot(\sqrt{f} \circ\xi_F), \quad (t,x,v) \in \mathbb{R}_+ \times \mathbb{R}^{2d},}
\end{align}
where $\rho=\int_{\mathbb{R}^d}fdv$. The stochastic forcing $\xi_F$ is white in time and correlated in both spatial and velocity variables. The symbol $\star_x$ denotes convolution in the spatial variable, and $\circ$ indicates Stratonovich stochastic integration. The interaction kernel $V : \mathbb{R}^d \rightarrow \mathbb{R}^d$ depends only on the spatial variable $x \in \mathbb{R}^d$. 

To the best of our knowledge, the original VFPDK equation was proposed by M\"uller, von Renesse, and the third author \cite{FMJ25}, and takes the form
\begin{equation}\label{SPDE-VFPDK-intro}
	{\color{black}\partial_t f_N = \frac{\kappa}{2} \Delta_v f_N - v \cdot \nabla_x f_N - \nabla_v \cdot (f_N V \star_x \rho_N) + \kappa \nabla_v \cdot (v f_N) - \sqrt{\frac{\kappa}{N}} \nabla_v \cdot (\sqrt{f_N} \, \xi),}
\end{equation}
for every $N \in \mathbb{N}_+$ and some $\kappa > 0$. Here, $\rho_N = \int_{\mathbb{R}^d} f_N \, dv$, and $\xi$ denotes space-velocity-time white noise. The sign '$-$' in front of the noise can also be replaced by a '$+$'; this does not change the quadratic variation structure used in the derivation. 
The model \eqref{SPDE-VFPDK-intro} is designed to describe the evolution of the empirical measure of the second-order mean-field particle system
\begin{align}\label{particles}
dX_i &= V_i \, dt, \notag \\
dV_i &= \frac{1}{N} \sum_{j=1}^N V(X_i - X_j) \, dt - \kappa V_i \, dt + \sqrt{\kappa} \, dB_i(t), \quad i = 1, \dots, N,
\end{align}
where $\{B_i(t)\}_{i=1}^N$ are independent Brownian motions. More precisely, applying It\^o's formula to the empirical measure $\pi_N = \frac{1}{N} \sum_{i=1}^N \delta_{(X_i, V_i)}$ yields, in a weak sense,
\begin{equation}\label{SPDE-empirical-intro}
	{\color{black}\partial_t \pi_N = \frac{\kappa}{2} \Delta_v \pi_N - v \cdot \nabla_x \pi_N - \nabla_v \cdot (\pi_N V \star_x \rho_{\pi_N}) + \kappa \nabla_v \cdot (v \pi_N) - \frac{\sqrt{\kappa}}{N} \nabla_v \cdot \left( \sum_{i=1}^N \delta_{(X_i, V_i)} \dot{B}_i \right),}
\end{equation}
{\color{black}Here $\rho_{\pi_N}:=(\operatorname{pr}_x)_\#\pi_N$ denotes the spatial marginal of the empirical measure; it is a measure rather than an $L^2$ function.} The deterministic component of \eqref{SPDE-empirical-intro} corresponds to the classical Vlasov-Fokker-Planck equation \cite{Risken84}. Following the methodology of \cite{D96}, and by comparing the covariance structure of the stochastic terms, the noise term $$\frac{\sqrt{\kappa}}{N} \nabla_v \cdot \left( \sum_{i=1}^N \delta_{(X_i, V_i)} \dot{B}_i \right)$$ can be statistically approximated by a Dean-Kawasaki-type noise of the form $\sqrt{\frac{\kappa}{N}} \nabla_v \cdot (\sqrt{\pi_N} \, \xi)$. This leads to the formal derivation of the VFPDK equation \eqref{SPDE-VFPDK-intro}. 

Due to the high singularity of the space-velocity-time white noise, the equation \eqref{SPDE-VFPDK-intro} falls into the supercritical regime in the sense of singular SPDE theory \cite{H14,GIP12}, making its well-posedness a significant challenge. Nonetheless, \cite{FMJ25}, extending the techniques developed in \cite{KLR19,KLR20}, provides a rigorous mathematical interpretation of \eqref{SPDE-VFPDK-intro}, showing that the empirical measure $\pi_N$ associated with the particle system \eqref{particles} uniquely (in law) solves \eqref{SPDE-VFPDK-intro} in the sense of a martingale problem. This demonstrates the close relationship between the VFPDK equation and the second-order mean-field particle system. 

The equation \eqref{SPDE-VFPDK-intro} can be interpreted as a fluctuating hydrodynamics equation. It combines the Vlasov-Fokker-Planck equation  with a stochastic perturbation of Dean-Kawasaki type. In the context of fluctuating hydrodynamics, such noise terms are typically governed by the fluctuation-dissipation relation, reflecting the interplay between microscopic fluctuations and macroscopic dissipation.

To place this in a broader physical context, we briefly recall the theoretical foundations of fluctuating hydrodynamics (FHD). Closely linked to this is Macroscopic Fluctuation Theory (MFT), which provides a unified framework for studying systems far from equilibrium, extending classical near-equilibrium linear response theory (see Bertini et al. \cite{BDGJL}, Derrida \cite{Derrida}). At its core, MFT is based on an ansatz concerning the large deviation principles for interacting particle systems. FHD, on the other hand, seeks to model microscopic fluctuations by incorporating stochastic perturbations into hydrodynamic equations, drawing on principles from statistical mechanics and non-equilibrium thermodynamics. In this setting, conservative stochastic partial differential equations (SPDEs) are postulated to encode the essential features of fluctuation dynamics in non-equilibrium systems (see \cite{LL87}, \cite{HS}). The foundational ansatz of MFT often emerges as the zero-noise large deviation limit of such SPDEs, providing a formal bridge between particle-level randomness and macroscopic fluctuation behavior.

In general, fluctuating hydrodynamics equations are intended to describe systems at the mesoscopic scale. Consequently, the driving noise should exhibit only short-range correlations in both space and velocity. Motivated by this, we replace the idealized space-velocity-time white noise $\xi$ with a spatially and velocity-correlated noise. Additionally, both physical and analytical considerations motivate further modifications to the noise structure. From a thermodynamic perspective, \"{O}ttinger \cite{Ottinger} introduced a stochastic formulation of the GENERIC framework driven by Klimontovich noise, with the aim of preserving the invariance of the dynamics under the Gibbs measure, albeit at a formal level. On the analytical side, Fehrman and Gess \cite{FG24} demonstrated that employing a Stratonovich noise of the Dean-Kawasaki equation yields a stochastic coercivity, which plays a crucial role in establishing well-posedness.

When comparing Klimontovich and Stratonovich noise, it is known that the associated It\^o corrections differ by a factor of $1/2$. However, this discrepancy is asymptotically negligible, having no impact on either the macroscopic limit or fluctuation behavior. We therefore adopt the Stratonovich noise and reformulate \eqref{SPDE-empirical-intro} as an SPDE driven by correlated Stratonovich noise. Furthermore, since the coefficient $\kappa$ in \eqref{SPDE-VFPDK-intro} does not play a substantive role in our analysis, we may, without loss of generality, set $\kappa=1$. These considerations lead us to the primary object of study in this work: the VFPDK equation \eqref{SPDE-00}.

    We now state the main result of this paper. We consider the equation 
    \begin{align}\label{SPDE-sigma-intro}
{\color{black}\partial_tf =\Delta_vf-v\cdot\nabla_xf-\nabla_v\cdot \big(fV\star_x\rho\big)+\nabla_v\cdot (vf)-\nabla_v\cdot(\sigma(f) \circ\xi_F), \quad (t,x,v) \in \mathbb{R}_+ \times \mathbb{R}^{2d}.} 
\end{align}

\begin{theorem}\label{thm-main-1}
{\color{black}Under the noise assumptions stated in Section~\ref{sec-2}, the following assertions hold.
\begin{enumerate}
\item[(i)] Suppose that $V\in L^\infty(\mathbb R_x^d)\cap L^2(\mathbb R_x^d)$, that $\sigma$ is a regular coefficient, and that the nonnegative initial datum satisfies
\begin{align*}
 {\color{black}\operatorname*{ess\,sup}_{\omega\in\Omega}}\|f_0(\omega)\|_{L^1}
 +\mathbb E\|f_0\|_{L^2}^2
 +\mathbb E\int_{\mathbb R^{2d}}(|x|^2+|v|^2)f_0(x,v)\,dxdv<\infty.
\end{align*}
{\color{black}Then \eqref{SPDE-sigma-intro} admits a unique probabilistically
strong weak solution in the sense of Definition~\ref{def-weaksolution}, and
this solution is a renormalized kinetic solution in the sense of
Definition~\ref{def-kineticsolution}.}
\item[(ii)] Suppose that $V\equiv0$, $\sigma(r)=\sqrt r$, and that the nonnegative initial datum satisfies
\begin{align*}
 {\color{black}\operatorname*{ess\,sup}_{\omega\in\Omega}}\|f_0(\omega)\|_{L^1}
 +\mathbb E\|f_0\|_{L^2}^2<\infty.
\end{align*}
Then \eqref{SPDE-sigma-intro} admits a probabilistically weak renormalized kinetic solution in the sense of Definition~\ref{def-prob-weak-square-root}. 
\end{enumerate}}
\end{theorem}

\subsection{New challenges: the VFPDK equation vs. the classical Dean-Kawasaki equation}
\ 

It is instructive to compare the VFPDK equation \eqref{SPDE-00} with the classical Dean--Kawasaki equation:
\begin{equation}\label{classical-dk}
	\partial_t\rho=\frac{1}{2}\Delta_x\rho-\nabla_x\cdot(\sqrt{\rho}\circ\xi_G), 
\end{equation}
where $\rho$ denotes the density depending only on the spatial and temporal variables, and $\xi_G$ denotes a noise that is white in time and spatially correlated. In order to establish the existence of solutions to \eqref{classical-dk}, \cite{FG24} developed both an energy-based compactness method and an entropy-based compactness method. We briefly illustrate the energy method: one can show that
\begin{align*}
\mathbb{E}\left(\int|\rho(s)|^2\,dx\right)\Big|^{s=t}_{s=0}
+\mathbb{E}\int^t_0\int|\nabla_x\rho(s)|^2\,dx\,ds
\leq C,
\end{align*}
for some finite constant $C>0$. This inequality provides control of the energy dissipation and, as a consequence, yields spatial regularity in the $x$-variable. In combination with suitable time-regularity estimates for an approximation sequence, one can then obtain compactness. 

Formally, one may attempt to apply an analogous argument to derive the following energy estimate for the VFPDK equation:
\begin{align*}
\mathbb{E}\left(\int | f(s)|^2\,dx\,dv\right)\Big|^{s=t}_{s=0}
+\mathbb{E}\int^t_0\int|\nabla_vf(s)|^2\,dx\,dv\,ds
\leq C.
\end{align*}
However, in this setting the energy dissipation term only provides regularity in the $v$-variable, without yielding any control in the $x$-variable. Consequently, this estimate alone is insufficient to apply the Aubin--Lions lemma and obtain strong compactness in Lebesgue spaces. To overcome this difficulty, we develop a new approach based on kinetic semigroup estimates, which yields the desired strong compactness and represents a significant departure from the methodology of \cite{FG24}.

{\color{black}A second difficulty is genuinely caused by the unbounded phase
space.  For a regular noise coefficient, uniform second spatial--velocity
moments prevent mass from escaping to infinity.  They upgrade local compactness
to strong convergence in global $L^1(\mathbb R^{2d})$, thereby allowing both
the conservation of total mass and the nonlocal interaction term to pass to
the limit.  The same moments control the negative part of the Boltzmann
entropy $f\log f-f$ and are therefore an essential ingredient in the entropy
dissipation estimate.  By contrast, the available estimates for smooth
approximations of the square-root coefficient are not uniform at the level of
spatial moments.  Nor can this difficulty be bypassed by demanding global
integrability of the logarithmic term appearing in the It\^o form of the
square-root equation.  Indeed, if $f\geq0$ belongs to
$L^1(\mathbb R^{2d})$, then the set $\{f>e^{-1}\}$ has finite Lebesgue
measure, whereas its complement has infinite measure; on the part of this
complement where $f>0$ one has $|\log f|\geq1$, and $\log f$ is singular on
$\{f=0\}$.  Thus $\log f$ cannot be globally integrable for a finite-mass
density on the whole space.  This explains why the whole-space relative
entropy theory of Fehrman and Gess~\cite{FG25} is formulated instead around a
strictly positive constant reference density $\gamma$.  Their relative
entropy is nonnegative, and, in the shifted class $L^1_\gamma$, it is
$\rho-\gamma$ rather than $\rho$ that is integrable; consequently, the
solution itself generally has infinite total mass.  This framework is well
suited to fluctuations around a homogeneous background, but it does not
describe the finite-mass density required here as an approximation of a
particle system.  We therefore work around the zero background and, for
regular coefficients, establish new moment estimates that preserve the
global $L^1$ structure.  At the square-root endpoint these uniform moment
estimates are unavailable: the compactness is only local in $L^1$, mass may
be lost at spatial infinity, and the global kinetic-measure estimates needed
for the standard doubling-of-variables uniqueness argument cannot be
justified.  These are the reasons why our square-root result is restricted to
the non-interacting case and gives probabilistically weak existence without
mass conservation or uniqueness.}

In addition, several further technical obstacles arise in this work. The damping term $\nabla_v\cdot(vf)$ in \eqref{SPDE-VFPDK-intro}, originating from the microscopic damping structure $-V_i$ in \eqref{particles}, introduces new challenges. Specifically, the function $\eta(v)=v$ appearing in this term lacks integrability and thus does not belong to any Besov space. Moreover, compared with \cite{FG24}, the present equation incorporates additional features such as transport, damping, and nonlocal interaction terms. These structural complexities necessitate substantially more delicate analysis in both the existence and uniqueness proofs. 

In the following, we outline the strategy employed to overcome these challenges.

\subsection{Key ideas and technical comments}\label{subsec-1-2}
{\color{black}We note that if one informally applies It\^o's formula to the entropy functional, the transport term $-v\cdot\nabla_x f$ vanishes. Consequently, the entropy dissipation estimate does not exploit any contribution of the transport term. However, a key observation is that a scaling relation is inherent in the kinetic operator $\Delta_v-v\cdot\nabla_x$. This suggests that by considering the combined effect of $\Delta_v$ and $v\cdot\nabla_x$, spatial regularity can be recovered through the kinetic semigroup generated by this hypoelliptic operator. We use the diffusion-one normalization throughout the analytical part of the paper. It is obtained from the physical diffusion-$1/2$ convention by a deterministic rescaling of time together with the corresponding drift and noise intensities; it is not a consequence of the Stratonovich--It\^o conversion.}

\ 

{\bf Scaling properties and anisotropic Besov spaces. } We first investigate the well-posedness of a regularized variant of the VFPDK equation with correlated noise:
\begin{align}\label{SPDE-0}
{\color{black}\partial_tf =\Delta_vf-v\cdot\nabla_xf-\nabla_v\cdot \big(fV\star_x\rho\big)+\nabla_v\cdot (vf)-\nabla_v\cdot(\sigma(f) \circ\xi_F), \quad (t,x,v) \in \mathbb{R}_+ \times \mathbb{R}^{2d}.} 
\end{align}
where the function $\sigma(\cdot)$ is a smooth approximation of a square-root-type coefficient. To utilize the kinetic semigroup estimates when studying this equation, we consider the model kinetic equation
\begin{equation}
{\color{black}\partial_t f=\big(\Delta_v-v\cdot\nabla_x\big)f.}	
\end{equation}
For every $\lambda>0$, let $f_\lambda(t,x,v):=f(\lambda^\alpha t, \lambda^\beta x, \lambda^\gamma v)$ be a rescaled version of $f$. A direct computation yields
\begin{align*}
\color{black}
\partial_t f_\lambda=\lambda^\alpha(\partial_tf)(\lambda^\alpha t,\lambda^\beta x, \lambda^\gamma v)&=\lambda^\alpha\left(\Delta_vf-\lambda^\gamma v\cdot\nabla_x f \right)(\lambda^\alpha t,\lambda^\beta x,\lambda^\gamma v) \\
&= \lambda^\alpha\left(\lambda^{-2\gamma}\Delta_v f_\lambda- \lambda^{\gamma-\beta}v\cdot\nabla_x f_\lambda \right).
\end{align*}
To ensure that the rescaled function $f_\lambda$ satisfies the same equation as $f$, we impose the scaling relations $\alpha-2\gamma=0$ and $\alpha+\gamma-\beta=0$. Choosing $\alpha=2$ leads to $\gamma=1$ and $\beta=3$. This scaling motivates the introduction of anisotropic Besov spaces with a $3:1$ scaling between the spatial and velocity variables,  see Definition \ref{def-ani-besov} for details. 

Informally, the operator $\Delta_v$ contributes a gain of two derivatives in the velocity direction. If this regularity is transferred to the time variable via the semigroup, one may gain at most one derivative in time. Similarly, transferring this regularity to the spatial variable yields at most a gain of $2/3$ of a derivative. In other words, the kinetic semigroup provides a total of two derivatives, distributed anisotropically among the variables. Let $r_t$, $r_x$, and $r_v$ denote the regularity gains in time, space, and velocity, respectively. Then we have the relation:
\begin{align*}
2r_t+3r_x+r_v=2.
\end{align*}
For instance, if $r_t=0$ and $r_v=1$, then the maximum possible spatial regularity gain is $r_x=\frac{1}{3}$. This heuristic underpins the choice of scaling parameters used in the scaling vector of the anisotropic Besov spaces. 

{\bf Frozen characteristic lines approach. }Based on the preceding analysis, our goal is to establish uniform regularity estimates for the approximation scheme of \eqref{SPDE-0} in anisotropic Besov spaces adapted to the space-velocity scaling. However, a new difficulty arises in this setting. We observe that the damping term $\nabla_v \cdot(vf)$ in \eqref{SPDE-0}, which originates from the damping structure $-V_i$ in the microscopic particle system \eqref{particles}, involves the function $\eta(v)=v$. This function does not belong to any Besov space due to its lack of integrability, which poses a technical obstruction to directly applying Besov space techniques. To address this issue, we developed a frozen characteristic line argument. 

The key idea is to exploit the rapid decay properties of the kernel associated with the kinetic semigroup in the velocity variable. Let $p$ denote the kernel of the kinetic semigroup. When analyzing the mild formulation of \eqref{SPDE-0}, one encounters integrals of the form
$$
\int_0^t \int_{\mathbb{R}^{2d}} p_{t-s}(x-y,v-w) \nabla_w \cdot (w f(s,y,w))dwdyds.
$$
Although the kernel $p$ decays rapidly in $v-w$, this is not sufficient to control the unboundedness of $w$, unless one takes $v = 0$. Motivated by this observation, we consider a shifted solution defined by
$$
f^{(x_0,v_0)}(t,x,v):=f(t,x+x_0,v+v_0), \text{ for every }(t,x,v)\in\mathbb{R}_+\times\mathbb{R}^{2d} \text{ and for some }(x_0,v_0)\in \mathbb{R}^{2d}, 
$$
and study the equation satisfied by $f^{(x_0, v_0)}$. However, the choice of the shift $(x_0, v_0)$ cannot be arbitrary. Due to the presence of the inhomogeneous coefficients $v$ in both the transport term $-v \cdot \nabla_x f$ and the damping term $\nabla_v \cdot (v f)$, the shifting operation introduces additional terms such as $-v_0 \cdot \nabla_x f$ and $\nabla_v \cdot (v_0 f)$. 

To retain the original structure of \eqref{SPDE-0} after shifting such that kinetic semigroup estimates can still be applied, we choose the shift $(x_0, v_0)$ to lie along the characteristic flow direction (see \eqref{shiftoperator}). By freezing the characteristic line and evaluating at $(x, v) = 0$, we obtain estimates for the Besov regularity of $f^{(x_0, v_0)}(t,0,0)$ with respect to $(x_0,v_0)$. Furthermore, comparing the Besov norms of $f^{(x_0, v_0)}(t,0,0)$ in the variables $(x_0,v_0)$ with those of $f(t, x, v)$ in the original variables $(x,v)$ (see Lemma \ref{lem-properties}), we are able to deduce the desired Besov regularity estimates.

This argument follows the prototype developed in \cite{HWZ20}, where the authors employed this approach to derive $L^{\infty}$ and H\"older regularity estimates for PDEs in the analysis of well-posedness for SDEs with inhomogeneous coefficients. In the present work, we extend this methodology to a general Besov framework and adapt it to our SPDE setting, with a refined analysis of the discrepancy between the Besov norms of the original variables and those of the transformed characteristic variables. These aspects will be discussed in detail in Lemma~\ref{lem-properties}, where we shall see that additional factors arise from the general $L^p$-norm and interact with the Littlewood--Paley block operators.

{\bf Two steps approximation scheme. }
{\color{black}To establish global-in-time existence, we introduce a two-step approximation scheme given by
\begin{align}\label{SPDE-iteration-intro}
\begin{split}
    df_{n,r}=\Big(&\Delta_v f_{n,r}+\frac{1}{n}\Delta_x f_{n,r}-\alpha_r^2(v)v\cdot\nabla_x f_{n,r}+\nabla_v\cdot(\alpha_r^2(v)v f_{n,r})-\nabla_v\cdot(f_{n,r}V\star_x\langle \alpha_r^2 f_{n,r}\rangle)\Big)dt\\
    &-\nabla_v\cdot(\sigma(f_{n,r})dW_F)+{\color{black}\frac{1}{2}\nabla_v\cdot\bigl(F_1\sigma'(f_{n,r})^2\nabla_v f_{n,r}\bigr)dt},
\end{split}
\end{align}
where $\alpha_r^2$, $r \geq 1$, denotes a sequence of velocity truncation functions. Owing to the spatial regularization $\frac{1}{n}\Delta_x$, the existence of solutions to \eqref{SPDE-iteration-intro} follows from a standard Galerkin approximation argument. This scheme serves as an approximation of \eqref{SPDE-0}.

Starting from \eqref{SPDE-iteration-intro}, we first establish the non-negativity of solutions, which in turn ensures the conservation of the $L^1(\mathbb{R}^{2d})$ mass and yields corresponding moment estimates. These properties allow us to prove strong tightness, relying crucially on the spatial regularization $\frac{1}{n}\Delta_x$, and to pass to the limit as $r \to \infty$. In this way, we obtain the existence of solutions to
\begin{align}\label{SPDE-L1-intro}
\begin{split}
    df_{n}=\Big(&\Delta_v f_{n}+\frac{1}{n}\Delta_x f_{n}-v\cdot\nabla_x f_{n}+\nabla_v\cdot(v f_{n})-\nabla_v\cdot(f_{n}V\star_x\langle f_{n}\rangle)\Big)dt\\
    &-\nabla_v\cdot(\sigma(f_{n})dW_F)+{\color{black}\frac{1}{2}\nabla_v\cdot\bigl(F_1\sigma'(f_{n})^2\nabla_v f_{n}\bigr)dt}.
\end{split}
\end{align}
Finally, we apply the aforementioned Besov regularity estimates, which are uniform in $n \geq 1$, to establish strong tightness in Lebesgue spaces and thereby conclude the existence of solutions to \eqref{SPDE-0}. 
}

{\color{black}{\bf Taming the square-root coefficient. }
Once we obtain the existence of \eqref{SPDE-0}, we are able to consider an approximation scheme
{\color{black}\begin{align}\label{SPDE-000}
\partial_tf_n
={\color{black}\Delta_vf_n}-v\cdot\nabla_xf_n
+\nabla_v\cdot(vf_n)
-\nabla_v\cdot\bigl(\sigma_n(f_n)\circ\xi_F\bigr),
\qquad (t,x,v)\in\mathbb R_+\times\mathbb R^{2d}.
\end{align}}
{ \color{black}where $(\sigma_n)_{n\geq1}$ is a sequence of smooth approximations of the square-root coefficient. In this part we set $V=0$. Following the zero-value truncation idea of \cite{FG24}, we combine a truncation near zero with kinetic semigroup estimates.}

In the context of VFPDK, we introduce a novel combination of this zero-value truncation technique with kinetic semigroup estimates. Specifically, we define a truncation $h_{\delta}(f_n)$ that cuts off values of $f_n$ in the interval $[0, \delta]$. We derive the equation satisfied by $u_n:=h_{\delta}(f_n)$ for every $n\geq1$:
{\color{black}\begin{align}\label{SPDE-un-intro}
	\partial_tu_n
 &={\color{black}\Delta_vu_n}-v\cdot\nabla_xu_n+\nabla_v\cdot(vu_n)\notag\\
 &\quad+g_{1,\delta,n}+g_{2,\delta,n}
 -h_{1,\delta,n}\xi_F{\color{black}-h_{2,\delta,n}\nabla_v\cdot\xi_F}.
\end{align}}
where 
\begin{align*}
	&{\color{black}g_{1,\delta,n} := \frac{1}{2} \nabla_v \cdot \bigl( h_{\delta}'(f_n) \sigma_n'(f_n)^2 \nabla_v f_n F_1 \bigr)}, \\
&g_{2,\delta,n} := - h_{\delta}''(f_n) |\nabla_v f_n|^2 + {\color{black}d h_{\delta}'(f_n) f_n - d h_{\delta}(f_n)} +\frac{1}{2}h_{\delta}''(f_n)\sigma_n(f_n)^2\, F_3,\\
    &h_{1,\delta,n}:=h_{\delta}'(f_n)\nabla_v\sigma_n(f_n),\quad h_{2,\delta,n}:=h_{\delta}'(f_n)\sigma_n(f_n). 
\end{align*}
{\color{black}We observe that \eqref{SPDE-un-intro} also takes the form of a generalized stochastic kinetic equation. However, in contrast to \eqref{SPDE-L1-intro}, the presence of the new forcing term $g_{2,\delta,n}$ in the space $L^1_{t,x,v}$ introduces additional analytical challenges.}

{\color{black}The resulting anisotropic Besov bounds are uniform in $n$ for each fixed $\delta>0$. They yield tightness on every bounded phase-space domain and, after a diagonal extraction, strong compactness only in
\begin{align*}
 L^1\bigl([0,T];L^1_{\rm loc}(\mathbb R^{2d})\bigr).
\end{align*}
No spatial moment estimate is used. Consequently, we neither obtain global $L^1$ convergence nor infer preservation of the total mass in the limit. The local compactness is sufficient to pass to the local terms of the equation with $V=0$ and to construct a probabilistically weak solution on a new stochastic basis.}

{\color{black}{\bf Scope of the square-root result. } In the absence of a valid spatial moment estimate, the square-root result is restricted to $V=0$. The limiting process is not obtained on the original stochastic basis, and no uniqueness or $C([0,T];L^1)$ modification is asserted.}

{\color{black}{\bf Renormalized kinetic solutions. } For regular coefficients, the global $L^1$ theory and the doubling-of-variables method yield uniqueness in the mass-preserving class. For the square-root equation, the available compactness is local and the limiting kinetic measure may have atoms in time. We therefore use a weaker, distribution-in-time kinetic formulation and make no claim of uniqueness.}
}

\subsection{Second-order mean-field systems and fluctuations} 
The second-order mean-field system \eqref{particles} serves as a microscopic model for a wide range of phenomena. In this system, $X$ denotes the position of a particle and $V$ denotes its velocity, so the first equation in \eqref{particles} reflects the basic kinematic relation $\dot{X}=V$. According to Newton's second law, the acceleration is determined by the total external forces acting on the particle. These forces comprise the interactions with other particles, frictional forces, and random fluctuations. The friction is modeled by a linear dissipative term $-\kappa V$ in \eqref{particles}. The stochastic forcing accounts for small-scale random perturbations that are not captured deterministically; to ensure consistency with the fluctuation-dissipation theorem, these random perturbations are introduced as energy-conserving noise with intensity $\sqrt{\kappa}$. By choosing different forms for the interaction term $\frac{1}{N}\sum_{j=1}^N V(X_i-X_j)$, the system \eqref{particles} can be adapted to model various physical and biological systems, including fluid dynamics, plasma physics, and collective behavior in biology; see  \cite{C43,ST08} for further references.   

The large-$N$ mean-field limits of the particle system, together with their associated fluctuations, have been the subject of extensive study in both the physics and mathematics communities. We refer the reader to \cite{Sznitman06,BGM10,M01,M03} for results concerning convergence to the mean-field limit, and to \cite{BD25,MS96,WZZ22,CG25} for works on fluctuation phenomena. We note that this body of literature addresses not only second-order particle systems, but also first-order particle systems. 

At the mesoscopic level, fluctuating hydrodynamics equations serve as analogues that capture similar fluctuation phenomena. Specifically, consider the zero-noise fluctuation regime of 
\begin{align}\label{SPDE-ep}
{\color{black}\partial_tf_{\varepsilon} =\Delta_vf_{\varepsilon}-v\cdot\nabla_xf_{\varepsilon}-\nabla_v\cdot \big(f_{\varepsilon}V\star_x\rho_{\varepsilon}\big)+\nabla_v\cdot (vf_{\varepsilon})-\varepsilon^{1/2}\nabla_v\cdot(\sqrt{f_{\varepsilon}} \circ\xi_{\delta(\varepsilon)}).} 
\end{align}
under singular scaling limits $(\varepsilon,\delta(\varepsilon)) \to (0, 0)$, where the space-velocity correlated noise converges to space-velocity-time white noise, $\xi_{\delta(\varepsilon)} \to \xi$. Informally, one expects the same limiting Gaussian SPDE to arise in the central limit regime for the empirical measure of \eqref{particles}, as well as the same large deviation rate function, see  \cite[(1.20)]{BD25} and \cite{DPZ13} for explanations. Therefore, a rigorous derivation of the zero-noise fluctuations in \eqref{SPDE-ep} would support the consistency between microscopic particle models and mesoscopic fluctuating hydrodynamic equations. Establishing this connection remains a subject of ongoing and future investigation.

\subsection{Comments on the literature} 

\ 

{\bf Vlasov-Fokker-Planck equations. }
{\color{black}According to Kac's program on the propagation of chaos, the limiting behavior of the empirical measure $\pi_N$ associated with the particle system \eqref{particles} is governed, for a general interaction force $V$, by the Vlasov--Fokker--Planck equation:}
\begin{align}\label{VPFP}
    {\color{black}\p_t f=\frac{\kappa}{2}\Delta_v f-v\cdot \nabla_x f-\nabla_v\cdot(f V\star_x\rho)+\kappa\nabla_v\cdot(v f).}
\end{align}
{\color{black}In three spatial dimensions, setting $\kappa=0$ corresponds to the inviscid case. If $V$ is taken to be the Newtonian or Coulomb force kernel, namely the negative gradient of the corresponding interaction potential, equation \eqref{VPFP} reduces to the classical Vlasov--Poisson equation.} This model describes the evolution of the particle distribution 
  and arises naturally in the contexts of galactic dynamics and plasma physics (see \cite[Section 1]{HRZ25} and references therein).

The global existence of smooth solutions to the Vlasov-Poisson equation under appropriate moment assumptions on the initial data was established in \cite{LP91, Sch91, Pf92}. The existence of classical solutions to the VPFP equation has been demonstrated in \cite{Bo93, CS97, OS00}, employing methodologies analogous to those in \cite{LP91}. Furthermore, when the Laplacian operator is replaced by its fractional counterpart, well-posedness results for the corresponding nonlocal VPFP-type equation have been recently obtained in \cite{HRZ25}.

We conclude with a brief overview of the literature concerning the long-time behavior of the VFP equation. This topic has been extensively studied, and the results reveal a clear dichotomy between convex and non-convex confining potentials. In the convex setting (and under sufficiently regular, weak interactions), one typically obtains uniqueness of the invariant measure together with exponential relaxation to equilibrium. Such results are established through hypocoercivity and spectral-gap techniques (see Villani \cite{Villani09}, H\'erau \cite{HF07}, Dolbeault-Mouhot-Schmeiser \cite{DMS15}). In addition, entropy dissipation methods based on log-Sobolev or Poincar\'e inequalities have been employed to derive explicit exponential convergence rates in Wasserstein or entropy metrics for weakly interacting kinetic Fokker-Planck equations (cf. \cite{BF17}). In contrast, in the non-convex setting, the presence of multiple wells may lead to several invariant measures and phase transitions. A common approach is to exploit a free-energy (Lyapunov) functional adapted to the kinetic framework, which ensures the existence of limit points and allows one to characterize invariant measures; however, such methods typically yield convergence (often only subsequential) without quantitative rates or selection principles among equilibria. Duong and Tugaut \cite{DT18}, combining free-energy techniques with tools from the McKean-Vlasov literature, provide further insight into this regime. We also refer the reader to \cite{BFL24,TG25,HR18} for additional developments.

{\bf Kinetic semigroup estimates. }

The kinetic operator $\p_t - \Delta_v + v \cdot \nabla_x$ plays a central role in various physical models, including the Vlasov-Poisson-Fokker-Planck and Landau equations. Consequently, the study of its associated semigroup has attracted significant attention. In \cite{Ko34}, Kolmogorov first identified the kinetic semigroup $P_t$ as being characterized by the density of the stochastic process $(\int_0^t B_s \dif s, B_t)$, as explicitly expressed in equations \eqref{k-semigroup} and \eqref{density} below. Later developments, such as in \cite{Me11}, provided heat kernel estimates for kinetic stochastic differential equations (SDEs) driven by multiplicative noise. For a comprehensive overview of the literature concerning heat kernel estimates, we refer the reader to \cite{RZ25} and the references therein.

In \cite{Bo02}, the following anisotropic regularity estimate was established:
\begin{align*}
    \|\Delta_x^{\frac13} f\|_{L^2_{t,x,v}} \lesssim \|\Delta_v f\|_{L^2_{t,x,v}}^{\frac13} \|(\p_t - v \cdot \nabla_x)f\|_{L^2_{t,x,v}}^{\frac23},
\end{align*}
which implies that the spatial regularity in the $x$-variable, quantified by $\Delta_x^{\frac13}$, can be inferred from the velocity regularity $\Delta_v f$ and the action of the kinetic transport operator $(\p_t - v \cdot \nabla_x)f$. This estimate rigorously confirms the scaling analysis presented in Section \ref{subsec-1-2}. Building upon this framework and using an anisotropic distance, anisotropic Besov norm estimates for the kinetic semigroup were developed in \cite{HZZZ24, HRZ25}, including the following key inequality:
\begin{align}\label{in:semi-est}
    \|P_t f\|_{\bB^{\alpha+\beta}_{{\bf p};\theta}} \lesssim t^{-\frac{\alpha}{2}} \|f\|_{\bB^{\beta}_{{\bf p};\theta}}, \quad \alpha \ge 0,\ \beta \in \mR,\ t > 0,
\end{align}
where the anisotropic Besov space $\bB^{\alpha+\beta}_{{\bf p};\theta}$ and the scaling index $\theta$ are defined in Definition \ref{def-ani-besov} and equation \eqref{parameters}, respectively.

Schauder-type estimates for (fractional) kinetic operators have also been developed in recent works such as \cite{HWZ20, CHM21, IS21}. In the present paper, we make use of a variant of Lemma 5.1 from \cite{HWZ20}; see Section \ref{sec-4} for details.

Notably, in \cite{ZZ24}, based on the semigroup estimate \eqref{in:semi-est}, the authors derived Schauder estimates for the following class of kinetic SPDEs:
\begin{align*}
    \dif f = (a_{ij}\p_{v_i}\p_{v_j} + v \cdot \nabla_x + b \cdot \nabla_v)f \dif t + g_1 \dif t + (\sigma \cdot \nabla_v f + g_2)\dif W_t,
\end{align*}
where $W_t$ is a finite-dimensional Brownian motion, $2a - \sigma^T\sigma$ is uniformly elliptic, and the coefficients $a, b, \sigma$ possess anisotropic H\"older regularity uniformly in $\omega \in \Omega$. However, their framework does not allow for nonlinear dependencies in $\sigma = \sigma(f)$, while certain nonlinearities such as $g_2 = g_2(f)$ are permissible, the structure $b = V \star_x \langle f \rangle$ is excluded. As a result, their results do not cover the nonlinear SPDEs given in equations \eqref{SPDE-00} and \eqref{SPDE-000}.

Finally, we also refer to \cite{KW24} and the references therein for discussions on the Harnack inequality in the context of kinetic equations.

{\bf Conservative SPDEs. }
Significant progress has been made in the analytical study of the Dean-Kawasaki equation and related stochastic PDEs through the development of sophisticated mathematical techniques. A foundational contribution was made by Debussche and Vovelle \cite{DV10}, who introduced the framework of kinetic solutions for addressing the Cauchy problem in stochastic conservation laws. This concept was later extended to more general parabolic-hyperbolic SPDEs with conservative noise by Gess and Souganidis \cite{GS17}, Fehrman and Gess \cite{FG19}, and Dareiotis and Gess \cite{DG20}. These approaches have proven instrumental in the study of Dean-Kawasaki-type equations.

For models with local interactions, Fehrman and Gess \cite{FG24} established the well-posedness of functional-valued solutions to the Dean-Kawasaki equation driven by correlated noise. Expanding on this work, they investigated the asymptotic behavior of the system under vanishing noise, establishing a large deviation principle in \cite{FG23}. Later on, they extend the results to the whole-space case \cite{FG25}, and Fehrman extends the well-posedness result to the case of non-stationary Stratonovich noise \cite{F25}. More recently, Clini and Fehrman {\color{black}\cite{CliniFehrman23}} derived a central limit theorem for the nonlinear Dean-Kawasaki equation with correlated noise, while Gess, the second author, and Zhang \cite{GWZ24} developed higher-order fluctuation expansions for this class of systems.

In the nonlocal interaction regime, Wang, the second author, and Zhang \cite{WWZ22}, as well as the second author and Zhang \cite{WZ24}, analyzed the well-posedness and large deviation behavior of Dean-Kawasaki equations with singular interaction potentials, with particular applications to the fluctuating Ising-Kac-Kawasaki model. For models with Dirichlet boundary conditions, Popat \cite{Shyam25,Shyam25-fluc} established well-posedness, large deviation, and fluctuation results. For discussions on weak error estimates, we refer the reader to \cite{DKP24,DJP25,CF23,CFIR23} and related works. Complementary to these results, Martini and Mayorcas \cite{AA25,AA24} established well-posedness and large deviation principles for an additive-noise approximation of the Keller-Segel-type Dean-Kawasaki equation. Further work concerning the regularized Dean-Kawasaki equation of second-order particle systems is considered in \cite{CS23,CSZ19,CSZ20}, and has been applied to computational chemistry \cite{JLR25}.

\subsection{Structure of the paper} 
{\color{black}This paper is organized as follows. In Section~\ref{sec-2}, we introduce the notation, assumptions, and solution concepts. Section~\ref{sec-3} proves uniqueness in the global mass-preserving solution class used for regular coefficients. Section~\ref{sec-4} develops kinetic semigroup estimates. In Section~\ref{sec-5}, we prove global probabilistically strong well-posedness for the regular-coefficient equation with bounded interaction. Finally, Section~\ref{sec-6} combines zero-value truncation with kinetic semigroup estimates to construct a probabilistically weak solution for the square-root equation with $V=0$, using only local $L^1$ compactness.}

\section{Preliminaries}\label{sec-2}
\subsection{Notations}
Throughout this work, we fix a time horizon $T > 0$. Let $(\Omega,\mathcal{F},\mathbb{P},\{\mathcal{F}_t\}_{t \in [0,T]})$ be a stochastic basis. {\color{black}We assume throughout that the filtration is complete and right-continuous. Let $\{B^k(t)\}_{t\in[0,T]},\ k\in\mathbb{N},$ be a sequence of independent $\{\mathcal{F}_t\}_{t \in [0,T]}$-adapted $d$-dimensional Wiener processes.} Expectations with respect to the probability measure $\mathbb{P}$ are denoted by $\mathbb{E}$. 

For every $p\in[1,\infty]$, we denote by $\|\cdot\|_{L^p(\mathbb{R}^{2d})}$ the norm in the Lebesgue space $L^p(\mathbb{R}^{2d})$ (or $L^p(\mathbb{R}^{2d}; \mathbb{R}^m)$ with $m\in\mN_+$ when vector-valued functions are considered). We further denote by $L^p_{\mathrm{loc}}(\mathbb{R}^{2d})$ the space of locally $p$-integrable functions. The inner product in $L^2(\mathbb{R}^{2d})$ is denoted by $\langle\cdot,\cdot\rangle$. We write $C^\infty(\mathbb{R}^{2d} \times (0,\infty))$ for the space of infinitely differentiable functions on $\mathbb{R}^{2d} \times (0,\infty)$, and $C_c^\infty(\mathbb{R}^{2d} \times (0,\infty))$ for its subspace of functions with compact support. Furthermore, for every $k \in \mathbb{N}$, we denote by $C^k(\mathbb{R}^{2d})$ the space of $k$-times continuously differentiable functions, and by $C^k_{\mathrm{loc}}(\mathbb{R}^{2d})$ the space of $k$-times continuously differentiable functions, endowed with locally $C^k$-convergence topology. Throughout this work, we consistently use $x$ to denote the spatial variable and $v$ the velocity variable. For convenience, we write $z = (x, v) \in \mathbb{R}^{2d}$ to represent the pair of spatial and velocity variables. 

Moreover, we define the following differential operators:
\begin{align*}
\nabla_x f(z) := \nabla f(\cdot, v)(x), \ \nabla_v f(z) := \nabla f(x, \cdot)(v), \
\Delta_x := \nabla_x\cdot\nabla_x, \ \Delta_v := \nabla_v\cdot\nabla_v,
\end{align*}
and for any $\bar{z} = (\bar{x}, \bar{v})$,
\begin{align*}
v \cdot \nabla_x f(\bar{z}) := \bar{v} \cdot \nabla_x f(\bar{z}).
\end{align*}
We emphasize that the symbols $(x, v)$ appearing in $\nabla_x$, $\nabla_v$, $\Delta_x$, $\Delta_v$, and $v \cdot \nabla_x$ are not variables but rather serve as notational placeholders to indicate which argument the differential operator acts upon.

\subsection{The correlated noise}
We consider a spatial-velocity correlated Brownian motion defined by
\begin{align*}
W_F(t,x,v) = \sum_{k\geq1} B^k(t) f_k(x,v), \quad (t,x,v) \in \mathbb{R}_+ \times \mathbb{R}^{2d},
\end{align*}
where $(B^k)_{k\geq1}$ is a sequence of $d$-dimensional Brownian motions and $(f_k)_{k \geq 1}$ is a sequence of smooth functions such that the following quantities are well-defined scalar or vector fields:
\begin{align*}
F_1 = \sum_{k \geq 1} f_k^2, \quad 
F_2 = \sum_{k \geq 1} f_k \nabla_v f_k, \quad 
F_3 = \sum_{k \geq 1} |\nabla_v f_k|^2, \quad 
F_4 = \sum_{k \geq 1} f_k \Delta_v f_k.
\end{align*}
Furthermore, we assume that $F_1,F_2\in L^{\infty}(\mathbb{R}^{2d})$, $F_3\in L^{\infty}(\mathbb{R}^{2d})\cap L^1(\mathbb{R}^{2d})$, and the following compatibility condition holds:
\begin{align*}
\nabla_v \cdot F_2 = F_3 + F_4 = 0.
\end{align*} 
{\color{black}As part of the noise hypothesis, we require the series above to define the displayed fields locally in $L^1$ and impose the distributional identities
\begin{align}\label{noise-field-identities}
 \nabla_vF_1=2F_2,
 \qquad
 \nabla_v\cdot F_2=F_3+F_4=0
 \quad\text{in }\mathcal D'(\mathbb R^{2d}).
\end{align}
These identities specify unambiguously how derivatives of the covariance fields are to be interpreted; no pointwise differentiation of a merely $L^\infty$ representative will be used.}
{\color{black}We record that the hypotheses just imposed actually imply
$F_2=0$. Indeed, in the sense of distributions,
\begin{align*}
 \Delta_vF_1=2\nabla_v\cdot F_2=0.
\end{align*}
Since $F_1\in L^\infty(\mathbb R^{2d})$, for almost every fixed $x$ the
function $F_1(x,\cdot)$ is a bounded distributional harmonic function on
$\mathbb R_v^d$ and is therefore constant by the Liouville theorem. Hence
$\nabla_vF_1=0$ and $F_2=0$. {\color{black}Accordingly, all subsequent
equations, estimates, and proofs are written under the standing consequences
$\nabla_vF_1=0$ and $F_2=0$; no vanishing $F_2$-term is retained. In
particular, none of the well-posedness results below is claimed for a noise
with $F_2\neq0$; such a result would require a genuinely different noise
assumption and a separate verification of the estimates.}}
\begin{remark}
As explained in the introduction, the study of fluctuation behavior and large deviation principles in the vanishing noise intensity and correlation limits is of particular interest for Dean-Kawasaki type equations. The main well-posedness results established in this paper provide a foundational step toward such analyses. In this context, we introduce the correlated noise structure above as a model for the truncated space-time white noise. For concrete examples of such noise structures on the torus, we refer the reader to \cite{FG24,DFG}, while examples on the whole Euclidean space can be found in \cite{FG25}.
\end{remark}

\subsection{Anisotropic Besov spaces} 
Let $n,N\in\mathbb{N}_+$, $m=(m_1,...,m_n)\in\mathbb{N}^n$ with $m_1+...+m_n=N$, and let $\theta=(\theta_1,...,\theta_n)\in[1,+\infty)^n$. We define $\theta\cdot m=\theta_1m_1+...+\theta_nm_n$. For every $x=(x_1,...,x_n), y=(y_1,...,y_n)\in\mathbb{R}^{m_1}\times...\times\mathbb{R}^{m_n}$, set 
\begin{align*}
|x-y|_{\theta}:=\sum_{i=1}^n|x_i-y_i|^{1/\theta_i}. 	
\end{align*}
For every $y\in\mathbb{R}^N$ and $r>0$, let $B^{\theta}_r(y)=\{x\in\mathbb{R}^N:|x-y|_{\theta}\leq r\}$ and $B^{\theta}_r:=B^{\theta}_r(0)$.
Let $\chi^{\theta}$ be a symmetric nonnegative $C^{\infty}(\mathbb{R}^N)$ function with 
\begin{align*}
\chi^{\theta}(\xi)=1, \ \text{for }\xi\in B^{\theta}_1,\ \chi^{\theta}(\xi)=0,\ \text{for }\xi\notin B^{\theta}_{4/3}.\end{align*}
{\color{black}Then for every $\xi=(\xi_1,...,\xi_n)\in\mathbb{R}^{m_1}\times...\times\mathbb{R}^{m_n}$ and every integer $j\geq-1$, we define}
\begin{align*}
\phi^{\theta}_j(\xi):=
\begin{cases}
\chi^{\theta}(\xi),\quad &j=-1,\\
\chi^{\theta}(2^{-\theta (j+1)}\xi)-\chi^{\theta}(2^{-\theta j}\xi), \quad &j\ge0,
\end{cases}
\end{align*}
where $2^{-\theta j}\xi:=(2^{-\theta_1j}\xi,...,2^{-\theta_nj}\xi_n)$. This definition implies that for every $j\in\mathbb{N}_0$, 
\begin{align*}
\phi^{\theta}_j(\xi)=\phi^{\theta}_0(2^{-\theta j}\xi)
\end{align*}
and 
\begin{align*}
{\rm supp}\ \phi^{\theta}_j\subset B^{\theta}_{2^{j+3}/3}\backslash B^{\theta}_{2^{j}},\ \ \sum_{j=-1}^n\phi^{\theta}_j(\xi)=\chi^{\theta}(2^{-(n+1)}\xi)\rightarrow1,\ \text{as }n\rightarrow\infty. 	
\end{align*}
In the following, we will introduce the Block operators with respect to the anisotropic dyadic decomposition. Let $\mathcal{S}$ be the space of Schwartz functions on $\mathbb{R}^N$, $\mathcal{S}'$ denotes the tempered distribution space. For every $j\geq-1$, we define 
\begin{align}\label{block-operator}
\mathcal{R}^{\theta}_jf:=\mathcal{F}^{-1}(\phi^{\theta}_j\mathcal{F}(f))=\mathcal{F}^{-1}(\phi^{\theta}_j)\ast f,\ \ \text{for every }f\in\mathcal{S}', 	
\end{align}
where $\mathcal{F}$ and $\mathcal{F}^{-1}$ denote the Fourier transform and its inverse, respectively. Precisely, for every $f\in\mathcal{S}$, 
$$
\mathcal{F}(f)(\xi):=\frac{1}{(2\pi)^{N/2}}\int_{\mathbb{R}^{N}}e^{-i\xi \cdot z}f(z)dz,\quad \xi\in\mathbb{R}^{N}.
$$
We also use $\hat{f}$ and $\check{f}$ to denote the Fourier transform $\mathcal{F}(f)$ and the inverse Fourier transform $\mathcal{F}^{-1}(f)$, respectively. Here, the Fourier transform of tempered distributions is defined by extension from the Fourier transform on $L^1(\mathbb{R}^{N})$ functions, and the convolution $\ast$ is understood in the sense of distributions. For further details on the definition of block operators on tempered distributions, we refer the reader to \cite{BC11}. In particular, we have 
\begin{align*}
\langle\mathcal{R}^{\theta}_jf,g\rangle=\langle f,\mathcal{R}^{\theta}_jg\rangle,\ \ \text{for every }f\in\mathcal{S}', g\in\mathcal{S}. 	
\end{align*}
For the sake of clarity and consistency in notation, we fix the model parameters as follows:
\begin{equation}\label{parameters}
	N = 2d,\quad n = 2,\quad m_1 = m_2 = d,\quad \theta = (3,1).
\end{equation}
In the following, we will introduce the anisotropic Besov spaces. 
\begin{definition}\label{def-ani-besov}
For every $s\in\mathbb{R}$, $p\in[1,+\infty]$, let 
\begin{align*}
\mathbf{B}^s_{p;\theta}:=\Big\{f\in\mathcal{S}':\|f\|_{\mathbf{B}^s_{p;\theta}}:=\sup_{j\geq-1}(2^{sj}\|\mathcal{R}^{\theta}_jf\|_{L^p(\mathbb{R}^{2d})})<\infty\Big\}. 	
\end{align*}
When $\theta=(1,...,1)$, we denote 
\begin{align*}
\mathbf{B}^s_p:=\mathbf{B}^s_{p;\theta},\ \mathcal{R}_j:=\mathcal{R}^{\theta}_j,\ \phi_j:=\phi^{\theta}_j.	
\end{align*}
\end{definition}
From the definition, we observe that for every $j\geq-1$, 
\begin{equation}\label{Besov-basic-inequality}
	\|\mathcal{R}^{\theta}_jf\|_{L^p(\mathbb{R}^{2d})}\leq 2^{-sj}\|f\|_{\mathbf{B}^s_{p;\theta}},\quad\text{for every }f\in\mathbf{B}^s_{p;\theta}. 
\end{equation}

In the following, we introduce a Bernstein type lemma \cite[Lemma 2.1]{HRZ25} for the anisotropic block operators. 
\begin{lemma}\label{Bernstein}
For every $k_1,k_2\in\mN_0$ and $1\le p\le q\le\infty$, there exists a constant $C=C(k_1,k_2,p,q,d)>0$ such that for every $j\geq-1$, 
\begin{align*}
\|\nabla^{k_1}_x\nabla^{k_2}_v\mathcal{R}^{\theta}_jf\|_{L^{q}(\mathbb{R}^{2d})}\leq C2^{j(3k_1+k_2+\frac{4d}{p}-\frac{4d}{q})}\|\mathcal{R}^{\theta}_jf\|_{L^{p}(\mathbb{R}^{2d})}. 	
\end{align*}	
\end{lemma}

\subsection{Assumptions: coefficients, initial data, interaction kernels} 
In this part, we present all the assumptions in this paper. 

{\bf Coefficients.} We first present assumptions for both regular and irregular coefficients. 
\begin{assumption}\label{Assump-sigma-n}
	\textbf{(Regular coefficients)} Assume that $\sigma(\cdot)\in\mathrm{C}_{\text{loc}}^{1}((0,\infty))\cap C_c((0,\infty))$. Furthermore, $\sigma$ has the following properties.
	 \begin{enumerate}
	 \item[(1)] $\sigma\in C(\left[0,\infty\right))\cap C^{\infty}((0,\infty))$ with $\sigma(0)=0$ and $\sigma'\in C^{\infty}_c(\left[0,\infty\right))$,
	 \item[(2)] there exists $c\in(0,\infty)$ such that for every $\zeta\in[0,\infty)$,
	 \begin{align}
	 	|\sigma(\zeta)|\le c\sqrt\zeta,\quad |\sigma(\zeta)\sigma'(\zeta)|\leq c,
	 	\quad\text{and}\quad |(\sigma(\zeta)\sigma'(\zeta))'|\leq c.
	 \end{align}
		{\color{black}In addition, for every $\zeta>0$,
		$|\sigma'(\zeta)|^2\leq c/\zeta$.}
	\end{enumerate}
\end{assumption}

\begin{assumption}\label{Assump-sigma}
	\textbf{(Irregular coefficients)} Assume that $\sigma\in C([0,\infty))$ satisfy that $\sigma\in C^{1}_{loc}((0,\infty))$ and $\sigma\sigma'\in C([0,\infty))$. Furthermore, there exists a constant $c>0$ such that $\limsup_{\zeta\rightarrow0+}\frac{\sigma(\zeta)^2}{\zeta}\leq c$, and there exists a constant $c'\in[1,\infty)$ such that 
	\begin{align}
	\sup_{\zeta'\in[0,\zeta]}\sigma(\zeta')^2\leq c'(1+\zeta+\sigma(\zeta)^2),	
	\end{align}
 for every $\zeta\in[0,\infty)$. 

\end{assumption}

{\color{black}We remark that the square-root coefficient $\sigma(\zeta)=\sqrt\zeta$ satisfies Assumption~\ref{Assump-sigma}. Assumption~\ref{Assump-sigma-n} is used both for the regular equation and for smooth approximations of the square-root coefficient. {\color{black}The uniqueness theorem in Section~\ref{sec-3} concerns regular weak solutions in the sense of Definition~\ref{def-weaksolution}; Proposition~\ref{equivalence} supplies their renormalized kinetic formulations.} The probabilistically weak solution constructed in Section~\ref{sec-6} is not asserted to belong to that class, and no uniqueness conclusion is drawn there. In accordance with \cite[Lemma 5.18]{FG24}, we use the following approximation of the square-root coefficient.}

\begin{lemma}\label{lem-sigma}
Let $\sigma(\zeta)=\sqrt{\zeta}$, for every $\zeta\in[0,\infty)$. Then $\sigma(\cdot)$ satisfies Assumption \ref{Assump-sigma}. Furthermore,  there exists a sequence of regular coefficients $(\sigma_n(\cdot))_{n\geq1}$ satisfying Assumption \ref{Assump-sigma-n}, such that $\sigma_n(\cdot)\rightarrow\sigma(\cdot)$ in $C^1_{loc}((0,\infty))$, as $n\rightarrow\infty$. 

{\color{black}Moreover, the sequence can be chosen so that there is a constant $C>0$, independent of $n$, for which, for every $n\geq1$ and $\zeta>0$,
\begin{align}\label{uniform-sigma-approximation-bounds}
 |\sigma_n(\zeta)|&\leq C\sqrt\zeta,
 &|\sigma_n(\zeta)\sigma_n'(\zeta)|&\leq C,
 &|\sigma_n'(\zeta)|^2&\leq \frac{C}{\zeta}.
\end{align}
The bound on $(\sigma_n\sigma_n')'$ required in
Assumption~\ref{Assump-sigma-n} is imposed separately for each $n$; its
constant is not uniform in $n$. Such a uniform derivative bound would be
incompatible with convergence to $\sqrt\zeta$ at the origin.}
	
\end{lemma}

{\bf Initial data.} We begin by introducing appropriate function spaces for the initial data. To describe regular initial data, we define
\begin{align}\label{Ini}
\operatorname{Ini}(\mathbb{R}^{2d}) := \left\{ f_0 \in L^1(\mathbb{R}^{2d}) \cap L^2(\mathbb{R}^{2d}) : f_0\geq0\right\},
\end{align}
equipped with the strong topology of $L^1(\mathbb{R}^{2d})$. 


{\color{black}For the regular-coefficient theory, we use the Boltzmann entropy density
\begin{align}\label{entropy-func}
 \Psi(\zeta):=
 \begin{cases}
  \zeta\log\zeta-\zeta,&\zeta>0,\\
  0,&\zeta=0.
 \end{cases}
\end{align}
The $L^2$ and second-moment conditions in Assumption~\ref{Assump-initialdata-regular}
imply that $\int_{\mathbb R^{2d}}\Psi(f_0)\,dz$ is well defined and integrable.
Indeed, its positive part is controlled by $f_0^2$, whereas its negative
part is controlled by the mass and the second spatial--velocity moment after
using a Gaussian weight. More precisely, every nonnegative measurable $h$
satisfies
\begin{align}\label{entropy-negative-part-bound}
 \int_{\mathbb R^{2d}}\bigl(\Psi(h(z))\bigr)^-\,dz
 \leq \int_{\mathbb R^{2d}}\bigl(1+|x|^2+|v|^2\bigr)h(z)\,dz
 +\frac1e\int_{\mathbb R^{2d}}e^{-(|x|^2+|v|^2)}\,dz.
\end{align}}

\begin{assumption}\label{Assump-initialdata-regular}
\textbf{(Initial data with finite moment estimates)} We assume that $f_0$ is an {\color{black}$\mathcal{F}_0$-measurable} random variable such that $f_0 \in \operatorname{Ini}(\mathbb{R}^{2d})$ almost surely. Moreover,
\begin{align*}
{\color{black}\operatorname*{ess\,sup}_{\omega\in\Omega}}\|f_0(\omega)\|_{L^1(\mathbb{R}^{2d})} 
+ \mathbb{E} \|f_0\|_{L^2(\mathbb{R}^{2d})}^2 
{\color{black}+\mathbb{E}\int_{\mathbb{R}^{2d}}(|v|^2+|x|^2)f_0\, dz} < \infty.
\end{align*}
\end{assumption}

\begin{assumption}\label{Assump-initialdata-entropy}
{\color{black}\textbf{($L^1\cap L^2$ initial data without moment assumptions)} We assume that $f_0$ is an $\mathcal{F}_0$-measurable random variable such that $f_0 \in \operatorname{Ini}(\mathbb{R}^{2d})$ almost surely. Furthermore,}
\begin{align*}
{\color{black}\operatorname*{ess\,sup}_{\omega\in\Omega}}\|f_0(\omega)\|_{L^1(\mathbb{R}^{2d})}+\mathbb{E} \|f_0\|_{L^2(\mathbb{R}^{2d})}^2 < \infty.
\end{align*}
\end{assumption}

{\color{black}Assumption~\ref{Assump-initialdata-regular} is imposed in the regular-coefficient well-posedness theory. Assumption~\ref{Assump-initialdata-entropy}, despite the legacy label, contains no entropy or moment condition and is used for the square-root existence result with $V=0$.}

{\bf Interaction kernels.} Finally, we specify the assumption on the interaction field.

\begin{assumption}\label{Assump-ker}
We assume that the interaction kernel $V\colon \mathbb{R}^d \to \mathbb{R}^d$ depends only on the spatial variable $x$, and satisfies $V \in L^{\infty}\cap L^2(\mathbb{R}^d)$.
\end{assumption}

%
%
\subsection{Definitions of solutions}

We begin by introducing the definition of renormalized kinetic solutions to \eqref{SPDE-0}. According to a straightforward computation, the equation \eqref{SPDE-0} can be rewritten in the It\^o form as follows:
\begin{align}\label{SPDE-0-ito}
df= \Big(\Delta_v f - v \cdot \nabla_x f - \nabla_v \cdot (f V \star_x &\rho) + \nabla_v \cdot (v f)\Big) dt - \nabla_v \cdot (\sigma(f) dW_F) \notag \\
&+ \frac{1}{2} \sum_{k\geq1}\nabla_v \cdot \left(\sigma'(f)f_k \nabla_v (f_k\sigma(f))\right)dt\notag\\
= \Big(\Delta_v f - v \cdot \nabla_x f - \nabla_v \cdot (f V \star_x &\rho) + \nabla_v \cdot (v f)\Big) dt - \nabla_v \cdot (\sigma(f) dW_F) \notag \\
&+{\color{black}\frac{1}{2}\nabla_v\cdot\left(
 F_1\sigma'(f)^2\nabla_vf
 \right)dt.}
\end{align}
Here we comment that we drop the factor $\frac{1}{2}$ in \eqref{SPDE-0}, as it does not affect the main well-posedness result. {\color{black}The last term in \eqref{SPDE-0-ito} is deliberately kept in divergence form. Indeed,
\begin{align*}
 \sum_{k\geq1}\sigma'(f)f_k\nabla_v\bigl(f_k\sigma(f)\bigr)
 ={\color{black}F_1\sigma'(f)^2\nabla_vf}.
\end{align*}
Whenever the expanded form is used for a smooth solution, \eqref{noise-field-identities} gives an equivalent expression.}

In the following, we present a formal computation that motivates the definition of renormalized kinetic solutions. Assuming that the solution $f$ is nonnegative, we define the kinetic function $\chi(x,v,\zeta,t)$ by 
$$
\chi(x,v,\zeta,t) = I_{\{0 < \zeta < f(x,v,t)\}} - I_{\{f(x,v,t) < \zeta < 0\}} = I_{\{0 < \zeta < f(x,v,t)\}},
$$
and observe that, in the sense of distributions, the following identities hold:
$$
\nabla_x \chi = \delta_{f=\zeta}\nabla_x f, \quad 
\nabla_v \chi = \delta_{f=\zeta}\nabla_v f, \quad 
\partial_{\zeta} \chi = \delta_0(\zeta) - \delta_{f=\zeta}.
$$

{\color{black}As a formal computation, let $s\in C_c^\infty((0,+\infty))$ and set
$S(\zeta):=\int_0^\zeta s(r)\,dr$. Thus $S(0)=0$, $S'=s$, and $S$ is smooth and bounded (although it need not have compact support).} The application of It\^o's formula to $S(f)$ yields
\begin{align*}
dS(f) ={}& S'(f)\left(\Delta_v f - v \cdot \nabla_x f - \nabla_v \cdot(f V \star_x \rho) + \nabla_v \cdot(v f)\right) dt - S'(f)\nabla_v \cdot(\sigma(f) dW_F) \notag \\
&+ {\color{black}\frac{1}{2} \nabla_v \cdot \left(F_1 \sigma'(f)^2 S'(f) \nabla_v f\right) dt} \\
&+ {\color{black}\frac{1}{2} S''(f) \sigma(f)^2 F_3\,dt} \\
={}& \left\langle S', \delta_{f=\zeta} \left(\Delta_v f - v \cdot \nabla_x f - \nabla_v \cdot(f V \star_x \rho) + \nabla_v \cdot(v f)\right) \right\rangle dt \\
&- \left\langle S', \delta_{f=\zeta} \nabla_v \cdot(\sigma(f) dW_F) \right\rangle \\
&+ {\color{black}\left\langle S', \frac{1}{2} \nabla_v \cdot \left(\delta_{f=\zeta}F_1 \sigma'(\zeta)^2 \nabla_v f \right) \right\rangle dt} \\
&- {\color{black}\left\langle S', \frac{1}{2} \partial_{\zeta} \left( \delta_{f=\zeta}\sigma(f)^2 F_3 \right) \right\rangle dt}.
\end{align*}

Using the chain rule, the following identities hold in the distributional sense:
\begin{align*}
\delta_{f = \zeta} \Delta_v f = \Delta_v \chi + \partial_{\zeta} \left( \delta_{f = \zeta} |\nabla_v f|^2 \right), \quad -\delta_{f = \zeta} \nabla_x f \cdot v = -v \cdot \nabla_x \chi. 
\end{align*}

Combining the above computations, we formally arrive at the following kinetic formulation:
\begin{align}\label{kineticformula-spde}
\partial_t \chi ={}& \Delta_v \chi - v \cdot \nabla_x \chi + \delta_{f = \zeta} \nabla_v \cdot (v f) - \delta_{f = \zeta} \nabla_v \cdot(f V \star_x \rho) \notag \\
& - \delta_{f = \zeta} \nabla_v \cdot(\sigma(f) \dot{W}_F) + {\color{black}\frac{1}{2} \nabla_v \cdot \left( \delta_{f = \zeta} F_1 \sigma'(\zeta)^2 \nabla_v f \right)} \notag \\
& - {\color{black}\frac{1}{2} \partial_{\zeta} \left( \delta_{f = \zeta}\sigma(f)^2 F_3 \right)} + \partial_{\zeta} p,
\end{align}
where the parabolic defect measure $p$ is given by
\begin{align}\label{kineticmeasure-regularity-1}
p = \delta_{f = \zeta} |\nabla_v f|^2.
\end{align}
This observation suggests defining solutions to equation \eqref{SPDE-0-ito} via a weak formulation based on the renormalized kinetic function. However, we remark that the quantity $p$ cannot always be characterized by an identity such as \eqref{kineticmeasure-regularity-1}. In particular, when proving the existence of solutions to \eqref{SPDE-0-ito}, the available \emph{a priori} estimates do not suffice to construct an approximation scheme $(f_n)_{n \geq 1}$ with strong convergence of $(\nabla_v f_n)_{n \geq 1}$ in $L^2(\mathbb{R}^{2d})$. Consequently, the limit of the sequence $(\delta_{f_n = \zeta} |\nabla_v f_n|^2)_{n \geq 1}$ yields a random measure that lacks an explicit characterization, and the identity \eqref{kineticmeasure-regularity-1} must be replaced by an inequality, thanks to the lower semi-continuity of the $L^2(\mathbb{R}^{2d})$-norm. This consideration motivates the introduction of the notion of a kinetic measure associated with renormalized kinetic solutions. We now proceed to give its definition.

\begin{definition}\label{def-kineticmeasure}
		 A \emph{kinetic measure} is a map $p$ from $\Omega$ to the space of nonnegative, locally finite measures on $\mathbb{R}^{2d}\times(0,\infty)\times[0,T]$ such that for every $\psi\in C^{\infty}_c(\mathbb{R}^{2d}\times (0,\infty))$, the map 
		\begin{align*}
		(\omega,t)\in \Omega\times[0,T]\to
		{\color{black}\int_{\mathbb R^{2d}\times(0,\infty)\times[0,t)}
		\psi(z,\zeta)\,p(dz,d\zeta,dr)}
		\end{align*}
		{\color{black}is $\mathcal{F}_t$-predictable.  The half-open interval is
		essential here: the cumulative process is then left-continuous and
		adapted.  In kinetic identities with an upper endpoint $t$, an integral
		written as $\int_0^t\partial_\zeta\psi\,dp$ means integration over
		$[0,t]$.  For the exact parabolic defect measures associated with regular
		weak solutions, the measure has no time atoms, so the two endpoint
		conventions coincide.}
\end{definition}
{\color{black}We fix the meaning of the stochastic term used in the kinetic
formulations below.  If $\psi\in
C_c^\infty(\mathbb R^{2d}\times(0,\infty))$ and $f$ belongs to the regular
solution class, then
\begin{align}\label{definition-kinetic-stochastic-integral}
 &\int_0^t\int_{\mathbb R^{2d}}
 \psi(z,f)\,\nabla_v\cdot\bigl(\sigma(f)\,dW_F\bigr)\,dz \notag\\
 &\quad:=\sum_{k\geq1}\sum_{\ell=1}^d\int_0^t
 \left(\int_{\mathbb R^{2d}}\psi(z,f)
 \partial_{v_\ell}\bigl(f_k\sigma(f)\bigr)\,dz\right)dB_s^{k,\ell}.
\end{align}
This scalar It\^o integral is well defined.  Indeed, if $K$ is the
$z$-projection of $\operatorname{supp}\psi$, spatial Cauchy--Schwarz gives
\begin{align*}
 \sum_{k\geq1}\sum_{\ell=1}^d
 \left|\int_K\psi(z,f)
 \partial_{v_\ell}\bigl(f_k\sigma(f)\bigr)\,dz\right|^2
 \leq C_\psi\int_K Q(f)\,dz,
\end{align*}
where $Q(f)=\sum_{k\geq1}|\nabla_v(f_k\sigma(f))|^2$ is locally
integrable by the regular-coefficient bounds and \eqref{L2-es}.}
Then we introduce the concept of renormalized kinetic solutions.
\begin{definition}\label{def-kineticsolution} 
Assume that $f_0,V$ and $\sigma(\cdot)$ satisfy Assumptions \ref{Assump-initialdata-regular}, \ref{Assump-ker} and \ref{Assump-sigma-n}, respectively. A nonnegative, $L^1(\mathbb{R}^{2d})$-valued $\mathcal{F}_t$-progressively measurable process $f$ is called a \emph{renormalized kinetic  solution} of \eqref{SPDE-0-ito} with initial data $f_0$ if  the following properties hold.
\begin{enumerate}
		{\color{black}\item Initial condition: $f(0)=f_0$ in
		$L^1(\mathbb R^{2d})$, almost surely.}
		\item Preservation of mass: almost surely for every $t\in[0,T]$,
		\begin{equation}\label{preservation-mass}
		\|f(t)\|_{L^{1}(\mathbb{R}^{2d})}=\|f_0\|_{L^{1}\left(\mathbb{R}^{2d}\right)}.
		\end{equation}
			{\color{black}\item Energy and Fisher-information bounds: one has
			\begin{align}\label{L2-es}
			 \sup_{t\in[0,T]}\mathbb E\|f(t)\|_{L^2(\mathbb R^{2d})}^2
			 +\mathbb E\int_0^T\!\int_{\mathbb R^{2d}}|\nabla_vf|^2\,dz\,dt
			 +\mathbb E\int_0^T\!\int_{\mathbb R^{2d}}|\nabla_v\sqrt f|^2\,dz\,dt<\infty.
			\end{align}
		Here and below,
		$\nabla_v\sqrt f:=\boldsymbol1_{\{f>0\}}\nabla_vf/(2\sqrt f)$.
		{\color{black}In particular,
		\begin{align}\label{Fisher-bound-regular}
		 \mathbb E\int_0^T\!\int_{\mathbb R^{2d}}
		 \frac{|\nabla_vf|^2}{f}\,dz\,dt
		 =4\mathbb E\int_0^T\!\int_{\mathbb R^{2d}}
		 |\nabla_v\sqrt f|^2\,dz\,dt<\infty,
		\end{align}
		where the quotient is defined to be zero on $\{f=0\}$.}
		In addition, the moment bound
		\begin{align}\label{moment-bound-solution-class}
		 \mathbb E\sup_{t\in[0,T]}
		 \int_{\mathbb R^{2d}}(|x|^2+|v|^2)f(t,x,v)\,dx\,dv<\infty
		\end{align}
		holds.}
		Furthermore, there exists a nonnegative kinetic measure $p$ satisfying the following properties.
		\item Regularity: in distributional sense, we have 
		\begin{align}\label{control KM}
		\delta_{f=\zeta}|\nabla_v f|^{2}\leq p\ \text{on}\ \mathbb{R}^{2d}\times(0,\infty)\times[0,T], \quad \mathbb{P}-a.s.
		\end{align}
		\item Vanishing at infinity:
		\begin{align}\label{control CKM VI}
		{\color{black}\liminf_{M\to \infty}\mathbb{E}\Big[p(\mathbb{R}^{2d}\times [M,M+1]\times [0,T])\Big]=0.}
		\end{align}
		\item The equation: almost surely for every $t\in[0,T]$ and $\psi\in\mathrm{C}_{c}^{\infty}\left(\mathbb{R}^{2d}\times(0,\infty)\right)$,
	\begin{align}\label{MC kenitic solution}
&\int_{\mathbb{R}}\int_{\mathbb{R}^{2d}} \chi(z, \zeta, t) \psi(z, \zeta)=\int_{\mathbb{R}} \int_{\mathbb{R}^{2d}} \chi(z, \zeta, 0)  \psi(z, \zeta)-\int_{0}^{t} \int_{\mathbb{R}^{2d}}\nabla_vf\cdot(\nabla_v\psi)(z,f)+\int^t_0\int_{\mathbb{R}}\int_{\mathbb{R}^{2d}}v\chi\cdot\nabla_x\psi\notag\\
&+\int^t_0\int_{\mathbb{R}^{2d}}\nabla_v\cdot(vf)\psi(z,f)-\int^t_0\int_{\mathbb{R}^{2d}}\psi(z,f)\nabla_v\cdot(fV\star_x\rho)-\int^t_0\int_{\mathbb{R}^{2d}}\psi(z,f)\nabla_v\cdot(\sigma(f)dW_F)\notag\\
&{\color{black}-\frac{1}{2}\int^t_0\int_{\mathbb{R}^{2d}}F_1\sigma'(f)^2\nabla_vf\cdot(\nabla_v\psi)(z,f)}\notag\\
&{\color{black}+\frac{1}{2}\int^t_0\int_{\mathbb{R}^{2d}}F_3\sigma^2(f)(\partial_{\zeta}\psi)(z,f)}-\int^t_0\int_{\mathbb{R}}\int_{\mathbb{R}^{2d}}\partial_{\zeta}\psi dp. 
\end{align}
\end{enumerate}			
\end{definition}


We now return to the regular setting, where the coefficient $\sigma(\cdot)$ is governed by Assumption \ref{Assump-sigma-n}. This regularity permits us to rigorously formulate the notion of weak solutions to \eqref{SPDE-0-ito}. 

\begin{definition}\label{def-weaksolution}
Assume that $f_0,V$ and $\sigma(\cdot)$ satisfy Assumptions \ref{Assump-initialdata-regular}, \ref{Assump-ker} and \ref{Assump-sigma-n}, respectively. A nonnegative, $L^1(\mathbb{R}^{2d})$-valued $\mathcal{F}_t$-progressively measurable process $f$ is called a \emph{weak solution} of \eqref{SPDE-0-ito} with initial data $f_0$ if the following properties are satisfied.
\begin{enumerate}
		{\color{black}\item Initial condition: $f(0)=f_0$ in
		$L^1(\mathbb R^{2d})$, almost surely.}
		\item Preservation of mass: almost surely for every $t\in[0,T]$,
		\begin{equation}\label{preservation-mass-weak}
		\|f(t)\|_{L^{1}(\mathbb{R}^{2d})}=\|f_0\|_{L^{1}\left(\mathbb{R}^{2d}\right)}.
		\end{equation}
			{\color{black}\item Energy and moment bounds: one has
			\begin{align}\label{weak-energy-bound-regular}
			 \sup_{t\in[0,T]}\mathbb E\|f(t)\|_{L^2(\mathbb R^{2d})}^2
			 +\mathbb E\int_0^T\!\int_{\mathbb R^{2d}}|\nabla_vf|^2\,dz\,dt<\infty,
			\end{align}
			and the moment estimate \eqref{moment-bound-solution-class} holds.  Entropy
			dissipation is not included in the definition of a weak solution: it is
			proved a posteriori in Proposition~\ref{prp-entropy-dissipation-regular}.}
       \item The equation: almost surely for every $\varphi\in C^{\infty}_c(\mathbb{R}^{2d})$ and every $t\in[0,T]$, 
	\begin{align*}
	\int_{\mathbb{R}^{2d}}f(t)\varphi=&\int_{\mathbb{R}^{2d}}f_0\varphi-\int^t_0\int_{\mathbb{R}^{2d}}\nabla_vf\cdot\nabla_v\varphi+\int^t_0\int_{\mathbb{R}^{2d}}vf\cdot\nabla_x\varphi\\
	&-\int^t_0\int_{\mathbb{R}^{2d}}vf\cdot\nabla_v\varphi+\int^t_0\int_{\mathbb{R}^{2d}}fV\star_x\rho\cdot\nabla_v\varphi+\int^t_0\int_{\mathbb{R}^{2d}}\nabla_v\varphi\cdot\sigma(f)dW_F\\
		&{\color{black}-\frac{1}{2}\int^t_0\int_{\mathbb{R}^{2d}}F_1\sigma'(f)^2\nabla_vf\cdot\nabla_v\varphi.}
	\end{align*} 
\end{enumerate}	
\end{definition}
{\color{black}In the following proposition, we prove the implication from weak solutions to renormalized kinetic solutions that is used later. The converse implication is not needed here.}

\begin{proposition}\label{equivalence} 
	Assume that $f_0,V$ and $\sigma(\cdot)$ satisfy Assumptions \ref{Assump-initialdata-regular}, \ref{Assump-ker} and \ref{Assump-sigma-n}, respectively. Let $f$ be a weak solution of \eqref{SPDE-0-ito} with initial data $f_0$. Let $\chi(x,v,\zeta,t)=I_{\{0<\zeta<f(x,v,t)\}}$ be the kinetic function. Then $f$ is a renormalized kinetic solution. In particular, $f$ satisfies the kinetic formulation: almost surely, for every {\color{black}$\psi \in C^\infty_c(\mathbb{R}^{2d}\times(0,\infty))$} and every $t \in [0, T]$,  
	\begin{align}\label{kinetic-formula-equivalence}
&\int_{\mathbb{R}}\int_{\mathbb{R}^{2d}} \chi(z, \zeta, t) \psi(z, \zeta)=\int_{\mathbb{R}} \int_{\mathbb{R}^{2d}} \chi(z, \zeta, 0)  \psi(z, \zeta)-\int_{0}^{t} \int_{\mathbb{R}^{2d}}\nabla_vf\cdot(\nabla_v\psi)(z,f)+\int^t_0\int_{\mathbb{R}}\int_{\mathbb{R}^{2d}}v\chi\cdot\nabla_x\psi\notag\\
&+\int^t_0\int_{\mathbb{R}^{2d}}\nabla_v\cdot(vf)\psi(z,f)-\int^t_0\int_{\mathbb{R}^{2d}}\psi(z,f)\nabla_v\cdot(fV\star_x\rho)-\int^t_0\int_{\mathbb{R}^{2d}}\psi(z,f)\nabla_v\cdot(\sigma(f)dW_F)\notag\\
&{\color{black}-\frac{1}{2}\int^t_0\int_{\mathbb{R}^{2d}}F_1\sigma'(f)^2\nabla_vf\cdot(\nabla_v\psi)(z,f)}\notag\\
&{\color{black}+\frac{1}{2}\int^t_0\int_{\mathbb{R}^{2d}}F_3\sigma^2(f)(\partial_{\zeta}\psi)(z,f)}-\int^t_0\int_{\mathbb{R}}\int_{\mathbb{R}^{2d}}\partial_{\zeta}\psi dp,  
\end{align}
with a precise kinetic measure 
\begin{align}\label{kineticmeasure-regularity}
	p=\delta_{f=\zeta}|\nabla_vf|^2. 
\end{align}
\end{proposition}
\begin{proof}
		{\color{black}Let $s\in C_c^\infty((0,\infty))$ and set
		$S(\zeta):=\int_0^\zeta s(r)\,dr$. Then $S$ is smooth and bounded,
		$S(0)=0$, and $S'=s$.  The It\^o formula used below is the
		generalized variational It\^o formula in the local Gelfand triple
		$H_v^1\subset L^2\subset H_v^{-1}$, after convolution in $x$ and
		localization by $\psi$.  All its derivatives of $S$ are bounded.}
	Due to the lack of spatial regularity of the weak solutions, we first {\color{black}perform} a spatial regularization. Let $(\eta_{\gamma})_{\gamma \in (0,1)}$ be a sequence of standard convolution kernels on $\mathbb{R}^d_x$. For every $\psi\in C^{\infty}_c(\mathbb{R}^{2d})$, applying It\^o's formula to $\int_{\mathbb{R}^{2d}} S(\eta_{\gamma} \ast f(t))\psi$, we obtain that almost surely, for every $t\in[0,T]$, 
	{\color{black}For brevity, put
	\begin{align*}
		 A(f)&:={\color{black}F_1\sigma'(f)^2\nabla_vf},\\
	 G_k(f)&:=\nabla_v\bigl(f_k\sigma(f)\bigr),\\
		 Q(f)&:={\color{black}F_1\sigma'(f)^2|\nabla_vf|^2+F_3\sigma(f)^2}.
	\end{align*}
	Then $Q(f)=\sum_{k\geq1}|G_k(f)|^2$.}
	\begin{align*}
	&\int_{\mathbb{R}^{2d}}S(\eta_{\gamma} \ast f(t))\psi\\
	=&\int_{\mathbb{R}^{2d}}S(\eta_{\gamma} \ast f_0)\psi-\int^t_0\int_{\mathbb{R}^{2d}}S'(\eta_{\gamma} \ast f)\nabla_v\eta_{\gamma} \ast f\cdot\nabla_v\psi\\
	&-\int^t_0\int_{\mathbb{R}^{2d}}S''(\eta_{\gamma} \ast f)|\nabla_v\eta_{\gamma} \ast f|^2\psi-\int^t_0\int_{\mathbb{R}^{2d}}S'(\eta_{\gamma} \ast f)v\cdot\nabla_x\eta_{\gamma} \ast f\psi\\
	&-\int^t_0\int_{\mathbb{R}^{2d}}S'(\eta_{\gamma} \ast f)\eta_{\gamma} \ast\nabla_v\cdot(fV\star_x\rho)\psi+\int^t_0\int_{\mathbb{R}^{2d}}S'(\eta_{\gamma} \ast f)\eta_{\gamma} \ast\nabla_v\cdot(vf)\psi\\
	&-\int^t_0\int_{\mathbb{R}^{2d}}\psi S'(\eta_{\gamma} \ast f)\eta_{\gamma} \ast\nabla_v\cdot(\sigma(f)dW_F)\\
	&{\color{black}-\frac{1}{2}\int^t_0\int_{\mathbb{R}^{2d}}S'(\eta_{\gamma} \ast f)(\eta_{\gamma} \ast A(f))\cdot\nabla_v\psi
	-\frac{1}{2}\int^t_0\int_{\mathbb{R}^{2d}}S''(\eta_{\gamma} \ast f)(\eta_{\gamma} \ast A(f))\cdot\nabla_v(\eta_{\gamma} \ast f)\psi}\\
	&{\color{black}\hspace{1.5em}+\frac{1}{2}\int^t_0\int_{\mathbb{R}^{2d}}S''(\eta_{\gamma} \ast f)
	\sum_{k\geq1}|\eta_\gamma\ast G_k(f)|^2\psi.}
	\end{align*}
Using the $L^2([0,T];L^2(\mathbb{R}^{2d}))$-integrability of $f$ and $\nabla_v f$, we obtain, along a subsequence,
\begin{align}\label{passtolimit-etagamma}
\eta_{\gamma}\ast f \to f, \quad 
\eta_{\gamma}\ast \nabla_v f \to \nabla_v f, \quad 
\eta_{\gamma}\ast \nabla_v\cdot(fV\star_x\rho) \to \nabla_v\cdot(fV\star_x\rho),
\end{align}
as $\gamma \to 0$, in $L^2([0,T];L^2_{loc}(\mathbb{R}^{2d}))$ and almost everywhere.  
By the dominated convergence theorem, we can pass to the limit:
\begin{align*}
	\int_{\mathbb{R}^{2d}} S(\eta_{\gamma} \ast f(t)) \psi 
	\to \int_{\mathbb{R}^{2d}} S(f(t)) \psi, 
	\quad 
	\int_{\mathbb{R}^{2d}} S(\eta_{\gamma} \ast f_0) \psi 
	\to \int_{\mathbb{R}^{2d}} S(f_0) \psi,
\end{align*}
as $\gamma \to 0$.  
Using \eqref{passtolimit-etagamma} and the boundedness of $S'$, we have
\begin{align*}
	-\int_0^t \!\! \int_{\mathbb{R}^{2d}} 
	S'(\eta_{\gamma} \ast f)\, \nabla_v(\eta_{\gamma} \ast f)\cdot\nabla_v\psi
	&= -\int_0^t \!\! \int_{\mathbb{R}^{2d}} S'(f)\, \nabla_v f \cdot \nabla_v\psi\\
	&\quad -\int_0^t \!\! \int_{\mathbb{R}^{2d}} S'(\eta_{\gamma} \ast f)
	(\nabla_v(\eta_{\gamma} \ast f) - \nabla_v f)\cdot\nabla_v\psi\\
	&\quad -\int_0^t \!\! \int_{\mathbb{R}^{2d}} 
	(S'(\eta_{\gamma} \ast f) - S'(f))\, \nabla_v f \cdot \nabla_v\psi\\
	&\to -\int_0^t \!\! \int_{\mathbb{R}^{2d}} 
	S'(f)\, \nabla_v f \cdot \nabla_v\psi,
\end{align*}
where in the last line we used H\"older's inequality and the dominated convergence theorem.  
Applying the same argument, we can also pass to the limit in
\begin{align*}
	&-\int_0^t \!\! \int_{\mathbb{R}^{2d}} 
	S'(\eta_{\gamma} \ast f)\, \eta_{\gamma} \ast \nabla_v\cdot(fV\star_x\rho)\, \psi
	+ \int_0^t \!\! \int_{\mathbb{R}^{2d}} 
	S'(\eta_{\gamma} \ast f)\, \eta_{\gamma} \ast \nabla_v\cdot(vf)\, \psi\\
	&{\color{black}-\frac{1}{2}\int_0^t \!\! \int_{\mathbb{R}^{2d}}
	S'(\eta_{\gamma} \ast f)(\eta_{\gamma}\ast A(f))\cdot\nabla_v\psi}\\
	&{\color{black}\hspace{1.5em}-\frac{1}{2}\int_0^t \!\! \int_{\mathbb{R}^{2d}}
	S''(\eta_{\gamma} \ast f)(\eta_{\gamma}\ast A(f))\cdot\nabla_v(\eta_{\gamma}\ast f)\psi}\\
	&{\color{black}\hspace{1.5em}+\frac{1}{2}\int_0^t \!\! \int_{\mathbb{R}^{2d}}
	S''(\eta_{\gamma} \ast f)\sum_{k\geq1}|\eta_\gamma\ast G_k(f)|^2\psi}\\
	&\to -\int_0^t \!\! \int_{\mathbb{R}^{2d}} 
	S'(f)\, \nabla_v\cdot(fV\star_x\rho)\, \psi
	+ \int_0^t \!\! \int_{\mathbb{R}^{2d}} 
	S'(f)\, \nabla_v\cdot(vf)\, \psi\\
	&{\color{black}\quad-\frac{1}{2}\int_0^t \!\! \int_{\mathbb{R}^{2d}}
	S'(f)A(f)\cdot\nabla_v\psi}\\
	&{\color{black}\hspace{1.5em}-\frac{1}{2}\int_0^t \!\! \int_{\mathbb{R}^{2d}}
	S''(f)A(f)\cdot\nabla_vf\,\psi
	+\frac{1}{2}\int_0^t \!\! \int_{\mathbb{R}^{2d}}S''(f)Q(f)\psi,}
\end{align*}
as $\gamma \to 0$.

{\color{black}For the quadratic-variation term used above, the regular-coefficient bounds and the assumptions on $F_1$ and $F_3$ imply
\begin{align*}
 (G_k(f))_{k\geq1}\in
 L^2\bigl(\Omega\times(0,T)\times K;\ell^2(\mathbb N;\mathbb R^d)\bigr)
\end{align*}
for every compact $K\Subset\mathbb R^{2d}$. Hence convolution is strongly continuous in this Hilbert-valued $L^2$ space, and
\begin{align*}
 \sum_{k\geq1}|\eta_\gamma\ast G_k(f)|^2\longrightarrow Q(f)
 \quad\text{in }L^1\bigl(\Omega\times(0,T)\times K\bigr).
\end{align*}}

Using the integration by parts formula and the dominated convergence theorem, we obtain
\begin{align*}
	-\int_0^t \int_{\mathbb{R}^{2d}} 
	S'(\eta_{\gamma} \ast f)\, v \cdot \nabla_x(\eta_{\gamma} \ast f)\, \psi
	= \int_0^t \int_{\mathbb{R}^{2d}} 
	S(\eta_{\gamma} \ast f)\, v \cdot \nabla_x\psi
	\;\longrightarrow\;
	\int_0^t \int_{\mathbb{R}^{2d}} 
	S(f)\, v \cdot \nabla_x\psi,
\end{align*}
as $\gamma \to 0$.

Furthermore, employing a similar decomposition argument, we obtain
\begin{align*}
	-\int_0^t \int_{\mathbb{R}^{2d}} 
	S''(\eta_{\gamma} \ast f)\, |\nabla_v(\eta_{\gamma} \ast f)|^2 \psi
	&= -\int_0^t \int_{\mathbb{R}^{2d}} 
	S''(f)\, |\nabla_v f|^2 \psi \\
	&\quad - \int_0^t \int_{\mathbb{R}^{2d}} 
	S''(\eta_{\gamma} \ast f)\, (|\nabla_v(\eta_{\gamma} \ast f)|^2 - |\nabla_v f|^2)\psi \\
	&\quad - \int_0^t \int_{\mathbb{R}^{2d}} 
	(S''(\eta_{\gamma} \ast f) - S''(f))\, |\nabla_v f|^2 \psi.
\end{align*}
Using the triangle inequality and the boundedness of $S''$, by the $L^2([0,T];L^2_{loc}(\mathbb{R}^{2d}))$-convergence of $\eta_{\gamma}\ast\nabla_vf\rightarrow\nabla_vf$, we have
\begin{align*}
	\left| \int_0^t \int_{\mathbb{R}^{2d}} 
	S''(\eta_{\gamma} \ast f)\, (|\nabla_v(\eta_{\gamma} \ast f)|^2 - |\nabla_v f|^2)\psi \right|
	\leq C(S'',\psi) \int_0^t \int_{\mathrm{supp}\,\psi} 
	||\eta_{\gamma} \ast \nabla_v f|^2 - |\nabla_v f|^2| 
	\;\rightarrow\; 0,
\end{align*}
as $\gamma \to 0$. 
By the dominated convergence theorem, the second term also converges to zero:
\begin{align*}
	-\int_0^t \int_{\mathbb{R}^{2d}} 
	(S''(\eta_{\gamma} \ast f) - S''(f))\, |\nabla_v f|^2 \psi
	\;\longrightarrow\; 0,
\end{align*}
as $\gamma \to 0$.

{\color{black}Notice also that
\begin{align}\label{renormalization-correction-cancellation}
 Q(f)-A(f)\cdot\nabla_vf
 ={\color{black}F_3\sigma(f)^2}.
\end{align}
Thus the last two terms in the preceding limit combine into precisely the
$\partial_\zeta$-term in \eqref{kinetic-formula-equivalence}.}

Finally, for the stochastic noise term, we decompose
\begin{align*}
	-\int_0^t \int_{\mathbb{R}^{2d}} 
	\psi\, S'(\eta_{\gamma} \ast f)\, 
	&\eta_{\gamma} \ast \nabla_v\cdot(\sigma(f)\, dW_F)
	= -\int_0^t \int_{\mathbb{R}^{2d}} 
	\psi\, S'(f)\, \nabla_v\cdot(\sigma(f)\, dW_F) \\
	&\quad - \int_0^t \int_{\mathbb{R}^{2d}} 
	\psi\, S'(\eta_{\gamma} \ast f)\,
	[\eta_{\gamma} \ast \nabla_v\cdot(\sigma(f)\, dW_F) 
	- \nabla_v\cdot(\sigma(f)\, dW_F)] \\
	&\quad - \int_0^t \int_{\mathbb{R}^{2d}} 
	\psi\, [S'(\eta_{\gamma} \ast f) - S'(f)]\,
	\nabla_v\cdot(\sigma(f)\, dW_F).
\end{align*}
{\color{black}We now justify the stochastic limit uniformly in time.  Let
$K$ be a compact set containing the support of $\psi$, put
$u_\gamma:=\eta_\gamma\ast f$, and recall
$G_k(f)=\nabla_v(f_k\sigma(f))$.  Spatial Cauchy--Schwarz, It\^o's
isometry, and Doob's maximal inequality give
\begin{align*}
&\mathbb E\sup_{\tau\in[0,T]}
\left|\int_0^\tau\!\int_{\mathbb R^{2d}}
 \psi S'(u_\gamma)
 \bigl(\eta_\gamma\ast\nabla_v\cdot(\sigma(f)dW_F)
       -\nabla_v\cdot(\sigma(f)dW_F)\bigr)\right|^2\\
&\qquad\leq C_{\psi,S}
 \mathbb E\int_0^T\!\int_K
 \sum_{k\geq1}|\eta_\gamma\ast G_k(f)-G_k(f)|^2\,dz\,ds
 \longrightarrow0.
\end{align*}
The same argument yields
\begin{align*}
&\mathbb E\sup_{\tau\in[0,T]}
\left|\int_0^\tau\!\int_{\mathbb R^{2d}}
 \psi\bigl(S'(u_\gamma)-S'(f)\bigr)
 \nabla_v\cdot(\sigma(f)dW_F)\right|^2\\
&\qquad\leq C_\psi
 \mathbb E\int_0^T\!\int_K
 |S'(u_\gamma)-S'(f)|^2
 \sum_{k\geq1}|G_k(f)|^2\,dz\,ds
 \longrightarrow0.
\end{align*}
Here the first convergence follows from the strong continuity of spatial
convolution in the Hilbert-valued $L^2$ space from the preceding paragraph;
the second follows from almost-everywhere convergence, the boundedness of
$S'$, and dominated convergence against the locally integrable function
$Q(f)=\sum_k|G_k(f)|^2$.  We may therefore choose a deterministic sequence
$\gamma_m\downarrow0$ along which both suprema converge to zero almost
surely.  On this common event,
\begin{align*}
&-\int_0^t\!\int_{\mathbb R^{2d}}
 \psi S'(\eta_{\gamma_m}\ast f)
 \eta_{\gamma_m}\ast\nabla_v\cdot(\sigma(f)dW_F)\\
&\qquad\longrightarrow
-\int_0^t\!\int_{\mathbb R^{2d}}
 \psi S'(f)\nabla_v\cdot(\sigma(f)dW_F)
\end{align*}
uniformly for $t\in[0,T]$.}

Combining all the above arguments, we conclude that, almost surely, for every $t \in [0,T]$, 
	\begin{align*}
	\int_{\mathbb{R}^{2d}}S(f(t))\psi=&\int_{\mathbb{R}^{2d}}S(f_0)\psi-\int^t_0\int_{\mathbb{R}^{2d}}S'(f)\nabla_vf\cdot\nabla_v\psi\\
	&-\int^t_0\int_{\mathbb{R}^{2d}}S''(f)|\nabla_vf|^2\psi+\int^t_0\int_{\mathbb{R}^{2d}}S(f)v\cdot\nabla_x\psi\\
	&-\int^t_0\int_{\mathbb{R}^{2d}}S'(f)\nabla_v\cdot(fV\star_x\rho)\psi+\int^t_0\int_{\mathbb{R}^{2d}}S'(f)\nabla_v\cdot(vf)\psi\\
	&-\int^t_0\int_{\mathbb{R}^{2d}}\psi S'(f)\nabla_v\cdot(\sigma(f)dW_F)\\
		&{\color{black}-\frac{1}{2}\int^t_0\int_{\mathbb{R}^{2d}}S'(f)A(f)\cdot\nabla_v\psi
		+\frac{1}{2}\int^t_0\int_{\mathbb{R}^{2d}}S''(f)F_3\sigma(f)^2\psi.}
	\end{align*}
	Let $\Psi_S(x,v,\zeta)=\psi(x,v)S'(\zeta)$, $(x,v,\zeta)\in\mathbb{R}^{2d+1}$. Then we have that almost surely for every $t\in[0,T]$, 
\begin{align*}
&\int_{\mathbb{R}}\int_{\mathbb{R}^{2d}} \chi(z, \zeta, t) \Psi_S(z, \zeta)=\int_{\mathbb{R}} \int_{\mathbb{R}^{2d}} \chi(z, \zeta, 0)  \Psi_S(z, \zeta)-\int_{0}^{t} \int_{\mathbb{R}^{2d}}\nabla_vf\cdot(\nabla_v\Psi_S)(z,f){\color{black}+\int^t_0\int_{\mathbb{R}}\int_{\mathbb{R}^{2d}}v\chi\cdot\nabla_x\Psi_S}\notag\\
&+\int^t_0\int_{\mathbb{R}^{2d}}\nabla_v\cdot(vf)\Psi_S(z,f){\color{black}-\int^t_0\int_{\mathbb{R}^{2d}}\Psi_S(z,f)\nabla_v\cdot(fV\star_x\rho)}-\int^t_0\int_{\mathbb{R}^{2d}}\Psi_S(z,f)\nabla_v\cdot(\sigma(f)dW_F)\notag\\
&{\color{black}-\frac{1}{2}\int^t_0\int_{\mathbb{R}^{2d}}F_1\sigma'(f)^2\nabla_vf\cdot(\nabla_v\Psi_S)(z,f)}\notag\\
&{\color{black}+\frac{1}{2}\int^t_0\int_{\mathbb{R}^{2d}}F_3\sigma^2(f)(\partial_{\zeta}\Psi_S)(z,f)}-\int^t_0\int_{\mathbb{R}}\int_{\mathbb{R}^{2d}}\partial_{\zeta}\Psi_S dp, 
\end{align*}
	with $p = \delta_{f=\zeta} |\nabla_v f|^2$. {\color{black}Finite linear combinations of tensor products $\psi(z)s(\zeta)$, with
$\psi\in C_c^\infty(\mathbb R^{2d})$ and
$s\in C_c^\infty((0,\infty))$, are dense in
	$C_c^\infty(\mathbb R^{2d}\times(0,\infty))$ in the usual test-function topology.  {\color{black}In removing the $x$-mollification, the deterministic terms converge in the variational dualities, while the stochastic integrands converge in
	$L^2(\Omega\times(0,T);\ell^2)$.  We first obtain the identity for a
	countable dense tensor-product family and rational times.  The continuous
	semimartingale representatives on both sides then extend it, on one
	common event, to every test function and every $t\in[0,T]$.}} It follows that the function $f$ satisfies the kinetic formulation \eqref{kinetic-formula-equivalence}. Moreover, invoking {\color{black}\eqref{weak-energy-bound-regular}}, we deduce that $p$ is almost surely a finite measure on {\color{black}$\mathbb R^{2d}\times(0,\infty)\times[0,T]$}. In particular, all the conditions listed in Definition~\ref{def-kineticmeasure} are fulfilled. {\color{black}Proposition~\ref{prp-entropy-dissipation-regular}, proved in Section~\ref{sec-5}, supplies the Fisher-information bound in \eqref{L2-es}--\eqref{Fisher-bound-regular}; it also proves the stronger entropy estimate \eqref{entropy-dissipation-regular}.} Consequently, $p$ qualifies as a kinetic measure, and we conclude that $f$ is indeed a renormalized kinetic solution. This completes the proof. 
\end{proof}

\begin{remark}
{\color{black}Proposition~\ref{equivalence} establishes only that every weak solution in the regular class is a renormalized kinetic solution. A converse statement would require a separate reconstruction argument and is neither asserted nor used in this paper. For a singular coefficient such as $\sigma(\zeta)=\sqrt\zeta$, the variational weak formulation above is in general unavailable, which is why Section~\ref{sec-6} uses the kinetic formulation directly.}
\end{remark}

\section{The uniqueness of the SPDE}\label{sec-3}
{\color{black}In this section, we prove pathwise uniqueness for the regular weak
solutions of Definition~\ref{def-weaksolution}.  Proposition~\ref{equivalence}
associates with such a solution the exact kinetic defect measure
$\delta_{f=\zeta}|\nabla_vf|^2\,dz\,d\zeta\,dt$; in particular, its
finite-variation part has no time atoms.  This is the class needed in
Section~\ref{sec-5}.  We do not assert uniqueness for an arbitrary kinetic
measure in Definition~\ref{def-kineticsolution}, or for the probabilistically
weak square-root solution of Section~\ref{sec-6}.}

In the following, we introduce suitable truncation functions. For every $\beta\in(0,1)$, let $\varphi_{\beta}:\mathbb{R}\to[0,1]$ be the unique nondecreasing piecewise linear function that satisfies
\begin{equation}
\varphi_{\beta}(\zeta)=1\  {\rm{if}}\  \zeta\ge\beta, \ \varphi_{\beta}(\zeta)=0\ {\rm{if}}\  \zeta\le\frac{\beta}{2},\  {\rm{and}}\  \varphi'_{\beta}=\frac{2}{\beta}I_{\{\frac{\beta}{2}<\zeta<\beta\}}.
\end{equation}
For every $M\in\mathbb{N}$, let $\zeta_M:\mathbb{R}\to[0,1]$ be the unique nonincreasing piecewise linear function satisfying
\begin{equation}
\zeta_M(\zeta)=0\ {\rm{if}}\ \zeta\ge M+1,\ \zeta_M(\zeta)=1\ {\rm{if}}\ \zeta\le M,\ {\rm{and}}\ \zeta'_M=-I_{\{M<\zeta<M+1\}}.
\end{equation}
For every $\varepsilon_x,\varepsilon_v,\delta\in(0,1)$, let $\kappa^{\varepsilon_x}_{d}:\mathbb{R}^{d}_x\to\left[0,\infty\right)$, $\kappa^{\varepsilon_v}_{d}:\mathbb{R}^{d}_v\to\left[0,\infty\right)$ and $\kappa^{\delta}_{1}:\mathbb{R}\to\left[0,\infty\right)$ be standard convolution kernels of scales $\varepsilon_x,\varepsilon_v$ and $\delta$ on $\mathbb{R}^{d}_x,\mathbb{R}^d_v$ and $\mathbb{R}$, respectively. Let $\kappa^{\varepsilon}(z,z')=\kappa^{\varepsilon_x}_d(x-x')\kappa^{\varepsilon_v}_d(v-v')$, for every $z=(x,v),\ z'=(x',v')\in\mathbb{R}^{2d}$, and we let $\kappa^{\varepsilon,\delta}$ be defined by
\begin{align}
\kappa^{\varepsilon,\delta}(z,z',\zeta,\eta)={\color{black}\kappa^{\varepsilon}(z,z')}\kappa^{\delta}_{1}(\zeta-\eta)\ \text{for\ every}\ (z,z',\zeta,\eta)\in(\mathbb{R}^{2d})^2\times\mathbb{R}^2.
\end{align}
For simplicity, we denote the joint limits $\lim_{\varepsilon_x\rightarrow0,\varepsilon_v\rightarrow0}$ and $\limsup_{\varepsilon_x\rightarrow0,\varepsilon_v\rightarrow0}$ by $\lim_{\varepsilon\rightarrow0}$ and $\limsup_{\varepsilon\rightarrow0}$, respectively. 

For every $R=(R_1,R_2)\in(0,+\infty)^2$, let $\alpha_{R_1}^1$ be a smooth cutoff function on $\mathbb{R}^{d}_x$ such that 
\begin{align}\label{truncation-property}
\alpha_{R_1}^1(x)=1,\ \ |x|\leq R_1,\ \ \alpha_{R_1}^1(x)=0,\ \ |x|>2R_1, 	
\end{align}
and 
\begin{align}\label{truncation-property-2}
	R_1|\nabla_x^2\alpha_{R_1}^1|+|\nabla_x\alpha_{R_1}^1|\leq\frac{c}{R_1}, 
\end{align}
for some $c\in(0,\infty)$ independent of $R_1$. Similarly, we define the same truncation function with respect to the velocity variable, denoted by $\alpha^2_{R_2}$. We set 
\begin{equation}\label{alpha-R}
	\alpha_R(z)=\alpha^1_{R_1}(x)\alpha^2_{R_2}(v),\ \text{for every }z=(x,v)\in\mathbb{R}^{2d}. 
\end{equation}
{\color{black}Throughout Section~\ref{sec-3}, limits involving the phase-space cutoff are iterated in the order
\begin{align}\label{ordered-phase-space-cutoff-limit}
 \lim_{R\to\infty}:=\lim_{R_2\to\infty}\lim_{R_1\to\infty},
 \qquad
 \limsup_{R\to\infty}:=\limsup_{R_2\to\infty}\limsup_{R_1\to\infty}.
\end{align}
Thus $R_2$ is kept fixed while the spatial cutoff is first removed. This order is essential on the whole space because the $x$-transport error is bounded by a constant times $R_2/R_1$.}

\begin{proposition}\label{vanish-0-kinetic}
{\color{black}Assume that $f_0$, $V$, and $\sigma(\cdot)$ satisfy
Assumptions~\ref{Assump-initialdata-regular}, \ref{Assump-ker}, and
\ref{Assump-sigma-n}, respectively. Let $f$ be a weak solution of
\eqref{SPDE-0-ito}, and let
$p=\delta_{f=\zeta}|\nabla_vf|^2\,dz\,d\zeta\,dt$ be the kinetic measure
provided by Proposition~\ref{equivalence}. Then, almost surely,}
\begin{align*}
{\color{black}\liminf_{\beta \downarrow0} \, \beta^{-1} p\big( \mathbb{R}^{2d} \times [\beta/2, \beta] \times [0, T] \big) = 0.}
\end{align*}
\end{proposition}
\begin{proof}
	{\color{black}For the exact defect measure supplied by
	Proposition~\ref{equivalence}, the Fisher-information identity gives, for
	every $\beta\in(0,1)$,
	\begin{align*}
	&\beta^{-1}p\bigl(\mathbb R^{2d}\times[\beta/2,\beta]
	                    \times[0,T]\bigr)\\
	&\quad=4\int_0^T\!\int_{\mathbb R^{2d}}
	 \frac{f}{\beta}\boldsymbol1_{\{\beta/2\leq f\leq\beta\}}
	 |\nabla_v\sqrt f|^2\,dz\,dt\\
	&\quad\leq4\int_0^T\!\int_{\mathbb R^{2d}}
	 \boldsymbol1_{\{0<f\leq\beta\}}|\nabla_v\sqrt f|^2\,dz\,dt.
	\end{align*}
	The integrand on the last line converges pointwise to zero as
	$\beta\downarrow0$ and is dominated by the integrable Fisher-information
	density. Hence dominated convergence and \eqref{Fisher-bound-regular} yield
	\begin{align}\label{vanish-zero-kinetic-L1}
	 \lim_{\beta\downarrow0}\mathbb E\left[
	 \beta^{-1}p\bigl(\mathbb R^{2d}\times[\beta/2,\beta]
	                    \times[0,T]\bigr)\right]=0.
	\end{align}
	Choose a deterministic sequence $\beta_j\downarrow0$ along which the
	expectations in \eqref{vanish-zero-kinetic-L1} are summable. Fatou's lemma
	then shows that the corresponding random variables converge to zero almost
	surely. Since the continuous-parameter liminf is bounded above by the
	liminf along this sequence, the claimed assertion follows.}
\end{proof}

{\color{black}In the following, we establish pathwise uniqueness for regular weak
solutions to \eqref{SPDE-0-ito}. Proposition~\ref{equivalence} first
associates with each such weak solution the renormalized kinetic formulation
used in the doubling argument.} Our approach adopts the strategy and proof structure of \cite[Theorem 4.6]{FG24}. Owing to the close structural analogy between the equations, many terms coincide with those treated in \cite{FG24}; for such terms, we carefully refer to \cite{FG24}. We focus instead on the technical aspects of the new contributions, in particular the kernel terms and the transport terms. 

\begin{theorem}\label{Uniqueness-spde}
{\color{black}Assume that $f_{1,0},f_{2,0}$ satisfy
Assumption~\ref{Assump-initialdata-regular} and that $V,\sigma(\cdot)$
satisfy Assumptions~\ref{Assump-ker}, \ref{Assump-sigma-n}, respectively.
Let $f_1,f_2$ be two weak solutions of \eqref{SPDE-0-ito}, defined on the
same filtered probability space and driven by the same correlated Wiener
process $W_F$, with initial data $f_{1,0},f_{2,0}$, respectively. If
$f_{1,0}=f_{2,0}$ in $L^1(\mathbb R^{2d})$ almost surely, then}
\begin{align*}
\sup_{t\in[0,T]}\|f_1(t)-f_2(t)\|_{L^1(\mathbb{R}^{2d})}=0
\qquad{\color{black}\text{almost surely}.} 	
\end{align*}
\end{theorem}
\begin{proof}
{\color{black}By Proposition~\ref{equivalence}, each $f_i$ has the kinetic
formulation with
\begin{align*}
 p_i(dw,d\zeta,ds)
 =\delta_{f_i(w,s)=\zeta}|\nabla_vf_i(w,s)|^2\,dw\,d\zeta\,ds.
\end{align*}
Moreover, \eqref{vanish-zero-kinetic-L1} allows us to fix one deterministic
sequence $\beta_j\downarrow0$, common to $i=1,2$, such that
\begin{align}\label{common-low-density-sequence}
 \sum_{j\geq1}\sum_{i=1}^2\mathbb E\left[
 \beta_j^{-1}p_i\bigl(\mathbb R^{2d}\times[\beta_j/2,\beta_j]
                       \times[0,T]\bigr)\right]<\infty.
\end{align}
Thus the mollified kinetic processes below have continuous semimartingale
representatives, and the ordinary It\^o product formula applies without a
jump-covariation term.}
Let $\chi_1,\chi_2$ be the renormalized kinetic functions corresponding to $f_1$ and $f_2$, respectively. For every $\varepsilon, \delta > 0$, recall that $\kappa^{\varepsilon,\delta}$ denotes the standard convolution kernel. For $i = 1,2$, define the regularized kinetic functions {\color{black}$\chi_i^{\varepsilon,\delta}(z,\eta)=\chi_i\ast \kappa^{\varepsilon,\delta}:=\int_{\mathbb{R}^{2d+1}}\chi_i(z',\zeta)\kappa^{\varepsilon,\delta}(z,z',\zeta,\eta)\,dz'\,d\zeta$.}

By the definition of kinetic functions and the elementary properties of indicator functions,  
\begin{align}\label{L1-difference}
\int_{\mathbb{R}^{2d}}|f_1(t)-f_2(t)|dz=&\int_{\mathbb{R}^{2d}}\int_{\mathbb{R}}|\chi_1(t)-\chi_2(t)|d\zeta dz=\int_{\mathbb{R}^{2d}}\int_{\mathbb{R}}|\chi_1(t)-\chi_2(t)|^2d\zeta dz\notag\\
=&\int_{\mathbb{R}^{2d+1}}\chi_1(t)^2+\chi_2(t)^2-2\chi_1(t)\chi_2(t)dzd\zeta\notag\\
=&\int_{\mathbb{R}^{2d+1}}\chi_1(t)+\chi_2(t)-2\chi_1(t)\chi_2(t)dzd\zeta\notag\\
&{\color{black}=\lim_{j\to\infty}\lim_{M\to\infty}\lim_{R\to\infty}
\lim_{\delta\to0}\lim_{\varepsilon\to0}
\int_{\mathbb{R}^{2d+1}}(\chi_1^{\varepsilon,\delta}(t)+\chi_2^{\varepsilon,\delta}(t)-2\chi_1^{\varepsilon,\delta}(t)\chi_2^{\varepsilon,\delta}(t))
\varphi_{\beta_j}\zeta_M\alpha_R\,dz\,d\zeta.}  
\end{align}
For $i=1,2$, $s\in[0,T]$, let {\color{black}$\bar{\kappa}^{\delta}_{s,i}(z,\eta)=\kappa_1^{\delta}(f_i(z,s)-\eta)$ and $\bar{\kappa}^{\varepsilon,\delta}_{s,i}(z,w,\eta)=\kappa^{\varepsilon,\delta}(z,w,f_i(w,s),\eta)$}, $(z,w,\eta)\in\mathbb{R}^{2d}\times\mathbb{R}^{2d}\times\mathbb{R}$. For $i=1,2$, applying the kinetic formula \eqref{MC kenitic solution} with test function $\kappa^{\varepsilon,\delta}$, we have that almost surely, for every $(z,\eta)\in\mathbb{R}^{2d}\times(\delta/2,+\infty)$, 
\begin{align*}
\chi^{\varepsilon,\delta}_{s,i}(z,\eta)|^t_{s=0}=&\nabla_v\cdot\Big(\int^t_0(\nabla_vf_i\ast\bar{\kappa}^{\varepsilon,\delta}_{s,i})(z,\eta)\Big)-\int^t_0((v\cdot\nabla_xf_i)\ast\bar{\kappa}^{\varepsilon,\delta}_{s,i})(z,\eta)+\int^t_0((v\cdot\nabla_vf_i)\ast\bar{\kappa}^{\varepsilon,\delta}_{s,i})(z,\eta)\\
&+\int^t_0({\color{black}(d f_i)\ast\bar{\kappa}^{\varepsilon,\delta}_{s,i}})(z,\eta)-\int^t_0(\nabla_v\cdot(f_iV\star_x\rho_i)\ast\bar{\kappa}^{\varepsilon,\delta}_{s,i})(z,\eta)\notag\\
&-\int^t_0(\bar{\kappa}^{\varepsilon,\delta}_{s,i}\ast\nabla_v\cdot(\sigma(f_i)dW_F))(z,\eta)+{\color{black}\nabla_v\cdot\Big(\frac{1}{2}\int^t_0(F_1\sigma'(f_i)^2\nabla_vf_i)\ast\bar{\kappa}^{\varepsilon,\delta}_{s,i}\Big)}\\
&{\color{black}-\partial_{\eta}\Big(\frac{1}{2}\int^t_0
F_3\sigma^2(f_i)
\ast\bar{\kappa}^{\varepsilon,\delta}_{s,i}\Big)
+\partial_\eta\!\left(\int_0^t
 P_i^{\varepsilon,\delta}(s,z,\eta)\,ds\right).}	
\end{align*}
{\color{black}For later use, denote the time density of the mollified kinetic
measure by
\begin{align*}
 P_i^{\varepsilon,\delta}(s,z,\eta)
 :=\int_{\mathbb R^{2d}}
 \kappa^{\varepsilon,\delta}(z,w,f_i(w,s),\eta)
 |\nabla_vf_i(w,s)|^2\,dw.
\end{align*}}
Here, we adopt the convolution notation for convenience. For instance, the first term on the right-hand side of the identity is written as $(\nabla_v f_i \ast \bar{\kappa}^{\varepsilon,\delta}_{s,i})(z,\eta) = \int_{\mathbb{R}^{2d}} \nabla_v f_i(w) \kappa^{\varepsilon,\delta}(z, w, f_i(w,s), \eta) \, dw$.

We begin by analyzing the first two terms on the right-hand side of \eqref{L1-difference}. We observe that, almost surely, for every $\varepsilon, \beta \in (0,1)$, $M \in \mathbb{N}$, {\color{black}$R=(R_1,R_2)\in (1,+\infty)^2$}, $\delta \in (0, \beta/4)$, $t \in [0,T]$, and $i = 1,2$, 
\begin{align*}
\int_{\mathbb{R}^{2d+1}}\chi^{\varepsilon,\delta}_{s,i}(z,\eta)\varphi_{\beta}(\eta)\zeta_M(\eta)\alpha_R(z)dzd\eta\Big|^t_{s=0}=I^{i,cut}_t+I^{i,trans-x}_t+I^{i,trans-v}_t+I^{i,mart}_t+I^{i,ker}_t, 	
\end{align*}
where 
\begin{align*}
I^{i,cut}_t=&\int_{\mathbb{R}^{2d+1}}\nabla_v\cdot\Big(\int^t_0(\nabla_vf_i\ast\bar{\kappa}^{\varepsilon,\delta}_{s,i})(z,\eta)\Big)\varphi_{\beta}(\eta)\zeta_M(\eta)\alpha_R(z)\\
&{\color{black}-\frac{1}{2}\int^t_0\int_{\mathbb{R}}\int_{\mathbb{R}^{2d}}(F_1\sigma'(f_i)^2\nabla_vf_i)\ast\bar{\kappa}^{\varepsilon,\delta}_{s,i}\cdot\nabla_v(\varphi_{\beta}\zeta_M\alpha_R)}\\
&{\color{black}+\frac{1}{2}\int^t_0\int_{\mathbb{R}}\int_{\mathbb{R}^{2d}}
{\color{black}F_3\sigma^2(f_i)}
\ast\bar{\kappa}^{\varepsilon,\delta}_{s,i}\partial_{\eta}(\varphi_{\beta}\zeta_M\alpha_R)}\\
&{\color{black}+\int_0^t\int_{\mathbb{R}^{2d+1}}\partial_{\eta}P_i^{\varepsilon,\delta}(s,z,\eta)\varphi_{\beta}(\eta)\zeta_M(\eta)\alpha_R(z)\,dz\,d\eta\,ds},  	
\end{align*}
and 
\begin{align*}
I^{i,trans-x}_t=&-\int^t_0\int_{\mathbb{R}^{2d+1}}((v\cdot\nabla_xf_i)\ast\bar{\kappa}^{\varepsilon,\delta}_{s,i})(z,\eta)\varphi_{\beta}(\eta)\zeta_M(\eta)\alpha_R(z),\\
I^{i,trans-v}_t=&\int^t_0\int_{\mathbb{R}^{2d+1}}((v\cdot\nabla_vf_i)\ast\bar{\kappa}^{\varepsilon,\delta}_{s,i})(z,\eta)\varphi_{\beta}(\eta)\zeta_M(\eta)\alpha_R(z)\\
&+\int^t_0\int_{\mathbb{R}^{2d+1}}({\color{black}(d f_i)\ast\bar{\kappa}^{\varepsilon,\delta}_{s,i}})(z,\eta)\varphi_{\beta}(\eta)\zeta_M(\eta)\alpha_R(z),\\
I^{i,mart}_t=&-\int^t_0\int_{\mathbb{R}^{2d+1}}(\bar{\kappa}^{\varepsilon,\delta}_{s,i}\ast\nabla_v\cdot(\sigma(f_i)dW_F))(z,\eta)\varphi_{\beta}(\eta)\zeta_M(\eta)\alpha_R(z),\\
I^{i,ker}_t=&-\int^t_0\int_{\mathbb{R}^{2d+1}}(\nabla_v\cdot(f_iV\star_x\rho_i)\ast\bar{\kappa}^{\varepsilon,\delta}_{s,i})(z,\eta)	\varphi_{\beta}(\eta)\zeta_M(\eta)\alpha_R(z). 
\end{align*}
For the mixed term in \eqref{L1-difference}, for every $\beta\in(0,1)$, $M\in\mathbb{N}$, {\color{black}$R=(R_1,R_2)\in(1,+\infty)^2$}, we test it against the truncation functions $\varphi_{\beta}\zeta_M\alpha_R$ and apply It\^o's product rule to see that, almost surely, for every $t\in[0,T]$, 
\begin{align*}
	\int_{\mathbb{R}^{2d+1}}\chi^{\varepsilon,\delta}_{s,1}(z,\eta)\chi^{\varepsilon,\delta}_{s,2}(z,\eta)\varphi_{\beta}(\eta)\zeta_M(\eta)\alpha_R(z)\Big|^t_{s=0}=&\int^t_0\int_{\mathbb{R}^{2d+1}}\chi^{\varepsilon,\delta}_{s,2}(z,\eta)d\chi^{\varepsilon,\delta}_{s,1}(z,\eta)\varphi_{\beta}(\eta)\zeta_M(\eta)\alpha_R(z)\\
	&+\int^t_0\int_{\mathbb{R}^{2d+1}}\chi^{\varepsilon,\delta}_{s,1}d\chi^{\varepsilon,\delta}_{s,2}(z,\eta)\varphi_{\beta}(\eta)\zeta_M(\eta)\alpha_R(z)\\
	&+\int^t_0\int_{\mathbb{R}^{2d+1}}d\langle\chi^{\varepsilon,\delta}_2,\chi^{\varepsilon,\delta}_1\rangle_s(z,\eta)\varphi_{\beta}(\eta)\zeta_M(\eta)\alpha_R(z). 
\end{align*}
We first expand the $\chi^{\varepsilon,\delta}_{s,2}d\chi^{\varepsilon,\delta}_{s,1}$-term. Applying the kinetic formula \eqref{MC kenitic solution} for $\chi^{\varepsilon,\delta}_{s,1}$, it follows that almost surely for every $\beta\in(0,1)$, $M\in\mathbb{N}$, {\color{black}$R=(R_1,R_2)\in(1,+\infty)^2$}, $t\in[0,T]$, 
\begin{align}\label{21-formula}
&\int^t_0\int_{\mathbb{R}^{2d+1}}\chi^{\varepsilon,\delta}_{s,2}(z,\eta)d\chi^{\varepsilon,\delta}_{s,1}(z,\eta)\varphi_{\beta}(\eta)\zeta_M(\eta)\alpha_R(z)\notag\\
=&I^{2,1,err}_t+I^{2,1,meas}_t+I^{2,1,cut}_t+I^{2,1,trans-x}_t+I^{2,1,trans-v}_t+I^{2,1,mart}_t+I^{2,1,ker}_t, 	
\end{align}
where the error term is 
\begin{align*}
I^{2,1,err}_t=&{\color{black}-\frac{1}{2}\int^t_0\int_{\mathbb{R}^{2d+1}}(F_1\sigma'(f_1)^2\nabla_vf_1)\ast\bar{\kappa}^{\varepsilon,\delta}_{s,1}\cdot\nabla_vf_2\ast\bar{\kappa}^{\varepsilon,\delta}_{s,2}\varphi_{\beta}(\eta)\zeta_M(\eta)\alpha_R(z)}\\
&{\color{black}-\frac{1}{2}\int^t_0\int_{\mathbb{R}^{4d+1}}[(F_3\sigma^2(f_1))\ast\bar{\kappa}^{\varepsilon,\delta}_{s,1}]\bar{\kappa}^{\varepsilon,\delta}_{s,2}\varphi_{\beta}(\eta)\zeta_M(\eta)\alpha_R(z)},  
\end{align*}
the measure term is defined by 
{\color{black}
\begin{align*}
I^{2,1,meas}_t={}&\int_0^t\int_{\mathbb R^{4d+1}}
 P_1^{\varepsilon,\delta}(s,z,\eta)
 \bar\kappa^{\varepsilon,\delta}_{s,2}(z,w,\eta)
 \varphi_\beta(\eta)\zeta_M(\eta)\alpha_R(z)
 \,dz\,dw\,d\eta\,ds\\
&-\int_0^t\int_{\mathbb R^{2d+1}}
 (\nabla_vf_1\ast\bar\kappa^{\varepsilon,\delta}_{s,1})(z,\eta)
 \cdot(\nabla_vf_2\ast\bar\kappa^{\varepsilon,\delta}_{s,2})(z,\eta)
 \varphi_\beta(\eta)\zeta_M(\eta)\alpha_R(z)
 \,dz\,d\eta\,ds.	
\end{align*}}
the cutoff term is defined by 
\begin{align*}
I^{2,1,cut}_t=&-\int^t_0\int_{\mathbb{R}^{2d+1}}(\kappa^{\varepsilon,\delta}\ast p_1)\chi^{\varepsilon,\delta}_{s,2}\partial_{\eta}(\varphi_{\beta}\zeta_{M})\alpha_R\\
&{\color{black}+\frac{1}{2}\int^t_0\int_{\mathbb{R}^{2d+1}}(F_3\sigma(f_1)^2)\ast\bar{\kappa}^{\varepsilon,\delta}_{s,1}\chi^{\varepsilon,\delta}_{s,2}\partial_{\eta}(\varphi_{\beta}\zeta_M)\alpha_R}\\
&-\int^t_0\int_{\mathbb{R}^{2d+1}}(\nabla_vf_1\ast\bar{\kappa}^{\varepsilon,\delta}_{s,1})\chi^{\varepsilon,\delta}_{s,2}\cdot\nabla_v\alpha_R\varphi_{\beta}\zeta_M\\
&{\color{black}-\frac{1}{2}\int^t_0\int_{\mathbb{R}^{2d+1}}[(F_1\sigma'(f_1)^2\nabla_vf_1)\ast\bar{\kappa}^{\varepsilon,\delta}_{s,1}]\chi^{\varepsilon,\delta}_{s,2}\cdot\nabla_v\alpha_R\varphi_{\beta}\zeta_M}, 	
\end{align*}
the martingale term is defined by 
\begin{align*}
I^{2,1,mart}_t=-\int^t_0\int_{\mathbb{R}^{2d+1}}\bar{\kappa}^{\varepsilon,\delta}_{s,1}\ast\nabla_v\cdot(\sigma(f_1)dW_F)\chi^{\varepsilon,\delta}_{s,2}\varphi_{\beta}\zeta_M\alpha_R, 
\end{align*}
the kernel term is defined by 
\begin{align*}
I^{2,1,ker}_t=-\int^t_0\int_{\mathbb{R}^{2d+1}}(\bar{\kappa}^{\varepsilon,\delta}_{s,1}\ast\nabla_v\cdot(f_1V\star_x\rho_1))\chi^{\varepsilon,\delta}_{s,2}\varphi_{\beta}\zeta_M\alpha_R,  	
\end{align*}
and the transport terms are defined by 
\begin{align*}
I^{2,1,trans-x}_t=-\int^t_0\int_{\mathbb{R}^{2d+1}}((v\cdot\nabla_xf_1)\ast\bar{\kappa}^{\varepsilon,\delta}_{s,1})\chi^{\varepsilon,\delta}_{s,2}\varphi_{\beta}\zeta_M\alpha_R, 	
\end{align*}
and 
\begin{align*}
I^{2,1,trans-v}_t=&\int^t_0\int_{\mathbb{R}^{2d+1}}((v\cdot\nabla_vf_1)\ast\bar{\kappa}^{\varepsilon,\delta}_{s,1})\chi^{\varepsilon,\delta}_{s,2}\varphi_{\beta}\zeta_M\alpha_R\\
&+\int^t_0\int_{\mathbb{R}^{2d+1}}({\color{black}(d f_1)\ast\bar{\kappa}^{\varepsilon,\delta}_{s,1}})\chi^{\varepsilon,\delta}_{s,2}\varphi_{\beta}\zeta_M\alpha_R. 	
\end{align*}
We remark that $p_1$ and $p_2$ denote the kinetic measures associated with $f_1$ and $f_2$, respectively, while $\rho_1$ and $\rho_2$ represent the corresponding marginal densities. {\color{black}Whenever no confusion can arise, convolution against the kinetic measure means its time density
\begin{align*}
 (\kappa^{\varepsilon,\delta}\ast p_i)(s,z,\eta)
 :=P_i^{\varepsilon,\delta}(s,z,\eta).
\end{align*}}

Analogously to \eqref{21-formula}, the term $\int_0^t\int_{\mathbb{R}^{2d+1}}\chi^{\varepsilon,\delta}_{s,1}d\chi^{\varepsilon,\delta}_{s,2}\varphi_{\beta}\zeta_M\alpha_R$ can be expressed using a similar representation. For brevity, we omit the details. In what follows, we expand the It\^o quadratic variation term:  
\begin{align*}
&\int^t_0\int_{\mathbb{R}^{2d+1}}d\langle\chi^{\varepsilon,\delta}_{s,1},\chi^{\varepsilon,\delta}_{s,2}\rangle_s(z,\eta)\varphi_{\beta}\zeta_M\alpha_R\\
=&\sum_{k\geq1}\int^t_0\int_{\mathbb{R}^{2d+1}}[(f_k\nabla_v\sigma(f_1))\ast\bar{\kappa}^{\varepsilon,\delta}_{s,1}]\cdot[(f_k\nabla_v\sigma(f_2))\ast\bar{\kappa}^{\varepsilon,\delta}_{s,2}]\varphi_{\beta}\zeta_M\alpha_R\\
&+\sum_{k\geq1}\int^t_0\int_{\mathbb{R}^{2d+1}}[(\sigma(f_1)\nabla_vf_k)\ast\bar{\kappa}^{\varepsilon,\delta}_{s,1}]\cdot[(\sigma(f_2)\nabla_vf_k)\ast\bar{\kappa}^{\varepsilon,\delta}_{s,2}]\varphi_{\beta}\zeta_M\alpha_R\\
&+\sum_{k\geq1}\int^t_0\int_{\mathbb{R}^{2d+1}}[(\nabla_v\sigma(f_1)f_k)\ast\bar{\kappa}^{\varepsilon,\delta}_{s,1}]\cdot[(\sigma(f_2)\nabla_vf_k)\ast\bar{\kappa}^{\varepsilon,\delta}_{s,2}]\varphi_{\beta}\zeta_M\alpha_R\\
&+\sum_{k\geq1}\int^t_0\int_{\mathbb{R}^{2d+1}}[(\sigma(f_1)\nabla_vf_k)\ast\bar{\kappa}^{\varepsilon,\delta}_{s,1}]\cdot[(f_k\nabla_v\sigma(f_2))\ast\bar{\kappa}^{\varepsilon,\delta}_{s,2}]\varphi_{\beta}\zeta_M\alpha_R. 	
\end{align*}
Combining all the above formulas, we conclude that 
\begin{align*}
\int_{\mathbb{R}^{2d+1}}\chi^{\varepsilon,\delta}_{s,1}\chi^{\varepsilon,\delta}_{s,2}\varphi_{\beta}\zeta_M\alpha_R\Big|^{t}_{s=0}=I^{err}_t&+I^{meas}_t+I^{cut}_t+I^{mart}_t\\
&+I^{ker}_t+I^{trans-x}_t+I^{trans-v}_t, 	
\end{align*}
where the error term is defined by 
\begin{align}\label{I-err}
&I^{err}_t=I^{2,1,err}_t+I^{1,2,err}_t+\int^t_0\int_{\mathbb{R}^{2d+1}}d\langle\chi^{\varepsilon,\delta}_{s,1},\chi^{\varepsilon,\delta}_{s,2}\rangle_s(z,\eta)\varphi_{\beta}\zeta_M\alpha_R\notag\\
=&-\frac{1}{2}\int^t_0\int_{\mathbb{R}^{2d+1}}\Big[(F_1\sigma'(f_1)^2\nabla_vf_1)\ast\bar{\kappa}^{\varepsilon,\delta}_{s,1}\cdot\nabla_vf_2\ast\bar{\kappa}^{\varepsilon,\delta}_{s,2}+(F_1\sigma'(f_2)^2\nabla_vf_2)\ast\bar{\kappa}^{\varepsilon,\delta}_{s,2}\cdot\nabla_vf_1\ast\bar{\kappa}^{\varepsilon,\delta}_{s,1}\Big]\varphi_{\beta}\zeta_M\alpha_R\notag\\
&{\color{black}-\frac{1}{2}\int^t_0\int_{\mathbb{R}^{4d+1}}[(F_3\sigma^2(f_1))\ast\bar{\kappa}^{\varepsilon,\delta}_{s,1}]\bar{\kappa}^{\varepsilon,\delta}_{s,2}\varphi_{\beta}\zeta_M\alpha_R}\notag\\
&{\color{black}-\frac{1}{2}\int^t_0\int_{\mathbb{R}^{4d+1}}[(F_3\sigma^2(f_2))\ast\bar{\kappa}^{\varepsilon,\delta}_{s,2}]\bar{\kappa}^{\varepsilon,\delta}_{s,1}\varphi_{\beta}\zeta_M\alpha_R}\notag\\
&+\sum_{k\geq1}\int^t_0\int_{\mathbb{R}^{2d+1}}[(f_k\nabla_v\sigma(f_1))\ast\bar{\kappa}^{\varepsilon,\delta}_{s,1}]\cdot[(f_k\nabla_v\sigma(f_2))\ast\bar{\kappa}^{\varepsilon,\delta}_{s,2}]\varphi_{\beta}\zeta_M\alpha_R\notag\\
&+\sum_{k\geq1}\int^t_0\int_{\mathbb{R}^{2d+1}}[(\sigma(f_1)\nabla_vf_k)\ast\bar{\kappa}^{\varepsilon,\delta}_{s,1}]\cdot[(\sigma(f_2)\nabla_vf_k)\ast\bar{\kappa}^{\varepsilon,\delta}_{s,2}]\varphi_{\beta}\zeta_M\alpha_R\notag\\
&+\sum_{k\geq1}\int^t_0\int_{\mathbb{R}^{2d+1}}[(\nabla_v\sigma(f_1)f_k)\ast\bar{\kappa}^{\varepsilon,\delta}_{s,1}]\cdot[(\sigma(f_2)\nabla_vf_k)\ast\bar{\kappa}^{\varepsilon,\delta}_{s,2}]\varphi_{\beta}\zeta_M\alpha_R\notag\\
&+\sum_{k\geq1}\int^t_0\int_{\mathbb{R}^{2d+1}}[(\sigma(f_1)\nabla_vf_k)\ast\bar{\kappa}^{\varepsilon,\delta}_{s,1}]\cdot[(f_k\nabla_v\sigma(f_2))\ast\bar{\kappa}^{\varepsilon,\delta}_{s,2}]\varphi_{\beta}\zeta_M\alpha_R, 
\end{align}
the measure term is defined by 
{\color{black}
\begin{align*}
I^{meas}_t={}&I^{2,1,meas}_t+I^{1,2,meas}_t\\
={}&\int_0^t\int_{\mathbb R^{4d+1}}
 P_1^{\varepsilon,\delta}(s,z,\eta)
 \bar\kappa^{\varepsilon,\delta}_{s,2}(z,w,\eta)
 \varphi_\beta(\eta)\zeta_M(\eta)\alpha_R(z)
 \,dz\,dw\,d\eta\,ds\\
&+\int_0^t\int_{\mathbb R^{4d+1}}
 P_2^{\varepsilon,\delta}(s,z,\eta)
 \bar\kappa^{\varepsilon,\delta}_{s,1}(z,w,\eta)
 \varphi_\beta(\eta)\zeta_M(\eta)\alpha_R(z)
 \,dz\,dw\,d\eta\,ds\\
&-2\int_0^t\int_{\mathbb R^{2d+1}}
 (\nabla_vf_1\ast\bar\kappa^{\varepsilon,\delta}_{s,1})(z,\eta)
 \cdot(\nabla_vf_2\ast\bar\kappa^{\varepsilon,\delta}_{s,2})(z,\eta)
 \varphi_\beta(\eta)\zeta_M(\eta)\alpha_R(z)
 \,dz\,d\eta\,ds.	
\end{align*}}
the cutoff term, martingale term, kernel term and transport terms are defined by 
\begin{align*}
I^{cut}_t=&I^{1,cut}_t+I^{2,cut}_t-2(I^{2,1,cut}_t+I^{1,2,cut}_t),\\
I^{mart}_t=&I^{1,mart}_t+I^{2,mart}_t-2(I^{2,1,mart}_t+I^{1,2,mart}_t),\\
I^{ker}_t=&I^{1,ker}_t+I^{2,ker}_t-2(I^{2,1,ker}_t+I^{1,2,ker}_t),\\
I^{trans-x}_t=&I^{1,trans-x}_t+I^{2,trans-x}_t-2(I^{2,1,trans-x}_t+I^{1,2,trans-x}_t),\\
I^{trans-v}_t=&I^{1,trans-v}_t+I^{2,trans-v}_t-2(I^{2,1,trans-v}_t+I^{1,2,trans-v}_t).  	
\end{align*}
Then it follows that  
\begin{align*}
&\int_{\mathbb{R}^{2d+1}}(\chi^{\varepsilon,\delta}_{s,1}+\chi^{\varepsilon,\delta}_{s,2}-2\chi^{\varepsilon,\delta}_{s,1}\chi^{\varepsilon,\delta}_{s,2})\varphi_{\beta}\zeta_M{\color{black}\alpha_R}\Big|^{t}_{s=0}\\
=&-2I^{err}_t-2I^{meas}_t+I^{mart}_t+I^{cut}_t+I^{ker}_t+I^{trans-x}_t+I^{trans-v}_t. 	
\end{align*}

{\bf The measure term. }Thanks to the {\color{black}kinetic-measure domination \eqref{control KM}}, applying H\"older's inequality and Young's inequality, we have that almost surely for every $t\in[0,T]$, 
\begin{align*}
I^{meas}_t\geq0. 	
\end{align*}

{\bf The error term. }
Returning to the identity \eqref{I-err}, by collecting the corresponding terms, we obtain the $F_1$-contribution involving the square of the difference between $\sigma'(f_1)$ and $\sigma'(f_2)$ and the analogous $F_3$-contribution. Invoking \cite[(4.16)]{FG24}, we can then pass to the limit as $\varepsilon \to 0$, which yields 
\begin{align*}
\lim_{\varepsilon \to 0} I^{\mathrm{err}}_t &= {\color{black}-\frac12\int_0^t \int_{\mathbb{R}^{2d+1}} \Big( F_1 (\sigma'(f_1) - \sigma'(f_2))^2 \nabla_v f_1 \cdot \nabla_v f_2 + F_3 (\sigma(f_1) - \sigma(f_2))^2 \Big) \bar{\kappa}^\delta_{s,1} \bar{\kappa}^\delta_{s,2} \varphi_\beta \zeta_M \alpha_R}.
\end{align*}

{\color{black}Under the regular-coefficient Assumption~\ref{Assump-sigma-n}, the following bound holds:}
\begin{align*}
\limsup_{\varepsilon\rightarrow0}|I^{err}_t|\leq c\int^T_0\int_{\mathbb{R}^{2d+1}}I_{\{0<|f_1-f_2|<c\delta\}}\Big(1+|\nabla_vf_1|^2+|\nabla_vf_2|^2\Big)(\delta\bar{\kappa}^{\delta}_{s,1})\bar{\kappa}^{\delta}_{s,2}\varphi_{\beta}\zeta_M\alpha_R. 	
\end{align*}
Thanks to the boundedness of $\delta \bar{\kappa}^\delta_{s,1}$ and the regularity estimate \eqref{L2-es}, and by applying the dominated convergence theorem, we conclude that almost surely, for every $t \in [0, T]$,
\begin{align*}
\limsup_{\delta \to 0} \limsup_{\varepsilon \to 0} |I^{\mathrm{err}}_t| = 0.
\end{align*}
{\bf The cutoff term. }
For the cutoff term, by the proof of \cite[(4.30)]{FG24}, 
{\color{black}\begin{align*}
&\limsup_{\delta\rightarrow0}\limsup_{\varepsilon\rightarrow0}I^{cut}_t\\
\leq{}&\int_{\mathbb R^{2d}\times(0,\infty)\times[0,t]}
 |\partial_\eta(\varphi_\beta\zeta_M)(\eta)|\alpha_R(z)
 \, (p_1+p_2)(dz,d\eta,ds)\\
&-\frac12\int_0^t\int_{\mathbb R^{2d}}
 \sgn(f_2-f_1)F_3
 \Big(\sigma(f_1)^2\partial_\eta(\varphi_\beta\zeta_M)(f_1)
 -\sigma(f_2)^2\partial_\eta(\varphi_\beta\zeta_M)(f_2)\Big)
 \alpha_R\,dz\,ds\\
&+\int_0^t\int_{\mathbb R^{2d}}
 \sgn(f_2-f_1)\Big(
 (\varphi_\beta\zeta_M)(f_1)\nabla_vf_1
 -(\varphi_\beta\zeta_M)(f_2)\nabla_vf_2\Big)
 \cdot\nabla_v\alpha_R\,dz\,ds\\
&+\int_0^t\int_{\mathbb R^{2d}}
 \sgn(f_2-f_1)\Big(
 F_1\sigma'(f_1)^2(\varphi_\beta\zeta_M)(f_1)\nabla_vf_1\\
&\hspace{42mm}
 -F_1\sigma'(f_2)^2(\varphi_\beta\zeta_M)(f_2)\nabla_vf_2\Big)
 \cdot\nabla_v\alpha_R\,dz\,ds.
\end{align*}}
{\color{black}By \eqref{common-low-density-sequence}, the two low-density
kinetic-measure terms tend to zero almost surely along the already fixed
sequence $(\beta_j)_{j\geq1}$. Applying the same argument as in
\cite[(4.30)]{FG24} along this sequence gives the upper bound needed below:}
\begin{align*}
{\color{black}\limsup_{j\rightarrow\infty}\limsup_{M\rightarrow\infty}
\limsup_{R\rightarrow\infty}\limsup_{\delta\rightarrow0}
\limsup_{\varepsilon\rightarrow0}
\sup_{t\in[0,T]}I^{cut}_t\leq0,
\qquad \beta=\beta_j.}	
\end{align*}
{\color{black}From this point until the end of the proof, every occurrence of
$\beta\to0$ is understood along the deterministic sequence
$(\beta_j)_{j\geq1}$ fixed in \eqref{common-low-density-sequence}.}

{\bf The martingale term. }
{\color{black}Applying the same strategy as in \cite[(4.19), (4.20)]{FG24},
we obtain, almost surely, along subsequences chosen first as
$\varepsilon\to0$ and then as $\delta\to0$,}
\begin{align*}
{\color{black}\lim_{\delta\rightarrow0}\lim_{\varepsilon\rightarrow0}}I^{mart}_t 
=& \int^t_0\int_{\mathbb{R}^{2d}} \sgn(f_2 - f_1) \varphi_{\beta}(f_1) \zeta_M(f_1) \alpha_R(z) \nabla_v\cdot(\sigma(f_1)\, dW_F) \\
& - \int^t_0\int_{\mathbb{R}^{2d}} \sgn(f_2 - f_1) \varphi_{\beta}(f_2) \zeta_M(f_2) \alpha_R(z) \nabla_v\cdot(\sigma(f_2)\, dW_F).
\end{align*}
{\color{black}For the martingale term, define the auxiliary function
$\Theta^{\sigma}_{\beta,M}(\zeta) := \int^{\zeta}_0
\varphi_{\beta}(\zeta') \zeta_M(\zeta') \sigma'(\zeta')\,d\zeta'$.}
Using the chain rule and the integration-by-parts formula, the martingale term can be rewritten as
\begin{align*}
{\color{black}\lim_{\delta\rightarrow0}\lim_{\varepsilon\rightarrow0}}I^{mart}_t 
=& {\color{black}\int^t_0\int_{\mathbb{R}^{2d}} \sgn(f_2 - f_1) \alpha_R(z) \nabla_v \cdot \Big( (\Theta^{\sigma}_{\beta,M}(f_1) - \Theta^{\sigma}_{\beta,M}(f_2)) dW_F \Big)} \\
& {\color{black}+ \int^t_0\int_{\mathbb{R}^{2d}} \sgn(f_2 - f_1) \alpha_R(z) \Big( \varphi_{\beta}(f_1)\zeta_M(f_1)\sigma(f_1) - \Theta^{\sigma}_{\beta,M}(f_1) \Big) \nabla_v\cdot dW_F} \\
& {\color{black}- \int^t_0\int_{\mathbb{R}^{2d}} \sgn(f_2 - f_1) \alpha_R(z) \Big( \varphi_{\beta}(f_2)\zeta_M(f_2)\sigma(f_2) - \Theta^{\sigma}_{\beta,M}(f_2) \Big) \nabla_v\cdot dW_F.}
\end{align*}
Finally, based on the estimates for {\color{black}$\Theta^{\sigma}_{\beta,M}$} and $\varphi_{\beta} \zeta_M \sigma$ given in the proof of \cite[(4.24)]{FG24} and \cite[(2.26)]{FG25}, we conclude that
\begin{align*}
{\color{black}\lim_{j\rightarrow\infty}\lim_{M\rightarrow\infty}
\lim_{R\rightarrow\infty}\lim_{\delta\rightarrow0}
\lim_{\varepsilon\rightarrow0}
\sup_{t\in[0,T]}|I^{mart}_t|=0,
\qquad \beta=\beta_j.}
\end{align*} 

{\bf The kernel term. } {\color{black}Set
\begin{align*}
 \mathsf M:=\max_{i=1,2}\operatorname*{ess\,sup}_{\omega\in\Omega}
 \|f_{i,0}(\omega)\|_{L^1(\mathbb R^{2d})}<\infty.
\end{align*}
For $i=1,2$, mass conservation and the Fisher-information estimate in
\eqref{L2-es} give, almost surely,}
\begin{align*}
\|\nabla_v \cdot (f_i V \star_x \rho_i)\|_{L^1([0,T];L^1(\mathbb{R}^{2d}))}
&= \int_0^T \int_{\mathbb{R}^{2d}} |\nabla_v f_i \, V \star_x \rho_i| \, dx\,dv\,dt \\
&{\color{black}\leq 2\|V\|_{L^\infty(\mathbb R^d)}\mathsf M^{3/2}T^{1/2}
 \left(\int_0^T\|\nabla_v\sqrt{f_i}\|_{L^2(\mathbb R^{2d})}^2\,ds\right)^{1/2}<\infty.}
\end{align*}
{\color{black}Indeed,
$\|V\star_x\rho_i\|_{L_x^\infty}\leq
\|V\|_{L_x^\infty}\|f_i\|_{L^1_{x,v}}$ and
$\|\nabla_vf_i\|_{L^1_{x,v}}{\color{black}\leq}
2\|f_i\|_{L^1_{x,v}}^{1/2}\|\nabla_v\sqrt{f_i}\|_{L^2_{x,v}}$.
This is the point at which the Fisher-information condition in Definition~\ref{def-kineticsolution} is essential on the whole velocity space.}
Therefore, due to the boundedness of $\bar{\kappa}^{\delta}_{s,i}$ and the $L^1(\mathbb{R}^{2d})$-integrability of $\nabla_v \cdot (f_i V \star_x \rho_i)$ for $i=1,2$, we can pass to the limit $\varepsilon \rightarrow 0$. More precisely, almost surely, for every $t\in[0,T]$, 
\begin{align} \label{eq-4.13}
{\color{black}\lim_{\varepsilon\rightarrow0}I^{ker}_t}
=& {\color{black}\int_0^t \int_{\mathbb{R}^{2d+1}} \bar{\kappa}^{\delta}_{s,1} (2\chi^{\delta}_{s,2} - 1) \varphi_{\beta} \zeta_M \alpha_R\nabla_v \cdot (f_1 V \star_x \rho_1) \, dx\,dv\,d\eta\,ds} \notag \\
& {\color{black}+ \int_0^t \int_{\mathbb{R}^{2d+1}} \bar{\kappa}^{\delta}_{s,2} (2\chi^{\delta}_{s,1} - 1) \varphi_{\beta} \zeta_M \alpha_R\nabla_v \cdot (f_2 V \star_x \rho_2) \, dx\,dv\,d\eta\,ds.}
\end{align}
Since for each $i = 1,2$, the function $\nabla_v \cdot (f_i V \star_x \rho_i)$ belongs to $L^1(\Omega \times [0,T] \times \mathbb{R}^{2d})$, we may pass to a subsequence as $\delta \to 0$. Then, almost surely, we have
\begin{equation} \label{eq-4.14}
{\color{black}\lim_{\delta \to 0} \left| \int_0^t \int_{\mathbb{R}^{2d+1}} \bar{\kappa}^{\delta}_{s,1} (2\chi^{\delta}_{s,2} - 1) \left( \varphi_{\beta}(\eta) \zeta_M(\eta) - \varphi_{\beta}(f_1) \zeta_M(f_1) \right) \alpha_R\nabla_v \cdot (f_1 V \star_x \rho_1)\,dx\,dv\,d\eta\,ds \right| = 0.}
\end{equation}
Similarly, for the second term on the right-hand side of \eqref{eq-4.13}, it holds that
\begin{equation} \label{kk-49}
{\color{black}\lim_{\delta \to 0} \left| \int_0^t \int_{\mathbb{R}^{2d+1}} \bar{\kappa}^{\delta}_{s,2} (2\chi^{\delta}_{s,1} - 1) \left( \varphi_{\beta}(\eta) \zeta_M(\eta) - \varphi_{\beta}(f_2) \zeta_M(f_2) \right) \alpha_R\nabla_v \cdot (f_2 V \star_x \rho_2)\,dx\,dv\,d\eta\,ds \right| = 0.}
\end{equation}
Moreover, referring to \cite[(4.22)]{FG24} and using the fact that $\varphi_{\beta}(0) = 0$, we deduce that {\color{black}the following identity holds almost surely for almost every $(s,x,v)\in[0,T]\times\mathbb R^{2d}$}:
\begin{equation} \label{eq-4.16}
{\color{black}\lim_{\delta \to 0} \left( \int_{\mathbb{R}} \bar{\kappa}^{\delta}_{s,1} (2\chi^{\delta}_{s,2} - 1) \, d\eta \right) \varphi_{\beta}(f_1)} 
= \left( I_{\{f_1 = f_2\}} + 2 I_{\{0 \le f_1 < f_2\}} - 1 \right) \varphi_{\beta}(f_1).
\end{equation}
Combining \eqref{eq-4.16} with the previous estimates and passing to a subsequence as $\delta \to 0$, we obtain almost surely:
\begin{align}
{\color{black}\lim_{\delta\to0}\lim_{\varepsilon\to0}I_t^{\mathrm{ker}}}
=& {\color{black}\int_0^t \int_{\mathbb{R}^{2d}} \left( I_{\{f_1 = f_2\}} + 2 I_{\{f_1 < f_2\}} - 1 \right) \varphi_{\beta}(f_1) \zeta_M(f_1) \alpha_R\nabla_v \cdot (f_1 V \star_x \rho_1)} \notag \\
&{\color{black}+ \int_0^t \int_{\mathbb{R}^{2d}} \left( I_{\{f_1 = f_2\}} + 2 I_{\{f_2 < f_1\}} - 1 \right) \varphi_{\beta}(f_2) \zeta_M(f_2) \alpha_R\nabla_v \cdot (f_2 V \star_x \rho_2).} \label{kk-50}
\end{align}

{\color{black}For the same sequence $(\beta_j)_{j\geq1}$ and the order of
limits fixed above, dominated convergence gives, for each $i=1,2$, almost
surely,}
\begin{equation} \label{eq-4.22}
{\color{black}\lim_{j\to\infty}\lim_{M\to\infty}
 \varphi_{\beta_j}(f_i)\zeta_M(f_i)
 \nabla_v\cdot(f_iV\star_x\rho_i)
=\nabla_v\cdot(f_iV\star_x\rho_i)
\quad\text{strongly in }L^1(\mathbb R^{2d}\times[0,T]).}
\end{equation}

{\color{black}After the $\delta$-limit, the bound $0\leq\alpha_R\leq1$ and
the $L^1$-integrability in \eqref{eq-4.22} permit removal of the phase-space
cutoff in the prescribed order $R_1\to\infty$ and then
$R_2\to\infty$.} Combining \eqref{kk-50} and \eqref{eq-4.22}, we conclude that
\begin{align}
&{\color{black}\lim_{j\to\infty}\lim_{M\to\infty}
 \lim_{R\to\infty}\lim_{\delta\to0}\lim_{\varepsilon\to0}
 I_t^{\mathrm{ker}}} \notag \\
=& \int_0^t \int_{\mathbb{R}^{2d}} \left( I_{\{f_1 = f_2\}} + 2 I_{\{f_1 < f_2\}} - 1 \right) \left( \nabla_v \cdot (f_1 V \star_x \rho_1) - \nabla_v \cdot (f_2 V \star_x \rho_2) \right) \notag \\
& + \int_0^t \int_{\mathbb{R}^{2d}} \left( I_{\{f_1 = f_2\}} + 2 I_{\{f_1 < f_2\}} - 1 + I_{\{f_1 = f_2\}} + 2 I_{\{f_2 < f_1\}} - 1 \right) \nabla_v \cdot (f_2 V \star_x \rho_2) \notag \\
=:& \, \tilde{J}_1 + \tilde{J}_2. \label{eq-4.23}
\end{align}

Since the identity 
$$
I_{\{f_2 < f_1\}} = 1 - I_{\{f_1 = f_2\}} - I_{\{f_1 < f_2\}}
$$
holds, and $\nabla_v \cdot (f_i V \star_x \rho_i) \in L^1(\Omega \times [0,T] \times \mathbb{R}^{2d})$ for each $i=1,2$, it follows that
\begin{equation} \label{eq-4.24}
\tilde{J}_2 = 0.
\end{equation}
	For the term $\tilde{J}_1$, using the chain rule and the identity
$$
\sgn(f_2 - f_1) = I_{\{f_1 = f_2\}} + 2I_{\{f_1 < f_2\}} - 1,
$$
we decompose:
\begin{align*}
\tilde{J}_1 := \tilde{J}_{11} + \tilde{J}_{12},
\end{align*}
where
\begin{align*}
\tilde{J}_{11} &= \int_0^t \int_{\mathbb{R}^{2d}} \sgn(f_2 - f_1)\nabla_v \cdot \left((f_1 - f_2)V \star_x \rho_1\right), \\
\tilde{J}_{12} &= \int_0^t \int_{\mathbb{R}^{2d}} \sgn(f_2 - f_1)\nabla_v \cdot \left(f_2 V \star_x (\rho_1 - \rho_2)\right).
\end{align*}

Define $\sgn^{\delta} := \sgn \ast \kappa_1^{\delta}$ for every $\delta \in (0,1)$. By the integration by parts formula, we obtain that almost surely for every $t \in [0,T]$,
\begin{align}
\tilde{J}_{11} 
&= \lim_{\delta \to 0} \int_0^t \int_{\mathbb{R}^{2d}} \sgn^{\delta}(f_2 - f_1)\nabla_v \cdot \left((f_1 - f_2)V \star_x \rho_1\right) \notag \\
&= -\lim_{\delta \to 0} \int_0^t \int_{\mathbb{R}^{2d}} (\sgn^{\delta})'(f_2 - f_1)(f_1 - f_2)\nabla_v(f_2 - f_1) \cdot V \star_x \rho_1. \notag
\end{align}

Due to the uniform boundedness of $(\delta \kappa_1^{\delta})$ for $\delta \in (0,\beta/4)$, there exists a constant $c > 0$, independent of $\delta$ but depending on the convolution kernel, such that for all $\delta \in (0,\beta/4)$,
$$
\left|(\sgn^{\delta})'(f_2 - f_1)(f_1 - f_2)\right| = 2\left| \kappa_1^{\delta}(f_2 - f_1)(f_1 - f_2) \right| 
\le c I_{\{0 < |f_1 - f_2| < \delta\}}.
$$

{\color{black}The preceding estimate also yields the
$L^1(\Omega\times[0,T]\times\mathbb R^{2d})$-integrability of
$\nabla_v(f_2-f_1)\cdot(V\star_x\rho_1)$.} Therefore, by the dominated convergence theorem, we deduce that, almost surely for every $t \in [0,T]$,
\begin{equation} \label{eq-4.27}
|\tilde{J}_{11}| {\color{black}\leq c\int^t_0\int_{\mathbb{R}^{2d}}\lim_{\delta\rightarrow0}I_{\{0<|f_1-f_2|<\delta\}}
|\nabla_v(f_1-f_2)|\,|V\star_x\rho_1|=0.}
\end{equation}

{\color{black}To estimate $\tilde J_{12}$, first note that
\begin{align*}
 \|V\star_x(\rho_1-\rho_2)\|_{L_x^\infty}
 \leq \|V\|_{L_x^\infty}\|f_1-f_2\|_{L^1_{x,v}}.
\end{align*}
Consequently, using the Fisher information rather than an invalid
$L^2(\mathbb R_v^d)$ estimate for a function independent of $v$, we obtain}
\begin{align*}
|\tilde{J}_{12}| 
&{\color{black}\leq 2\|V\|_{L^{\infty}(\mathbb{R}^{d})}\mathsf M^{1/2}
 \int_0^t \|\nabla_v \sqrt{f_2(s)}\|_{L^2(\mathbb{R}^{2d})}
 \|f_1(s) - f_2(s)\|_{L^1(\mathbb{R}^{2d})}\,ds.}
\end{align*}

Consequently, we conclude that
\begin{align*}
{\color{black}\lim_{j\to\infty}\lim_{M\to\infty}\lim_{R\to\infty}
\lim_{\delta\to0}\lim_{\varepsilon\to0}I_t^{\mathrm{ker}}}
{\color{black}\leq 2\|V\|_{L^{\infty}(\mathbb{R}^{d})}\mathsf M^{1/2}
 \int_0^t \|\nabla_v \sqrt{f_2(s)}\|_{L^2(\mathbb{R}^{2d})}
 \|f_1(s) - f_2(s)\|_{L^1(\mathbb{R}^{2d})}\,ds.}
\end{align*}	

{\bf The $x$-transport term. }
We first apply the integration by parts formula to the term $I^{2,1,\mathrm{trans\text{-}x}}_t$, yielding
\begin{align*}
I^{2,1,\mathrm{trans\text{-}x}}_t
&= -\int_0^t \int_{\mathbb{R}^{2d+1}} (v \cdot \nabla_x \chi^{\delta}_{s,1}) \ast \kappa^{\varepsilon} \, \chi^{\varepsilon,\delta}_{s,2} \, \varphi_{\beta} \zeta_M \alpha_R \\
&= \int_0^t \int_{\mathbb{R}^{2d+1}} (v \chi^{\delta}_{s,1}) \ast \kappa^{\varepsilon} \, \nabla_x \chi^{\varepsilon,\delta}_{s,2} \, \varphi_{\beta} \zeta_M \alpha_R  + \int_0^t \int_{\mathbb{R}^{2d+1}} (v \chi^{\delta}_{s,1}) \ast \kappa^{\varepsilon} \, \chi^{\varepsilon,\delta}_{s,2} \, \varphi_{\beta} \zeta_M \nabla_x \alpha_R.
\end{align*}
Recall also that
$$
I^{1,2,\mathrm{trans\text{-}x}}_t = -\int_0^t \int_{\mathbb{R}^{2d+1}} (v \cdot \nabla_x \chi^{\delta}_{s,2}) \ast \kappa^{\varepsilon} \, \chi^{\varepsilon,\delta}_{s,1} \, \varphi_{\beta} \zeta_M \alpha_R.
$$

We now observe that, after formally passing to the limit $\varepsilon_v \to 0$, the first term on the right-hand side of the expression for $I^{2,1,\mathrm{trans\text{-}x}}_t$ cancels with $I^{1,2,\mathrm{trans\text{-}x}}_t$. 

To make this rigorous, we next pass to the limit $\varepsilon_v \to 0$. Thanks to the dominated convergence theorem, this limit is valid, and we obtain:
\begin{align*}
\lim_{\varepsilon_v \to 0} I^{2,1,\mathrm{trans\text{-}x}}_t 
&= \int_0^t \int_{\mathbb{R}^{2d+1}} v \chi^{\varepsilon_x,\delta}_{s,1} \nabla_x \chi^{\varepsilon_x,\delta}_{s,2} \, \varphi_{\beta} \zeta_M \alpha_R  + \int_0^t \int_{\mathbb{R}^{2d+1}} v \chi^{\varepsilon_x,\delta}_{s,1} \chi^{\varepsilon_x,\delta}_{s,2} \, \varphi_{\beta} \zeta_M \nabla_x \alpha_R,
\end{align*}
and
$$
\lim_{\varepsilon_v \to 0} I^{1,2,\mathrm{trans\text{-}x}}_t 
= -\int_0^t \int_{\mathbb{R}^{2d+1}} v \chi^{\varepsilon_x,\delta}_{s,1} \nabla_x \chi^{\varepsilon_x,\delta}_{s,2} \, \varphi_{\beta} \zeta_M \alpha_R.
$$
Here we adopt the notation $\chi^{\varepsilon_x,\delta}_{s,i}$ for $i = 1, 2$ to denote $\chi^{\varepsilon_x,\delta}_{s,i} := \chi_{s,i} \ast \left( \kappa^{\varepsilon_x}_d \, \kappa^{\delta}_1 \right), \quad i = 1, 2$. Consequently, we deduce that almost surely, for every $t\in[0,T]$, 
\begin{align} \label{epxto0}
\lim_{\varepsilon_v \to 0} \left( I^{1,2,\mathrm{trans\text{-}x}}_t + I^{2,1,\mathrm{trans\text{-}x}}_t \right) 
= \int_0^t \int_{\mathbb{R}^{2d+1}} v \chi^{\varepsilon_x,\delta}_{s,1} \chi^{\varepsilon_x,\delta}_{s,2} \, \varphi_{\beta} \zeta_M \nabla_x \alpha_R.
\end{align}
Based on the above analysis, a direct computation yields
\begin{align*}
I^{\mathrm{trans\text{-}x}}_t 
=&{\color{black}\ 2\int_0^t \int_{\mathbb{R}^{2d+1}} (v \cdot \nabla_x \chi^\delta_{s,1}) \ast \kappa^\varepsilon \, \chi^{\varepsilon,\delta}_{s,2} \, \varphi_\beta \zeta_M \alpha_R}\\
&{\color{black}+2\int_0^t \int_{\mathbb{R}^{2d+1}} (v \cdot \nabla_x \chi^\delta_{s,2}) \ast \kappa^\varepsilon \, \chi^{\varepsilon,\delta}_{s,1} \, \varphi_\beta \zeta_M \alpha_R}\\
&{\color{black}-\int_0^t \int_{\mathbb{R}^{2d+1}} (v \cdot \nabla_x \chi^\delta_{s,1}) \ast \kappa^\varepsilon \, \varphi_\beta \zeta_M \alpha_R
-\int_0^t \int_{\mathbb{R}^{2d+1}} (v \cdot \nabla_x \chi^\delta_{s,2}) \ast \kappa^\varepsilon \, \varphi_\beta \zeta_M \alpha_R.}
\end{align*}

Using the result from \eqref{epxto0}, we obtain
\begin{align*}
\lim_{\varepsilon_v \to 0} I^{\mathrm{trans\text{-}x}}_t 
=&{\color{black}\int_0^t \int_{\mathbb{R}^{2d+1}}
 v\bigl(\chi^{\varepsilon_x,\delta}_{s,1}
 +\chi^{\varepsilon_x,\delta}_{s,2}
 -2\chi^{\varepsilon_x,\delta}_{s,1}\chi^{\varepsilon_x,\delta}_{s,2}\bigr)
 \varphi_\beta \zeta_M \cdot \nabla_x \alpha_R.}
\end{align*}

Due to the choice of the truncation function $\alpha_R$, we have the estimate
$$
|\nabla_x \alpha_R| \leq \frac{c}{R_1} \alpha_{R_2}^2.
$$
Hence, by the dominated convergence theorem, we may pass to the limits
{\color{black}$\varepsilon_v\to0$, then $\varepsilon_x\to0$, and finally
$\delta\to0$}, leading to
\begin{align*}
{\color{black}\limsup_{\delta \to 0}  \limsup_{\varepsilon_x \to 0}\limsup_{\varepsilon_v \to 0}
 |I^{\mathrm{trans\text{-}x}}_t|}
\leq&{\color{black}\ \frac{c}{R_1} \int_0^t \int_{\mathbb{R}^{2d+1}}
 |v|\bigl(\chi_{s,1}+\chi_{s,2}+2\chi_{s,1}\chi_{s,2}\bigr)
 \varphi_\beta \zeta_M \alpha_{R_2}^2}\\
\leq&{\color{black}\ \frac{C R_2}{R_1}\int_0^t
 \bigl(\|f_1(s)\|_{L^1}+\|f_2(s)\|_{L^1}\bigr)\,ds.}
\end{align*}

Since the kinetic functions $\chi_{s,1}$ and $\chi_{s,2}$ are uniformly bounded, and $\alpha_{R_2}^2$ serves as a velocity cutoff, the integrals on the right-hand side are well-defined. Therefore we conclude that 
$$
{\color{black}\lim_{R_1 \to \infty} \limsup_{\delta \to 0}
\limsup_{\varepsilon_x \to 0}\limsup_{\varepsilon_v \to 0}
|I^{\mathrm{trans\text{-}x}}_t| = 0,}
$$
and hence, almost surely for every $t \in [0,T]$,
$$
{\color{black}\lim_{j\to\infty}\lim_{M\to\infty}\lim_{R\to\infty}
\limsup_{\delta\to0}\limsup_{\varepsilon_x\to0}
\limsup_{\varepsilon_v\to0} I^{\mathrm{trans\text{-}x}}_t=0,
\qquad\beta=\beta_j.}
$$

{\bf The $v$-transport term. } By the definitions of $I^{1,\mathrm{trans\text{-}v}}_t$, $I^{2,\mathrm{trans\text{-}v}}_t$, $I^{2,1,\mathrm{trans\text{-}v}}_t$, and $I^{1,2,\mathrm{trans\text{-}v}}_t$, we have
\begin{align*}
I^{\mathrm{trans\text{-}v}}_t 
=&\ I^{1,\mathrm{trans\text{-}v}}_t + I^{2,\mathrm{trans\text{-}v}}_t - 2\big(I^{2,1,\mathrm{trans\text{-}v}}_t + I^{1,2,\mathrm{trans\text{-}v}}_t\big) \\
=&\ \int_0^t \int_{\mathbb{R}^{2d+1}} \big(\nabla_v \cdot (v f_1) \ast \bar{\kappa}^{\varepsilon,\delta}_{s,1} \big) \big(1 - 2 \chi^{\varepsilon,\delta}_{s,2}\big) \varphi_\beta \zeta_M \alpha_R \\
&+ \int_0^t \int_{\mathbb{R}^{2d+1}} \big(\nabla_v \cdot (v f_2) \ast \bar{\kappa}^{\varepsilon,\delta}_{s,2} \big) \big(1 - 2 \chi^{\varepsilon,\delta}_{s,1}\big) \varphi_\beta \zeta_M \alpha_R.
\end{align*}

Using the same argument as for the martingale term, we can pass first to
{\color{black}$\varepsilon\to0$ and then to $\delta\to0$}; applying the chain
rule, we deduce
\begin{align*}
{\color{black}\lim_{\delta\to0}\lim_{\varepsilon\to0}}
I^{\mathrm{trans\text{-}v}}_t 
=&\ \int_0^t \int_{\mathbb{R}^{2d}} \sgn(f_2 - f_1) \varphi_\beta(f_1) \zeta_M(f_1) \alpha_R(z) \nabla_v \cdot (v f_1) \\
&- \int_0^t \int_{\mathbb{R}^{2d}} \sgn(f_2 - f_1) \varphi_\beta(f_2) \zeta_M(f_2) \alpha_R(z) \nabla_v \cdot (v f_2) \\
=&\ \int_0^t \int_{\mathbb{R}^{2d}} \sgn(f_2 - f_1) \varphi_\beta(f_1) \zeta_M(f_1) \alpha_R(z) \big(v \cdot \nabla_v f_1 + {\color{black}d f_1} \big) \\
&- \int_0^t \int_{\mathbb{R}^{2d}} \sgn(f_2 - f_1) \varphi_\beta(f_2) \zeta_M(f_2) \alpha_R(z) \big(v \cdot \nabla_v f_2 + {\color{black}d f_2} \big).
\end{align*}

{\color{black}To avoid confusion with the noise primitive introduced above, define}
$$
{\color{black}H_{\beta,M}(\zeta) := \int_0^\zeta \varphi_\beta(\zeta') \zeta_M(\zeta')\,d\zeta'.}
$$
Then, by the chain rule, we obtain
\begin{align*}
{\color{black}\lim_{\delta\to0}\lim_{\varepsilon\to0}}
I^{\mathrm{trans\text{-}v}}_t 
=&\ {\color{black}\int_0^t \int_{\mathbb{R}^{2d}} \sgn(f_2 - f_1) \big( \nabla_v H_{\beta,M}(f_1) - \nabla_v H_{\beta,M}(f_2) \big) \cdot \alpha_R(z) v} \\
&+ {\color{black}d\int_0^t \int_{\mathbb{R}^{2d}} \sgn(f_2 - f_1) \big( \varphi_\beta(f_1) \zeta_M(f_1) f_1 - \varphi_\beta(f_2) \zeta_M(f_2) f_2 \big)} \\
=:\ &I_1 + I_2.
\end{align*}

We first consider $I_2$. By the dominated convergence theorem, we can pass to the limits and conclude
$$
{\color{black}\lim_{j\to\infty}\lim_{M\to\infty}I_2
=d\int_0^t \int_{\mathbb{R}^{2d}}
\sgn(f_2-f_1)(f_1-f_2)
=-d\int_0^t\|f_1-f_2\|_{L^1(\mathbb R^{2d})},
\qquad\beta=\beta_j.}
$$

For $I_1$, again by the chain rule and dominated convergence, we obtain
\begin{align*}
I_1 
=&\ {\color{black}\lim_{\delta \to 0} \int_0^t \int_{\mathbb{R}^{2d}} \sgn^\delta(f_2 - f_1) \nabla_v \big( H_{\beta,M}(f_1) - H_{\beta,M}(f_2) \big) \cdot \alpha_R(z) v} \\
=&\ {\color{black}\lim_{\delta \to 0} \int_0^t \int_{\mathbb{R}^{2d}} \sgn^\delta(f_2 - f_1) \nabla_v \cdot \Big( \alpha_R(z) v \big( H_{\beta,M}(f_1) - H_{\beta,M}(f_2) \big) \Big)} \\
&\quad {\color{black}- \lim_{\delta \to 0} \int_0^t \int_{\mathbb{R}^{2d}} \sgn^\delta(f_2 - f_1) \big( H_{\beta,M}(f_1) - H_{\beta,M}(f_2) \big) \nabla_v \cdot ( \alpha_R(z) v )} \\
=:\ &I_{11} + I_{12}.
\end{align*}

For $I_{11}$, the boundedness of $\varphi_{\beta} \zeta_M$ implies that {\color{black}$H_{\beta,M}$} is a Lipschitz function. Applying integration by parts and using {\color{black}the definition of $\sgn^\delta$}, thanks to the dominated convergence theorem, we estimate 
\begin{align*}
|I_{11}| 
\leq&\ \lim_{\delta \to 0} \int_0^t \int_{\mathbb{R}^{2d}} I_{\{0 \leq |f_1 - f_2| < \delta\}} |\nabla_v(f_1 - f_2)| \, |\alpha_R(z) v| \\
\leq&\ {\color{black}2R_2} \int_0^t \int_{\mathbb{R}^{2d}} \lim_{\delta \to 0} I_{\{0 \leq |f_1 - f_2| < \delta\}} |\nabla_v(f_1 - f_2)| = 0.
\end{align*}
{\color{black}In the last equality we used the Sobolev property
$\nabla_v(f_1-f_2)=0$ almost everywhere on $\{f_1=f_2\}$, together with dominated convergence on the support of $\alpha_R$.}

For $I_{12}$, by the properties of the truncation function $\alpha_R$, we have
\begin{align*}
|I_{12}|
\leq&\ \int_0^t \int_{\mathbb{R}^{2d}} \Big| \int_{f_2}^{f_1} \varphi_\beta(\zeta') \zeta_M(\zeta') \, d\zeta' \Big| \Big| \nabla_v \alpha_R(z) \cdot v + {\color{black}d\alpha_R(z)} \Big| \\
\leq&\ \int_0^t \int_{\mathbb{R}^{2d}} |f_1 - f_2| \left( \frac{c}{R_2} I_{\{|v| \leq 2R_2\}} |v| + 1 \right) \leq\ C \int_0^t \|f_1 - f_2\|_{L^1(\mathbb{R}^{2d})}.
\end{align*}

Combining all the above estimates, we conclude that
$$
{\color{black}\lim_{j\to\infty}\lim_{M\to\infty}\lim_{R\to\infty}
\limsup_{\delta\to0}\limsup_{\varepsilon\to0}
I^{\mathrm{trans\text{-}v}}_t}
\leq C \int_0^t \|f_1 - f_2\|_{L^1(\mathbb{R}^{2d})}.
$$

{\bf Conclusion. } Combining all the estimates derived above, we now conclude the proof of uniqueness. Suppose that the initial data coincide, i.e., $f_{1,0} = f_{2,0} =: f_0$ almost surely, for almost every $z \in \mathbb{R}^{2d}$. Then, almost surely, for every $t \in [0, T]$, we obtain the estimate
\begin{align*}
&\|f_1(t) - f_2(t)\|_{L^1(\mathbb{R}^{2d})}\\
 &{\color{black}= \lim_{j\to\infty}\lim_{M\to\infty}\lim_{R\to\infty}
 \lim_{\delta\to0}\lim_{\varepsilon_x\to0}\lim_{\varepsilon_v\to0}
 \int_{\mathbb{R}^{2d+1}} \left( \chi_1^{\varepsilon,\delta}(t) + \chi_2^{\varepsilon,\delta}(t) - 2 \chi_1^{\varepsilon,\delta}(t)\chi_2^{\varepsilon,\delta}(t) \right) \varphi_{\beta_j} \zeta_M \alpha_R \, dz\,d\zeta} \\
{\color{black}\leq}&{\color{black}\ C\int_0^t
 \left(1+2\|V\|_{L^\infty(\mathbb R^d)}\mathsf M^{1/2}
 \|\nabla_v\sqrt{f_2(s)}\|_{L^2(\mathbb R^{2d})}\right)
 \|f_1(s)-f_2(s)\|_{L^1(\mathbb R^{2d})}\,ds.}
\end{align*}
{\color{black}The coefficient on the right-hand side belongs to $L^1(0,T)$ almost surely by \eqref{L2-es}. Applying Gronwall's inequality yields}
\begin{align*}
\sup_{t \in [0, T]} \|f_1(t) - f_2(t)\|_{L^1(\mathbb{R}^{2d})} = 0,
\end{align*}
which implies that almost surely, $f_1 \equiv f_2$ almost everywhere in $[0, T] \times \mathbb{R}^{2d}$. This completes the proof. 
\end{proof}

\section{Kinetic semigroup estimates}\label{sec-4}
In this section, we derive estimates for the kinetic semigroup. 
We begin by introducing the definition of the kinetic semigroup generated by the operator $\Delta_v - v \cdot \nabla_x$. 
This operator consists of a Laplacian acting only in the $v$-variable, together with a transport operator in the $x$-variable with an inhomogeneous coefficient $v$. Because of this inhomogeneity, the transport structure does not commute with the action of $\Delta_v$, which makes the construction of the corresponding semigroup nontrivial. To address this, we introduce a stochastic process whose transition density plays a central role in the definition of the kinetic semigroup.

Let $B(t)$ denote a $d$-dimensional Brownian motion, and define the stochastic processes
$$
(X_t, V_t) := \left(-\sqrt{2} \int_0^t B_s\,ds, \sqrt{2} B_t\right).
$$
The kinetic semigroup $(P_t)_{t \geq 0}$ associated with this process is defined by
\begin{align}\label{k-semigroup}
(P_t f)(x,v) = \mathbb{E} f(x-tv+X_t, v+V_t) = (\Gamma_t p_t) \ast (\Gamma_t f)(x,v),
\end{align}
where $p_t$ is the density of the time marginal law of the process $(X_t,V_t)$, and the translation operator $\Gamma_t$ is given by
\begin{align}\label{t-semigroup}
(\Gamma_t f)(x,v) = f(x-tv, v).
\end{align}
An explicit formula for the transition density $p_t$ is provided in \cite{Ko34}; see also \cite{HZ24} for analogous heat kernel estimates in the case of $\alpha$-stable processes. The density takes the form
{\color{black}
\begin{align}\label{density}
p_t(x,v) = \left(\frac{\sqrt3}{2\pi t^2}\right)^d
 \exp\left\{-\frac{3|x|^2}{t^3}
                  -\frac{3x\cdot v}{t^2}-\frac{|v|^2}{t}\right\}.
\end{align}}
For any smooth function $u_0$ and smooth solution $u$ to the following kinetic PDE:
\begin{align*}
    \p_t u=(\Delta_v-v\cdot \nabla_x)u,\quad u(0)=u_0,
\end{align*}
we apply It\^o's formula to $s\to u(t-s,x-sv+X_s,v+V_s)$ and have:
\begin{align*}
    \mE u(0,x-tv+X_t,v+V_t)=u(t,x,v),
\end{align*}
which is 
\begin{align*}
    u(t)=P_t u_0.
\end{align*}
Moreover, for any test function $\varphi \in C^\infty_b(\mR^{2d})$, applying It\^o's formula to the process $t \mapsto \varphi(X_t, V_t)$ yields
\begin{align*}
    \mE \varphi(X_t, V_t) = \varphi(0) + \mE \int_0^t (\Delta_v - V_s \cdot \nabla_x) \varphi(X_s, V_s) \, \dif s,
\end{align*}
which, upon noting that $\mE \varphi(X_t, V_t) = \int_{\mR^{2d}} \varphi(x, v) p_t(x, v) \, dx dv$, implies that the following Fokker-Planck equation holds:
\begin{align}\label{S4:FPE}
    \p_t p_t = (\Delta_v + v \cdot \nabla_x)p_t, \quad \lim_{t \to 0} p_t = \delta_0 \ \text{in the distributional sense.}
\end{align}
{\color{black} As for the approximation employed in Section~\ref{sec-5}, for any $n\in\mN$, we also consider the elliptic operator
\begin{align}\label{op:elliptic}
    \Delta_v+\frac1n\Delta_x-v\cdot \nabla_x,
\end{align}
which is associated with the process
\begin{align*}
    (\tilde{X}^n_t, \tilde{V}^n_t) := \Bigl(-\sqrt{2} \int_0^t B_s\,ds+\sqrt{\frac2n}\,\tilde{B}_t,\; \sqrt{2} B_t\Bigr)
    =(X_t+\sqrt{\tfrac2n}\,\tilde{B}_t,\; V_t),
\end{align*}
where $\tilde{B}$ is a standard $d$-dimensional Brownian motion independent of $B$. 
In the sequel, we use
\begin{align*}
    H^n_t f(x,v):=(p^n_t * f(\cdot,v))(x):=\mE f\bigl(x+\sqrt{\tfrac2n}\,\tilde{B}_t,\;v\bigr),
\end{align*}
{\color{black}with
\begin{align*}
 p_t^n(x)=\left(\frac{n}{4\pi t}\right)^{d/2}
          e^{-n|x|^2/(4t)}.
\end{align*}}
Note that
\begin{align*}
    \Gamma_t H^n_s=H^n_s \Gamma_t, \quad \forall s,t\in \mR_+,
\end{align*}
which implies
\begin{align}\label{ex:kinetic}
    P_tH^n_s=H^n_s P_t,\quad \forall s,t\in \mR_+.
\end{align}
Consequently, the semigroup associated with \eqref{op:elliptic} can be expressed as
\begin{align*}
     \mE f(x-tv+\tilde{X}^n_t,\,v+\tilde{V}^n_t)
     &=\mE\Bigl( \mE f(x-tv+{X}_t+y,\,v+{V}_t)\big|_{y=\sqrt{\frac2n}\,\tilde{B}_t}\Bigr)\\
     &=H^n_t P_t f=P_t H^n_t f.
\end{align*}
}

We now consider the following model kinetic SPDE:
\begin{align}\label{s4:SPDE}
    \dif u = \big[(\Delta_v{\color{black}+\frac{1}{n}\Delta_x} - v \cdot \nabla_x)u + f\big] \, \dif t + h \nabla_v \cdot \big(g \, \dif W_F(t)\big),
\end{align}
{\color{black}Here $f,h$, and $g$ are scalar progressively measurable fields on
$\Omega\times\mathbb R_+\times\mathbb R^{2d}$.  In particular,
$g\,dW_F=\sum_{k\geq1}gf_k\,dB^k$ is $\mathbb R^d$-valued and its
$v$-divergence is well defined.} We further denote by $(H_t^{n})_{t \in [0,T]}$ the heat semi-group generated by $\frac{1}{n}\Delta_x$.  

The following Duhamel formula is well known and frequently used in the literature. However, as we could not locate a clear proof in existing references, we include one here for the reader's convenience.

\bl[Duhamel's formula]\label{Duhamel}
{\color{black}Assume that $u_0\in L^2(\mathbb R^{2d})$,
$f\in L^1(\Omega\times[0,T];L^2(\mathbb R^{2d}))$, and
\begin{align*}
 \mathbb E\int_0^T\!\int_{\mathbb R^{2d}}
 \left(F_1|h\nabla_vg|^2+F_3|hg|^2\right)\,dz\,dt<\infty.
\end{align*}
Let $u$ be an adapted continuous $L^2$-valued variational solution of
\eqref{s4:SPDE}; equivalently, $u$ is the corresponding
distribution-valued continuous semimartingale.} Let $(P_t)_{t \in [0,T]}$ denote the kinetic semigroup defined by \eqref{k-semigroup} {\color{black}and let $(H^n_t)_{t\in[0,T]}$ be the heat semi-group generated by $\frac{1}{n}\Delta_x$}. {\color{black}Then, for every $t\in[0,T]$, almost surely, the following identity holds in the variational sense (and in $L^2$ whenever the three terms belong to $L^2$):}
\begin{align}\label{eq:Duhamel}
    u_t = P_t {\color{black}H^n_t}u_0 + \int_0^t P_{t-s} {\color{black}H^n_{t-s}}f_s \, \dif s + \int_0^t P_{t-s} {\color{black}H^n_{t-s}}\big[h_s \nabla_v \cdot (g_s \, \dif W_F(s))\big].
\end{align}
\el
\begin{proof}
{\color{black}The calculation is first performed for spatially smooth
approximations.  It then extends to the stated variational solution by the
preceding integrability assumptions and the continuity of deterministic and
stochastic convolution.  In particular, no time differentiability of $u$ is
assumed.}
Set
$$
    \sL_n:=\Delta_v+\frac1n\Delta_x-v\cdot\nabla_x,
    \qquad S_t^n:=P_tH_t^n .
$$
By the construction above,
$$
    S_t^n\phi(x,v)
    =\mE \phi\bigl(x-tv+\tilde X_t^n,\, v+\tilde V_t^n\bigr).
$$
Let
$$
    Z_t^{x,v}:=(X_t^{x,v},V_t^v)
    :=\bigl(x-tv+\tilde X_t^n,\ v+\tilde V_t^n\bigr).
$$
Then for every smooth $\phi$, It\^o's formula gives
\begin{align*}
    \dif \phi(Z_t^{x,v})
    &= (\sL_n\phi)(Z_t^{x,v})\,\dif t
      +\sqrt{\frac2n}\,\nabla_x\phi(Z_t^{x,v})\cdot\dif\tilde B_t
      +\sqrt2\,\nabla_v\phi(Z_t^{x,v})\cdot\dif B_t .
\end{align*}
Taking expectations yields
$$
    \partial_t S_t^n\phi=S_t^n\sL_n\phi .
$$

Fix $t\in(0,T]$ and define, for $0\le s<t$,
$$
    Y_s:=S_{t-s}^n u_s .
$$
Since $s\mapsto S_{t-s}^n$ is deterministic and of finite variation, the product rule gives no quadratic variation term. Using the preceding identity and the equation
$$
    \dif u_s=(\sL_nu_s+f_s)\,\dif s
        +h_s\nabla_v\cdot\bigl(g_s\,\dif W_F(s)\bigr),
$$
we obtain
$$
\begin{aligned}
    \dif Y_s
    &= -S_{t-s}^n\sL_nu_s\,\dif s+S_{t-s}^n\,\dif u_s  \\
    &= -S_{t-s}^n\sL_nu_s\,\dif s
       +S_{t-s}^n(\sL_nu_s+f_s)\,\dif s
       +S_{t-s}^n\bigl[h_s\nabla_v\cdot(g_s\,\dif W_F(s))\bigr] \\
    &= S_{t-s}^n f_s\,\dif s
       +S_{t-s}^n\bigl[h_s\nabla_v\cdot(g_s\,\dif W_F(s))\bigr].
\end{aligned}
$$
Integrating over $[0,t]$ gives
$$
\begin{aligned}
    u_t
    =S_t^n u_0+\int_0^t S_{t-s}^n f_s\,\dif s+\int_0^t
        S_{t-s}^n\bigl[h_s\nabla_v\cdot(g_s\,\dif W_F(s))\bigr].
\end{aligned}
$$
Finally, since $S_r^n=P_rH_r^n$, this is precisely \eqref{eq:Duhamel} and completes the proof.
\end{proof}

We now state the following lemma, which provides the main estimate for the kinetic semigroup and is structurally analogous to \cite[Lemma 3.1]{HWZ20}.

\begin{lemma}\label{lem-kineticestimates}
	{\color{black}For every $\alpha,\beta\geq0$, $m,n\in\mN_0$ and $\ell\ge0$, there is a constant $C=C(T,d,\alpha,\beta,m,n,\ell)>0$ such that for all $t\in(0,T]$ and $j\ge-1$,} 
	\begin{align*}	\int_{\mathbb{R}^{2d}}|x|^{\alpha}|v|^{\beta}|\nabla_x^m\nabla_v^n\mathcal{R}^{\theta}_j\Gamma_tp_t(x,v)|dz\leq C2^{{[3(m-\alpha)+(n-\beta)]}j}(2^{-2j}t^{-1})^\ell.	
	\end{align*}

\end{lemma}
\begin{proof}
    For notational convenience, we define $q_t := \Gamma_t p_t$ and $\phi_j := \check{\phi}_j^\theta$ for all $t \in [0,T]$ and $j \geq -1$. By definition, we then have
\begin{align*}
    I_j(t) &:= \int_{\mathbb{R}^{2d}} |x|^{\alpha} |v|^{\beta} \left| \nabla_x^m \nabla_v^n \mathcal{R}^\theta_j \Gamma_t p_t(x,v) \right| \, dz \\
    &= \int_{\mathbb{R}^{2d}} |x|^{\alpha} |v|^{\beta} \left| \int_{\mathbb{R}^{2d}} \nabla_x^m \nabla_v^n \phi_j(x-y, v-w) q_t(y,w) \, dy \, dw \right| dx \, dv,
\end{align*}
where, for $j \geq 0$ and any $\ell \in \mathbb{N}$,
\begin{align*}
    &\quad \int_{\mathbb{R}^{2d}} \nabla_x^m \nabla_v^n \phi_j(x-y, v-w) q_t(y,w) \, dy \, dw \\
    &= 2^{3mj + nj} 2^{4dj} \int_{\mathbb{R}^{2d}} \left(\nabla_x^m \nabla_v^n \phi_0\right)(2^{3j}(x-y), 2^j(v-w)) q_t(y,w) \, dy \, dw.
\end{align*}
By applying a change of variables, we deduce that for any $j \geq 0$,
\begin{align}
    I_j(t)
    &= 2^{[3(m - \alpha) + (n - \beta)]j} \int_{\mathbb{R}^{2d}} |x|^{\alpha} |v|^{\beta} \left| \int_{\mathbb{R}^{2d}} \left(\nabla_x^m \nabla_v^n \phi_0\right)(x - 2^{3j}y, v - 2^j w) q_t(y,w) \, dy \, dw \right| dx \, dv \notag \\
    &= 2^{[3(m - \alpha) + (n - \beta)]j} 2^{-4dj} \int_{\mathbb{R}^{2d}} |x|^{\alpha} |v|^{\beta} \left| \int_{\mathbb{R}^{2d}} \left(\nabla_x^m \nabla_v^n \phi_0\right)(x - y, v - w) q_t(2^{-3j} y, 2^{-j} w) \, dy \, dw \right| dx \, dv \notag \\
    &=: 2^{[3(m - \alpha) + (n - \beta)]j} J_j(t). \label{0219:00}
\end{align}
{\color{black}Fix an integer $L\geq1$, to be chosen below. We observe that
\begin{align*}
    &\quad \int_{\mathbb{R}^{2d}} \left(\nabla_x^m \nabla_v^n \phi_0\right)(x - y, v - w) q_t(2^{-3j} y, 2^{-j} w) \, dy \, dw \\
    &= \int_{\mathbb{R}^{2d}} \left(\nabla_x^m \nabla_v^n \phi_0\right)(x - y, v - w) (-\Delta_x - \Delta_v)^{-L} (-\Delta_x - \Delta_v)^{L} q_t(2^{-3j} \cdot, 2^{-j} \cdot)(y,w) \, dy \, dw \\
    &= \int_{\mathbb{R}^{2d}} \left((-\Delta_x - \Delta_v)^{-L} \nabla_x^m \nabla_v^n \phi_0\right)(x - y, v - w) (-\Delta_x - \Delta_v)^{L} q_t(2^{-3j} \cdot, 2^{-j} \cdot)(y,w) \, dy \, dw,
\end{align*}
where $(-\Delta_x-\Delta_v)^{-L}$ is the multiplier operator defined by
\begin{align*}
    \sF\left((-\Delta_x-\Delta_v)^{-L}\phi_0\right)(\xi,\eta):=(|\xi|^2+|\eta|^2)^{-L}\hat{\phi_0}(\xi,\eta) \in \sS(\mR^{2d}),
\end{align*}
since the support of $\hat{\phi_0}$ is contained in the annulus.
Define $\psi := (-\Delta_x-\Delta_v)^{-L} \nabla_x^m \nabla_v^n \phi_0 \in \mathscr{S}(\mathbb{R}^{2d})$ and denote $\Delta := \Delta_x + \Delta_v$. Then it follows that
\begin{align*}
    J_j(t)
    &= 2^{-4dj} \int_{\mathbb{R}^{2d}} |x|^{\alpha} |v|^{\beta} \left| \int_{\mathbb{R}^{2d}} \psi(x - y, v - w) \Delta^L q_t(2^{-3j} \cdot, 2^{-j} \cdot)(y,w) \, dy \, dw \right| dx \, dv \\
    &\leq 2^{-4dj} \int_{\mathbb{R}^{2d}} \int_{\mathbb{R}^{2d}} (1 + |x - y|)^{\alpha} (1 + |v - w|)^{\beta} (1 + |y|)^{\alpha} (1 + |w|)^{\beta} \\
    &\qquad\qquad\qquad\quad \times |\psi(x - y, v - w)| \left| \Delta^L q_t(2^{-3j} \cdot, 2^{-j} \cdot)(y,w) \right| \, dy \, dw \, dx \, dv,
\end{align*}
where we used the inequality $|x| \leq (1 + |x - y|)(1 + |y|)$ (and similarly for $|v|$). Define
$$
\Phi(x,v) := (1 + |x|)^{\alpha} (1 + |v|)^{\beta} |\psi(x,v)|.
$$
Then by convolutional Young's inequality, we obtain 
\begin{align*}
    J_j(t)
    &\leq 2^{-4dj} \left( \int_{\mathbb{R}^{2d}} \Phi(x,v) \, dx \, dv \right) \left( \int_{\mathbb{R}^{2d}} (1 + |y|)^{\alpha} (1 + |w|)^{\beta} \left| \Delta^L q_t(2^{-3j} \cdot, 2^{-j} \cdot)(y,w) \right| \, dy \, dw \right) \\
    &\lesssim 2^{-4dj} \int_{\mathbb{R}^{2d}} (1 + |y|)^{\alpha} (1 + |w|)^{\beta} \left| \Delta^L q_t(2^{-3j} \cdot, 2^{-j} \cdot)(y,w) \right| \, dy \, dw. 
\end{align*}
We observe that the scaling property of $q_t$ reads
$$
q_t(x,v) = t^{-2d} q_1(t^{-3/2} x, t^{-1/2} v).
$$
Set $ h := (t^{1/2} 2^j)^{-1} $, then we compute
\begin{align*}
    \Delta^L q_t(2^{-3j} \cdot, 2^{-j} \cdot)(y,w)
    &= t^{-2d}\sum_{k=0}^L\binom{L}{k}
       h^{6k+2(L-k)}
       (\Delta_x^k\Delta_v^{L-k}q_1)(h^3y,hw).
\end{align*}
Inserting this into the expression for $J_j(t)$, we obtain
\begin{align*}
    J_j(t)
    &\lesssim \sum_{k=0}^L h^{2L+4k}
       (1\vee h^{-3\alpha})(1\vee h^{-\beta})
       \int_{\mathbb R^{2d}}(1+|y|)^\alpha(1+|w|)^\beta
       |\Delta_x^k\Delta_v^{L-k}q_1(y,w)|\,dy\,dw \\
    &\lesssim \sum_{k=0}^L h^{2L+4k}
       (1\vee h^{-3\alpha})(1\vee h^{-\beta}).
\end{align*}
If $h\leq1$, choose $L$ so large that
$2L-3\alpha-\beta\geq2\ell$. Then
\begin{align}\label{0219:01}
 J_j(t)\lesssim h^{2L-3\alpha-\beta}\leq h^{2\ell}.
\end{align}
If $h\geq1$, the same argument with $L=0$ gives
$J_j(t)\lesssim1\leq h^{2\ell}$. Since
$h^{2\ell}=(2^{-2j}t^{-1})^\ell$, combining this with
\eqref{0219:00} proves the asserted estimate for every $j\geq0$.}

\medskip

\noindent
Finally, for the case $j = -1$, we appeal to Young's inequality. Noting that $\phi_{-1} \in \mathscr{S}(\mathbb{R}^{2d})$ and $q_t \in \mathscr{S}'(\mathbb{R}^{2d})$ with sufficient moment decay, we estimate
\begin{align*}
    I_{-1}(t)
    &\leq \int_{\mathbb{R}^{2d}} (1 + |x|)^\alpha (1 + |v|)^\beta |\nabla_x^m \nabla_v^n \phi_{-1}(x,v)| \, dx \, dv \\
    &\qquad\qquad \times \int_{\mathbb{R}^{2d}} (1 + |y|)^\alpha (1 + |w|)^\beta |q_t(y,w)| \, dy \, dw \\
    &\leq C(\alpha,\beta,m,n,T) \leq {\color{black}C(\alpha,\beta,m,n,T) T^\ell (2^2 t^{-1})^\ell}, \quad {\color{black}\forall t \in (0,T]},
\end{align*}
where we used the fact that
\begin{align*}
   \int_{\mathbb{R}^{2d}} (1 + |y|)^\alpha (1 + |w|)^\beta |q_t(y,w)| \, dy \, dw &=\int_{\mathbb{R}^{2d}} (1 + |t^{3/2}y|)^\alpha (1 + |t^{1/2}w|)^\beta |q_1(y,w)| \, dy \, dw \\
   &\lesssim_T\int_{\mathbb{R}^{2d}} (1 + |y|)^\alpha (1 + |w|)^\beta |q_1(y,w)| \, dy \, dw.
\end{align*}
Therefore, we complete the proof.\end{proof}

At the end of this section, we introduce the following result,  which is from \cite[Lemma 6.7]{HWZ20}.
\bl\label{s4:lem-01} 
For $t\geq 0$ and $j\in\mN$, define
\begin{align*}
\Theta^{t}_j:=\Big\{\ell\in\mN_0: 2^{\ell}\leq 2^4 (2^j+t2^{3j}),\ 
2^{j}\leq 2^4 (2^\ell+t2^{3\ell})\Big\}.
\end{align*}
\begin{enumerate}
\item Let $t\ge0$ and $j\in\mN$. For any $\ell\notin\Theta^{t}_j$, it holds that
\begin{align*}
\<\cR^\theta_jf,\Gamma_{t}\cR^\theta_\ell g\>=\int_{\mR^{2d}}\cR^\theta_jf(x,v)\cdot\Gamma_{t}
\cR^\theta_\ell g(x,v) dxdv=0.
\end{align*}
\item For any $\beta>0$, there is a constant $C=C(\beta)>0$ such that for all $j\in\mN$ and $t\geq 0$,
\begin{align*}
\sum_{\ell\in\Theta^t_j}2^{-\beta\ell}\leq C\Big(2^{-j}+t2^{j}\Big)^{\beta},\ \ \sum_{\ell\in\Theta^t_j}2^{\beta\ell}\leq 
C\Big(2^j+t2^{3j}\Big)^\beta.
\end{align*}
\end{enumerate}
\el
\begin{remark}
We remark that although the operators $\Gamma_t$ and $\mathcal{R}_j^\theta$ do not commute, i.e., $\Gamma_t \mathcal{R}_j^\theta \ne \mathcal{R}_j^\theta \Gamma_t$, it is therefore not straightforward to approximate the expression $\langle \mathcal{R}_j^\theta f, \Gamma_t g \rangle$ by $\langle \mathcal{R}_j^\theta f, \Gamma_t \mathcal{R}_j^\theta g \rangle$. Nevertheless, as justified by Lemma~\ref{s4:lem-01}~(1), we observe that
\begin{equation}\label{Thetat-property}
\langle \mathcal{R}_j^\theta f, \Gamma_t g \rangle = \sum_{l \in \Theta^t_j} \langle \mathcal{R}_j^\theta f, \Gamma_t \mathcal{R}_l^\theta g \rangle.
\end{equation}
Here, with the help of Lemma~\ref{s4:lem-01}~(2), the summation over $\Theta^t_j$ can be effectively approximated by the single term $\langle \mathcal{R}_j^\theta f, \Gamma_t \mathcal{R}_j^\theta g \rangle$. Although $\Gamma_t$ and $\mathcal{R}_j^\theta$ do not commute, Lemma~\ref{s4:lem-01} provides a framework that allows us to derive meaningful estimates for the kinetic semigroup. Further estimates have been developed in~\cite{HWZ20}, and we will make use of them in Section~\ref{sec-5}.
\end{remark}

\section{The existence of the regularized SPDE}\label{sec-5}

In this section, we investigate the existence of a regularized equation   
\begin{align}\label{SPDE-6}
df=\Big(\Delta_vf-v\cdot\nabla_xf&-\nabla_v\cdot(fV\star_x\rho)+\nabla_v\cdot(vf)\Big)dt-\nabla_v\cdot(\sigma(f)dW_F)\notag\\
 &{\color{black}+\frac{1}{2}\nabla_v\cdot\bigl(F_1\sigma'(f)^2\nabla_vf\bigr)dt}.    
\end{align}
The main result can be stated as follows. 
\begin{theorem}\label{thm-existence}
	Assume that $f_0,V$ and $\sigma(\cdot)$ satisfy Assumptions \ref{Assump-initialdata-regular}, \ref{Assump-ker} and \ref{Assump-sigma-n}, respectively. Then \eqref{SPDE-6} admits a unique weak solution. 
\end{theorem}
{\color{black}Throughout this section, we use the data size
\begin{align}\label{regular-data-size}
 M_*&:=\operatorname*{ess\,sup}_{\omega\in\Omega}
   \|f_0(\omega)\|_{L^1(\mathbb R^{2d})},\notag\\
 \mathcal K_{\rm reg}&:=
 1+\|f_0\|_{L^2(\Omega;L^2(\mathbb R^{2d}))}
 +M_*.
\end{align}
Constants may also depend on the fixed quantities
$T,d,V,\sigma,F_1,F_3$ and, in localized estimates, on the cutoff.}
{\color{black}We organize the proof of Theorem~\ref{thm-existence} into four steps.}
\begin{figure}[htbp]
\centering
\resizebox{\textwidth}{!}{%
\begin{tikzpicture}[
    >=Latex,
    node distance=5mm,
    every node/.style={font=\scriptsize},
    arrow/.style={-{Latex[length=2mm]}, thick},
    box/.style={
        rectangle,
        rounded corners=2pt,
        draw,
        align=center,
        inner sep=3pt,
        text width=30mm,
        minimum height=9mm
    },
    sidebox/.style={box, fill=gray!8, draw=gray!60!black},
    bluebox/.style={box, fill=blue!8, draw=blue!55!black},
    orangebox/.style={box, fill=orange!10, draw=orange!65!black},
    greenbox/.style={box, fill=green!8, draw=green!50!black},
    purplebox/.style={box, fill=purple!8, draw=purple!55!black}
]

\node[sidebox] (ass) at (0,0)
    {Assumptions\\ \ref{Assump-initialdata-regular}, \ref{Assump-ker}, \ref{Assump-sigma-n}\\regular $f_0$};

\node[sidebox, below=8mm of ass] (goal)
    {Goal\\weak well-posedness\\of \eqref{SPDE-6}};

\node[bluebox] (s1a) at (3.7,0)
    {Step 1\\linear model \eqref{lkSPDE}};

\node[bluebox, below=of s1a] (s1b)
    {frozen\\characteristic shift\\$\tau_t^{z_0}$};

\node[bluebox, below=of s1b] (s1c)
    {Duhamel formula\\kinetic semigroup\\$P_t$};

\node[bluebox, below=of s1c] (s1d)
    {Lemma \ref{lem-properties}\\shift estimates};

\node[bluebox, below=of s1d] (s1e)
    {Theorem \ref{thm:Besov}\\anisotropic\\Besov estimate};

\node[orangebox] (s2a) at (7.4,0)
    {Step 2\\approximation \eqref{SPDE-iteration}};

\node[orangebox, below=of s2a] (s2b)
    {parameters $n,r$\\$\frac1n\Delta_x$ and\\$\alpha_r^2(v)v$};

\node[orangebox, below=of s2b] (s2c)
    {Galerkin\\construction};

\node[orangebox, below=of s2c] (s2d)
    {Proposition~\ref{prp-uniforml2}\\$L^2$ energy};

\node[orangebox, below=of s2d] (s2e)
    {Lemma \ref{lem-L1-es}\\Prop. \ref{moment-es}\\Lemma \ref{time-besov-regularity-fn}\\positivity, mass,\\moments\\time regularity};

\node[greenbox] (r1) at (11.1,0)
    {Step 3\\Fix $n$, pass $r\to\infty$};

\node[greenbox, below=of r1] (r2)
    {Theorem \ref{tightness}\\tightness};

\node[greenbox, below=of r2] (r3)
    {Skorokhod / Krylov\\diagonal argument};

\node[greenbox, below=of r3] (r4)
    {identify nonlinear\\term and martingale\\term};

\node[greenbox, below=of r4] (r5)
    {Theorem~\ref{thm:nonb}\\viscous solution\\$f_n$ of \eqref{SPDE-L1}};

\node[greenbox] (n1) at (14.8,0)
    {Step 4\\pass $n\to\infty$};

\node[greenbox, below=of n1] (n2)
    {Theorem~\ref{besov-regularity-fn}\\Besov bound\\uniform in $n$\\(uses Thm. \ref{thm:Besov})};

\node[greenbox, below=of n2] (n3)
    {Lemma \ref{App:cpt}\\compact embedding};

\node[greenbox, below=of n3] (n4)
    {Corollary~\ref{co-tightness}\\tightness of\\$(f_n,\nabla_v f_n)$};

\node[greenbox, below=of n4] (n5)
    {$\frac1n\Delta_x f_n\to0$\\limit solves\\\eqref{SPDE-6}};

\node[greenbox, below=of n5] (n6)
    {Existence\\weak solution\\of \eqref{SPDE-6}};

\node[purplebox] (u1) at (18.5,-3.3)
    {Proposition~\ref{equivalence}\\weak $\Rightarrow$\\renormalized kinetic};

\node[purplebox, below=of u1] (u2)
    {Theorem~\ref{Uniqueness-spde}\\$L^1$ contraction\\uniqueness};

\node[purplebox, below=of u2] (u3)
    {Theorem~\ref{thm-existence}\\unique weak\\solution};

\draw[arrow] (s1a.south) -- (s1b.north);
\draw[arrow] (s1b.south) -- (s1c.north);
\draw[arrow] (s1c.south) -- (s1d.north);
\draw[arrow] (s1d.south) -- (s1e.north);

\draw[arrow] (s2a.south) -- (s2b.north);
\draw[arrow] (s2b.south) -- (s2c.north);
\draw[arrow] (s2c.south) -- (s2d.north);
\draw[arrow] (s2d.south) -- (s2e.north);

\draw[arrow] (r1.south) -- (r2.north);
\draw[arrow] (r2.south) -- (r3.north);
\draw[arrow] (r3.south) -- (r4.north);
\draw[arrow] (r4.south) -- (r5.north);

\draw[arrow] (n1.south) -- (n2.north);
\draw[arrow] (n2.south) -- (n3.north);
\draw[arrow] (n3.south) -- (n4.north);
\draw[arrow] (n4.south) -- (n5.north);
\draw[arrow] (n5.south) -- (n6.north);

\draw[arrow] (u1.south) -- (u2.north);
\draw[arrow] (u2.south) -- (u3.north);

\draw[arrow] (ass.east) -- (s1a.west);
\draw[arrow] (goal.east) -- (s1b.west);
\draw[arrow] (s1a.east) -- (s2a.west);
\draw[arrow] (s2a.east) -- (r1.west);
\draw[arrow] (r1.east) -- (n1.west);
\draw[arrow] (n6.east) -- (u2.west);

\end{tikzpicture}%
}
\caption{Proof diagram for the weak well-posedness of \eqref{SPDE-6}.}
\label{fig:proof-weak-wellposedness}
\end{figure}
\begin{itemize}
  \item \textbf{Subsection~\ref{subsec-5-1}.} We study the model kinetic SPDE \eqref{lkSPDE}, from which we derive uniform Besov regularity estimates. The frozen characteristic line method plays a key role.

  \item \textbf{Subsection~\ref{subsec-5-2}.} We introduce an approximation scheme (see \eqref{SPDE-iteration}). We then prove the existence of solutions to this scheme and establish non-negativity and moment estimates.

  \item \textbf{Subsection~\ref{subsec-5-3}.} The approximation scheme incorporates a velocity truncation, associated with a parameter $r \geq 1$, and a spatial regularization term $\frac{1}{n}\Delta_x$. For fixed $n \geq 1$, we first pass to the limit as $r \to \infty$ to remove the truncation. In this step, the regularization term $\frac{1}{n}\Delta_x$ provides strong compactness.

  \item \textbf{Subsection~\ref{subsec-5-4}.} In this subsection, we exploit the Besov regularity estimates from Subsection~\ref{subsec-5-1} to transfer regularity from $v$ to $x$ and establish strong compactness. Finally, we complete the proof of Theorem~\ref{thm-existence} by letting $n \to \infty$.
\end{itemize}

\subsection{The model equation and Besov regularity estimates}\label{subsec-5-1}
In this section, we begin by analyzing the following linear kinetic SPDE:
\begin{align}\label{lkSPDE}
    \dif u=\left[(\Delta_v{\color{black}+\frac{1}{n}\Delta_x} - v\cdot\nabla_x)u + \nabla_v\cdot(vu) - \nabla_v\cdot(bu) + [g_1 + g_2]\right]\dif t + h_1 \dif W_F+h_2 \nabla_v \cdot \dif W_F,
\end{align}
where $b=(b_1,...,b_d)$, $h_1=(h_{1,1},...,h_{1,d})$, $h_2$, $g_2$ are function-valued progressively measurable processes, while $g_1$ is a distribution-valued progressively measurable process. To guarantee the well-definedness of the stochastic integration, we impose the assumption 
$$
\int_0^T |h_1(t, z)|^2+|h_2(t, z)|^2 \, \dif t < \infty,
$$ almost surely for all $z \in \mR^{2d}$. Further assumptions on the regularity and integrability properties of the coefficients $b$, $g_1$, $g_2$, and $h_1$, $h_2$ will be specified below. We now state the main results concerning equation \eqref{lkSPDE}. 
 \begin{theorem}\label{thm:Besov}
{\color{black}Let $b\in L^\infty(\Omega\times[0,T]\times
\mathbb R^{2d})$.  The following two estimates hold, with constants
independent of the artificial viscosity parameter $n$.

\smallskip
\noindent\emph{(H) The Hilbert case.}  Suppose that $p=q=2$ and
\begin{align*}
 \mathcal K_H:=&\ \|u_0\|_{L^2(\Omega;L^2)}
 +\|g_1\|_{L^2(\Omega\times[0,T];\mathbf B^{-1}_{2;\theta})}
 +\|g_2\|_{L^2(\Omega\times[0,T]\times\mathbb R^{2d})}\\
 &+\|h_1\|_{L^2(\Omega\times[0,T];L^2)}
 +\|h_2\|_{L^2(\Omega\times[0,T];L^2)}<\infty.
\end{align*}
If $u$ is a mild solution of \eqref{lkSPDE} in
$L^2(\Omega\times[0,T]\times\mathbb R^{2d})$, then, for every
$\beta\in(0,1)$,
\begin{align}\label{apriori-besov-es}
 \|u\|_{L^2(\Omega\times[0,T];\mathbf B^\beta_{2;\theta})}
 \leq C\mathcal K_H.                                      \tag{H}
\end{align}

\smallskip
\noindent\emph{(S) The first-moment smoothing case.}  Let
\begin{align*}
 1<p<\frac{2d+1}{2d},\qquad q=1,\qquad
 \beta_0:=\frac{4d+2}{p}-4d>0,
\end{align*}
and assume
\begin{align*}
 \mathcal K_{S,p}:={}&\ \mathbb E\|u_0\|_{L^1}
 +\|g_1\|_{L^2(\Omega\times[0,T];\mathbf B^{-1}_{p;\theta})}
 +\mathbb E\|g_2\|_{L^1([0,T]\times\mathbb R^{2d})}\\
 &+\|h_1\|_{L^2(\Omega\times[0,T];L^p)}
 +\|h_2\|_{L^2(\Omega\times[0,T];L^p)}<\infty.
\end{align*}
If $u$ is a mild solution of \eqref{lkSPDE} such that
$\mathbb E\|u\|_{L^p([0,T]\times\mathbb R^{2d})}<\infty$, then, for
every $\beta\in(0,\min\{\beta_0,1\})$,
\begin{align}\label{apriori-besov-es-first-moment}
 \mathbb E\|u\|_{L^p([0,T];\mathbf B^\beta_{p;\theta})}
 \leq C\mathcal K_{S,p}.                                  \tag{S}
\end{align}
In both cases $C$ may depend on
$d,T,p,\beta,\|b\|_\infty,\|F_1\|_\infty$ and
$\|F_3\|_\infty$, but not on $n$.  The outer probability exponent is
therefore never improved by the kinetic semigroup: it is $2$ in (H) and
$1$ in (S).}
\end{theorem}

We emphasize that the results of Theorem \ref{thm:Besov} are formulated as estimates in Besov spaces. However, the lack of integrability of the inhomogeneous coefficient $v$ in the term $\nabla_v \cdot (v f)$ poses a significant challenge for obtaining such estimates. To overcome this difficulty, we introduce a novel argument based on the method of frozen characteristic lines. 

As outlined in the introduction, our strategy for handling the term 
$\nabla_v \cdot (v f)$ is to shift the solution by introducing a new variable $z_0=(x_0,v_0)$. After applying Duhamel's formula (Lemma \ref{Duhamel}) and evaluating the solution at $z=(x,v)=(0,0)$, we estimate the Besov norm with respect to $z_0$. In Lemma~\ref{lem-properties}, we establish several properties of this shifting operator and clarify the relationship between the Besov norm of the original solution (in the original variables) and that of the shifted solution (in the variable $z_0$). These properties allow us to recover the Besov regularity of the original solution. To apply Duhamel's formula and take advantage of the available estimates for the kinetic semigroup, however, the shift direction $z_0$ cannot be chosen arbitrarily; it must preserve the kinetic structure of the operator in the original equation. For this reason, we perform the shift along the characteristic curves. Precisely, for each $z_0 = (x_0, v_0) \in \mathbb{R}^{2d}$, consider the flow $\theta^{z_0}_t = (\theta^{z_0}_t(1), \theta^{z_0}_t(2))$ governed by the system
\begin{align*}
    \dot{\theta}^{z_0}_t(1) = \theta^{z_0}_t(2), \quad
    \dot{\theta}^{z_0}_t(2) = - \theta^{z_0}_t(2), \quad
    \theta^{z_0}_0 = z_0.
\end{align*}
A straightforward computation shows that
\begin{align}\label{flow}
    \theta^{z_0}_t = \left(x_0 + (1 - e^{-t}) v_0,\ e^{-t} v_0\right).
\end{align}

For each $t \in [0,T]$, define the shift operator $\tau^{z_0}_t$ acting on functions $h$ by
\begin{align}\label{shiftoperator}
    \tau^{z_0}_t h(x, v) := h\left(x + \theta^{z_0}_t(1),\ v + \theta^{z_0}_t(2)\right).
\end{align}
By applying the chain rule, it follows that
\begin{align*}
&\partial_t\tau^{z_0}_tu(t,x,v)=\tau^{z_0}_t	(\partial_tu)+\tau^{z_0}_t(\nabla_xu)\cdot\dot{\theta}^{z_0}_t(1)+\tau^{z_0}_t(\nabla_vu)\cdot\dot{\theta}^{z_0}_t(2)\\
=&\tau^{z_0}_t\Big[(\Delta_v{\color{black}+\frac{1}{n}\Delta_x}-v\cdot\nabla_x)u-\nabla_v\cdot(bu)+\nabla_v\cdot(vu)+g_1+g_2+h_1\dot{W}_F+h_2\nabla_v\cdot\dot{W}_F\Big]\\
&+\tau^{z_0}_t(\nabla_xu)\cdot\theta^{z_0}_t(2)-\tau^{z_0}_t(\nabla_vu)\cdot\theta^{z_0}_t(2)\\
=&(\Delta_v{\color{black}+\frac{1}{n}\Delta_x}-v\cdot\nabla_x)\tau^{z_0}_tu-\nabla_v\cdot(\tau^{z_0}_tu\tau^{z_0}_tb)+\nabla_v\cdot(v\tau^{z_0}_tu)+\tau^{z_0}_t[g_1+g_2]+\tau^{z_0}_t(h_1\dot{W}_F+h_2\nabla_v\cdot\dot{W}_F)\\
	&+\tau^{z_0}_t(\nabla_xu)\cdot\theta^{z_0}_t(2)-\tau^{z_0}_t(\nabla_vu)\cdot\theta^{z_0}_t(2)-\tau^{z_0}_t(\nabla_xu)\cdot\theta^{z_0}_t(2)+\tau^{z_0}_t(\nabla_vu)\cdot\theta^{z_0}_t(2)\\
=&(\Delta_v{\color{black}+\frac{1}{n}\Delta_x}-v\cdot\nabla_x)\tau^{z_0}_tu-\nabla_v\cdot(\tau^{z_0}_tu\tau^{z_0}_tb)+\nabla_v\cdot(v\tau^{z_0}_tu)+\tau^{z_0}_t[g_1+g_2]+\tau^{z_0}_t (h_1\dot{W}_F+h_2\nabla_v\cdot\dot{W}_F). 
\end{align*}
One then observes that, owing to the choice of the characteristic curve, the additional terms arising from the chain rule applied to the time derivative cancel precisely with the additional terms generated by the action of the shift operator on the terms involving the inhomogeneous coefficient $v$. For the reader's convenience, we summarize the resulting evolution equation for $\tau^{z_0}_t u$ as follows:
\begin{align}\label{SPDE-7}
\begin{split}
    \dif(\tau^{z_0}_t u) = \Big[&(\Delta_v{\color{black}+\frac{1}{n}\Delta_x} - v \cdot \nabla_x)\tau^{z_0}_t u 
    - \nabla_v \cdot (\tau^{z_0}_t u \cdot \tau^{z_0}_t b)
    + \nabla_v \cdot (v \tau^{z_0}_t u) \\
    &+ \tau^{z_0}_t [g_1 + g_2] \Big] \dif t 
    + \tau^{z_0}_t [h_1\dif W_F+h_2\nabla_v\cdot \dif W_F].
\end{split}
\end{align}

For each $j \ge -1$, recall that $\mathcal{R}^{\theta}_j$ denotes the anisotropic block operator defined in \eqref{block-operator} and $(P_t)_{t\in[0,T]}$ denotes the kinetic semigroup defined in \eqref{k-semigroup},  {\color{black} and $(H^n_t)_{t\in[0,T]}$ denotes the heat semi-group generated by $\frac{1}{n}\Delta_x$. We then apply the Duhamel formula with respect to the kinetic semigroup (see Lemma \ref{Duhamel}), which yields}
\begin{align*}
    \mathcal{R}^{\theta}_j \tau^{z_0}_t u(t)
    =&\ {\color{black}\mathcal{R}^{\theta}_j P_t H^n_{t}\tau^{z_0}_0 u_0} 
    - \int_0^t \mathcal{R}^{\theta}_j P_{t-s}{\color{black}H^n_{t-s}} \nabla_v \cdot (\tau^{z_0}_s u \cdot \tau^{z_0}_s b) \, \dif s 
    + \int_0^t \mathcal{R}^{\theta}_j P_{t-s}{\color{black}H^n_{t-s}} \nabla_v \cdot (v \tau^{z_0}_s u) \, \dif s \\
    &+ \int_0^t \mathcal{R}^{\theta}_j P_{t-s}{\color{black}H^n_{t-s}} \left[ \tau^{z_0}_s (h_1\dif W_F+h_2\nabla_v\cdot \dif W_F) \right] 
    + \int_0^t \mathcal{R}^{\theta}_j P_{t-s}{\color{black}H^n_{t-s}} \tau^{z_0}_s [g_1 + g_2] \, \dif s.
\end{align*}

Noting that $\mathcal{R}^{\theta}_j \tau^{z_0}_t = \tau^{z_0}_t \mathcal{R}^{\theta}_j$, and evaluating at $z = (x, v) = 0$, we obtain
{\color{black}
\begin{align}\label{tauR}
	\tau^{z_0}_t\mathcal{R}^{\theta}_ju(t,0)=&(\mathcal{R}^{\theta}_jP_tH^n_t\tau^{z_0}_0u_0)(0)-\int^t_0\mathcal{R}^{\theta}_jP_{t-s}H^n_{t-s}\nabla_v\cdot(\tau^{z_0}_su\tau^{z_0}_sb)(0)ds+\int^t_0\mathcal{R}^{\theta}_jP_{t-s}H^n_{t-s}\nabla_v\cdot(v\tau^{z_0}_su)(0)ds\notag\\
	&+\int^t_0\mathcal{R}^{\theta}_jP_{t-s}H^n_{t-s}\left[ \tau^{z_0}_s (h_1\dif W_F+h_2\nabla_v\cdot \dif W_F) \right] (0)+\int^t_0\mathcal{R}^{\theta}_jP_{t-s}H^n_{t-s}\tau^{z_0}_s[g_1+g_2]\dif s\notag\\
    :=&(\mathcal{R}^{\theta}_jP_tH^n_t\tau^{z_0}_0u_0)(0)+Z_{1,j}+Z_{2,j}+Z_{3,j}+Z_{4,j}.
\end{align}}

We next state several properties of the shift operator $\tau^{z_0}_t$ that will be used in our analysis.

\begin{lemma}\label{lem-properties}
Let $0 \le s \le t \le T$ and recall that $\Gamma_t$ denotes the transport operator defined in \eqref{t-semigroup}. Let $f:\mathbb{R}^{2d}\rightarrow\mathbb{R}$ be a suitable regular function, the following properties hold for each $j\ge -1$ :

\begin{enumerate}
  \item[(i)] 
  The identity 
  $$
    \mathcal{R}^{\theta}_jP_t\bigl(\nabla_v\!\cdot f\bigr)(0)
    = -\int_{\mathbb{R}^{2d}}
      \Bigl(\nabla_v\mathcal{R}^{\theta}_j\Gamma_t p_t 
      + t\,\nabla_x\mathcal{R}^{\theta}_j\Gamma_t p_t\Bigr)\bigl(x + t v,\,v\bigr)\,
      f(x,v)\,\dif z.
  $$

  \item[(ii)] 
  For every $p\ge1$,
  $$
    \sup_{(x,v)\in\mathbb{R}^{2d}}
      \bigl\|\tau_s^{\cdot} f(x,v)\bigr\|_{L^p(\mathbb{R}^{2d})}
    \;\le\;
    e^{\tfrac{s d}{p}}\,\|f\|_{L^p(\mathbb{R}^{2d})}.
  $$

  \item[(iii)] 
  {\color{black}For every $m\ge0$, $\beta\in\mathbb{R}$, and $p\ge q\ge1$
  such that
  $\beta+4d/p-4d/q\neq0$, there exists}
  {\color{black}$C=C(T,d,m,\beta,p,q)>0$} such that
  \begin{align}\label{0219:02}
    \bigl\|\mathcal{R}^{\theta}_jP_t(\tau^{\cdot}_s f)(0)\bigr\|_{L^p(\mathbb{R}^{2d})}
    &\le 
    C\,2^{-\bigl(\beta+\frac{4d}{p}-\frac{4d}{q}\bigr)j}
      \bigl(t^{-1}2^{-2j}\bigr)^m
      e^{\tfrac{s d}{p}}\,
      \|f\|_{\mathbf{B}^{\beta}_{q;\theta}}.
  \end{align}

  {\color{black}\item[(iv)] For every $m\geq0$ and $p\in[1,\infty)$,
  there is $C=C(T,d,m,p)>0$ such that, for $j\geq-1$,
  \begin{align}\label{critical-kernel-shift-estimate}
   \bigl\|\mathcal R_j^\theta P_t(\tau_s^{\cdot}f)(0)
   \bigr\|_{L^p(\mathbb R^{2d})}
   \leq C e^{sd/p}
   \left[1\wedge(2^{-2j}t^{-1})^m\right]
   \|f\|_{L^p(\mathbb R^{2d})}.
  \end{align}
  This is the critical estimate used below when $p=q$ and $\beta=0$;
  it does not assert the generally false endpoint bound with
  $\mathbf B^0_{p;\theta}$ in place of $L^p$.}
\end{enumerate}
\end{lemma}
\begin{proof}
	(i) By using the definition of the kinetic semigroup and the change variables formula, 
	\begin{align*}
		\mathcal{R}^{\theta}_jP_t(\nabla_v\cdot f)(0)=&\mathcal{R}^{\theta}_j(\Gamma_tp_t)\ast(\Gamma_t\nabla_v\cdot f)(0)\\
		=&\int_{\mathbb{R}^{2d}}\mathcal{R}^{\theta}_j(\Gamma_tp_t)(-x,-v)(\Gamma_t\nabla_v\cdot f)(x,v)dz\\
		=&\int_{\mathbb{R}^{2d}}\mathcal{R}^{\theta}_j(\Gamma_tp_t)(-x,-v)(\nabla_v\cdot f)(x-tv,v)dz\\
		=&\int_{\mathbb{R}^{2d}}\mathcal{R}^{\theta}_j(\Gamma_tp_t)(-x-tv,-v)\nabla_v\cdot f(x,v)dz\\
		=&-\int_{\mathbb{R}^{2d}}\nabla_v[\mathcal{R}^{\theta}_j(\Gamma_tp_t)(-x-tv,-v)]\cdot f(x,v)dz\\
		=&\int_{\mathbb{R}^{2d}}(\nabla_v\mathcal{R}^{\theta}_j\Gamma_tp_t+t\nabla_x\mathcal{R}^{\theta}_j\Gamma_tp_t)(-x-tv,-v)f(x,v)dz\\
		=&-\int_{\mathbb{R}^{2d}}(\nabla_v\mathcal{R}^{\theta}_j\Gamma_tp_t+t\nabla_x\mathcal{R}^{\theta}_j\Gamma_tp_t)(x+tv,v)f(x,v)dz.   
	\end{align*}
To avoid ambiguity, we stress that the last step does not follow from the change of variables formula.  
Instead, we use the fact that $\Gamma_t p_t$ is an even function in $(x,v) \in \mathbb{R}^{2d}$.  
Since the Littlewood--Paley block multipliers $(\phi^{\theta}_j)_{j \geq -1}$ are also even, it follows that $\mathcal{R}^{\theta}_j \Gamma_t p_t$ is even as well.  
Differentiating, we conclude that $\nabla_v \mathcal{R}^{\theta}_j \Gamma_t p_t + t \nabla_x \mathcal{R}^{\theta}_j \Gamma_t p_t$
is an odd function in $(x,v) \in \mathbb{R}^{2d}$.

	(ii) Using the change variables formula, we can see that for any $s>0$ and $(x,v)\in\mathbb{R}^{2d}$, 
	\begin{align*}
	\|\tau^{\cdot}_sf(x,v)\|_{L^p(\mathbb{R}^{2d})}^p=\int_{\mathbb{R}^{2d}}|f(x+x_0+(1-e^{-s})v_0,v+e^{-s}v_0)|^pdx_0dv_0=e^{sd}\|f\|_{L^p(\mathbb{R}^{2d})}^p. 	
	\end{align*}
	
(iii) {\color{black}Put
$\beta_{p,q}:=\beta+\frac{4d}{p}-\frac{4d}{q}\neq0$ and first let
$j\geq0$.}
Based on Lemma \ref{s4:lem-01} and \eqref{Thetat-property}, noting that $\cR_\ell^\theta \tau^{z_0}_s=\tau^{z_0}_s\cR_\ell^\theta$ holds for every $\ell\geq-1$ and $z_0\in\mathbb{R}^{2d}$, we have
\begin{align*}
    \mathcal{R}^{\theta}_jP_t(\tau^{z_0}_sf)(0)=\int_{\mR^{2d}} \mathcal{R}^{\theta}_j\Gamma_{t}p_t(x,v)\sum_{\ell\in\Theta_j^t}\Gamma_t\tau^{z_0}_s\cR_\ell^\theta f(x,v)dxdv,
\end{align*}
which by (ii) and \eqref{Besov-basic-inequality} yield that
\begin{align*}
\|\mathcal{R}^{\theta}_jP_t(\tau^{\cdot}_sf)(0)\|_{L^p(\mathbb{R}^{2d})}&\leq\int_{\mR^{2d}} |\mathcal{R}^{\theta}_j\Gamma_{t}p_t(x,v)|\sum_{\ell\in\Theta_j^t}\|\tau^{\cdot}_s\cR_\ell^\theta f(x,v)\|_{L^p(\mathbb{R}^{2d})}dxdv\\
&\le e^{\frac{sd}{p}}\int_{\mR^{2d}} |\mathcal{R}^{\theta}_j\Gamma_{t}p_t(x,v)|\sum_{\ell\in\Theta_j^t}\|\cR_\ell^\theta f\|_{L^p(\mathbb{R}^{2d})}dxdv\\
&\lesssim {\color{black}e^{\frac{sd}{p}}}\int_{\mR^{2d}} |\mathcal{R}^{\theta}_j\Gamma_{t}p_t(x,v)|\sum_{\ell\in\Theta_j^t}2^{-\ell(\frac{4d}{p}-\frac{4d}{q})}\|\cR_\ell^\theta f\|_{L^q(\mathbb{R}^{2d})}dxdv\\
&\lesssim {\color{black}e^{\frac{sd}{p}}}\int_{\mR^{2d}} |\mathcal{R}^{\theta}_j\Gamma_{t}p_t(x,v)|\sum_{\ell\in\Theta_j^t}2^{-(\beta+\frac{4d}{p}-\frac{4d}{q}) \ell}dxdv \|f\|_{\bB^\beta_{q;\theta}}.
\end{align*}
Then by Lemma \ref{s4:lem-01} again, we have
\begin{align*}
\|\mathcal{R}^{\theta}_jP_t(\tau^{\cdot}_sf)(0)\|_{L^p(\mathbb{R}^{2d})}
&\lesssim e^{\frac{sd}{p}}\|f\|_{\bB^\beta_{q;\theta}} 2^{-\beta_{p,q} j}(1+t2^{2j})^{|\beta_{p,q}|}\int_{\mR^{2d}} |\mathcal{R}^{\theta}_j\Gamma_{t}p_t(x,v)|dxdv\\
&\lesssim_{p,q} e^{\frac{sd}{p}}\|f\|_{\bB^\beta_{q;\theta}} 2^{-\beta_{p,q} j}[1+(t2^{2j})^{|\beta_{p,q}|}]\int_{\mR^{2d}} |\mathcal{R}^{\theta}_j\Gamma_{t}p_t(x,v)|dxdv. 
\end{align*}
In the last step, we decompose $(1+t2^{2j})^{|\beta_{p,q}|}$ into two terms, $1$ and $(t2^{2j})^{|\beta_{p,q}|}$, which will be treated separately. For each $m \ge 0$, applying Lemma \ref{lem-kineticestimates} with $l = m$ to the term with factor $1$ and with $l = |\beta_{p,q}| + m$ to the term with factor $(t2^{2j})^{|\beta_{p,q}|}$, we obtain that, for any $m \ge 0$,  
\begin{align*}
\|\mathcal{R}^{\theta}_jP_t(\tau^{\cdot}_sf)(0)\|_{L^p(\mathbb{R}^{2d})}
&\lesssim_{\beta,m,p,q} e^{\frac{sd}{p}}\|f\|_{\bB^\beta_{q;\theta}} 2^{-\beta_{p,q} j}\left([2^{-2j}t^{-1}]^m+(t2^{2j})^{|\beta_{p,q}|}[2^{-2j}t^{-1}]^{|\beta_{p,q}|+m}\right)\\
&\lesssim_{\beta,m,p,q} e^{\frac{sd}{p}}\|f\|_{\bB^\beta_{q;\theta}} 2^{-\beta_{p,q} j}[2^{-2j}t^{-1}]^m.
\end{align*}
{\color{black}For $j=-1$, the low-frequency kernel is treated directly.
Young's inequality, the change of variables in part (ii), and the bounded
moments of $q_t=\Gamma_tp_t$ on $0<t\leq T$ give
\begin{align*}
 \|\mathcal R_{-1}^\theta P_t(\tau_s^{\cdot}f)(0)\|_{L^p}
 \leq C_Te^{sd/p}\|f\|_{\mathbf B^\beta_{q;\theta}},
\end{align*}
after using the finite low-frequency overlap and Bernstein's inequality
block by block.
Since $2^{-2(-1)}t^{-1}=4/t$ is bounded below by $4/T$, enlarging $C_T$
proves \eqref{0219:02} also for $j=-1$.

(iv) For $j\geq0$, write the left-hand side as an integral against
$\mathcal R_j^\theta\Gamma_tp_t$.  Minkowski's inequality and part (ii)
give
\begin{align*}
 \|\mathcal R_j^\theta P_t(\tau_s^{\cdot}f)(0)\|_{L^p}
 \leq e^{sd/p}\|\mathcal R_j^\theta\Gamma_tp_t\|_{L^1}\|f\|_{L^p}.
\end{align*}
Lemma~\ref{lem-kineticestimates}, first with $\ell=0$ and then with
$\ell=m$, bounds the kernel norm by
$C[1\wedge(2^{-2j}t^{-1})^m]$.  The same low-frequency argument as above
handles $j=-1$.  This proves \eqref{critical-kernel-shift-estimate} and
completes the proof.}
\end{proof}

\begin{proof}[Proof of Theorem \ref{thm:Besov}]
{\color{black}
We prove (H) and (S) in parallel, but keep the probability norm outside
the spatial--time estimates.  In (H) we put $p=q=2$ and
$\beta_*=1$; in (S) we use the displayed $p$, put $q=1$, and set
$\beta_*=\beta_0$.  All estimates up to the stochastic convolution are
pathwise.  At the end we take the $L^2(\Omega)$ norm in (H), whereas in
(S) we take expectation of the pathwise $L^p_t$ norm.  In particular, no
spatial smoothing is used to improve integrability in $\omega$.

For $j\geq-1$, $a\geq0$, and $r>0$, set
\begin{align*}
 k_{j,a}(r):=1\wedge(2^{-2j}r^{-1})^a.
\end{align*}
The heat semigroup $H_t^n$ commutes with $P_t$ and the dyadic blocks and
is an $L^p$-contraction.  For the initial term, estimate
\eqref{critical-kernel-shift-estimate} is used in (H).  In (S), the two
estimates \eqref{0219:02} with $\beta=0$ and $m=0,2$ are combined.  Thus,
with constants also covering $j=-1$,
\begin{align}\label{Besov-initial-two-cases}
 \bigl\|(\mathcal R_j^\theta P_\cdot H^n_\cdot
       \tau_0^{\cdot}u_0)(0)\bigr\|_{L^p(0,T;L^p_{z_0})}
 \leq
 \begin{cases}
  C2^{-j}\|u_0\|_{L^2},&\text{in (H)},\\
  C2^{-\beta_0j}\|u_0\|_{L^1},&\text{in (S)}.
 \end{cases}
\end{align}
Indeed, in (S) the pointwise bound is
\begin{align*}
 C2^{4d(1-1/p)j}k_{j,2}(t)\|u_0\|_{L^1},
\end{align*}
and its $L^p(0,T)$ norm is bounded by
$C2^{-\beta_0j}\|u_0\|_{L^1}$ after the change of variables
$t=2^{-2j}r$.  The same calculation without the Bernstein factor gives
$2^{-j}$ in (H).  Notice that the critical case $p=q=2$, $\beta=0$ is
handled by Lemma~\ref{lem-properties}(iv), not by part (iii).

We next estimate the two terms containing $u$.  Lemma
\ref{lem-properties}(iii), with $\beta=-1$, $p=q$, and
$m=1-\varepsilon$, gives, for $0<\varepsilon<1/2$,
\begin{align}\label{Besov-Z1-correct}
 \|Z_{1,j}\|_{L^p(0,T;L^p_{z_0})}
 &\leq C e^{Td/p}\|b\|_\infty\,2^j
  \|k_{j,1-\varepsilon}*\|u(\cdot)\|_{L^p}\|_{L^p(0,T)}\notag\\
 &\leq C_\varepsilon e^{Td/p}\|b\|_\infty
 T^\varepsilon2^{-(1-2\varepsilon)j}
 \|u\|_{L^p(0,T;L^p)}.
\end{align}
For $Z_{2,j}$, Lemma~\ref{lem-properties}(i) and
Lemma~\ref{lem-kineticestimates} show that its kernel satisfies
\begin{align*}
 \|q_{2,j}^\theta(r)\|_{L^1}
 \leq C2^{-j}r^{-1/2}.
\end{align*}
Consequently,
\begin{align}\label{Besov-Z2-correct}
 \|Z_{2,j}\|_{L^p(0,T;L^p_{z_0})}
 \leq Ce^{Td/p}T^{1/2}2^{-j}
 \|u\|_{L^p(0,T;L^p)}.
\end{align}
The factors $e^{Td/p}$ in \eqref{Besov-Z1-correct} and
\eqref{Besov-Z2-correct} follow from Lemma~\ref{lem-properties}(ii).

We now make the square-function estimate for the stochastic convolution
explicit.  Componentwise in the $d$ Brownian coordinates, define
\begin{align*}
 \mathscr T^{k,1}_{j,r,s}(z_0)
 &:=[\mathcal R_j^\theta P_rH_r^n
       \tau_s^{z_0}(h_1(s)f_k)](0),\\
 \mathscr T^{k,2}_{j,r,s}(z_0)
 &:=[\mathcal R_j^\theta P_rH_r^n
       \tau_s^{z_0}(h_2(s)\nabla_vf_k)](0).
\end{align*}
Minkowski's inequality in $\ell^2$, the absolute kernel bound in
Lemma~\ref{lem-properties}(iv), and the $L^p$-contraction of $H_r^n$
give
\begin{align}\label{gamma-kernel-bound-Besov}
 &\bigl\|(\mathscr T^{k,1}_{j,r,s})_k
          \bigr\|_{L^p_{z_0}(\ell^2)}
 +\bigl\|(\mathscr T^{k,2}_{j,r,s})_k
          \bigr\|_{L^p_{z_0}(\ell^2)}\notag\\
 &\qquad\leq Ce^{sd/p}k_{j,1}(r)
 \left(\|F_1\|_\infty^{1/2}\|h_1(s)\|_{L^p}
      +\|F_3\|_\infty^{1/2}\|h_2(s)\|_{L^p}\right).
\end{align}
Here we used, before applying the convolution kernel,
\begin{align*}
 \left(\sum_k|h_1f_k|^2\right)^{1/2}
 \leq\|F_1\|_\infty^{1/2}|h_1|,
 \qquad
 \left(\sum_k|h_2\nabla_vf_k|^2\right)^{1/2}
 \leq\|F_3\|_\infty^{1/2}|h_2|.
\end{align*}
Thus the factors $f_k$ and $\nabla_vf_k$ are retained throughout; they
are not commuted through multiplication by $h_1$ or $h_2$.

Since $L^p(\mathbb R^{2d})$ is a UMD space for $1<p<\infty$, its
$\gamma$-radonifying norm is equivalent to the displayed
$L^p(\ell^2)$ square-function norm.  The
Burkholder--Davis--Gundy inequality and \eqref{gamma-kernel-bound-Besov}
therefore imply, for the present range $1<p\leq2$,
\begin{align}\label{Besov-Z3-correct}
 \left(\mathbb E\|Z_{3,j}\|_{L^p(0,T;L^p_{z_0})}^p\right)^{1/p}
 \leq C2^{-j}\Big(&\|F_1\|_\infty^{1/2}
       \|h_1\|_{L^2(\Omega\times[0,T];L^p)}\notag\\
 &+\|F_3\|_\infty^{1/2}
       \|h_2\|_{L^2(\Omega\times[0,T];L^p)}\Big).
\end{align}
For completeness, this follows by first applying the stochastic inequality
at each time, then using
\begin{align*}
 \int_0^T A(t)^{p/2}\,dt
 \leq T^{1-p/2}\left(\int_0^T A(t)\,dt\right)^{p/2},
 \qquad
 \int_0^T k_{j,1}(r)^2\,dr\leq C_T2^{-2j},
\end{align*}
and finally Jensen's inequality in $\omega$.  This also covers $p=2$.

It remains to estimate $Z_{4,j}$.  The $g_1$ term is treated by
Lemma~\ref{lem-properties}(iii) with $\beta=-1$ and $p=q$.  For $g_2$,
we use part (iv) in (H), and part (iii) with $q=1$, $\beta=0$, and
$m=0,2$ in (S).  Young's inequality in time yields
\begin{align}\label{Besov-Z4-two-cases}
 \|Z_{4,j}\|_{L^p(0,T;L^p_{z_0})}
 \leq C2^{-j}\|g_1\|_{L^p(0,T;\mathbf B^{-1}_{p;\theta})}
 +\begin{cases}
  C2^{-2j}\|g_2\|_{L^2(0,T;L^2)},&\text{in (H)},\\
  C2^{-\beta_0j}\|g_2\|_{L^1(0,T;L^1)},&\text{in (S)}.
 \end{cases}
\end{align}

We now take the probability norms.  Lemma~\ref{lem-properties}(ii) gives
\begin{align*}
 \|\mathcal R_j^\theta u(t)\|_{L^p}
 =e^{-td/p}\|\tau_t^{\cdot}\mathcal R_j^\theta u(t,0)\|_{L^p_{z_0}}.
\end{align*}
In (H), take the $L^2(\Omega)$ norm in
\eqref{Besov-initial-two-cases}--\eqref{Besov-Z4-two-cases}.  If
$0<\varepsilon<(1-\beta)/2$, this gives
\begin{align}\label{Besov-block-H}
 \|\mathcal R_j^\theta u\|_{L^2(\Omega\times[0,T];L^2)}
 \leq C2^{-(1-2\varepsilon)j}
 \left(\mathcal K_H+
 T^{\varepsilon\wedge1/2}\|u\|_{L^2(\Omega\times[0,T];L^2)}\right).
\end{align}
In (S), take expectation of every pathwise $L^p_t$ estimate and use
\begin{align*}
 \mathbb E\|g_1\|_{L^p_t\mathbf B^{-1}_{p;\theta}}
 &\leq C_T\|g_1\|_{L^2(\Omega\times[0,T];
                 \mathbf B^{-1}_{p;\theta})},\\
 \mathbb E\|Z_{3,j}\|_{L^p_tL^p}
 &\leq\left(\mathbb E\|Z_{3,j}\|_{L^p_tL^p}^p\right)^{1/p}.
\end{align*}
With $\bar\beta:=\min\{\beta_0,1\}$ and
$0<\varepsilon<(\bar\beta-\beta)/2$, we obtain
\begin{align}\label{Besov-block-S}
 \mathbb E\|\mathcal R_j^\theta u\|_{L^p(0,T;L^p)}
 \leq C2^{-(\bar\beta-2\varepsilon)j}
 \left(\mathcal K_{S,p}+
 T^{\varepsilon\wedge1/2}\mathbb E\|u\|_{L^p(0,T;L^p)}\right).
\end{align}
This is precisely where the outer probability exponent in (S) remains
one.

To pass from the blocks to the Besov norm in (H), we use
\begin{align*}
 \|u\|_{L^2(\Omega\times[0,T];\mathbf B^\beta_{2;\theta})}
 \leq\sum_{j\geq-1}2^{\beta j}
 \|\mathcal R_j^\theta u\|_{L^2(\Omega\times[0,T];L^2)}.
\end{align*}
In (S), Minkowski's inequality is applied pathwise before expectation:
\begin{align*}
 \mathbb E\|u\|_{L^p(0,T;\mathbf B^\beta_{p;\theta})}
 \leq\sum_{j\geq-1}2^{\beta j}
 \mathbb E\|\mathcal R_j^\theta u\|_{L^p(0,T;L^p)}.
\end{align*}
The respective series converge by \eqref{Besov-block-H} and
\eqref{Besov-block-S}.

Finally, we remove the terms involving the unregularized norm of $u$ by
an exponential time weight.  Multiplying each Duhamel term by
$e^{-\lambda t}$ replaces every Volterra kernel $K(t-s)$ by
$e^{-\lambda(t-s)}K(t-s)$ and each input at time $s$ by its weighted
version.  In particular,
\begin{align*}
 \int_0^T e^{-\lambda r}r^{-1+\varepsilon}\,dr
 \leq C_\varepsilon\lambda^{-\varepsilon},
 \qquad
 \int_0^T e^{-\lambda r}r^{-1/2}\,dr
 \leq C\lambda^{-1/2}.
\end{align*}
Hence \eqref{Besov-block-H} and \eqref{Besov-block-S}, followed by the
same dyadic summation, give respectively
\begin{align*}
 \|e^{-\lambda\cdot}u\|_{L^2(\Omega\times[0,T];
       \mathbf B^\beta_{2;\theta})}
 &\leq C\mathcal K_H+C(\lambda^{-\varepsilon}+\lambda^{-1/2})
 \|e^{-\lambda\cdot}u\|_{L^2(\Omega\times[0,T];L^2)},\\
 \mathbb E\|e^{-\lambda\cdot}u\|_{L^p(0,T;
       \mathbf B^\beta_{p;\theta})}
 &\leq C\mathcal K_{S,p}+C(\lambda^{-\varepsilon}+\lambda^{-1/2})
 \mathbb E\|e^{-\lambda\cdot}u\|_{L^p(0,T;L^p)}.
\end{align*}
For $\beta>0$, $\mathbf B^\beta_{p;\theta}$ continuously embeds into
$L^p$.  Choosing $\lambda$ large absorbs the last term in the appropriate
left-hand side.  Since $e^{-\lambda T}\leq e^{-\lambda t}\leq1$, removing
the weight proves \eqref{apriori-besov-es} and
\eqref{apriori-besov-es-first-moment}.}
\end{proof}

{\color{black}
\subsection{An approximation scheme}\label{subsec-5-2}
Recall that $\alpha_r^2$ denote the smooth truncation function on $\mathbb{R}^d_v$, for every $r>1$. We introduce the following approximation scheme: 
\begin{align}\label{SPDE-iteration}
\begin{split}
    df_{n,r}=\Big(&\Delta_vf_{n,r}+\frac{1}{n}\Delta_xf_{n,r}-\alpha^2_r(v)v\cdot\nabla_xf_{n,r}+\nabla_v\cdot(\alpha^2_r(v)vf_{n,r})-\nabla_v\cdot(f_{n,r}V\star_x\langle\alpha^2_rf_{n,r}\rangle)\Big)dt\\
    &-\nabla_v\cdot(\sigma(f_{n,r})dW_F)
    +{\color{black}\frac12\nabla_v\cdot\left(
       F_1\sigma'(f_{n,r})^2\nabla_vf_{n,r}\right)dt.}  
\end{split}
\end{align}
{\color{black}Throughout this subsection, a \emph{weak solution} of an
approximating equation means an adapted $L^2$-valued variational solution
with the indicated Sobolev regularity, satisfying the identity obtained by
testing that equation against $C_c^\infty(\mathbb R^{2d})$.  This terminology
does not include nonnegativity, mass preservation, or moment bounds; those
properties are established below.}

In the following, we derive an a priori  $L^2(\mathbb{R}^{2d})$-estimate for the approximation scheme \eqref{SPDE-iteration} for fixed $n,r\geq1$. 

\begin{proposition}\label{prp-uniforml2}
{\color{black}Assume that $f_0$, $V$, and $\sigma(\cdot)$ satisfy
Assumptions~\ref{Assump-initialdata-regular}, \ref{Assump-ker}, and
\ref{Assump-sigma-n}, respectively.} Then for every $n,r\geq1$, there exists {\color{black}at least one} weak solution $f_{n,r}$ to equation~\eqref{SPDE-iteration} with initial data $f_0$.

Moreover, there exists a constant {\color{black}$C=C(T,d,\sigma,
\|F_3\|_{L^\infty(\mathbb R^{2d})},
\|F_1\|_{L^\infty(\mathbb R^{2d})})>0$} such that
\begin{align}\label{0425:00}
    \sup_{n,r\geq1}\left(\sup_{t \in [0,T]} \mathbb{E} \|f_{n,r}(t)\|_{L^2(\mathbb{R}^{2d})}^2 
    + \mathbb{E} \int_0^T \|\nabla_v f_{n,r}(s)\|_{L^2(\mathbb{R}^{2d})}^2 \, \dif s +\frac{1}{n}\mathbb{E} \int_0^T \|\nabla_x f_{n,r}(s)\|_{L^2(\mathbb{R}^{2d})}^2 \, \dif s\right)\le C\mathbb{E} \|f_0\|_{L^2(\mathbb{R}^{2d})}^2.
\end{align}
\end{proposition}
\begin{proof}

We first derive an a priori estimate for \eqref{SPDE-iteration}. For each $R=(R_1,R_2)\in\mathbb{R}_+^2$, let $\alpha_R$ be the cut-off function defined in \eqref{alpha-R}. Assume that $f_{n,r}$ is a weak solution of \eqref{SPDE-iteration} with $\nabla_xf_{n,r}\in L^2([0,T]\times\mathbb{R}^{2d})$ almost surely, applying It\^o's formula (\cite{Krylov}), we obtain for every $t \in [0,T]$,
\begin{align*}
\frac{1}{2}\int_{\mathbb{R}^{2d}} |f_{n,r}(t)|^2 \alpha_R 
=& \ \frac{1}{2} \int_{\mathbb{R}^{2d}} |f_0|^2 \alpha_R 
- \int_0^t \int_{\mathbb{R}^{2d}} \nabla_v(f_{n,r} \alpha_R) \cdot \nabla_v f_{n,r} - \frac{1}{n}\int_0^t \int_{\mathbb{R}^{2d}} \nabla_x(f_{n,r} \alpha_R) \cdot \nabla_x f_{n,r}\\
&-\int^t_0\int_{\mathbb{R}^{2d}}\alpha_Rf_{n,r}\alpha^2_r(v)v\cdot\nabla_xf_{n,r} + \int_0^t \int_{\mathbb{R}^{2d}} \nabla_v(f_{n,r} \alpha_R)\cdot f_{n,r} V\star_x\langle\alpha_r^2f_{n,r}\rangle \\ 
&+ \int_0^t \int_{\mathbb{R}^{2d}} \alpha_R f_{n,r} \nabla_v \cdot (\alpha_r^2v f_{n,r})  + \int_0^t\int_{\mathbb{R}^{2d}} \nabla_v(f_{n,r} \alpha_R)\cdot \sigma(f_{n,r}) \, dW_F \\
& - \frac{1}{2} \int_0^t \int_{\mathbb{R}^{2d}} \nabla_v(f_{n,r} \alpha_RF_1)\cdot \sigma'(f_{n,r})^2 \nabla_v f_{n,r} \\
& + \frac{1}{2} \int_0^t \int_{\mathbb{R}^{2d}} \alpha_R |\nabla_v \sigma(f_{n,r})|^2 F_1 
 + \frac{1}{2} \int_0^t \int_{\mathbb{R}^{2d}} \alpha_R \sigma(f_{n,r})^2 F_3 \\
=: & \frac{1}{2} \int_{\mathbb{R}^{2d}} |f_0|^2 \alpha_R - \int_0^t \int_{\mathbb{R}^{2d}} \alpha_R |\nabla_v f_{n,r}|^2- \frac{1}{n}\int_0^t \int_{\mathbb{R}^{2d}} \alpha_R |\nabla_x f_{n,r}|^2\\
&- \int_0^t \int_{\mathbb{R}^{2d}} f_{n,r} \nabla_v f_{n,r} \cdot \nabla_v \alpha_R-\frac{1}{n} \int_0^t \int_{\mathbb{R}^{2d}} f_{n,r} \nabla_x f_{n,r} \cdot \nabla_x \alpha_R\\
&-\int^t_0\int_{\mathbb{R}^{2d}}f_{n,r}\alpha^2_rv\cdot\nabla_xf_{n,r}\alpha_R+ {\color{black}\sum_{i=1}^5 I_i}.
\end{align*}

Thanks to the assumed regularity $\nabla_v f_{n,r},\nabla_xf_{n,r} \in L^2([0,T] \times \mathbb{R}^{2d})$ almost surely, and the properties of the cut-off function $\alpha_R$, we estimate
\begin{align*}
	- \int_0^t \int_{\mathbb{R}^{2d}} f_{n,r} \nabla_v f_{n,r} \cdot \nabla_v \alpha_R 
	\lesssim \frac{1}{R_2} \int_0^t \|f_{n,r}\|_{L^2(\mathbb{R}^{2d})} \|\nabla_v f_{n,r}\|_{L^2(\mathbb{R}^{2d})}  
	\rightarrow 0, \quad \text{as } R_2 \rightarrow +\infty,
\end{align*}
almost surely. Similarly, we obtain that 
\begin{align*}
- \int_0^t \int_{\mathbb{R}^{2d}} f_{n,r} \nabla_x f_{n,r} \cdot \nabla_x \alpha_R 
	\lesssim \frac{1}{R_1} \int_0^t \|f_{n,r}\|_{L^2(\mathbb{R}^{2d})} \|\nabla_x f_{n,r}\|_{L^2(\mathbb{R}^{2d})}  
	\rightarrow 0, \quad \text{as } R_1 \rightarrow +\infty,
\end{align*}
and
\begin{align*}
	\frac{1}{2}\int^t_0\int_{\mathbb{R}^{2d}}\nabla_x\alpha_R\cdot \alpha_r^2v|f_{n,r}|^2
	\lesssim \frac{R_2}{R_1} \int_0^t \|f_{n,r}\|_{L^2(\mathbb{R}^{2d})}^2 
	\rightarrow 0, \quad \text{as } R_1 \rightarrow +\infty.
\end{align*}

We now estimate the terms $I_1$ through $I_8$. For $I_1$, since $\nabla_v \cdot V\star_x\langle\alpha^2_rf_{n,r}\rangle = 0$, we integrate by parts and apply Young's convolution inequality to obtain
\begin{align*}
I_1 
= \frac{1}{2} \int_0^t \int_{\mathbb{R}^{2d}} |f_{n,r}|^2 \nabla_v \alpha_R \cdot V\star_x\langle\alpha^2_rf_{n,r}\rangle \, ds 
\lesssim \frac{1}{R_2} \|f_{n,r}\|_{L^1([0,T]; L^2(\mathbb{R}^{2d}))}^2 \|V\star_x\langle\alpha^2_rf_{n,r}\rangle\|_{L^\infty([0,T] \times \mathbb{R}^d)} 
\rightarrow 0,
\end{align*}
as $R \to +\infty$, almost surely.

For $I_2$, on the one hand, we use the chain rule to see that
\begin{align*}
\int_0^t \int_{\mathbb{R}^{2d}} \alpha_R f_{n,r} \nabla_v \cdot (\alpha_r^2v f_{n,r})   
&= \int_0^t \int_{\mathbb{R}^{2d}} (d\alpha^2_r+\nabla_v\alpha^2_r\cdot v)\alpha_R f_{n,r} f_{n,r}   
+ \int_0^t \int_{\mathbb{R}^{2d}} \alpha_R f_{n,r} \alpha^2_rv \cdot \nabla_v f_{n,r}.
\end{align*}
On the other hand, by applying the integration by parts formula and the chain rule, it follows that  
\begin{align*}
\int_0^t \int_{\mathbb{R}^{2d}} \alpha_R f_{n,r} \nabla_v \cdot (\alpha^2_rv f_{n,r})   
&= - \int_0^t\int_{\mathbb{R}^{2d}} \nabla_v (\alpha_R f_{n,r})\cdot \alpha^2_rv f_{n,r}  \\
&= - \int_0^t \int_{\mathbb{R}^{2d}} \alpha_R \nabla_v f_{n,r}\cdot \alpha^2_rv f_{n,r} 
- \int_0^t\int_{\mathbb{R}^{2d}} \nabla_v \alpha_R f_{n,r}\cdot \alpha^2_rv f_{n,r}.
\end{align*}
Combining the above two identities, we deduce 
\begin{align*}
\int_0^t\int_{\mathbb{R}^{2d}} \alpha_Rf_{n,r} \nabla_v \cdot (\alpha^2_rv f_{n,r})   
= \int_0^t \int_{\mathbb{R}^{2d}} \frac{1}{2}(d\alpha^2_r+\nabla_v\alpha^2_r\cdot v)\alpha_R f_{n,r} f_{n,r}  
- \frac{1}{2} \int_0^t \int_{\mathbb{R}^{2d}} \nabla_v \alpha_R f_{n,r}\cdot \alpha^2_rv f_{n,r}.
\end{align*}
By the support properties of $\nabla_v \alpha_R$, it follows that
\begin{align*}
\frac{1}{2} \int_0^t\int_{\mathbb{R}^{2d}} \nabla_v \alpha_R f_{n,r}\cdot v f_{n,r} 
\leq \int_0^t \int_{\mathbb{R}^{2d}} I_{\{R_2 \leq |v| \leq 2R_2\}} |f_{n,r}|^2 
\rightarrow 0,
\end{align*}
as $R_2 \to +\infty$, almost surely. Hence, by the dominated convergence theorem, we obtain
\begin{align*}
\int_0^t\int_{\mathbb{R}^{2d}} f_{n,r} \nabla_v \cdot (\alpha^2_rv f_{n,r})  
\leq  \int_0^t \|f_{n,r}\|_{L^2(\mathbb{R}^{2d})}^2,
\end{align*}
almost surely.

For $I_3$, the stochastic integral has zero expectation due to the martingale property, i.e.,
$$
\mathbb{E} I_3 = 0.
$$

For $I_4$, we apply the chain rule to obtain
$$
	{\color{black}I_4 = -\frac{1}{2}\int_0^t\int_{\mathbb{R}^{2d}} \nabla_v \alpha_R \cdot f_{n,r} \sigma'(f_{n,r})^2 \nabla_v f_{n,r} F_1},
$$
{\color{black}where the interior part of the It\^o correction cancels exactly with the $F_1$-part of the quadratic variation.} Since $f_{n,r}, \nabla_v f_{n,r} \in L^2([0,T]; L^2(\mathbb{R}^{2d}))$ almost surely, and using the properties of $\alpha_R$, we deduce that
$$
\frac{1}{2} \int_0^t \int_{\mathbb{R}^{2d}} \nabla_v \alpha_R \cdot f_{n,r} \sigma'(f_{n,r})^2 \nabla_v f_{n,r} F_1 
\leq \frac{C(\sigma')}{R_2} \int_0^t \|f_{n,r}\|_{L^2(\mathbb{R}^{2d})} \|\nabla_v f_{n,r}\|_{L^2(\mathbb{R}^{2d})}  \to 0,
$$
as $R_2 \to +\infty$, almost surely.

For ${\color{black}I_5}$, the linear growth of $\sigma(\cdot)$ yields
$$
{\color{black}I_5 \leq C(\sigma,F_3)\int^t_0\int_{\mathbb{R}^{2d}}|f_{n,r}|^2}.
$$
Then applying Gronwall's inequality, we complete the proof of the $L^2(\mathbb{R}^{2d})$-a priori estimate.  
 
 {\color{black}We now construct a solution; the results cited above for the torus
or for bounded domains do not apply verbatim on $\mathbb R^{2d}$.  Fix
$n,r\geq1$.  In the construction below, $f_{n,r,m}$ denotes the $m$th
Galerkin approximation, and $f_{n,r}$ denotes a subsequential limit, which
will be the weak solution asserted in the proposition.  The symbols $g$
and $h$ are used only as generic elements of $L^2(\mathbb R^{2d})$ in the
coefficient estimates.  We use the Gelfand triple
\begin{align*}
 \mathcal V_n:=H^1(\mathbb R^{2d})\hookrightarrow L^2(\mathbb R^{2d})
 \hookrightarrow\mathcal V_n^*,\qquad
 \|g\|_{\mathcal V_n}^2
 :=\|g\|_2^2+\|\nabla_vg\|_2^2+n^{-1}\|\nabla_xg\|_2^2.
\end{align*}
Because $\operatorname{supp}\sigma\Subset(0,\infty)$ in
Assumption~\ref{Assump-sigma-n}, the zero extension
\begin{align}\label{regular-coefficient-extension}
 \widehat\sigma(a):=\begin{cases}\sigma(a),&a>0,\\0,&a\leq0,\end{cases}
\end{align}
is smooth on $\mathbb R$, with bounded derivatives and
$|\widehat\sigma(a)|\leq\|\widehat\sigma'\|_\infty|a|$.

We record explicitly the nonlinear estimates used below.  Set
\begin{align*}
 \rho_r(g)(x):=\int_{\mathbb R^d}\alpha_r^2(v)g(x,v)\,\dif v,
 \qquad b_r(g):=V\star_x\rho_r(g).
\end{align*}
Cauchy--Schwarz in $v$ and Young's convolution inequality give
\begin{align}\label{fixed-r-interaction-bounds}
 \|\rho_r(g)\|_{L_x^2}
 &\leq\|\alpha_r^2\|_{L_v^2}\|g\|_2,\qquad
 \|b_r(g)\|_{L_x^\infty}
 \leq\|V\|_{L_x^2}\|\rho_r(g)\|_{L_x^2},\notag\\
 \|g b_r(g)-h b_r(h)\|_2
 &\leq C_r\|V\|_{L_x^2}
       (\|g\|_2+\|h\|_2)\|g-h\|_2.
\end{align}
Thus $g\mapsto-\nabla_v\cdot(gb_r(g))$ is locally Lipschitz from
$L^2$ to $\mathcal V_n^*$.  If
$B(g)e_{k,i}:=-\partial_{v_i}(\widehat\sigma(g)f_k)$, then the definition
of $F_1$ also yields
\begin{align}\label{fixed-r-noise-lipschitz}
 \|B(g)-B(h)\|_{L_2(\ell^2;\mathcal V_n^*)}^2
 \leq\|F_1\|_\infty\|\widehat\sigma'\|_\infty^2\|g-h\|_2^2.
\end{align}

Choose an $L^2$-orthonormal family
$(e_j)_{j\geq1}\subset C_c^\infty(\mathbb R^{2d})$ whose span is dense
in $\mathcal V_n$, and let $P_m$ be the corresponding projection.
Project the equation with $\sigma$ replaced by $\widehat\sigma$, initial
datum $P_mf_0$, and first multiply the interaction drift by a smooth
cutoff of $\|f_{n,r,m}\|_2$.  The resulting finite-dimensional drift and
Hilbert--Schmidt diffusion are locally Lipschitz by
\eqref{fixed-r-interaction-bounds}--\eqref{fixed-r-noise-lipschitz};
therefore the projected equation has a unique maximal solution.
For its exit time
$\tau_{L,m}:=\inf\{t:\|f_{n,r,m}(t)\|_2\geq L\}$, It\^o's
formula, $\langle P_mG,f_{n,r,m}\rangle
=\langle G,f_{n,r,m}\rangle$, and the
Burkholder--Davis--Gundy inequality imply
\begin{align}\label{fixed-n-r-Galerkin-estimate}
 \mathbb E\sup_{t\leq T\wedge\tau_{L,m}}\|f_{n,r,m}(t)\|_2^2
 +\mathbb E\int_0^{T\wedge\tau_{L,m}}
  \left(\|\nabla_vf_{n,r,m}\|_2^2
  +n^{-1}\|\nabla_xf_{n,r,m}\|_2^2\right)\dif t
 \leq C_{n,r,T}\mathbb E\|f_0\|_2^2.
\end{align}
Indeed, the interaction term cancels after integration by parts in $v$,
because $b_r(f_{n,r,m})$ is independent of $v$, and all other lower-order
coefficients are bounded for fixed $r$.  Hence
$\mathbb P(\tau_{L,m}<T)\leq
C_{n,r,T}L^{-2}\mathbb E\|f_0\|_2^2$; first the norm cutoff and then
$L$ can be removed.  Repeating the sharper expectation computation
already given above, without the time supremum, yields \eqref{0425:00}
uniformly in $m,n,r$.

For completeness, we describe the limit
$f_{n,r,m}\to f_{n,r}$ explicitly.  The drift is
bounded in $L^1(\Omega\times(0,T);\mathcal V_n^*)$, and
\eqref{fixed-r-noise-lipschitz} gives Aldous' increment estimate for the
martingale part in $H^{-s}_{\rm loc}$ for sufficiently large $s$.
Local Rellich compactness, the Aubin--Lions lemma, and Aldous' criterion
give tightness in
\begin{align*}
 L^2(0,T;L^2_{\rm loc}(\mathbb R^{2d}))
 \cap D([0,T];H^{-s}_{\rm loc}(\mathbb R^{2d})),
\end{align*}
together with weak compactness in $L^2(0,T;\mathcal V_n)$.  On a
Skorokhod representation, $f_{n,r,m}\to f_{n,r}$ strongly in the first
space and weakly in the last.  Thus
$\widehat\sigma'(f_{n,r,m})^2\nabla_vf_{n,r,m}
\rightharpoonup\widehat\sigma'(f_{n,r})^2\nabla_vf_{n,r}$ locally.
Moreover, $\rho_r(f_{n,r,m})\rightharpoonup\rho_r(f_{n,r})$ in
$L^2(0,T;L_x^2)$, while \eqref{fixed-r-interaction-bounds} bounds
$b_r(f_{n,r,m})$ in $L^2(0,T;L_x^\infty)$.  Testing the convolution
against compactly supported functions identifies every weak local limit
as $b_r(f_{n,r})$.  Splitting
\begin{align*}
 f_{n,r,m}b_r(f_{n,r,m})-f_{n,r}b_r(f_{n,r})
 &=(f_{n,r,m}-f_{n,r})b_r(f_{n,r,m})
 +f_{n,r}\bigl(b_r(f_{n,r,m})-b_r(f_{n,r})\bigr)
\end{align*}
then shows convergence in distributions: the first term tends to zero
by strong local $L^2$ convergence and the
$L_t^2L_x^\infty$ bound, and the second by the just identified weak
convergence.  This suffices to pass to the interaction drift.
The strong local convergence of $\widehat\sigma(f_{n,r,m})$ identifies the
quadratic and cross variations of the limiting martingales, and hence
the limiting stochastic integral.  We obtain $f_{n,r}$ as a martingale
solution with coefficient $\widehat\sigma$.

Finally, apply the generalized It\^o formula to smooth convex
approximations of $S(a)=\frac12(a^-)^2$.  Since
$\widehat\sigma=\widehat\sigma'=0$ on $(-\infty,0]$, the stochastic
term, its quadratic variation, and its correction vanish on the
negative set.  The $x$-transport and interaction contributions cancel,
the Laplacians are dissipative, and the damping term is bounded by
$C_r\|f_{n,r}^-\|_2^2$.  Passing to the convex limit gives
\begin{align*}
 \mathbb E\|f_{n,r}^-(t)\|_2^2
 +2\mathbb E\int_0^t
  \left(\|\nabla_vf_{n,r}^-\|_2^2
  +n^{-1}\|\nabla_xf_{n,r}^-\|_2^2\right)\dif s
 \leq C_r\int_0^t\mathbb E\|f_{n,r}^-(s)\|_2^2\,\dif s.
\end{align*}
Gronwall's inequality shows that $f_{n,r}\geq0$, so
$\widehat\sigma(f_{n,r})=\sigma(f_{n,r})$ and the auxiliary extension has
been removed.  This gives the asserted martingale solution $f_{n,r}$.
The passage back to the
original stochastic basis, and uniqueness, are completed without
circularity in Corollary~\ref{fixed-n-r-pathwise-uniqueness} below, after
the estimates needed for the kinetic uniqueness argument have been
proved for every weak solution.}

\end{proof}

In the following, we show the non-negativity of the solution to
\eqref{SPDE-iteration} and establish an $L^1(\mathbb{R}^{2d})$ estimate.
Fix $n,r\geq1$, and let $f_{n,r}$ be a weak solution furnished by
Proposition~\ref{prp-uniforml2}.  Only within the proof below, we abbreviate
$f_{n,r}$ as $f$.

\begin{lemma}\label{lem-L1-es}
	{\color{black}Assume that $f_0$, $V$, and $\sigma(\cdot)$ satisfy
 Assumptions~\ref{Assump-initialdata-regular}, \ref{Assump-ker}, and
 \ref{Assump-sigma-n}, respectively.} Let $f_{n,r}$ be a weak solution to
 equation~\eqref{SPDE-iteration} with initial condition $f_0$. Then, for
 every $t \in [0,T]$, the function $f_{n,r}(t)$ is almost surely
 nonnegative. Moreover, the total mass is almost surely preserved in time,
 that is,
\begin{align}\label{preservationofmean}
    \int_{\mathbb{R}^{2d}} f_{n,r}(t,z)\,\dif z
    = \int_{\mathbb{R}^{2d}} f_0(z)\,\dif z,
    \quad \text{for all } t \in [0,T].
\end{align}

\end{lemma}
\begin{proof}
	For the rest of this proof, write $f:=f_{n,r}$.
	We begin by establishing the integrability estimate
\begin{align}\label{est:L^1_finite}
    \mathbb{E} \left[\sup_{t \in [0,T]} \int_{\mathbb{R}^{2d}} |f(t,z)|\,\mathrm{d}z \right] < \infty,
\end{align}
and verifying that the mass conservation identity \eqref{preservationofmean} holds.

For each $R > 1$, recall $\alpha_R(x,v) := \alpha^1_{R^2}(x)\alpha^2_{R}(v)$ defined as \eqref{alpha-R} with $(R_1,R_2)=(R^2,R)$. We emphasize that, with a slight abuse of notation, $R$ here denotes a real number rather than a vector. This differs from the parameter $R$ introduced in Section \ref{sec-3}. It is easy to see that for any  $(x,v)\in \text{supp}\alpha_R$, $\alpha_{2R}(x,v)=1$.
Moreover, to approximate the absolute value function $|\cdot|$, we introduce the smooth function $a_\delta(r) := \sqrt{r^2 + \delta^2} - \delta$, for every $\delta\in(0,1)$. $a_\delta$ satisfies the bounds
\begin{align}\label{0422:-1}
    |a_\delta(r)| \leq |r|, \quad |a_\delta'(r)| \leq 1, \quad |a_\delta''(r)| \leq |r|^{-1} I_{\{|r| \leq \sqrt{\delta}\}} + \sqrt{\delta} I_{\{|r| > \sqrt{\delta}\}}.
\end{align}

Let $(\eta_{\gamma})_{\gamma \in (0,1)}$ be a sequence of standard convolution kernels on $\mathbb{R}^d_x$. Applying It\^o's formula to $\int_{\mathbb{R}^{2d}} a_\delta(\eta_{\gamma} \ast f(t))\,\alpha_R$, and then using the integration by parts formula and passing to the limit $\gamma \to 0$, we obtain the following identity in the distributional sense:

\begin{align*}
    \int_{\mathbb{R}^{2d}} a_\delta(f(t)) \alpha_R =:\; & \int_{\mathbb{R}^{2d}} a_\delta(f_0) \alpha_R +I_{b,\delta}(f,t)+I_{W,\delta}(f,t)+I_{1,\delta}(f,t)+I_{2,3,\delta}(f,t).
\end{align*}
where 
\begin{align*}
	I_{b,\delta}(f,t)=& - \int_0^t \int_{\mathbb{R}^{2d}} a_\delta'(f) \nabla_v f \cdot \nabla_v \alpha_R
    - \int_0^t \int_{\mathbb{R}^{2d}} a_\delta''(f) |\nabla_v f|^2 \alpha_R \notag\\
   & - \frac{1}{n}\int_0^t \int_{\mathbb{R}^{2d}} a_\delta'(f) \nabla_x f \cdot \nabla_x \alpha_R
    - \frac{1}{n}\int_0^t \int_{\mathbb{R}^{2d}} a_\delta''(f) |\nabla_x f|^2 \alpha_R \notag\\
    & + \int_0^t \int_{\mathbb{R}^{2d}} a_\delta(f)\alpha^2_r v \cdot \nabla_x\alpha_R
	    + {\color{black}\int_0^t \int_{\mathbb{R}^{2d}} a_\delta'(f) \nabla_v \cdot ([\alpha^2_rv-V\star_x\langle\alpha^2_rf\rangle] f) \alpha_R,}\\
    I_{W,\delta}(f,t)=&-\int_0^t \int_{\mathbb{R}^{2d}}\alpha_Ra_\delta'(f) \nabla_v \cdot (\sigma(f)\,\mathrm{d}W_F),\\
    I_{1,\delta}(f,t)=&-\frac{1}{2}\int_0^t \int_{\mathbb{R}^{2d}}\nabla_v(\alpha_Ra_\delta'(f)F_1) (\sigma'(f)^2 \nabla_v f)+\frac{1}{2}\int_0^t \int_{\mathbb{R}^{2d}}\alpha_R a_\delta''(f)  |\nabla_v \sigma(f)|^2 F_1,
\end{align*}
and
\begin{align*}
	{\color{black}I_{2,3,\delta}(f,t)=\frac12\int_0^t\int_{\mathbb{R}^{2d}} \alpha_R a_\delta''(f)|\sigma(f)|^2 F_3}. 
\end{align*}

First, by the integral by parts and observing that $a_{\delta}'' \ge 0$, one can see that
\begin{align*}
 I_{b,\delta}(f,t)
    \le& - \int_0^t \int_{\mathbb{R}^{2d}} a_\delta'(f) \nabla_v f \cdot \nabla_v \alpha_R - \frac{1}{n}\int_0^t \int_{\mathbb{R}^{2d}} a_\delta'(f) \nabla_x f \cdot \nabla_x \alpha_R\\
    &+ \int_0^t \int_{\mathbb{R}^{2d}} a_\delta(f) \alpha^2_rv \cdot \nabla_x\alpha_R+ {\color{black}\int_0^t \int_{\mathbb{R}^{2d}} a_\delta'(f) \nabla_v \cdot ([\alpha^2_rv-V\star_x\langle\alpha^2_rf\rangle] f) \alpha_R}
\end{align*}
Applying the chain rule, we have
\begin{align*}
&\quad- \int_0^t \int_{\mathbb{R}^{2d}} a_\delta'(f) \nabla_v f \cdot \nabla_v \alpha_R- \frac{1}{n}\int_0^t \int_{\mathbb{R}^{2d}} a_\delta'(f) \nabla_x f \cdot \nabla_x \alpha_R+ \int_0^t \int_{\mathbb{R}^{2d}} a_\delta(f) \alpha^2_rv \cdot \nabla_x\alpha_R\notag \\
&
=- \int_0^t \int_{\mathbb{R}^{2d}} \nabla_v a_\delta(f) \cdot \nabla_v \alpha_R- \frac{1}{n}\int_0^t \int_{\mathbb{R}^{2d}} \nabla_xa_\delta(f) \cdot \nabla_x \alpha_R + \int_0^t \int_{\mathbb{R}^{2d}} a_\delta(f) \alpha^2_rv \cdot \nabla_x\alpha_R \notag\\
&= \int_0^t \int_{\mathbb{R}^{2d}} a_{\delta}(f)[\Delta_v \alpha_R+\frac{1}{n}\Delta_x\alpha_R+ \alpha^2_rv \cdot \nabla_x \alpha_R] 
  {\color{black}\lesssim \left(\frac{1}{R^2}+\frac{1}{R^4}+ \frac{1}{R}\right) \int_0^t \int_{\mathbb{R}^{2d}} |f| \alpha_{2R}.} 
\end{align*}
{\color{black}The estimate \eqref{0425:00} controls
$\sup_{t\leq T}\mathbb E\|f(t)\|_2^2$, but it does not provide a
deterministic pathwise bound for $\|f(t)\|_2$.  Before mass preservation
has been proved, the estimate
$\|V\star_x\langle\alpha_r^2f\rangle\|_\infty
\leq C_r\|V\|_2\|f\|_2$ would therefore produce an uncontrolled mixed
moment when inserted directly into the $L^1$ cutoff argument.  To avoid
using mass preservation circularly, fix $N\geq1$ and introduce the
$L^2$-stopping time
\begin{align}\label{L2-stopping-L1-proof}
 \vartheta_N:=\inf\left\{t\in[0,T]:\|f(t)\|_{L^2}>N\right\}\wedge T.
\end{align}
The fixed-$r$ estimate \eqref{fixed-r-interaction-bounds} gives, on
$[0,\vartheta_N]$,
\begin{align*}
 \|V\star_x\langle\alpha_r^2f\rangle\|_{L_x^\infty}
 \leq C_r\|V\|_{L_x^2}\|f\|_{L^2}\leq C_{r,V}N.
\end{align*}
Moreover, for fixed $n,r$, the estimate
\eqref{fixed-n-r-Galerkin-estimate}, passed to the limit by lower
semicontinuity, gives
\begin{align*}
 \mathbb E\sup_{t\leq T}\|f(t)\|_2^2
 \leq C_{n,r,T}\mathbb E\|f_0\|_2^2.
\end{align*}
Consequently,
$\mathbb P(\vartheta_N<T)\leq
N^{-2}\mathbb E\sup_{t\leq T}\|f(t)\|_2^2\to0$, and hence
$\vartheta_N\uparrow T$ almost surely.  The stopping is nevertheless
needed because it turns the random coefficient $C_r\|f(t)\|_2$ into the
deterministic bound $C_rN$ required by the $L^1$ cutoff iteration.  Until
the removal of this stopping time below, every identity is applied at
$t\wedge\vartheta_N$; to avoid overloading the notation we continue to
write $t$.  In particular, all
time integrals and suprema in \eqref{0422:05}--\eqref{0422:07} are
understood in this stopped sense.

For the drift term, put
$c_r(f):=\alpha_r^2v-V\star_x\langle\alpha_r^2f\rangle$.
Using \eqref{0422:-1} and integration by parts, we compute}
\begin{align*}
&{\color{black}\left|\int_0^t \int_{\mathbb{R}^{2d}} a_{\delta}'(f)
 \nabla_v \cdot (c_r(f) f) \alpha_R\right|}\notag\\
&{\color{black}\leq \int_0^t \int_{\mathbb{R}^{2d}} a_{\delta}''(f)
 |\nabla_v f|\,|c_r(f)| f \alpha_R
 + \int_0^t \int_{\mathbb{R}^{2d}} a_{\delta}'(f) f
 |c_r(f)| |\nabla_v \alpha_R|}\notag\\
&{\color{black}\leq (C_{r,V}N+R) \int_0^t \int_{\mathbb{R}^{2d}} (1_{\{0<|f| \le \sqrt{\delta}\}} + \sqrt{\delta} |f|)|\nabla_v f|\alpha_R} \notag
\\
&{\color{black}\quad+  \int_0^t \int_{\mathbb{R}^{2d}} |f| \alpha_{2R}\left(\frac{C_{r,V}N}{R}+CI_{\{R\le|v|\le 2R\}}\right),} 
\end{align*}
where we used the fact $|v \cdot \nabla_v \alpha_R| \lesssim 1_{\{R \le |v| \le 2R\}}$. Thus, taking expectations, thanks to the $L^2(\mathbb{R}^{2d})$-estimates, we have
\begin{align}\label{0422:05}
     \mathbb{E}I_{b,\delta}(f,t)\lesssim (R+C_{r,V}N)\mathbb{E} \int_0^t \int_{\mathbb{R}^{2d}} (1_{\{0<|f| \le \sqrt{\delta}\}} + \sqrt{\delta} |f|)|\nabla_v f|\alpha_R
+  \mathbb{E}\int_0^t \int_{\mathbb{R}^{2d}} |f| \alpha_{2R}\left(\frac{C_{r,V}N}{R}+I_{\{R\le|v|\le 2R\}}\right).
\end{align}

For the martingale term $I_W(f)$, by applying the chain rule and integration by parts, we obtain
\begin{align*}
I_{W,\delta}(f,t)=& - \int_0^t \int_{\mathbb{R}^{2d}} \alpha_R a_{\delta}'(f) \nabla_v \sigma(f) \cdot dW_F 
- \int_0^t \int_{\mathbb{R}^{2d}} \alpha_R a_{\delta}'(f) \sigma(f) \nabla_v \cdot dW_F \\
=& - \int_0^t \int_{\mathbb{R}^{2d}} \alpha_R \nabla_v \left( \int_0^f a_{\delta}'(\zeta)\sigma'(\zeta)\,d\zeta \right) \cdot dW_F 
- \int_0^t \int_{\mathbb{R}^{2d}} \alpha_R a_{\delta}'(f)\sigma(f)\nabla_v \cdot dW_F \\
=& \int_0^t \int_{\mathbb{R}^{2d}} \nabla_v \alpha_R \left( \int_0^f a_{\delta}'(\zeta)\sigma'(\zeta)\,d\zeta \right) \cdot dW_F \\
&+ \int_0^t \int_{\mathbb{R}^{2d}} \alpha_R \left( \int_0^f a_{\delta}'(\zeta)\sigma'(\zeta)\,d\zeta - a_{\delta}'(f)\sigma(f) \right)\nabla_v \cdot dW_F.
\end{align*}
{\color{black}Set
\begin{align*}
 J_\delta(a):=\int_0^a a_\delta'(\zeta)\sigma'(\zeta)\,d\zeta,
 \qquad K_\delta(a):=J_\delta(a)-a_\delta'(a)\sigma(a).
\end{align*}
The quadratic variation contains squares of the spatial integrals.  By
Minkowski's inequality in $\ell^2$, followed by Cauchy--Schwarz for the
second spatial integral,
\begin{align*}
&\mathbb E\sup_{t\in[0,T]}|I_{W,\delta}(f,t)|\\
&\quad\leq C\|F_1\|_\infty^{1/2}
 \mathbb E\!\left[\left\{\int_0^T
 \left(\int_{\mathbb R^{2d}}|\nabla_v\alpha_R|\,
 |J_\delta(f)|\,dz\right)^2dt\right\}^{1/2}\right]\notag\\
&\qquad+C\|F_3\|_{L^1}^{1/2}
 \left(\mathbb E\int_0^T\!\int_{\mathbb R^{2d}}
 \alpha_R^2|K_\delta(f)|^2\,dz\,dt\right)^{1/2}.
\end{align*}
Here the expectation in the first term remains outside the square root;
this avoids replacing the square of a spatial integral by an integral of
pointwise squares.  Since $|J_\delta(a)|+|K_\delta(a)|\leq C|a|$, the
first term is at most
\begin{align*}
 \frac{C\sqrt T}{R}\,
 \mathbb E\sup_{t\in[0,T]}
 \int_{\mathbb R^{2d}}|f(t)|\alpha_{2R},
\end{align*}
whereas the second term tends to zero as $\delta\downarrow0$ by dominated
convergence and the $L^2$ estimate.}
Since
$$
\lim_{\delta \to 0} \left( \int_0^f a_{\delta}'(\zeta)\sigma'(\zeta)\,d\zeta - a_{\delta}'(f)\sigma(f) \right) 
= \int_0^f \sgn(\zeta)\sigma'(\zeta)\,d\zeta - \sgn(f)\sigma(f) = 0,
$$
we obtain the convergence estimate
\begin{align}
{\color{black}\limsup_{\delta\to0}\mathbb E\sup_{t\in[0,T]}
 |I_{W,\delta}(f,t)|
 \leq \frac{C\sqrt T}{R}\,
 \mathbb E\sup_{t\in[0,T]}
 \int_{\mathbb R^{2d}}|f(t)|\alpha_{2R}.}\label{0422:07}
\end{align}
We now estimate the term $I_{1,\delta}(f,t)$. By integration by parts and the standing identity ${\color{black}\nabla_vF_1=0}$, one sees that
\begin{align*}
    I_{1,\delta}(f,t)=&\frac12\left(\int_0^t \int_{\mR^{2d}}-(\sigma'(f)^2\nabla_v f)\nabla_v(a'_\delta(f)F_1\alpha_R)+a''_\delta(f)|\nabla_v\sigma(f)|^2F_1\alpha_R\right)\\
    =&\frac12\left(\int_0^t \int_{\mR^{2d}}-\sigma'(f)^2|\nabla_v f|^2 a''_\delta(f)F_1\alpha_R+a_\delta''(f)|\nabla_v\sigma(f)|^2F_1\alpha_R\right)\\
    &-\frac12\int_0^t \int_{\mR^{2d}}(\sigma'(f)^2\nabla_v f)\cdot \nabla_vF_1 a'_\delta(f)\alpha_R-\frac12\int_0^t \int_{\mR^{2d}}(\sigma'(f)^2\nabla_v f)\cdot \nabla_v\alpha_R a'_\delta(f)F_1\\
    =&{\color{black}-\frac12\int_0^t \int_{\mR^{2d}}(\sigma'(f)^2\nabla_v f)\cdot \nabla_v\alpha_R a'_\delta(f)F_1},
\end{align*}
where we noticed that $\sigma'(f)^2|\nabla_v f|^2 =|\nabla_v\sigma(f)|^2 $. Since
$$
a'_\delta(f)\sigma'(f)^2\nabla_v f=\nabla_v\left(\int_0^f a'_\delta(\zeta)\sigma'(\zeta)^2\dif \zeta\right), 
$$ we have
\begin{align}
    I_{1,\delta}(f,t)
    &={\color{black}\frac12\int_0^t \int_{\mR^{2d}}\left(\int_0^f a'_\delta(\zeta)\sigma'(\zeta)^2\dif \zeta\right) F_1\Delta_v\alpha_R}\notag\\
    &{\color{black}\le C(\sigma,F_1)R^{-1}\int_0^t \int_{\mR^{2d}}|f|\alpha_{2R}.} \label{est:I_1}
\end{align}
For the term $I_{2,3,\delta}(f,t)$, the assumptions on $\sigma$ and \eqref{0422:-1} give

{\color{black}\begin{align}\label{est:I_23}
 I_{2,3,\delta}(f,t)
 \leq {\color{black}C(\sigma,F_3)\int_0^t\!\int_{\mathbb R^{2d}}
 \bigl(1_{\{0<|f|\leq\sqrt\delta\}}+\sqrt\delta|f|\bigr)
 \alpha_R\,dz\,ds}.
\end{align}}
{\color{black}Letting $\delta\to0$ in \eqref{0422:05} and
\eqref{est:I_23}, and using \eqref{est:I_1}, we obtain, for fixed $N$
and $R>1$,
\begin{align*}
 &\limsup_{\delta\to0}\mathbb E\left[
 \sup_{t\in[0,T]}\bigl(
 I_{b,\delta}(f,t)+I_{1,\delta}(f,t)+I_{2,3,\delta}(f,t)
 \bigr)\right]\\
 &\qquad\leq C_N\mathbb E\int_0^{T\wedge\vartheta_N}
 \!\int_{\mathbb R^{2d}}|f(s,z)|\alpha_{2R}(z)
 \left(R^{-1}+1_{\{R\leq|v|\leq2R\}}\right)dz\,ds,
\end{align*}
where $C_N$ may depend on $d,r,V,\sigma,F_1,F_3$ and $N$, but not
on $R$ or $\delta$.  Together with \eqref{0422:07}, this gives
\begin{align}\label{0422:12}
 J_{R,N}(t):={}&\mathbb E\left[\sup_{s\in[0,t\wedge\vartheta_N]}
 \int_{\mathbb R^{2d}}|f(s)|\alpha_R\right]\notag\\
 \leq{}& I_0+D_{R,N}(t)+\frac{C\sqrt T}{R}J_{2R,N}(t),
\end{align}
where
\begin{align*}
 I_{R,N}(t)&:=\mathbb E\int_{\mathbb R^{2d}}
 |f(t\wedge\vartheta_N)|\alpha_R,\qquad
 I_0:=\mathbb E\int_{\mathbb R^{2d}}|f_0|,\\
 D_{R,N}(t)&:=C_N\mathbb E\int_0^{t\wedge\vartheta_N}
 \!\int_{\mathbb R^{2d}}|f|\alpha_{2R}
 \left(R^{-1}+I_{\{R\leq|v|\leq2R\}}\right).
\end{align*}
In particular,
$D_{R,N}(t)\leq C_N\int_0^t I_{2R,N}(s)\,ds$.  Repeating the calculation
at a fixed stopped time makes the martingale mean zero, and the same
Volterra iteration gives
\begin{align*}
 \sup_{R>1}\sup_{t\in[0,T]}I_{R,N}(t)
 \leq C_{T,N}\left(I_0+
 \|f\|_{L^2(\Omega\times[0,T];L^2)}\right).
\end{align*}
Substitution in \eqref{0422:12} yields
\begin{align*}
 J_{R,N}(T)\leq C_{T,N}+\frac{C_T}{R}J_{2R,N}(T).
\end{align*}
On $\{\vartheta_N>0\}$, Cauchy--Schwarz on the support of $\alpha_R$
gives the bound $C N(2^mR)^{3d/2}$, while on
$\{\vartheta_N=0\}$ the stopped value is the initial datum and is bounded
by $M_*$.  Hence
$J_{2^mR,N}(T)\leq C N(2^mR)^{3d/2}+M_*$.  Iteration leaves the remainder
\begin{align*}
 \frac{C_T^m}{R^m2^{m(m-1)/2}}J_{2^mR,N}(T),
\end{align*}
which tends to zero as $m\to\infty$.  Hence
$\sup_{R>1}J_{R,N}(T)<\infty$.  Monotone convergence first gives the
stopped version of \eqref{est:L^1_finite}; it then implies
$D_{R,N}(t)\to0$ by dominated convergence.  Letting $R\to\infty$ in
\eqref{0422:12}, we find
\begin{align*}
 \mathbb E\left[\sup_{s\in[0,t\wedge\vartheta_N]}
 \int_{\mathbb R^{2d}}|f(s)|\right]\leq
 \mathbb E\int_{\mathbb R^{2d}}|f_0|.
\end{align*}
Finally, $\vartheta_N\uparrow T$ almost surely.  Fatou's lemma, followed
by $N\to\infty$, proves both \eqref{est:L^1_finite} and
\begin{align}\label{0422:13}
 \mathbb E\left[\sup_{s\in[0,t]}
 \int_{\mathbb R^{2d}}|f(s)|\right]\leq
 \mathbb E\int_{\mathbb R^{2d}}|f_0|.
\end{align}}
Now we establish \eqref{preservationofmean}. By the definition of a weak solution, taking $\alpha_R$ as a test function, we obtain
{\color{black}
\begin{align*}
			\int_{\mathbb{R}^{2d}}f(t)\alpha_R =& \int_{\mathbb{R}^{2d}}f_0\alpha_R
 - \int_0^t \int_{\mathbb{R}^{2d}} \nabla_v f \cdot \nabla_v \alpha_R
 -\frac1n\int_0^t\int_{\mathbb R^{2d}}\nabla_xf\cdot\nabla_x\alpha_R\\
			&+\int_0^t\int_{\mathbb R^{2d}}\alpha_r^2v f\cdot\nabla_x\alpha_R
			-\int_0^t\int_{\mathbb R^{2d}}\alpha_r^2v f\cdot\nabla_v\alpha_R\\
			&+\int_0^t\int_{\mathbb R^{2d}}f
			\bigl(V\star_x\langle\alpha_r^2f\rangle\bigr)\cdot\nabla_v\alpha_R \\
		&+ \int_0^t \int_{\mathbb{R}^{2d}} \nabla_v \alpha_R \cdot (\sigma(f) \, dW_F)\\
			&\quad{\color{black}-\frac12\int_0^t\int_{\mathbb R^{2d}}
			 F_1\sigma'(f)^2\nabla_vf\cdot\nabla_v\alpha_R}.
\end{align*}}

{\color{black}Using estimate \eqref{0422:13} and spatial dominated
convergence, on one event of probability one the endpoint terms satisfy,
for every fixed $t\in[0,T]$,}
$$
		\int_{\mathbb{R}^{2d}} f(t) \alpha_R \to \int_{\mathbb{R}^{2d}} f(t), \quad \text{and} \quad \int_{\mathbb{R}^{2d}} f_0 \alpha_R \to \int_{\mathbb{R}^{2d}} f_0.
$$
{\color{black}Moreover, we integrate the two diffusion terms by parts once
more.  Estimate \eqref{est:L^1_finite}, the cutoff bounds, and dominated
convergence then give}
{\color{black}
\begin{align*}
		\int_0^t\int_{\mathbb R^{2d}}f
		\left(\Delta_v\alpha_R+\frac1n\Delta_x\alpha_R\right)\to0,\quad
		\int_0^t\int_{\mathbb R^{2d}}\alpha_r^2v f\cdot\nabla_x\alpha_R\to0,\quad
		\int_0^t\int_{\mathbb R^{2d}}\alpha_r^2v f\cdot\nabla_v\alpha_R\to0,\\
		\int_0^t\int_{\mathbb R^{2d}}f
		\bigl(V\star_x\langle\alpha_r^2f\rangle\bigr)\cdot\nabla_v\alpha_R\to0.
\end{align*}
uniformly with respect to $t\in[0,T]$ as $R\to\infty$, because every
time integral is bounded by the corresponding absolute integral over
$[0,T]$.}

{\color{black}For the martingale term, the Burkholder--Davis--Gundy
inequality gives control of the supremum in time.  Since
$|\sigma(r)|\leq\|\sigma'\|_\infty |r|$, define
\begin{align*}
 \tau_N:=\inf\{s\in[0,T]:\|f(s)\|_{L^1}>N\}\wedge T.
\end{align*}
Then
\begin{align*}
&\mathbb E\sup_{t\in[0,T]}\left|
 \int_0^{t\wedge\tau_N}\int_{\mathbb R^{2d}}\nabla_v\alpha_R\cdot
 \bigl(\sigma(f)\,dW_F\bigr)\right|^2\\
&\quad\leq C\mathbb E\int_0^{T\wedge\tau_N}\sum_{k\geq1}
 \left|\int_{\mathbb R^{2d}}\nabla_v\alpha_R(z)
 \sigma(f(s,z))f_k(z)\,dz\right|^2ds\\
&\quad\leq \|F_1\|_\infty\mathbb E\int_0^{T\wedge\tau_N}
 \left(\int_{\mathbb R^{2d}}|\nabla_v\alpha_R|\,|\sigma(f)|\,dz\right)^2ds
\leq \frac{CTN^2}{R^2}.
\end{align*}
Choose a deterministic sequence $R_m\uparrow\infty$ such that
$\sum_mR_m^{-2}<\infty$.  For each integer $N$, Chebyshev's inequality
and Borel--Cantelli imply
\begin{align*}
 \sup_{t\in[0,T]}\left|
 \int_0^{t\wedge\tau_N}\int_{\mathbb R^{2d}}\nabla_v\alpha_{R_m}\cdot
 \bigl(\sigma(f)\,dW_F\bigr)\right|\longrightarrow0
 \quad\text{almost surely}.
\end{align*}
By \eqref{est:L^1_finite}, for almost every $\omega$ there is an integer
$N=N(\omega)$ for which $\tau_N(\omega)=T$.  Therefore the unstopped
martingale converges to zero almost surely in $C([0,T])$ along $(R_m)$.}

{\color{black}It remains to treat the It\^o correction. Define
\begin{align*}
 A(r):=\int_0^r\sigma'(s)^2\,ds.
\end{align*}
Then $|A(r)|\leq C r$ for $r\geq0$. Using $\nabla_vF_1=0$, its contribution equals
\begin{align*}
&-\frac12\int_0^t\int F_1\nabla_vA(f)\cdot\nabla_v\alpha_R
={\color{black}\frac12\int_0^t\int A(f)F_1\Delta_v\alpha_R}
 \longrightarrow0
 \quad\text{uniformly for }t\in[0,T]\text{ along }(R_m),
\end{align*}
where all unlabelled integrals are over $\mathbb R^{2d}$.  The last
limit follows from the $L^1$ bound and the cutoff estimates.  The weak
identity holds for all $t$ on a common event, the deterministic and
stochastic remainders converge uniformly in $t$, and the spatial endpoint
limits hold for each $t$ along the same sequence.  Hence
\eqref{preservationofmean} holds on one event of probability one,
simultaneously for every $t\in[0,T]$.}

To establish the non-negativity of $f$, define the negative part $f^- := -(f \wedge 0)$, and observe the identity $f^- = \frac{1}{2}(|f| - f)$. Then, using \eqref{preservationofmean}, we compute
$$
	\int_{\mathbb{R}^{2d}} f^-(t) = \frac{1}{2} \int_{\mathbb{R}^{2d}} |f(t)| - \frac{1}{2} \int_{\mathbb{R}^{2d}} f(t) = \frac{1}{2} \int_{\mathbb{R}^{2d}} |f(t)| - \frac{1}{2} \int_{\mathbb{R}^{2d}} f_0.
$$
Hence, applying \eqref{0422:13}, we find
$$
	\mathbb{E} \left( \sup_{t \in [0,T]} \int_{\mathbb{R}^{2d}} f^-(t) \right)
	= \frac{1}{2} \mathbb{E} \left( \sup_{t \in [0,T]} \int_{\mathbb{R}^{2d}} |f(t)| - \int_{\mathbb{R}^{2d}} f_0 \right)
	\leq \frac{1}{2} \left( \mathbb{E} \int_{\mathbb{R}^{2d}} f_0 - \mathbb{E} \int_{\mathbb{R}^{2d}} f_0 \right) = 0.
$$
This implies that $f$ is almost surely nonnegative for all $t \in [0,T]$.
\end{proof}

In the following, we establish a moment estimate that is uniform in $n, r \geq 1$. 

\begin{proposition}\label{moment-es}
{\color{black}Assume that $f_0$, $V$ and $\sigma(\cdot)$ satisfy Assumptions~\ref{Assump-initialdata-regular}, \ref{Assump-ker}, and \ref{Assump-sigma-n}, respectively, and let $M_*$ be defined in \eqref{regular-data-size}. For every $n \in \mathbb{N}$, let $f_{n,r}$ be the weak solution to \eqref{SPDE-iteration} with initial data $f_0$. Then the following uniform bound holds:}
\begin{align*}
	\sup_{n,r\geq1} 
	 \mathbb{E}\left(\sup_{t \in [0,T]} \int_{\mathbb{R}^{2d}}(|v|^2+|x|^2)f_{n,r}(t)\right) 
		{\color{black}\leq C(T,\sigma,V,F_1,F_3,M_*)
 \left(\mathbb{E}\int_{\mathbb{R}^{2d}}(|v|^2+|x|^2)f_0+M_*\right).}
\end{align*}
\end{proposition}
\begin{proof}
For each $R > 1$, recall that $\alpha_R(x,v) := \alpha^1_{R^3}(x)\alpha^2_{R}(v)$, as defined in \eqref{alpha-R} with $(R_1,R_2)=(R^3,R)$. Taking $|v|^2\alpha_R$ as a test function, we deduce that, almost surely, for every $t\in[0,T]$,
\begin{align*}
\int_{\mathbb{R}^{2d}}|v|^2\alpha_R f_{n,r}(t)
=&\int_{\mathbb{R}^{2d}}|v|^2\alpha_R f_0
-\int^t_0\int_{\mathbb{R}^{2d}}\nabla_v(|v|^2\alpha_R)\cdot\nabla_v f_{n,r}-\frac{1}{n}\int^t_0\int_{\mathbb{R}^{2d}}\nabla_x(|v|^2\alpha_R)\cdot\nabla_x f_{n,r}\\
&+\int^t_0\int_{\mathbb{R}^{2d}}\nabla_x\alpha_R\cdot |v|^2\alpha^2_rv f_{n,r}
-\int^t_0\int_{\mathbb{R}^{2d}}\nabla_v(|v|^2\alpha_R)\cdot \alpha^2_rv f_{n,r}\\
&+\int^t_0\int_{\mathbb{R}^{2d}}\nabla_v(|v|^2\alpha_R)\cdot f_{n,r} V\star_x\langle\alpha^2_rf_{n,r}\rangle\\
&+\int^t_0\int_{\mathbb{R}^{2d}}\nabla_v(|v|^2\alpha_R)\cdot \sigma(f_{n,r})\, dW_F\\
&-\frac{1}{2}\int^t_0\int_{\mathbb{R}^{2d}}\nabla_v(|v|^2\alpha_R)\cdot \sigma'(f_{n,r})^2 \nabla_v f_{n,r} F_1\\
=:&\int_{\mathbb{R}^{2d}}|v|^2\alpha_R f_0+{\color{black}\sum_{i=1}^7 I_i}.
\end{align*}

In the following, we estimate $(I_i)_{1\leq i\leq 7}$ term by term.

For $I_1$, by the integration by parts formula and the chain rule,
\begin{align*}
I_1=\int^t_0\int_{\mathbb{R}^{2d}}(2d\alpha_R+4v\cdot\nabla_v\alpha_R+|v|^2\Delta_v\alpha_R)f_{n,r}
{\color{black}\leq C_dTM_*.}
\end{align*}
Similarly, for $I_2$, we have that 
\begin{align*}
I_2=\frac{1}{n}\int^t_0\int_{\mathbb{R}^{2d}}(|v|^2\Delta_x\alpha_R)f_{n,r}
{\color{black}\lesssim\frac{M_*}{R^4}.}
\end{align*}
\noindent
For $I_3$, we obtain
\begin{align*}
I_3\leq \int^t_0\int_{\mathbb{R}^{2d}} f_{n,r}(s)
{\color{black}\leq TM_*.}
\end{align*}

\noindent
For $I_4$, using the chain rule,
\begin{align*}
I_4=-\int^t_0\int_{\mathbb{R}^{2d}}\alpha^2_r(2|v|^2\alpha_R+|v|^2\nabla_v\alpha_R\cdot v)f_{n,r}
\leq \int^t_0\int_{\mathbb{R}^{2d}}|v|^2\alpha_{2R}f_{n,r}.
\end{align*}

\noindent
For $I_5$, by the chain rule, H\"older's inequality, and Young's inequality,
\begin{align*}
I_5
=&\int^t_0\int_{\mathbb{R}^{2d}}(2v\alpha_R+|v|^2\nabla_v\alpha_R)\cdot f_{n,r}V\star_x\langle\alpha^2_rf_{n,r}\rangle\\
\leq&{\color{black}\|V\|_{L^{\infty}(\mathbb{R}^d)}M_*}
\left(\int^t_0\int_{\mathbb{R}^{2d}}|v|^2\alpha_{2R}f_{n,r}\right)^{1/2}
\left(\int^t_0\int_{\mathbb{R}^{2d}}\alpha_{2R}f_{n,r}\right)^{1/2}\\
\leq&\int^t_0\int_{\mathbb{R}^{2d}}|v|^2\alpha_{2R}f_{n,r}+{\color{black}C(V,M_*)}.
\end{align*}

\noindent
{\color{black}For the martingale term $I_6$, set
$\varphi_R(v,x)=|v|^2\alpha_R(x,v)$. Its quadratic variation is
\begin{align*}
 \langle I_6\rangle_t
 =\int_0^t\sum_{k\geq1}
 \left|\int_{\mathbb R^{2d}}\nabla_v\varphi_R(z)
 \sigma(f_{n,r}(s,z))f_k(z)\,dz\right|^2ds.
\end{align*}
Thus Minkowski's inequality in $\ell^2$, the bound
$|\sigma(a)|\leq\|\sigma'\|_\infty a$ for $a\geq0$, and mass
conservation imply
\begin{align*}
 \sum_{k\geq1}\left|\int\nabla_v\varphi_R\,\sigma(f_{n,r})f_k\right|^2
 &\leq \|F_1\|_\infty
 \left(\int|\nabla_v\varphi_R|\,|\sigma(f_{n,r})|\right)^2\\
 &\leq C\left(\int |v|\alpha_{2R}f_{n,r}\right)^2\\
 &{\color{black}\leq CM_*}
 \int |v|^2\alpha_{2R}f_{n,r}.
\end{align*}
The Burkholder--Davis--Gundy inequality and Young's inequality therefore
give
\begin{align*}
 \mathbb E\sup_{t\in[0,T]}|I_6(t)|
 {\color{black}\leq CM_*^{1/2}}
 \left(\mathbb E\int_0^T\!\int
 |v|^2\alpha_{2R}f_{n,r}\right)^{1/2}
 {\color{black}\leq C(M_*)\int_0^T\mathbb E\int|v|^2\alpha_{2R}f_{n,r}(s)\,ds+C(M_*).}
\end{align*}}

\noindent
For $I_7$, using the integration by parts formula, the chain rule, H\"older's inequality, and Young's inequality,
\begin{align*}
I_7
=&\frac{1}{2}\int^t_0\int_{\mathbb{R}^{2d}}\nabla_v\cdot(2v\alpha_R+|v|^2\nabla_v\alpha_R)
\left(\int^{f_{n,r}}_0\sigma'(\zeta)^2 d\zeta\right)F_1\\
\leq& C(\sigma)\int^t_0\int_{\mathbb{R}^{2d}}(2d\alpha_R+4v\cdot\nabla_v\alpha_R+|v|^2\Delta_v\alpha_R)f_{n,r}F_1\\
\leq& {\color{black}C(\sigma,F_1)\int^t_0\int_{\mathbb{R}^{2d}}|v|^2\alpha_{2R}f_{n,r}+C(M_*)}.
\end{align*}

\medskip

Combining the above estimates, we conclude that
\begin{align*}
\mathbb{E}\left(\sup_{t\in[0,T]}\int_{\mathbb{R}^{2d}}|v|^2\alpha_R f_{n,r}(t)\right)
\leq& \mathbb{E}\int_{\mathbb{R}^{2d}}|v|^2 f_0\\
&+{\color{black}C(M_*,F_1,F_3,\sigma)}\mathbb{E}\int^T_0\int_{\mathbb{R}^{2d}}|v|^2\alpha_R f_{n,r}(t)
+{\color{black}C(M_*,V)}.
\end{align*}

Let
\begin{align*}
I_R(t):=\mathbb{E}\int_{\mathbb{R}^{2d}}|v|^2\alpha_R f_{n,r}(t),
\qquad
I_0:=\mathbb{E}\int_{\mathbb{R}^{2d}}|v|^2 f_0.
\end{align*}
Then
\begin{align*}
I_R(t)\leq I_0+{\color{black}C(M_*,F_1,F_3)}\int^t_0 I_{2R}(s)\,ds + {\color{black}C(M_*,V)}.
\end{align*}

By iterating this estimate as in the proof of Lemma \ref{lem-L1-es}, we obtain
\begin{align*}
\mathbb{E}\left(\sup_{t\in[0,T]}\int_{\mathbb{R}^{2d}}|v|^2 f_{n,r}(t)\right)
{\color{black}\leq C(T,\sigma,V,F_1,F_3,M_*)}
\Big(\mathbb{E}\int_{\mathbb{R}^{2d}}(|v|^2+|x|^2)f_0
+{\color{black}M_*}\Big).
\end{align*}

Finally, we take $\alpha_R(x,v)=\alpha_{R_1}^1(x)\alpha_{R_2}^2(v)$ with independent parameters $R_1,R_2>1$. Using {\color{black}$|x|^2\alpha_R$} as a test function and proceeding as above, then letting $R_2\to\infty$, we deduce that
\begin{align*}
\mathbb{E}\left(\int_{\mathbb{R}^{2d}}\alpha^1_{R_1}|x|^2f_{n,r}(t)\right)
\leq& \mathbb{E}\int_{\mathbb{R}^{2d}}|x|^2f_0
+\frac{1}{n}\mathbb{E}\int^t_0\int_{\mathbb{R}^{2d}}\Delta_x(\alpha^1_{R_1}|x|^2)f_{n,r}+\mathbb{E}\int^t_0\int_{\mathbb{R}^{2d}}\nabla_x(|x|^2\alpha_{R_1}^1)\cdot \alpha^2_rv f_{n,r}\\
\leq& \mathbb{E}\int_{\mathbb{R}^{2d}}|x|^2f_0
+\frac{1}{n}\mathbb{E}\int^t_0\int_{\mathbb{R}^{2d}}f_{n,r}+\mathbb{E}\int^t_0\int_{\mathbb{R}^{2d}}\left(x\alpha_{R_1}^1+|x|^2\nabla_x\alpha_{R_1}^1\right)\cdot \alpha^2_rv f_{n,r}\\
\leq& \mathbb{E}\int_{\mathbb{R}^{2d}}|x|^2f_0
+\mathbb{E}\int^t_0\int_{\mathbb{R}^{2d}}\alpha^1_{2R_1}|x||v| f_{n,r}+{\color{black}M_*T}\\
\lesssim& \mathbb{E}\int_{\mathbb{R}^{2d}}|x|^2f_0
+\mathbb{E}\int^t_0\int_{\mathbb{R}^{2d}}|v|^2 f_{n,r}+\mathbb{E}\int^t_0\int_{\mathbb{R}^{2d}}\alpha^1_{2R_1}|x|^2 f_{n,r}
+{\color{black}M_*T}.
\end{align*}
By iterating this estimate as in the proof of Lemma \ref{lem-L1-es}, we obtain
\begin{align*}
\mathbb{E}\left(\sup_{t\in[0,T]}\int_{\mathbb{R}^{2d}}|x|^2 f_{n,r}(t)\right)
{\color{black}\leq C(T,\sigma,V,F_1,F_3,M_*)}
\Big(\mathbb{E}\int_{\mathbb{R}^{2d}}(|v|^2+|x|^2)f_0
+{\color{black}M_*}\Big).
\end{align*}
This completes the proof.
\end{proof}

{\color{black}
\begin{proposition}[Uniform entropy estimate for the approximation]
\label{approximate-entropy-dissipation}
Under the assumptions of Proposition~\ref{prp-uniforml2}, every weak
solution $f_{n,r}$ of \eqref{SPDE-iteration} satisfies
\begin{align}\label{approximate-entropy-bound}
 \sup_{n,r\geq1}\Bigg\{
 &\mathbb E\sup_{t\in[0,T]}
   \left|\int_{\mathbb R^{2d}}\Psi(f_{n,r}(t,z))\,dz\right|
 +\mathbb E\int_0^T\!\int_{\mathbb R^{2d}}
       \frac{|\nabla_vf_{n,r}|^2}{f_{n,r}}\,dz\,dt\notag\\
 &+\frac1n\mathbb E\int_0^T\!\int_{\mathbb R^{2d}}
       \frac{|\nabla_xf_{n,r}|^2}{f_{n,r}}\,dz\,dt\Bigg\}<\infty .
\end{align}
The quotients are defined to be zero on $\{f_{n,r}=0\}$.  The constant
is uniform in the artificial parameters $n,r$, but it depends on the
fixed regular coefficient $\sigma$.  In particular, this proposition
does not give a uniform entropy estimate for the sequence
$(\sigma_n)_{n\geq1}$ approximating the square root in
Section~\ref{sec-6}.
\end{proposition}

\begin{proof}
We give a term-by-term entropy computation for the present approximation.
Fix $n,r\geq1$ and write $f:=f_{n,r}$.  For
$\delta\in(0,1)$, set
\begin{align*}
 \Psi_\delta(a)&:=\int_0^a\log(s+\delta)\,ds,
 &\psi_\delta(a)&:=\Psi_\delta'(a)=\log(a+\delta),\\
 H_\delta(a)&:=a\psi_\delta(a)-\Psi_\delta(a)
 =a-\delta\log\left(1+\frac a\delta\right),
 &S_\delta(a)&:=\int_0^a\frac{\sigma(s)}{s+\delta}\,ds.
\end{align*}
Thus $\Psi_\delta''(a)=(a+\delta)^{-1}$ and
$0\leq H_\delta(a)\leq a$ for $a\geq0$.

We first make the It\^o calculation rigorous by spatial mollification, as
in the argument preceding \eqref{passtolimit-etagamma}.  Let
$(\eta_\gamma)_{\gamma\in(0,1)}$ be standard convolution kernels on
$\mathbb R_x^d$ and set
\begin{align*}
 f^\gamma:=\eta_\gamma\ast_x f,
 \qquad
 \alpha_R(x,v):=\alpha^1_{R^3}(x)\alpha^2_R(v).
\end{align*}
We convolve \eqref{SPDE-iteration} in $x$ and apply the generalized
variational It\^o formula to
\begin{align*}
 \int_{\mathbb R^{2d}}\Psi_\delta(f^\gamma(t,z))\alpha_R(z)\,dz.
\end{align*}
This does not require a truncation of the values of $f$: for fixed
$\delta>0$, $\Psi_\delta''(a)=(a+\delta)^{-1}$ is bounded on
$[0,\infty)$ and $\Psi_\delta'$ has at most linear growth.  When applying
the formula, we extend $\Psi_\delta$ smoothly to $\mathbb R$ so that its
second derivative remains bounded.  The choice of extension is immaterial,
since $f^\gamma\geq0$ because $f\geq0$ and the mollifier is nonnegative.

We then let $\gamma\downarrow0$.  For fixed $n,r,\delta$, the energy
bounds imply, locally in phase space and along a subsequence,
\begin{align*}
 f^\gamma&\longrightarrow f,&
 \nabla_vf^\gamma&\longrightarrow\nabla_vf,&
 \nabla_xf^\gamma&\longrightarrow\nabla_xf
 \quad\text{in }L^2(\Omega\times(0,T);L^2_{\rm loc})
 \quad\text{and a.e.}
\end{align*}
Consequently all deterministic terms converge in their variational
dualities.  For the stochastic terms, the regular-coefficient assumptions
give
$(\nabla_v(f_k\sigma(f)))_{k\geq1}\in
L^2(\Omega\times(0,T)\times K;\ell^2)$ for every compact
$K\Subset\mathbb R^{2d}$.  Strong continuity of convolution in this
Hilbert-valued space yields
\begin{align*}
 \sum_{k\geq1}\left|\eta_\gamma\ast_x
 \nabla_v\bigl(f_k\sigma(f)\bigr)\right|^2
 \longrightarrow \sum_{k\geq1}
 \left|\nabla_v\bigl(f_k\sigma(f)\bigr)\right|^2
 \quad\text{in }L^1(\Omega\times(0,T)\times K).
\end{align*}
It\^o's isometry and the Burkholder--Davis--Gundy inequality likewise
give convergence of the stochastic integrals in
$L^1(\Omega;C([0,T]))$.  Thus, exactly as in the passage to
\eqref{passtolimit-etagamma}, the limit $\gamma\downarrow0$ gives the
localized It\^o identity below directly with $\Psi_\delta$ and
$\psi_\delta$.  Namely, almost surely, for every $t\in[0,T]$,
\begin{align}
 &\int_{\mathbb R^{2d}}\Psi_\delta(f(t))\alpha_R
 +\int_0^t\!\int_{\mathbb R^{2d}}
       \frac{|\nabla_vf|^2}{f+\delta}\alpha_R
 +\frac1n\int_0^t\!\int_{\mathbb R^{2d}}
       \frac{|\nabla_xf|^2}{f+\delta}\alpha_R\notag\\
 &={}
 \int_{\mathbb R^{2d}}\Psi_\delta(f_0)\alpha_R
 +\int_0^t\!\int_{\mathbb R^{2d}}
       \Psi_\delta(f)\left(\Delta_v\alpha_R
       +\frac1n\Delta_x\alpha_R\right)\notag\\
 &\quad+\int_0^t\!\int_{\mathbb R^{2d}}
       \Psi_\delta(f)\alpha_r^2v\cdot\nabla_x\alpha_R\notag\\
 &\quad+\int_0^t\!\int_{\mathbb R^{2d}}
       H_\delta(f)\nabla_v\!\cdot(\alpha_r^2v)\alpha_R
 -\int_0^t\!\int_{\mathbb R^{2d}}
       \Psi_\delta(f)\alpha_r^2v\cdot\nabla_v\alpha_R\notag\\
 &\quad+\int_0^t\!\int_{\mathbb R^{2d}}
       \Psi_\delta(f)
       \bigl(V\star_x\langle\alpha_r^2f\rangle\bigr)
       \cdot\nabla_v\alpha_R\notag\\
 &\quad-\frac12\int_0^t\!\int_{\mathbb R^{2d}}
       \psi_\delta(f)F_1\sigma'(f)^2\nabla_vf\cdot\nabla_v\alpha_R\notag\\
 &\quad-\frac12\int_0^t\!\int_{\mathbb R^{2d}}
       \frac{F_1\sigma'(f)^2|\nabla_vf|^2}{f+\delta}\alpha_R
 +\frac12\int_0^t\!\int_{\mathbb R^{2d}}
       \frac{F_1\sigma'(f)^2|\nabla_vf|^2
       +F_3\sigma(f)^2}{f+\delta}\alpha_R\notag\\
 &\quad-\int_0^t\!\int_{\mathbb R^{2d}}
       \psi_\delta(f)\alpha_R
       \nabla_v\cdot\bigl(\sigma(f)\,dW_F\bigr).
 \label{approximate-entropy-localized-Ito}
\end{align}
Here and below, all unlabelled space--time integrals are over
$[0,t]\times\mathbb R^{2d}$.  The penultimate line deliberately keeps the
interior $F_1$ contribution of the divergence-form It\^o correction
separate from the $F_1$ contribution of the quadratic variation; the two
terms cancel exactly.

The first two positive terms on the left-hand side are the entropy
dissipation produced by $\Delta_vf$ and $n^{-1}\Delta_xf$; their cutoff
contribution is the second term on the right-hand side.  Since
$\alpha_r^2v$ is independent of $x$, the truncated kinetic transport gives
the next term in \eqref{approximate-entropy-localized-Ito}.

For the velocity damping, the identity
\begin{align*}
 \psi_\delta(f)\nabla_v\cdot(\alpha_r^2vf)
 =\nabla_v\cdot\bigl(\alpha_r^2v\Psi_\delta(f)\bigr)
  +H_\delta(f)\nabla_v\cdot(\alpha_r^2v)
\end{align*}
gives, after multiplication by $\alpha_R$ and integration by parts, the
two corresponding terms in
\eqref{approximate-entropy-localized-Ito}.  Likewise,
$V\star_x\langle\alpha_r^2f\rangle$ is independent of $v$, and therefore
\begin{align*}
 -\psi_\delta(f)\nabla_v\cdot
 \bigl(fV\star_x\langle\alpha_r^2f\rangle\bigr)
 =-\nabla_v\cdot\left(
 \Psi_\delta(f)V\star_x\langle\alpha_r^2f\rangle\right).
\end{align*}
This yields the interaction cutoff term in
\eqref{approximate-entropy-localized-Ito}.

It remains to combine the It\^o correction with the quadratic variation
of the conservative noise.  The latter is
\begin{align*}
 \frac12\int_0^t\!\int_{\mathbb R^{2d}}
 \frac{\alpha_R}{f+\delta}
 {\color{black}\bigl(F_1\sigma'(f)^2|\nabla_vf|^2
 +F_3\sigma(f)^2\bigr)}\,dz\,ds.
\end{align*}
The $F_1$-term cancels exactly with the corresponding interior term from
the divergence-form It\^o correction in \eqref{SPDE-iteration}.  Hence only
the $F_3$-term and the $F_1$ cutoff term displayed in
\eqref{approximate-entropy-localized-Ito} remain.  To handle the latter,
set
\begin{align*}
 C_\delta(a):=\int_0^a\psi_\delta(s)\sigma'(s)^2\,ds.
\end{align*}
Then $\psi_\delta(f)\sigma'(f)^2\nabla_vf=\nabla_vC_\delta(f)$.
Since $\nabla_vF_1=0$, one more integration by parts gives
\begin{align*}
 -\frac12\int_0^t\!\int_{\mathbb R^{2d}}
 \psi_\delta(f)F_1\sigma'(f)^2\nabla_vf\cdot\nabla_v\alpha_R
 =\frac12\int_0^t\!\int_{\mathbb R^{2d}}
 C_\delta(f)F_1\Delta_v\alpha_R.
\end{align*}

The martingale in the cutoff identity is
\begin{align*}
 M_t^{n,r,\delta,R}
 &:=-\int_0^t\!\int_{\mathbb R^{2d}}
       \psi_\delta(f)\alpha_R\nabla_v\cdot(\sigma(f)\,dW_F).
\end{align*}
Integrating by parts and using $S_\delta'=\sigma/(\,\cdot+\delta)$ gives
\begin{align*}
 M_t^{n,r,\delta,R}
 =\sum_{k\geq1}\int_0^t\Bigg[&-\int_{\mathbb R^{2d}}
       \alpha_RS_\delta(f)\nabla_vf_k\,dz\\
 &+\int_{\mathbb R^{2d}}
       \bigl(\psi_\delta(f)\sigma(f)-S_\delta(f)\bigr)
       f_k\nabla_v\alpha_R\,dz\Bigg]\cdot dB_s^k.
\end{align*}

We next remove the spatial cutoff.  We first fix $(n,r,\delta)$.  For
fixed $\delta$, one has
\begin{align*}
 |\Psi_\delta(a)|+|C_\delta(a)|
 \leq C_\delta(a+a^2),\qquad a\geq0.
\end{align*}
The $L^1\cap L^2$ estimates, Proposition~\ref{moment-es}, mass
preservation, and the cutoff bounds therefore show that every deterministic
term in \eqref{approximate-entropy-localized-Ito} containing a derivative
of $\alpha_R$ converges to zero as $R\to\infty$.  For the transport term,
the velocity factor $\alpha_R^2(v)$ in the phase-space cutoff restricts
$|v|$ to $2R$, and hence
$|v||\nabla_x\alpha_R|\lesssim R^{-2}$.  We also use
\begin{align*}
 \sup_{r\geq1}\|\nabla_v\cdot(\alpha_r^2v)\|_\infty<\infty,
 \qquad
 \left\|V\star_x\langle\alpha_r^2f\rangle\right\|_\infty
 \leq\|V\|_\infty M_*.
\end{align*}
The velocity-damping cutoff term vanishes once $R$ is sufficiently large
relative to the fixed $r$.  After the cutoff has been removed, the two
bounds above show that all constants in the resulting whole-space identity
are independent of $n$ and $r$.

To justify the martingale limit, observe directly that the function
$a\mapsto\psi_\delta(a)\sigma(a)-S_\delta(a)$ vanishes at zero and has
derivative $\psi_\delta(a)\sigma'(a)$.  The regularity of the fixed
coefficient therefore gives
$|\psi_\delta(a)\sigma(a)-S_\delta(a)|\leq C_\delta a$.  Consequently,
\begin{align*}
 \left(\sum_{k\geq1}\left|\int
       \bigl(\psi_\delta(f)\sigma(f)-S_\delta(f)\bigr)f_k
       \nabla_v\alpha_R\,dz\right|^2\right)^{1/2}
 &\leq \|F_1\|_\infty^{1/2}
       \int \left|\psi_\delta(f)\sigma(f)-S_\delta(f)\right|
       |\nabla_v\alpha_R|\,dz\\
 &\leq C_\delta R^{-1}M_*.
\end{align*}
By the Burkholder--Davis--Gundy inequality, this boundary part of
$M^{n,r,\delta,R}$ converges to zero in $L^1(\Omega;C([0,T]))$.  For the
other part, Cauchy--Schwarz gives
\begin{align*}
 &\sum_{k\geq1}\left|\int(1-\alpha_R)S_\delta(f)
       \nabla_vf_k\,dz\right|^2\\
 &\qquad\leq
 \left(\int(1-\alpha_R)S_\delta(f)^2\,dz\right)
 \left(\int(1-\alpha_R)F_3\,dz\right)\longrightarrow0.
\end{align*}
Thus the first part converges, again by the
Burkholder--Davis--Gundy inequality, to
\begin{align*}
 M_t^{n,r,\delta}
 :=-\sum_{k\geq1}\int_0^t
 \left(\int_{\mathbb R^{2d}}S_\delta(f(s,z))\nabla_vf_k(z)\,dz\right)
 \cdot dB_s^k.
\end{align*}
We have thus proved the whole-space identity
\begin{align}\label{approximate-entropy-delta-identity}
 &\int_{\mathbb R^{2d}}\Psi_\delta(f(t))\,dz
 +\int_0^t\!\int_{\mathbb R^{2d}}
       \frac{|\nabla_vf|^2}{f+\delta}\,dz\,ds
 +\frac1n\int_0^t\!\int_{\mathbb R^{2d}}
       \frac{|\nabla_xf|^2}{f+\delta}\,dz\,ds\notag\\
 &=\int_{\mathbb R^{2d}}\Psi_\delta(f_0)\,dz
 +M_t^{n,r,\delta}
 +\frac12\int_0^t\!\int_{\mathbb R^{2d}}
       \frac{\sigma(f)^2}{f+\delta}F_3\,dz\,ds\notag\\
 &\quad+\int_0^t\!\int_{\mathbb R^{2d}}
       H_\delta(f)\nabla_v\cdot(\alpha_r^2v)\,dz\,ds.
\end{align}

The fixed structural bound $|\sigma(a)|\leq c\sqrt a$ implies, uniformly
in $\delta\in(0,1)$,
\begin{align*}
 |S_\delta(a)|
 \leq c\int_0^a\frac{\sqrt s}{s+\delta}\,ds
 \leq2c\sqrt a,
 \qquad
 0\leq\frac{\sigma(a)^2}{a+\delta}\leq c^2.
\end{align*}
By Cauchy--Schwarz and mass preservation,
\begin{align*}
 \frac d{dt}\langle M^{n,r,\delta}\rangle_t
 &=\sum_{k\geq1}\left|
 \int_{\mathbb R^{2d}}S_\delta(f)\nabla_vf_k\,dz\right|^2\\
 &\leq\left(\int_{\mathbb R^{2d}}S_\delta(f)^2\,dz\right)
       \left(\int_{\mathbb R^{2d}}F_3\,dz\right)
 \leq C M_*\|F_3\|_{L^1}.
\end{align*}
The Burkholder--Davis--Gundy inequality therefore gives
\begin{align*}
 \sup_{n,r\geq1}\sup_{\delta\in(0,1)}
 \mathbb E\sup_{t\in[0,T]}|M_t^{n,r,\delta}|<\infty.
\end{align*}
The $F_3$-term in \eqref{approximate-entropy-delta-identity} is bounded
by $CT\|F_3\|_{L^1}$, and the damping term is bounded by $CTM_*$ because
$0\leq H_\delta(f)\leq f$.

Finally, $\Psi_\delta\geq\Psi$, and hence
$(\Psi_\delta)^-\leq\Psi^-$.  Thus
\eqref{entropy-negative-part-bound} and Proposition~\ref{moment-es} give
\begin{align*}
 \sup_{n,r\geq1}\sup_{\delta\in(0,1)}
 \mathbb E\sup_{t\in[0,T]}
 \int_{\mathbb R^{2d}}\Psi_\delta(f_{n,r}(t))^-\,dz<\infty.
\end{align*}
The positive part of $\Psi_\delta(f_0)$ is bounded by a constant times
$f_0+f_0^2$.  Taking the supremum in time in
\eqref{approximate-entropy-delta-identity}, and then expectations, yields
\begin{align*}
 \sup_{n,r\geq1}\sup_{\delta\in(0,1)}\Bigg[
 &\mathbb E\sup_{t\in[0,T]}
 \left|\int_{\mathbb R^{2d}}\Psi_\delta(f_{n,r}(t))\,dz\right|\\
 &+\mathbb E\int_0^T\!\int_{\mathbb R^{2d}}
       \frac{|\nabla_vf_{n,r}|^2}{f_{n,r}+\delta}\,dz\,dt\\
 &+\frac1n\mathbb E\int_0^T\!\int_{\mathbb R^{2d}}
       \frac{|\nabla_xf_{n,r}|^2}{f_{n,r}+\delta}\,dz\,dt\Bigg]<\infty.
\end{align*}
Since $\Psi_\delta(a)\to\Psi(a)$,
$(\Psi_\delta)^-\leq\Psi^-$, and
$(\Psi_\delta)^+\leq C(a+a^2)$, the $L^1\cap L^2$ bounds give, for every
$t$, convergence of the spatial entropy integrals.  Therefore
\begin{align*}
 \sup_{t\in[0,T]}\left|\int\Psi(f(t))\,dz\right|
 \leq\liminf_{\delta\downarrow0}
 \sup_{t\in[0,T]}\left|\int\Psi_\delta(f(t))\,dz\right|.
\end{align*}
The two dissipation integrands increase as $\delta\downarrow0$.
Fatou's lemma and monotone convergence now give
\eqref{approximate-entropy-bound}.  For a nonnegative Sobolev
function, its weak gradient vanishes almost everywhere on its zero set,
so the stated convention for the quotients is consistent.  This completes
the proof.

\end{proof}}

{\color{black}
\begin{corollary}[Pathwise uniqueness for the approximation]
\label{fixed-n-r-pathwise-uniqueness}
For every fixed $n,r\geq1$, the weak solution of
\eqref{SPDE-iteration} is pathwise unique.  Consequently, the martingale
solution constructed in Proposition~\ref{prp-uniforml2} is realized as
the unique adapted weak solution on the original stochastic basis.
\end{corollary}

\begin{proof}
Let $f^1$ and $f^2$ be two weak solutions driven by the same Brownian
motions and with the same initial datum.  Lemma~\ref{lem-L1-es},
Proposition~\ref{moment-es}, and
Proposition~\ref{approximate-entropy-dissipation} were proved from the
weak equation and \eqref{0425:00}, and therefore apply separately to
both solutions.  In particular, both have the mass and Fisher-information
control required in the renormalized kinetic argument.

The proof of Proposition~\ref{equivalence} applies to
\eqref{SPDE-iteration}; the two extra elliptic terms contribute the
nonnegative kinetic defects
\begin{align*}
 \delta_{f^i=\zeta}|\nabla_vf^i|^2
 +n^{-1}\delta_{f^i=\zeta}|\nabla_xf^i|^2.
\end{align*}
We may therefore repeat the doubling argument in
Theorem~\ref{Uniqueness-spde}.  The $x$-Laplacian and its defect have
the favorable sign.  The fixed fields $\alpha_r^2v$ and
$\nabla_v(\alpha_r^2v)$ are bounded, so the corresponding transport and
damping commutators vanish with the mollification and cutoff parameters.
For the nonlocal term, $0\leq\alpha_r^2\leq1$ gives, exactly as in that
theorem,
\begin{align*}
 \|V\star_x(\rho_r(f^1)-\rho_r(f^2))\|_{L_x^\infty}
 \leq\|V\|_{L_x^\infty}\|f^1-f^2\|_{L^1_{x,v}},
\end{align*}
and the remaining factor is controlled by mass conservation and the
Fisher information.  After the same successive limits, one obtains,
almost surely,
\begin{align*}
 \|f^1(t)-f^2(t)\|_{L^1}
 \leq C_{n,r}\int_0^t
 \left(1+\|\nabla_v\sqrt{f^1(s)}\|_2
          +\|\nabla_v\sqrt{f^2(s)}\|_2\right)
 \|f^1(s)-f^2(s)\|_{L^1}\,\dif s.
\end{align*}
The coefficient is time integrable by
Proposition~\ref{approximate-entropy-dissipation}, so Gronwall's
inequality proves pathwise uniqueness.  Apply this conclusion
to every joint subsequential limit of two Galerkin approximations and
the common Brownian motions.  The Gy\"ongy--Krylov criterion then yields
convergence in probability on the original stochastic basis, whose
limit is the asserted adapted weak solution.
\end{proof}}

In the following, we provide a time-regularity estimate. 
\begin{lemma}\label{time-besov-regularity-fn}
Assume that $f_0$, $V$, and $\sigma(\cdot)$ satisfy
Assumptions~\ref{Assump-initialdata-regular}, \ref{Assump-ker}, and
\ref{Assump-sigma-n}, respectively. For every $n,r\geq1$, let $f_{n,r}$
be a weak solution of \eqref{SPDE-iteration} with initial datum $f_0$.
Then, for every cutoff $\chi\in C_c^\infty(\mathbb R^{2d})$ and every
$\gamma\in(0,1/2)$,
\begin{align}\label{time-holder-regularity-fnr}
 \sup_{n,r\geq1}\mathbb E
 \|\chi f_{n,r}\|_{C^\gamma([0,T];H^{-1}(\mathbb R^{2d}))}^2
 \leq C_{\gamma,T,d,V,\sigma,F_1,F_3,\chi}
          \mathcal K_{\rm reg}^2.
\end{align}
In particular, since
$H^{-1}(\mathbb R^{2d})\hookrightarrow\bB^{-6}_{2;\theta}$, for every
$\alpha\in(0,1/2)$,
\begin{align}\label{0521:05-1}
 \sup_{n,r\geq1}
 \|\chi f_{n,r}\|_{L^2(\Omega;
 W^{\alpha,2}([0,T];\bB^{-6}_{2;\theta}))}
 \leq C_{\alpha,T,d,V,\sigma,F_1,F_3,\chi}
          \mathcal K_{\rm reg}.
\end{align}
\end{lemma}
\begin{proof}
Set $f=f_{n,r}$ and write the localized equation in
$H^{-1}(\mathbb R^{2d})$ as
\begin{align}\label{localized-drift-martingale-decomposition}
 \chi f(t)=\chi f_0+\int_0^tA_{n,r}^\chi(s)\,\dif s
                 +M_{n,r}^\chi(t),
\end{align}
where $A_{n,r}^\chi$ contains all deterministic terms and
\begin{align*}
 M_{n,r}^\chi(t)
 =-\sum_{k\geq1}\sum_{\ell=1}^d
   \int_0^t\chi\partial_{v_\ell}
       \bigl(\sigma(f(s))f_k\bigr)\,\dif B_s^{k,\ell}.
\end{align*}
All derivatives in $A_{n,r}^\chi$ are kept in divergence form. For
example,
\begin{align*}
 -\chi\alpha_r^2v\cdot\nabla_xf
 &=-\nabla_x\cdot(\chi\alpha_r^2vf)
     +\alpha_r^2vf\cdot\nabla_x\chi,\\
 \frac1n\chi\Delta_xf
 &=\frac1n\nabla_x\cdot(\chi\nabla_xf)
     -\frac1n\nabla_x\chi\cdot\nabla_xf,
\end{align*}
and the analogous identities are used for the $v$-divergences. On the
support of $\chi$, $|v|$ is bounded uniformly in $r$, while mass
conservation gives
\begin{align*}
 \|V\star_x\langle\alpha_r^2f(t)\rangle\|_{L_x^\infty}
 \leq\|V\|_{L^\infty}M_*.
\end{align*}
Moreover, $\sigma'$ is bounded and $(\sigma\sigma')(0)=0$ with
$(\sigma\sigma')'$ bounded. Consequently, \eqref{0425:00} and
\begin{align*}
 \frac1{n^2}\mathbb E\int_0^T\|\nabla_xf\|_2^2\,\dif t
 \leq\frac{C}{n}\mathbb E\|f_0\|_2^2
\end{align*}
imply the schematic bound
\begin{align*}
 \|A_{n,r}^\chi(t)\|_{H^{-1}}
 \leq C_\chi\left((1+M_*)\|f(t)\|_2
       +\|\nabla_vf(t)\|_2+\frac1n\|\nabla_xf(t)\|_2\right).
\end{align*}
Here the interaction flux is controlled by the preceding bound on
$V\star_x\langle\alpha_r^2f\rangle$, whereas the two fluxes in the
It\^o correction are controlled by boundedness of $\sigma'$ and the
linear growth of $\sigma\sigma'$. We therefore obtain
\begin{align}\label{localized-drift-Hminus1}
 \sup_{n,r\geq1}\mathbb E\int_0^T
 \|A_{n,r}^\chi(t)\|_{H^{-1}}^2\,\dif t
 \leq C_\chi\mathcal K_{\rm reg}^2.
\end{align}
Thus the deterministic primitive in
\eqref{localized-drift-martingale-decomposition} is $1/2$-H\"older in
time with values in $H^{-1}$, with its squared H\"older norm bounded in
expectation.

Let $\mathfrak U=\ell^2(\mathbb N\times\{1,\ldots,d\})$ and denote the
integrand in $M_{n,r}^\chi$ by $G_{n,r}^\chi$. Multiplication by $\chi$
and one distributional derivative map $L^2$ continuously into $H^{-1}$.
Hence, using $F_1=\sum_k f_k^2$, $\sigma(a)^2\leq Ca$, and mass
conservation, we have pathwise
\begin{align}\label{localized-noise-HS-bound}
 \|G_{n,r}^\chi(t)\|_{L_2(\mathfrak U;H^{-1})}^2
 &\leq C_\chi\sum_{k\geq1}\|f_k\sigma(f(t))\|_2^2\notag\\
 &=C_\chi\int_{\mathbb R^{2d}}F_1\sigma(f(t))^2\,dz
 \leq C_\chi M_*.
\end{align}
The Hilbert-space Burkholder--Davis--Gundy inequality therefore yields,
for every $q\geq2$,
\begin{align*}
 \mathbb E\|M_{n,r}^\chi(t)-M_{n,r}^\chi(s)\|_{H^{-1}}^q
 \leq C_{q,\chi}|t-s|^{q/2}.
\end{align*}
Choose $q>2/(1-2\gamma)$. The Kolmogorov--Chentsov theorem (equivalently,
the Garsia--Rodemich--Rumsey inequality) gives a version satisfying
\begin{align*}
 \sup_{n,r\geq1}\mathbb E
 \|M_{n,r}^\chi\|_{C^\gamma([0,T];H^{-1})}^2
 \leq C_{\gamma,\chi}.
\end{align*}
Together with \eqref{localized-drift-Hminus1} and
$\chi f_0\in L^2(\Omega;H^{-1})$, this proves
\eqref{time-holder-regularity-fnr}. Finally, choose
$\gamma\in(\alpha,1/2)$ and use
$C^\gamma([0,T];H^{-1})\hookrightarrow
W^{\alpha,2}([0,T];\bB^{-6}_{2;\theta})$ to obtain \eqref{0521:05-1}.
\end{proof}

Combining the above estimates, we conclude the following tightness result.
\begin{theorem}\label{tightness}
Under the assumptions of Proposition~\ref{prp-uniforml2}, for fixed
$n\geq1$, the sequence
$\{(f_{n,r},\nabla_v f_{n,r})\}_{r\geq1}$ is tight in
\begin{align}\label{fixed-n-tightness-space}
 \mathbb X_n:=
 \Big(L^2([0,T];L^2_{\rm loc}(\mathbb R^{2d}))
      \cap C([0,T];H^{-2}_{\rm loc}(\mathbb R^{2d}))\Big)
 \times\big(L^2([0,T];L^2(\mathbb R^{2d})),w\big).
\end{align}
Here the local $C([0,T];H^{-2})$ topology is induced by a countable
family of compactly supported cutoffs. In particular, the first
components are also tight locally in
$C([0,T];\bB^{-6}_{2;\theta})$.
\end{theorem}
\begin{proof}
For fixed $n$, Proposition~\ref{prp-uniforml2} gives the local spatial
$H^1$ bound, and Lemma~\ref{time-besov-regularity-fn} gives both the
time regularity needed in the Aubin--Lions lemma and tightness in
$C^\gamma([0,T];H^{-1})$ after localization. This proves compactness in
$L^2([0,T];L^2_{\rm loc})$. On each fixed compact set,
$H^{-1}\Subset H^{-2}$; hence the latter H\"older bound and the
Arzel\`a--Ascoli theorem give compactness in
$C([0,T];H^{-2})$. A diagonal argument over a countable exhaustion of
$\mathbb R^{2d}$ proves the assertion. The final claim follows from the
continuous embedding $H^{-2}\hookrightarrow\bB^{-6}_{2;\theta}$.
\end{proof}

\subsection{Passage to the limit for $r \to \infty$}\label{subsec-5-3} 
In this subsection, we consider the following nonlinear SPDE for any fixed $n\geq1$: 
\begin{align}\label{SPDE-L1}
df_{n}=\Big(&\Delta_vf_{n}+\frac{1}{n}\Delta_xf_{n}-v\cdot\nabla_xf_{n}+\nabla_v\cdot(vf_{n})-\nabla_v\cdot(f_{n}V\star_x\langle f_{n}\rangle)\Big)dt\notag\\
    &-\nabla_v\cdot(\sigma(f_{n})dW_F)
    +{\color{black}\frac12\nabla_v\cdot\left(
       F_1\sigma'(f_n)^2\nabla_vf_n\right)dt.}
\end{align}
{\color{black}A weak solution of \eqref{SPDE-L1} is understood in the same
variational sense specified after \eqref{SPDE-iteration}, with the
coefficients in the displayed equation.  In Theorem~\ref{thm:nonb}, existence
and uniqueness are asserted in the regular mass-preserving class: the
solution is nonnegative and satisfies the energy, second-moment, and
Fisher-information estimates obtained below.  The cutoff, moment, and
entropy computations apply to every variational solution in this class, so
these conditions are a priori properties rather than extra choices made in
the compactness argument.}
To establish the existence of solutions, we will make use of the following version of Gy\"ongy-Krylov diagonal argument.

{\color{black}
\begin{lemma}[Gy\"ongy--Krylov criterion on a Jakubowski space]
\label{lem-diagonal}
Let $E$ be a Hausdorff topological space admitting a countable family of
continuous functions which separates points, and assume that every compact
subset of $E$ is metrizable.  Let $(Z_n)_{n\geq1}$ be a tight sequence of
$E$-valued random variables.  Suppose that for every two subsequences
$(Z_{n_k})_k$ and $(Z_{m_k})_k$ there are further subsequences, not relabelled,
such that the joint laws of $(Z_{n_k},Z_{m_k})$ converge weakly to a probability
measure supported on the diagonal
\begin{align*}
 \Delta_E:=\{(z,z):z\in E\}\subset E\times E.
\end{align*}
Then $(Z_n)_n$ converges in probability to an $E$-valued random variable, in
the following intrinsic sense: for every compact $K\subset E$, every metric
$d_K$ generating the relative topology on $K$, and every $\varepsilon>0$,
\begin{align*}
 \lim_{n,m\to\infty}\mathbb P\bigl(
 Z_n,Z_m\in K,\ d_K(Z_n,Z_m)>\varepsilon\bigr)=0,
\end{align*}
and tightness permits the compact set $K$ to be chosen with arbitrarily large
probability.  If $E$ is Polish, this is the usual convergence in probability.
The conclusion also applies to finite or countable products of such spaces.
\end{lemma}

\begin{proof}
This is the Jakubowski-space form of the Gy\"ongy--Krylov criterion; see
\cite{GK21} and the quasi-Polish formulation in
\cite[Theorem~1.1]{holden2022global}.  For completeness, tightness reduces the
argument, up to an event of arbitrarily small probability, to a compact set
$K$.  By assumption $K$ is metrizable.  If the displayed Cauchy-in-probability
property failed on $K$, two subsequences would have a joint weak limit giving
positive mass to the complement of every sufficiently small neighbourhood of
$\Delta_E$, contradicting the diagonal-support hypothesis.  The usual
Cauchy-in-probability argument on the metrizable compact sets, followed by
tightness, produces the asserted $E$-valued limit.
\end{proof}

In the applications below, the strong local Lebesgue and local path spaces
are Polish.  A weak $L^2$ factor has a countable separating family of
continuous linear functionals, and its norm-bounded balls are weakly compact
and metrizable because $L^2$ is separable.  The energy estimates make these
balls a compact exhaustion in probability.  Consequently, the product spaces
$\widehat{\mathbb X}_n$, $\mathbb X$, and their countable path-coordinate
extensions used below satisfy the hypotheses of Lemma~\ref{lem-diagonal}.}

\bt\label{thm:nonb}
{\color{black}Assume that $f_0$, $V$, and $\sigma(\cdot)$ satisfy Assumptions~\ref{Assump-initialdata-regular}, \ref{Assump-ker}, and~\ref{Assump-sigma-n}, respectively.} Let $n\geq1$, then there exists a unique weak solution $f$ to the nonlinear equation \eqref{SPDE-L1}. Moreover, the solution $f$ satisfies 
\begin{align}\label{0425:04}
    {\color{black}\|f\|_{L^\infty([0,T];L^2(\Omega;L^2(\mathbb{R}^{2d})))} 
    + \|\nabla_v f\|_{L^2([0,T] \times \Omega; L^2(\mathbb{R}^{2d}))}+\frac{1}{n}\|\nabla_x f\|_{L^2([0,T] \times \Omega; L^2(\mathbb{R}^{2d}))} 
    {\color{black}\leq C_{T,d,V,\sigma,F_1,F_3}\,\mathcal K_{\rm reg}.}} 
\end{align}
\et
\begin{proof}
Fix $n\geq1$ and suppress it from the notation.  Let $f_r=f_{n,r}$ be
the solution of \eqref{SPDE-iteration}, and set
\begin{align*}
 b_r(g):=V\star_x\langle\alpha_r^2g\rangle .
\end{align*}
We retain the structure of the compactness argument above, but incorporate
the time-continuous topology supplied by
Lemma~\ref{time-besov-regularity-fn}.  Choose a countable dense family
$(\varphi_\ell)_{\ell\geq1}\subset C_c^\infty(\mathbb R^{2d})$.  For a
sufficiently regular $g$, define
\begin{align*}
 \mathcal D_r^\ell(g):={}&
 \int_{\mathbb R^{2d}}g\left(
 \Delta_v\varphi_\ell+\frac1n\Delta_x\varphi_\ell
 +\alpha_r^2v\cdot\nabla_x\varphi_\ell
 -\alpha_r^2v\cdot\nabla_v\varphi_\ell
 +b_r(g)\cdot\nabla_v\varphi_\ell\right)\,dz\\
 &-\frac12\int_{\mathbb R^{2d}}
 {\color{black}F_1\sigma'(g)^2\nabla_vg}
       \cdot\nabla_v\varphi_\ell\,dz
\end{align*}
and
\begin{align}\label{fixed-n-martingale-coefficient}
 G_{r,k}^\ell(t)
 :=\int_{\mathbb R^{2d}}
       \sigma(f_r(t,z))f_k(z)\nabla_v\varphi_\ell(z)\,dz
 \in\mathbb R^d.
\end{align}
The weak formulation of \eqref{SPDE-iteration} reads
\begin{align}\label{fixed-n-residual-martingale}
 M_r^\ell(t)
 &:=
 \langle f_r(t),\varphi_\ell\rangle
 -\langle f_0,\varphi_\ell\rangle
 -\int_0^t\mathcal D_r^\ell(f_r(s))\,\dif s\notag\\
 &=\sum_{k\geq1}\int_0^tG_{r,k}^\ell(s)\cdot\dif B_s^k.
\end{align}
Spatial Cauchy--Schwarz, $F_1\in L^\infty$, the bound
$\sigma(a)^2\leq Ca$, and mass conservation give
\begin{align}\label{fixed-n-bracket-uniform}
 \sum_{k\geq1}|G_{r,k}^\ell(t)|^2
 \leq C_\ell\int_{\operatorname{supp}\varphi_\ell}
       F_1\sigma(f_r(t))^2\,dz
 \leq C_\ell M_*
\end{align}
almost surely.  The Burkholder--Davis--Gundy inequality and
Kolmogorov's criterion therefore imply tightness of
$(M_r^\ell)_r$ in $C([0,T])$, simultaneously for every $\ell$.

We now follow the two-subsequence argument in the original proof.  Let
$(r_j)$ and $(r_j')$ be arbitrary subsequences and put
\begin{align*}
 Y_r:=\bigl(f_r,\nabla_vf_r,\nabla_xf_r,
             (M_r^\ell)_{\ell\geq1}\bigr).
\end{align*}
Since $n$ is fixed, \eqref{0425:00} also makes the third coordinate tight
in the weak topology of $L^2(0,T;L^2)$. Set
\begin{align*}
 \widehat{\mathbb X}_n
 :=\mathbb X_n\times\big(L^2(0,T;L^2(\mathbb R^{2d})),w\big).
\end{align*}
By Theorem~\ref{tightness}, \eqref{fixed-n-bracket-uniform}, and the
tightness of the common initial datum and Brownian motions, the laws of
\begin{align*}
 \bigl(Y_{r_j},Y_{r_j'},(B^k)_{k\geq1},f_0\bigr)
\end{align*}
are tight on two copies of
\begin{align*}
 \mathbb Y_n:=\widehat{\mathbb X}_n\times C([0,T])^{\mathbb N},
\end{align*}
together with
$C([0,T];\mathbb R^d)^{\mathbb N}\times L^1(\mathbb R^{2d})$.
Each weak $L^2$ factor has a countable family of continuous linear
functionals separating its points, and the uniform energy bounds give
tightness on weakly compact balls. Thus the full product satisfies the
hypotheses of the Jakubowski--Skorokhod representation theorem. On a new probability
space, along a further subsequence, we consequently obtain copies
\begin{align*}
 \bigl(\bar Y_j,\bar Y_j',(\bar B_j^k)_{k\geq1},\bar f_{0,j}\bigr)
 \longrightarrow
 \bigl(\bar Y,\bar Y',(\bar B^k)_{k\geq1},\bar f_0\bigr)
\end{align*}
almost surely in the preceding product topology, with the same laws as
the original variables.  Suppressing the prime for the moment, this
means
\begin{align}\label{fixed-n-Skorokhod-convergences}
 \bar f_j&\longrightarrow\bar f
 &&\text{in }L^2(0,T;L^2_{\rm loc})
   \cap C([0,T];H^{-2}_{\rm loc}),\notag\\
 \nabla_v\bar f_j&\rightharpoonup\bar g
 &&\text{in }L^2(0,T;L^2),\notag\\
 \nabla_x\bar f_j&\rightharpoonup\bar h
 &&\text{in }L^2(0,T;L^2),\notag\\
 \bar M_j^\ell&\longrightarrow\bar M^\ell
 &&\text{in }C([0,T])\quad\text{for every }\ell,
\end{align}
and the Brownian coordinates converge uniformly. Testing the second and
third lines against compactly supported smooth vector fields and integrating
by parts shows that $\bar g=\nabla_v\bar f$ and
$\bar h=\nabla_x\bar f$. The same conclusions hold
for the primed sequence.  Passing to a further subsequence if necessary,
the local strong convergence also holds almost everywhere with respect
to $\bar{\mathbb P}\otimes dt\otimes dz$.

The additional time topology is used at precisely the point at which an
$L^2$-in-time limit is insufficient.  If
$\varphi\in C_c^\infty(\mathbb R^{2d})$, then
\begin{align}\label{fixed-n-all-time-endpoint}
 \sup_{t\in[0,T]}
 |\langle\bar f_j(t)-\bar f(t),\varphi\rangle|
 \longrightarrow0
\end{align}
almost surely.  In particular, since
$\bar f_j(0)=\bar f_{0,j}$, the convergence of the initial coordinate in
$L^1$ gives $\bar f(0)=\bar f_0$ in distributions and hence in $L^1$.
Thus the continuous path obtained in the Skorokhod representation is
the continuous distribution-valued modification of the
$L^2(0,T;L^2_{\rm loc})$ limit.  Since every $\bar f_j(t)$ is
nonnegative, \eqref{fixed-n-all-time-endpoint} also shows that
$\bar f(t)$ is a nonnegative distribution, and hence a nonnegative
Radon measure, for every $t\in[0,T]$.

We next pass the moment estimate to the limit.  Let
$(\omega_R)_{R\geq1}\subset C_c^\infty(\mathbb R^{2d})$ be nonnegative
and increasing, with
$\omega_R(z)\uparrow|x|^2+|v|^2$.  By
\eqref{fixed-n-all-time-endpoint}, equality in law, Fatou's lemma, and
Proposition~\ref{moment-es},
\begin{align*}
 \bar{\mathbb E}\sup_{t\leq T}\langle\bar f(t),\omega_R\rangle
 &\leq\liminf_{j\to\infty}
 \bar{\mathbb E}\sup_{t\leq T}\langle\bar f_j(t),\omega_R\rangle\\
 &\leq C_{\rm mom},
\end{align*}
where
\begin{align*}
 C_{\rm mom}:=
 C(T,\sigma,V,F_1,F_3,M_*)
 \left(\mathbb E\int_{\mathbb R^{2d}}
 (|x|^2+|v|^2)f_0\,dz+M_*\right).
\end{align*}
Letting $R\to\infty$ and using monotone convergence therefore yields
\begin{align}\label{fixed-n-limit-moment}
 \bar{\mathbb E}\sup_{t\leq T}\int_{\mathbb R^{2d}}
 (|x|^2+|v|^2)\bar f(t,z)\,dz\leq C_{\rm mom}.
\end{align}

As in the original proof, this moment estimate upgrades local convergence
to the global convergence needed for the interaction term.  Indeed,
for every $R>1$,
\begin{align*}
 \int_0^T\int_{\{|x|+|v|>R\}}(\bar f_j+\bar f)\,dz\,dt
 \leq\frac{C}{R^2}\int_0^T\int_{\mathbb R^{2d}}
 (|x|^2+|v|^2)(\bar f_j+\bar f)\,dz\,dt .
\end{align*}
The expectations of the right-hand side are bounded uniformly in $j$.
Local $L^2$ convergence and Markov's inequality first give convergence
in probability in $L^1(0,T;L^1)$; after extracting a further subsequence
we have
\begin{align}\label{fixed-n-global-L1-convergence}
 \bar f_j\longrightarrow\bar f
 \quad\text{in }L^1(0,T;L^1(\mathbb R^{2d}))
 \quad\text{almost surely}.
\end{align}
The mass bound in fact holds simultaneously for every time. Let
$0\leq\eta_R\in C_c^\infty$ increase to $1$. For each $t$, the
time-continuous convergence gives
\begin{align*}
 \langle\bar f(t),\eta_R\rangle
 =\lim_{j\to\infty}\langle\bar f_j(t),\eta_R\rangle\leq M_*.
\end{align*}
Since $\bar f(t)$ is a nonnegative Radon measure, monotone convergence
therefore yields $\|\bar f(t)\|_1\leq M_*$ for every $t\in[0,T]$, on the
same event.

We now pass term by term in the weak equation, as in the original
argument, but keep all time primitives explicit.  The diffusion,
transport, and damping terms converge uniformly in their upper time
limit by \eqref{fixed-n-Skorokhod-convergences}, local strong
convergence, and the locally uniform convergence $\alpha_{r_j}^2\to1$.
The $x$-viscosity is written as
\begin{align*}
 \frac1n\int_0^t\int\bar f_j\Delta_x\varphi,
\end{align*}
so no separate convergence of $\nabla_x\bar f_j$ is needed.  For a weakly
convergent scalar sequence $h_j\rightharpoonup h$ in $L^2(0,T)$, the
primitives $\int_0^\cdot h_j$ converge uniformly: they converge first
pointwise by weak convergence and are uniformly $1/2$-H\"older.
This observation applies to every remaining term containing
$\nabla_v\bar f_j$.

For completeness, consider explicitly the nonlocal term.  With
$\bar b_j=b_{r_j}(\bar f_j)$ and
$\bar b=V\star_x\langle\bar f\rangle$, we write
\begin{align*}
 \bar f_j\bar b_j-\bar f\bar b
 =(\bar f_j-\bar f)\bar b_j+\bar f(\bar b_j-\bar b).
\end{align*}
Mass conservation for the approximations gives
$\|\bar b_j\|_\infty\leq\|V\|_\infty M_*$.  Thus the first term
converges locally in $L^1(0,T;L^1)$ by the local strong convergence.
Moreover,
\begin{align*}
 \|\bar b_j-\bar b\|_{L_x^\infty}
 \leq\|V\|_\infty
 \|\alpha_{r_j}^2\bar f_j-\bar f\|_{L^1_{x,v}},
\end{align*}
and hence the second term converges by
\eqref{fixed-n-global-L1-convergence}.  Both time primitives consequently
converge uniformly on $[0,T]$.

It remains to consider the It\^o correction.  We preserve its divergence
form:
\begin{align*}
 {\color{black}F_1\sigma'(\bar f_j)^2\nabla_v\bar f_j}.
\end{align*}
Since $\sigma'^2$ is bounded and Lipschitz in the regular setting, the coefficient converges strongly
locally, while the gradients converge weakly. More explicitly, with
$a_j=\sigma'(\bar f_j)^2$ and $a=\sigma'(\bar f)^2$, we use
\begin{align*}
 a_j\nabla_v\bar f_j-a\nabla_v\bar f
 =(a_j-a)\nabla_v\bar f_j
   +a(\nabla_v\bar f_j-\nabla_v\bar f).
\end{align*}
The time primitive of the first term converges uniformly by strong
$L^2$ convergence of $a_j-a$ and the uniform $L^2$ gradient bound. After
pairing the second term with $F_1\nabla_v\varphi$, its scalar integrand
converges weakly in $L^2(0,T)$; its primitives converge pointwise and
are uniformly $1/2$-H\"older, hence converge uniformly. We have therefore
proved
\begin{align}\label{fixed-n-deterministic-residual-convergence}
 \sup_{t\leq T}\left|
 \int_0^t\mathcal D_{r_j}^\ell(\bar f_j(s))\,\dif s
 -\int_0^t\mathcal D^\ell(\bar f(s))\,\dif s
 \right|\longrightarrow0
\end{align}
for every $\ell$, where $\mathcal D^\ell$ is obtained from
$\mathcal D_r^\ell$ by replacing $\alpha_r^2$ by $1$ and $b_r$ by
$V\star_x\langle\bar f\rangle$.  Combining
\eqref{fixed-n-all-time-endpoint} and
\eqref{fixed-n-deterministic-residual-convergence} identifies the
uniform limit $\bar M^\ell$ as the residual of the limiting equation for
every $t\in[0,T]$.

We finally identify the limiting noise, retaining the martingale
characterization in the original proof.  Define
\begin{align*}
 G_k^\ell(t):=\int_{\mathbb R^{2d}}
 \sigma(\bar f(t,z))f_k(z)\nabla_v\varphi_\ell(z)\,dz.
\end{align*}
We write $\bar G_{j,k}^\ell$ for \eqref{fixed-n-martingale-coefficient}
with $\bar f_j$ in place of $f_r$.
The Lipschitz continuity of $\sigma$ and $F_1\in L^\infty$ imply
\begin{align}\label{fixed-n-noise-coefficient-convergence}
 \int_0^T\sum_{k\geq1}
 |\bar G_{j,k}^\ell-G_k^\ell|^2\,\dif t
 \leq C_\ell\int_0^T
 \|\bar f_j-\bar f\|_{L^2(\operatorname{supp}\varphi_\ell)}^2\,\dif t
 \longrightarrow0
\end{align}
almost surely.  The deterministic bracket bound
\eqref{fixed-n-bracket-uniform} gives uniform integrability of all
martingales and products below.

Let $\bar{\mathcal G}_t$ be the augmented canonical filtration generated
up to time $t$ by both limiting solution coordinates, their residual
martingales, the common initial datum, and $(\bar B^k)_k$.  Testing the
Brownian increment, quadratic-variation, and cross-variation identities
against bounded continuous functions of the canonical past and passing
to the limit shows, by a monotone-class argument, that the
$(\bar B^k)_k$ are independent $\mathbb R^d$-valued Brownian motions
relative to $(\bar{\mathcal G}_t)$.  The same procedure, using
\eqref{fixed-n-noise-coefficient-convergence}, gives that
\begin{align*}
 \bar M^\ell,\qquad
 \bar M^\ell\bar B^{k,a}
   -\int_0^\cdot G_{k,a}^\ell(s)\,\dif s,\qquad
 (\bar M^\ell)^2
   -\int_0^\cdot\sum_{k\geq1}|G_k^\ell(s)|^2\,\dif s
\end{align*}
are $(\bar{\mathcal G}_t)$-martingales for every $k$, $a$, and $\ell$.
Consequently,
\begin{align}\label{fixed-n-limit-stochastic-integral}
 \bar M^\ell(t)
 =\sum_{k\geq1}\int_0^tG_k^\ell(s)\cdot\dif\bar B_s^k,
 \qquad t\in[0,T],
\end{align}
because the difference of the two sides is a continuous martingale with
zero quadratic variation.  We have thus obtained the weak formulation
of \eqref{SPDE-L1} first for the countable dense family and then for all
$\varphi\in C_c^\infty(\mathbb R^{2d})$, on one event and for every
$t\in[0,T]$.

It remains to verify the a priori properties needed for uniqueness.
The preceding argument already gives nonnegativity for every time.
Applying the cutoff argument of
Lemma~\ref{lem-L1-es} to the limiting weak equation yields mass
conservation on one event for all $t\in[0,T]$.  Lower semicontinuity in
\eqref{0425:00} gives the energy estimate.  Finally, the convex
lower-semicontinuity of
\begin{align*}
 (a,\xi)\longmapsto
 \begin{cases}|\xi|^2/a,&a>0,\\0,&(a,\xi)=(0,0),\\+\infty,&\text{otherwise},
 \end{cases}
\end{align*}
combined with Proposition~\ref{approximate-entropy-dissipation}, gives
the required Fisher-information bound for $\bar f$ (and the analogous
$n^{-1}$-weighted $x$-dissipation). Thus the proof of
Proposition~\ref{equivalence} applies to the limiting equation. The
additional viscosity contributes $n^{-1}\Delta_x\chi$ to the kinetic
equation and the nonnegative kinetic defect
$n^{-1}\delta_{\bar f=\zeta}|\nabla_x\bar f|^2$.

The primed limit satisfies the same equation, has the same initial datum,
and is driven by the same Brownian motions.  In the doubling argument of
Theorem~\ref{Uniqueness-spde}, the additional $x$-viscosity has the
favorable sign, so the two limits coincide.  Their
$C([0,T];H^{-2}_{\rm loc})$ modifications then coincide for every time.
Therefore every joint subsequential limit is supported on the diagonal
of $\widehat{\mathbb X}_n$. The Gy\"ongy--Krylov criterion gives convergence in
probability on the original stochastic basis and produces the unique
adapted weak solution of \eqref{SPDE-L1}. For each fixed time,
convergence in $H^{-2}_{\rm loc}$ together with the uniform global
$L^2$ bound implies weak $L^2$ convergence along a subsequence; hence
the time-slice norm is lower semicontinuous. The gradient bounds pass by
the weak convergences recorded above. Taking the supremum over time and
using \eqref{0425:00} proves \eqref{0425:04}.
\end{proof}

\subsection{Passage to the limit for $n \to \infty$}\label{subsec-5-4}
For every $n \geq 1$, let $f_n$ be the solution to \eqref{SPDE-L1}. In this subsection, we pass to the limit as $n \to \infty$ to prove Theorem \ref{thm-existence}. 

{\color{black}With the help of Proposition \ref{prp-uniforml2}, we apply the
Hilbert case (H) of Theorem~\ref{thm:Besov}, so throughout this application
$p=q=2$.  Put $\rho_n:=\langle f_n\rangle$ and
\begin{align*}
 b_n&:=V\star_x\rho_n,\\
 g_1&:={\color{black}\frac12\nabla_v\cdot\left(
 F_1\sigma'(f_n)^2\nabla_vf_n\right)},\\
 g_2&:=0,\qquad
 h_1:=-\sigma'(f_n)\nabla_vf_n,\qquad
 h_2:=-\sigma(f_n).
\end{align*}
This is the full divergence form of the It\^o correction in
\eqref{SPDE-0-ito}; in particular, no derivative of the merely bounded
field $F_1$ is taken separately.  Mass conservation and
\eqref{regular-data-size} give
\begin{align*}
 \|b_n\|_{L^\infty(\Omega\times[0,T]\times\mathbb R^{2d})}
 \leq\|V\|_{L^\infty}M_*.
\end{align*}
Since divergence maps $L^2$ continuously into
$\mathbf B^{-1}_{2;\theta}$, $\sigma(0)=0$, and $\sigma'$ is bounded,
\begin{align*}
 \|g_1\|_{L^2(\Omega\times[0,T];\mathbf B^{-1}_{2;\theta})}
 &\leq C_{F_1,\sigma}
 \left(\|\nabla_vf_n\|_{L^2(\Omega\times[0,T]\times\mathbb R^{2d})}
       +\|f_n\|_{L^2(\Omega\times[0,T]\times\mathbb R^{2d})}\right),\\
 \|h_1\|_{L^2(\Omega\times[0,T];L^2)}
 +\|h_2\|_{L^2(\Omega\times[0,T];L^2)}
 &\leq C_\sigma
 \left(\|\nabla_vf_n\|_{L^2(\Omega\times[0,T]\times\mathbb R^{2d})}
       +\|f_n\|_{L^2(\Omega\times[0,T]\times\mathbb R^{2d})}\right).
\end{align*}
Therefore, by \eqref{0425:00}, the data norm in (H) satisfies
\begin{align}\label{0425:06}
 \sup_{n\geq1}\mathcal K_{H,n}
 \leq C_{T,d,V,\sigma,F_1,F_3,M_*}\mathcal K_{\rm reg},
\end{align}
where $\mathcal K_{H,n}$ denotes $\mathcal K_H$ with the above
coefficients and initial datum $f_0$.}

As a consequence of Theorem \ref{thm:Besov} and Lemma \ref{time-besov-regularity-fn}, we have the following result. 

\begin{theorem}\label{besov-regularity-fn}
{\color{black}Assume that $f_0$, $V$, and $\sigma(\cdot)$ satisfy
Assumptions~\ref{Assump-initialdata-regular}, \ref{Assump-ker}, and
\ref{Assump-sigma-n}, respectively.} For every $n \in \mathbb{N}$, let $f_{n}$ be the solution to \eqref{SPDE-L1} with initial data $f_0$. Then for any $\beta \in (0, 1)$, we have the uniform anisotropic Besov estimate
\begin{align}\label{0425:020}
    {\color{black}\sup_{n \geq 1} \|f_{n}\|_{L^2([0, T] \times \Omega; \bB^\beta_{2; \theta})}
    \leq C_{\beta,T,d,V,\sigma,F_1,F_3,M_*}\,\mathcal K_{\rm reg}.}
\end{align}
Furthermore, for every cutoff $\chi\in C_c^\infty(\mathbb R^{2d})$ and
every $\gamma\in(0,1/2)$,
\begin{align}\label{time-holder-regularity-fn-limit}
 \sup_{n\geq1}\mathbb E
 \|\chi f_n\|_{C^\gamma([0,T];H^{-1}(\mathbb R^{2d}))}^2
 \leq C_{\gamma,T,d,V,\sigma,F_1,F_3,M_*,\chi}
       \mathcal K_{\rm reg}^2.
\end{align}
Consequently, for every $\alpha\in(0,1/2)$,
\begin{align}\label{time-besov-regularity-fn-limit}
    {\color{black}\sup_{n \geq 1} \|f_{n} \chi\|_{L^2(\Omega; W^{\alpha, 2}([0, T]; \bB^{-6}_{2; \theta}))}
    \leq C_{\alpha,T,d,V,\sigma,F_1,F_3,M_*,\chi}\,\mathcal K_{\rm reg}.}
\end{align}
\end{theorem}
\begin{proof}
The Besov estimate is the Hilbert case of Theorem~\ref{thm:Besov}, with
the coefficients estimated in \eqref{0425:06}. We prove the new time
estimate directly. Set $f=f_n$ and localize \eqref{SPDE-L1} by $\chi$.
Exactly as in the proof of Lemma~\ref{time-besov-regularity-fn}, in
$H^{-1}(\mathbb R^{2d})$ we have
\begin{align}\label{localized-fn-drift-martingale-decomposition}
 \chi f(t)=\chi f_0+\int_0^tA_n^\chi(s)\,\dif s+M_n^\chi(t),
\end{align}
where all deterministic derivatives in $A_n^\chi$ are retained in
divergence form and
\begin{align*}
 M_n^\chi(t)
 =-\sum_{k\geq1}\sum_{\ell=1}^d
   \int_0^t\chi\partial_{v_\ell}
       \bigl(\sigma(f(s))f_k\bigr)\,\dif B_s^{k,\ell}.
\end{align*}
Mass conservation, the lower-semicontinuous form of \eqref{0425:00},
boundedness of $\sigma'$, and the
linear growth of $\sigma\sigma'$ give
\begin{align*}
 \|A_n^\chi(t)\|_{H^{-1}}
 \leq C_\chi\left((1+M_*)\|f(t)\|_2
       +\|\nabla_vf(t)\|_2+\frac1n\|\nabla_xf(t)\|_2\right).
\end{align*}
Since
\begin{align*}
 \frac1{n^2}\mathbb E\int_0^T\|\nabla_xf_n\|_2^2\,\dif t
 \leq\frac{C}{n}\mathbb E\|f_0\|_2^2,
\end{align*}
we obtain
\begin{align}\label{localized-fn-drift-Hminus1}
 \sup_{n\geq1}\mathbb E\int_0^T
 \|A_n^\chi(t)\|_{H^{-1}}^2\,\dif t
 \leq C_\chi\mathcal K_{\rm reg}^2.
\end{align}

Let $G_n^\chi$ be the stochastic integrand in $M_n^\chi$. The product
rule for $\chi\partial_v$, together with $F_1=\sum_kf_k^2$,
$\sigma(a)^2\leq Ca$, and mass conservation, yields pathwise
\begin{align*}
 \|G_n^\chi(t)\|_{L_2(\ell^2(\mathbb N\times\{1,\ldots,d\});H^{-1})}^2
 \leq C_\chi\int F_1\sigma(f_n(t))^2\,dz
 \leq C_\chi M_*.
\end{align*}
Thus, for every $q\geq2$, the Hilbert-space
Burkholder--Davis--Gundy inequality gives
\begin{align*}
 \mathbb E\|M_n^\chi(t)-M_n^\chi(s)\|_{H^{-1}}^q
 \leq C_{q,\chi}|t-s|^{q/2}.
\end{align*}
Given $\gamma<1/2$, choose $q>2/(1-2\gamma)$ and apply the
Kolmogorov--Chentsov theorem, or equivalently the
Garsia--Rodemich--Rumsey inequality. Combining the resulting martingale
estimate with \eqref{localized-fn-drift-Hminus1} proves
\eqref{time-holder-regularity-fn-limit}. Finally,
$H^{-1}\hookrightarrow\bB^{-6}_{2;\theta}$ and
$C^\gamma_t\hookrightarrow W^{\alpha,2}_t$ for
$\alpha<\gamma<1/2$, which proves
\eqref{time-besov-regularity-fn-limit}.
\end{proof}
Before establishing the tightness of the approximation sequence, we introduce a technical lemma regarding compact embeddings.

\begin{lemma}\label{App:cpt}
{\color{black}Let $\alpha \in (0,1)$ and $1 \leq p < \infty$.} Then the embedding
$$
\bB^\alpha_{p;\theta} \hookrightarrow L^p
$$
is compact on any bounded domain. More precisely, for any uniformly bounded sequence $(f_n)_{n\geq1} \subset \bB^\alpha_{p;\theta}$, the sequence $(f_n)_{n\geq1}$ is relatively compact in $L^p(D)$, for any bounded domain $D \subset \mathbb{R}^{2d}$. 
\end{lemma}

\begin{proof}
By \cite[Theorem 2.7]{HZZZ24}, there exists a constant $C = C(d, \alpha, p) > 0$ such that for all $h \in \mathbb{R}^{2d}$ and $n \in \mathbb{N}$,
$$
\| f_n(\cdot - h) - f_n(\cdot) \|_{L^p(D)} \leq {\color{black}\| f_n(\cdot - h) - f_n(\cdot) \|_{L^p(\mathbb{R}^{2d})}} \leq C |h|_\theta^\alpha \| f_n \|_{\bB^\alpha_{p;\theta}}.
$$
Since the sequence $\{f_n\}$ is uniformly bounded in $\bB^\alpha_{p;\theta}$, the above estimate shows uniform equicontinuity in the $L^p$-norm. By the Kolmogorov-Riesz compactness theorem, the embedding is compact.
\end{proof}

As a consequence, the embedding $\bB^\beta_{2;\theta} \hookrightarrow L^2$ is compact on any bounded domain, which allows us to obtain the following corollary.

\begin{corollary}\label{co-tightness}
Under the assumptions of Theorem~\ref{besov-regularity-fn}, the sequence
$(f_n,\nabla_vf_n)_{n\geq1}$ is tight in
\begin{align}\label{n-tightness-space}
 \mathbb X:=
 \Big(L^2(0,T;L^2_{\rm loc}(\mathbb R^{2d}))
      \cap C([0,T];H^{-2}_{\rm loc}(\mathbb R^{2d}))\Big)
 \times\big(L^2(0,T;L^2(\mathbb R^{2d})),w\big).
\end{align}
The local continuous topology is induced by a countable family of
compactly supported cutoffs. In particular, the first components are
also tight locally in $C([0,T];\bB^{-6}_{2;\theta})$.

\end{corollary}

\begin{proof}
The Besov estimate, Lemma~\ref{App:cpt}, and the time estimate
\eqref{time-holder-regularity-fn-limit} give compactness in
$L^2(0,T;L^2_{\rm loc})$ by the Aubin--Lions lemma. For each localized
family with support in a fixed compact set,
$H^{-1}\Subset H^{-2}$; hence the H\"older estimate and the
Arzel\`a--Ascoli theorem also give compactness in
$C([0,T];H^{-2})$. A diagonal argument over a countable exhaustion gives
the first factor in \eqref{n-tightness-space}. The energy estimate gives
weak-$L^2$ tightness of $\nabla_vf_n$, and the continuous embedding
$H^{-2}\hookrightarrow\bB^{-6}_{2;\theta}$ proves the final assertion.
\end{proof}

With the help of the above tightness result, we are ready to prove Theorem \ref{thm-existence}

\begin{proof}[Proof of Theorem~\ref{thm-existence}]
For each $n\geq1$, let $f_n$ be the unique weak solution of
\eqref{SPDE-L1}, and put
\begin{align*}
 \rho_n:=\langle f_n\rangle,
 \qquad b(g):=V\star_x\langle g\rangle.
\end{align*}
First, lower semicontinuity in the passage $r\to\infty$, together with
Proposition~\ref{approximate-entropy-dissipation}, gives
\begin{align}\label{uniform-n-Fisher}
 \sup_{n\geq1}\mathbb E\int_0^T\int_{\mathbb R^{2d}}
       \frac{|\nabla_vf_n|^2}{f_n}\,dz\,dt<\infty.
\end{align}

We retain the structure of the preceding compactness proof, now using
the time-continuous topology in Corollary~\ref{co-tightness}. Choose a
countable dense family
$(\varphi_\ell)_{\ell\geq1}\subset C_c^\infty(\mathbb R^{2d})$ and set
\begin{align*}
 \mathcal D_n^\ell(g):={}&
 \int_{\mathbb R^{2d}}g\left(
 \Delta_v\varphi_\ell+\frac1n\Delta_x\varphi_\ell
 +v\cdot\nabla_x\varphi_\ell-v\cdot\nabla_v\varphi_\ell
 +b(g)\cdot\nabla_v\varphi_\ell\right)\,dz\\
 &-\frac12\int_{\mathbb R^{2d}}
 {\color{black}F_1\sigma'(g)^2\nabla_vg}
       \cdot\nabla_v\varphi_\ell\,dz,
\end{align*}
and
\begin{align}\label{n-main-noise-coefficient}
 G_{n,k}^\ell(t):=
 \int_{\mathbb R^{2d}}
       \sigma(f_n(t,z))f_k(z)\nabla_v\varphi_\ell(z)\,dz
 \in\mathbb R^d.
\end{align}
The weak formulation is equivalently
\begin{align}\label{n-main-residual-martingale}
 M_n^\ell(t)
 &:={}
 \langle f_n(t),\varphi_\ell\rangle
 -\langle f_0,\varphi_\ell\rangle
 -\int_0^t\mathcal D_n^\ell(f_n(s))\,\dif s\notag\\
 &=\sum_{k\geq1}\int_0^tG_{n,k}^\ell(s)\cdot\dif B_s^k.
\end{align}
Spatial Cauchy--Schwarz, $F_1\in L^\infty$,
$\sigma(a)^2\leq Ca$, and mass conservation give the deterministic
bracket bound
\begin{align}\label{n-main-bracket-bound}
 \sum_{k\geq1}|G_{n,k}^\ell(t)|^2
 \leq C_\ell\int_{\operatorname{supp}\varphi_\ell}
       F_1\sigma(f_n(t))^2\,dz
 \leq C_\ell M_*.
\end{align}
Consequently $(M_n^\ell)_n$ is tight in $C([0,T])$ for every $\ell$.

Let $(n_j)$ and $(n_j')$ be arbitrary subsequences and define
\begin{align*}
 Y_n:=\bigl(f_n,\nabla_vf_n,(M_n^\ell)_{\ell\geq1}\bigr).
\end{align*}
By Corollary~\ref{co-tightness}, \eqref{n-main-bracket-bound}, and
tightness of the common Brownian motions and initial datum, the laws of
\begin{align*}
 \bigl(Y_{n_j},Y_{n_j'},(B^k)_{k\geq1},f_0\bigr)
\end{align*}
are tight on two copies of
\begin{align*}
 \mathbb Y:=\mathbb X\times C([0,T])^{\mathbb N},
\end{align*}
together with
$C([0,T];\mathbb R^d)^{\mathbb N}\times L^1(\mathbb R^{2d})$.
Each weak $L^2$ factor has a countable separating family of continuous
linear functionals, and the energy estimate gives tightness on weakly
compact balls. The Jakubowski--Skorokhod theorem therefore supplies a
new probability space on which, along a further subsequence,
\begin{align*}
 \bigl(\bar Y_j,\bar Y_j',(\bar B_j^k)_{k\geq1},\bar f_{0,j}\bigr)
 \longrightarrow
 \bigl(\bar Y,\bar Y',(\bar B^k)_{k\geq1},\bar f_0\bigr)
\end{align*}
almost surely, with the same laws as the original variables. Suppressing
the prime while identifying either limit, we have
\begin{align}\label{n-main-Skorokhod-convergences}
 \bar f_j&\longrightarrow\bar f
 &&\text{in }L^2(0,T;L^2_{\rm loc})
   \cap C([0,T];H^{-2}_{\rm loc}),\notag\\
 \nabla_v\bar f_j&\rightharpoonup\bar g
 &&\text{in }L^2(0,T;L^2),\notag\\
 \bar M_j^\ell&\longrightarrow\bar M^\ell
 &&\text{in }C([0,T])\quad\text{for every }\ell,
\end{align}
and the Brownian coordinates converge uniformly. Distributional
integration by parts identifies $\bar g=\nabla_v\bar f$. The same holds
for the primed variables. After a further subsequence, the local strong
convergence also holds almost everywhere with respect to
$\bar{\mathbb P}\otimes dt\otimes dz$.

The new time topology gives, for every
$\varphi\in C_c^\infty(\mathbb R^{2d})$,
\begin{align}\label{n-main-all-time-endpoint}
 \sup_{t\in[0,T]}
 |\langle\bar f_j(t)-\bar f(t),\varphi\rangle|
 \longrightarrow0
\end{align}
almost surely. Hence $\bar f(0)=\bar f_0$ in distributions. Moreover,
each time slice of $\bar f$ is a nonnegative distribution and therefore
a nonnegative Radon measure.

We next pass the moment estimate to the limit. Let
$0\leq\omega_R\in C_c^\infty(\mathbb R^{2d})$ increase pointwise to
$|x|^2+|v|^2$, and let $C_{\rm mom}$ denote the uniform constant in
Proposition~\ref{moment-es}. By \eqref{n-main-all-time-endpoint}, equality in law,
Fatou's lemma, and Proposition~\ref{moment-es},
\begin{align*}
 \bar{\mathbb E}\sup_{t\leq T}
 \langle\bar f(t),\omega_R\rangle
 \leq\liminf_{j\to\infty}
 \bar{\mathbb E}\sup_{t\leq T}
 \langle\bar f_j(t),\omega_R\rangle
 \leq C_{\rm mom}.
\end{align*}
Letting $R\to\infty$ gives
\begin{align}\label{n-main-limit-moment}
 \bar{\mathbb E}\sup_{t\leq T}
 \int_{\mathbb R^{2d}}(|x|^2+|v|^2)\bar f(t,z)\,dz
 \leq C_{\rm mom}.
\end{align}
Thus, for $R>1$,
\begin{align*}
 \int_0^T\int_{\{|x|+|v|>R\}}(\bar f_j+\bar f)\,dz\,dt
 \leq\frac{C}{R^2}\int_0^T\int_{
 \mathbb R^{2d}}(|x|^2+|v|^2)(\bar f_j+\bar f)\,dz\,dt.
\end{align*}
Local $L^2$ convergence and Markov's inequality first give global
$L^1(0,T;L^1)$ convergence in probability. Passing to a further
subsequence, we may assume
\begin{align}\label{n-main-global-L1}
 \bar f_j\longrightarrow\bar f
 \quad\text{in }L^1(0,T;L^1(\mathbb R^{2d}))
 \quad\text{almost surely}.
\end{align}
The mass bound also holds simultaneously for every time. Indeed, if
$0\leq\eta_R\in C_c^\infty$ increases to $1$, then
\begin{align*}
 \langle\bar f(t),\eta_R\rangle
 =\lim_{j\to\infty}\langle\bar f_j(t),\eta_R\rangle\leq M_*
\end{align*}
for every $t$. Monotone convergence yields
$\|\bar f(t)\|_1\leq M_*$ for all $t$ on the same event. In particular,
the initial identity above is an equality of $L^1$ functions.

We now pass term by term in the weak equation, keeping the upper time
limit uniform. The artificial viscosity vanishes because
\begin{align}\label{n-main-vanishing-viscosity}
 \sup_{t\leq T}\left|
 \frac1{n_j}\int_0^t\int\bar f_j\Delta_x\varphi_\ell\,dz\,ds
 \right|
 \leq\frac{C_{T,\ell}}{n_j}
 \|\bar f_j\|_{L^2(0,T;L^2(\operatorname{supp}\varphi_\ell))}
 \longrightarrow0.
\end{align}
The diffusion, transport, and damping terms converge uniformly in their
upper time limit by \eqref{n-main-Skorokhod-convergences} and local
strong convergence. Here, for every scalar sequence
$h_j\rightharpoonup h$ in $L^2(0,T)$, the primitives
$\int_0^\cdot h_j$ converge uniformly: weak convergence gives pointwise
convergence, while the uniform $L^2$ bound gives a common
$1/2$-H\"older modulus. This observation handles the diffusion term
containing $\nabla_v\bar f_j$. For the interaction, write
\begin{align*}
 \bar f_jb(\bar f_j)-\bar fb(\bar f)
 =(\bar f_j-\bar f)b(\bar f_j)
   +\bar f\,b(\bar f_j-\bar f).
\end{align*}
The first term converges locally in $L^1$ because
$\|b(\bar f_j)\|_\infty\leq\|V\|_\infty M_*$, and the second converges
by \eqref{n-main-global-L1}. Hence both primitives converge uniformly.

For the It\^o correction we retain the divergence form
\begin{align*}
 {\color{black}F_1\sigma'(\bar f_j)^2\nabla_v\bar f_j}.
\end{align*}
Set
$a_j=\sigma'(\bar f_j)^2$ and $a=\sigma'(\bar f)^2$ and use
\begin{align*}
 a_j\nabla_v\bar f_j-a\nabla_v\bar f
 =(a_j-a)\nabla_v\bar f_j
   +a(\nabla_v\bar f_j-\nabla_v\bar f).
\end{align*}
The first term has a primitive converging uniformly by strong local
$L^2$ convergence and the gradient bound. After pairing the second with
$F_1\nabla_v\varphi_\ell$, its scalar integrand converges weakly in
$L^2(0,T)$; its primitives converge pointwise and are uniformly
$1/2$-H\"older, and hence converge uniformly. We have proved
\begin{align}\label{n-main-deterministic-residual-convergence}
 \sup_{t\leq T}\left|
 \int_0^t\mathcal D_{n_j}^\ell(\bar f_j(s))\,\dif s
 -\int_0^t\mathcal D^\ell(\bar f(s))\,\dif s
 \right|\longrightarrow0,
\end{align}
where $\mathcal D^\ell$ is obtained by deleting the artificial-viscosity
term from $\mathcal D_n^\ell$. Together with
\eqref{n-main-all-time-endpoint}, this identifies $\bar M^\ell$ as the
residual of the limiting equation for every time.

It remains to identify the stochastic term, following the martingale
characterization in the original proof. Define
\begin{align*}
 G_k^\ell(t):=\int_{\mathbb R^{2d}}
 \sigma(\bar f(t,z))f_k(z)\nabla_v\varphi_\ell(z)\,dz,
\end{align*}
and let $\bar G_{j,k}^\ell$ denote \eqref{n-main-noise-coefficient} with
$\bar f_j$ in place of $f_n$. The Lipschitz continuity of $\sigma$ and
$F_1\in L^\infty$ imply
\begin{align}\label{n-main-noise-coefficient-convergence}
 \int_0^T\sum_{k\geq1}|\bar G_{j,k}^\ell-G_k^\ell|^2\,\dif t
 \leq C_\ell\int_0^T
 \|\bar f_j-\bar f\|_{L^2(\operatorname{supp}\varphi_\ell)}^2\,\dif t
 \longrightarrow0
\end{align}
almost surely. Bound \eqref{n-main-bracket-bound} gives the required
uniform integrability. Let $\bar{\mathcal G}_t$ be the augmented
canonical filtration generated up to time $t$ by both limiting solution
coordinates, their residuals, the common initial datum, and the Brownian
coordinates. Passing the Brownian increment, quadratic-variation, and
cross-variation identities against bounded continuous functions of the
canonical past shows that $(\bar B^k)_k$ remain independent
$\mathbb R^d$-valued Brownian motions for this filtration. The same
argument gives that
\begin{align*}
 \bar M^\ell,\qquad
 \bar M^\ell\bar B^{k,a}-\int_0^\cdot G_{k,a}^\ell(s)\,\dif s,
 \qquad
 (\bar M^\ell)^2-
 \int_0^\cdot\sum_{k\geq1}|G_k^\ell(s)|^2\,\dif s
\end{align*}
are martingales. Therefore
\begin{align*}
 \bar M^\ell(t)=\sum_{k\geq1}
 \int_0^tG_k^\ell(s)\cdot\dif\bar B_s^k,
 \qquad t\in[0,T].
\end{align*}
We have obtained the weak formulation of \eqref{SPDE-6}, on one event
and for every time, first for the countable family and then for every
compactly supported smooth test function.

The preceding argument already gives nonnegativity and the initial trace.
Applying the cutoff argument of Lemma~\ref{lem-L1-es} to the limiting
equation yields mass conservation for every time. The energy and moment
estimates pass by lower semicontinuity. More precisely, convergence in
$H^{-2}_{\rm loc}$ at each fixed time makes the extended global $L^2$
norm lower semicontinuous, while the gradient estimate passes through
the weak convergence in \eqref{n-main-Skorokhod-convergences}. The convex lower-semicontinuity
of
\begin{align*}
 (a,\xi)\longmapsto
 \begin{cases}|\xi|^2/a,&a>0,\\0,&(a,\xi)=(0,0),\\
 +\infty,&\text{otherwise},\end{cases}
\end{align*}
together with \eqref{uniform-n-Fisher}, gives the Fisher-information
bound for $\bar f$. Proposition~\ref{equivalence} consequently makes
both the primed and unprimed limits renormalized kinetic solutions. They
have the same initial datum and are driven by the same Brownian motions,
so Theorem~\ref{Uniqueness-spde} shows that they coincide. Their
$C([0,T];H^{-2}_{\rm loc})$ modifications therefore coincide at every
time. Every joint subsequential limit is supported on the diagonal of
$\mathbb X$, and the Gy\"ongy--Krylov criterion gives convergence in
probability on the original stochastic basis. This proves existence and
uniqueness of the weak solution asserted in Theorem~\ref{thm-existence}.
Indeed, an almost surely convergent subsequence on the original basis has
an adapted $C([0,T];H^{-2}_{\rm loc})$ limit; its continuous modification
is progressively measurable, so the limiting process has all the
measurability properties required in Definition~\ref{def-weaksolution}.
\end{proof}

{\color{black}
We conclude this section with the entropy estimate needed in the
$L^1$-contraction argument. This estimate is asserted only for the fixed
regular coefficient of Theorem~\ref{thm-existence}; its constant is not
uniform over the square-root approximations used in Section~\ref{sec-6}.

\begin{proposition}[Entropy dissipation]\label{prp-entropy-dissipation-regular}
Assume that $f_0$, $V$, and $\sigma$ satisfy Assumptions
\ref{Assump-initialdata-regular}, \ref{Assump-ker}, and
\ref{Assump-sigma-n}, respectively. Then every weak solution $f$ of
\eqref{SPDE-6} in the sense of Definition~\ref{def-weaksolution} satisfies
\begin{align}\label{entropy-dissipation-regular}
 &\mathbb E\sup_{t\in[0,T]}
 \left|\int_{\mathbb R^{2d}}\Psi(f(t,z))\,dz\right|
 +\mathbb E\int_0^T\int_{\mathbb R^{2d}}
 \frac{|\nabla_vf|^2}{f}\,dz\,dt
 \le C,
\end{align}
where the quotient is defined to be zero on $\{f=0\}$ and
$C$ depends on $T,d,f_0,V,\sigma,F_1$, and $F_3$. In particular,
\eqref{Fisher-bound-regular} holds.
\end{proposition}

\begin{proof}
For $\delta\in(0,1)$, define
\begin{align*}
 \Psi_\delta(r)&:=\int_0^r\log(s+\delta)\,ds,
 &S_\delta(r)&:=\int_0^r\frac{\sigma(s)}{s+\delta}\,ds,
 &A_\delta(r)&:=\int_0^r
       \frac{\sigma(s)\sigma'(s)}{s+\delta}\,ds.
\end{align*}
Then $\Psi_\delta''(r)=(r+\delta)^{-1}$ and
$\Psi_\delta(r)\to\Psi(r)$ as $\delta\downarrow0$.

Before applying It\^o's formula, choose a family increasing pointwise in $L$,
$\vartheta_L\in C_c^\infty([0,\infty))$ such that
$0\leq\vartheta_L\leq1$, $\vartheta_L=1$ on $[0,L]$, and
$\vartheta_L=0$ on $[2L,\infty)$. Define the value-truncated entropy by
\begin{align*}
 \Psi_{\delta,L}'(r)
 &:=\log\delta+\int_0^r\frac{\vartheta_L(s)}{s+\delta}\,ds,
 &
 \Psi_{\delta,L}(r)&:=\int_0^r\Psi_{\delta,L}'(s)\,ds.
\end{align*}
Then $\Psi_{\delta,L}'$ and $\Psi_{\delta,L}''$ are bounded and
$\Psi_{\delta,L}=\Psi_\delta$ on $[0,L]$.  We apply It\^o's formula first
with $\Psi_{\delta,L}$ and, for fixed $\delta$ and spatial cutoff, let
$L\to\infty$.  The $L^1\cap L^2$ bounds and dominated convergence justify
this limit in every term below.

Let $\alpha_R(x,v)=\alpha^1_{R^3}(x)\alpha^2_R(v)$ be the cutoff used
above. {\color{black}We apply the generalized It\^o formula in the variational
Gelfand triple to
$\int_{\mathbb R^{2d}}\Psi_{\delta,L}(f(t))\alpha_R$.  More precisely, we
first mollify the equation only in the spatial variables and then remove
the mollification, using the same commutator convergence as in the display
preceding \cite[(2.26)]{FG25} and the convergence in
\cite[(2.26)]{FG25}.  No time convolution is used, since it would destroy
adaptedness.}
After the preceding $L\to\infty$ limit, the resulting terms contain
$\Psi_\delta$.
The deterministic diffusion contributes
\begin{align*}
 -\int_0^t\int_{\mathbb R^{2d}}
 \frac{|\nabla_vf|^2}{f+\delta}\alpha_R\,dz\,ds
 -\int_0^t\int_{\mathbb R^{2d}}
 \log(f+\delta)\nabla_vf\cdot\nabla_v\alpha_R\,dz\,ds.
\end{align*}
The transport and interaction terms are integrated by parts. Their
interior contributions vanish because $v$ is independent of $x$ and
$V\ast_x\rho$ is independent of $v$.

We next record the cancellation generated by the conservative noise. We
use the It\^o correction in the divergence form
\begin{align*}
 {\color{black}\frac12\nabla_v\cdot\bigl(
 F_1\sigma'(f)^2\nabla_vf\bigr)}.
\end{align*}
The quadratic variation in It\^o's formula is
\begin{align*}
 \frac12\int_0^t\int_{\mathbb R^{2d}}
 \Psi_\delta''(f)\alpha_R
 {\color{black}\bigl(F_1\sigma'(f)^2|\nabla_vf|^2
 +F_3\sigma(f)^2\bigr)}\,dz\,ds.
\end{align*}
The term containing $F_1\sigma'(f)^2|\nabla_vf|^2$ cancels exactly with
the corresponding term obtained from the It\^o correction. In treating
the remaining cutoff contribution, we use the standing identity
${\color{black}\nabla_vF_1=0}$.

For fixed $\delta$ and the fixed regular coefficient $\sigma$, all the
primitives occurring in the cutoff terms have at most linear or quadratic
growth. The $L^1\cap L^2$ bounds, the moment estimate, and the properties
of $\alpha_R$ therefore allow us to let $R\to\infty$. In particular, the
terms produced by $-v\cdot\nabla_xf$ and
$-\nabla_v\cdot(fV\ast_x\rho)$ vanish. The damping term satisfies
\begin{align*}
 f\Psi_\delta'(f)-\Psi_\delta(f)
 =\int_0^f\frac{r}{r+\delta}\,dr
 =f-\delta\log\left(1+\frac f\delta\right).
\end{align*}
Consequently,
\begin{align}\label{regular-entropy-delta-identity}
 &\int_{\mathbb R^{2d}}\Psi_\delta(f(t))\,dz
 +\int_0^t\int_{\mathbb R^{2d}}
       \frac{|\nabla_vf|^2}{f+\delta}\,dz\,ds                 \notag\\
 &\quad=\int_{\mathbb R^{2d}}\Psi_\delta(f_0)\,dz
 +M_t^\delta
 +\frac12\int_0^t\int_{\mathbb R^{2d}}
       \frac{\sigma(f)^2}{f+\delta}F_3\,dz\,ds                \\
 &\qquad\quad
 +d\int_0^t\int_{\mathbb R^{2d}}
 \left(f-\delta\log\left(1+\frac f\delta\right)\right)
 \,dz\,ds.                                                    \notag
\end{align}
Here
\begin{align*}
 M_t^\delta
 :=-\sum_{k\geq1}\int_0^t
 \left(\int_{\mathbb R^{2d}}
 S_\delta(f(s,z))\nabla_vf_k(z)\,dz\right)\cdot dB_s^k.
\end{align*}

The structural bound $|\sigma(r)|\leq c\sqrt r$ gives, uniformly in
$\delta\in(0,1)$,
\begin{align*}
 |S_\delta(r)|
 \leq c\int_0^r\frac{\sqrt s}{s+\delta}\,ds
 \leq2c\sqrt r.
\end{align*}
It follows from the Cauchy--Schwarz inequality and mass conservation that
\begin{align*}
 \frac{d}{dt}\langle M^\delta\rangle_t
 &\leq
 \left(\int_{\mathbb R^{2d}}S_\delta(f)^2\,dz\right)
 \left(\int_{\mathbb R^{2d}}F_3\,dz\right)\\
 &\leq C\|f_0\|_{L^1}\|F_3\|_{L^1}.
\end{align*}
The Burkholder--Davis--Gundy inequality thus yields
\begin{align}\label{regular-entropy-martingale-bound}
 \mathbb E\sup_{t\in[0,T]}|M_t^\delta|\leq C,
\end{align}
with a constant independent of $\delta$. Moreover,
\begin{align*}
 0\leq\frac{\sigma(f)^2}{f+\delta}\leq c^2,
 \qquad
 0\leq f-\delta\log\left(1+\frac f\delta\right)\leq f.
\end{align*}
The last two terms in \eqref{regular-entropy-delta-identity} are therefore
bounded by $C T\|F_3\|_{L^1}$ and
$dT\|f_0\|_{L^1}$, respectively.

It remains to control the negative part of the entropy on the whole
space. Since $\Psi_\delta\geq\Psi$, estimate
\eqref{entropy-negative-part-bound} and the moment estimate imply
\begin{align*}
 \sup_{\delta\in(0,1)}
 \mathbb E\sup_{t\in[0,T]}
 \int_{\mathbb R^{2d}}\Psi_\delta(f(t))^-\,dz\leq C.
\end{align*}
The positive part of $\Psi_\delta(f_0)$ is bounded by a constant times
$f_0+f_0^2$. Combining these facts with
\eqref{regular-entropy-delta-identity} and
\eqref{regular-entropy-martingale-bound}, we obtain
\begin{align*}
 \sup_{\delta\in(0,1)}\Bigg[
 \mathbb E\sup_{t\in[0,T]}
 \left|\int_{\mathbb R^{2d}}\Psi_\delta(f(t))\,dz\right|
 +\mathbb E\int_0^T\int_{\mathbb R^{2d}}
 \frac{|\nabla_vf|^2}{f+\delta}\,dz\,dt\Bigg]\leq C.
\end{align*}
Finally, $\Psi_\delta(f)\to\Psi(f)$, whereas
$(f+\delta)^{-1}|\nabla_vf|^2$ increases to
$f^{-1}|\nabla_vf|^2$. Fatou's lemma and the monotone convergence theorem
prove \eqref{entropy-dissipation-regular}. The Sobolev chain rule for
nonnegative functions gives
\begin{align*}
 \frac{|\nabla_vf|^2}{f}=4|\nabla_v\sqrt f|^2
\end{align*}
almost everywhere, with both sides defined to be zero on $\{f=0\}$.
\end{proof}
}

}

\section{Taming the square-root coefficient}\label{sec-6}
\begingroup\color{black}
{\color{black}In this section, we construct a probabilistically weak solution of the non-interacting Vlasov--Fokker--Planck--Dean--Kawasaki equation with square-root coefficient:}
{\color{black}
\begin{equation}\label{SPDE-squareroot}
df=\Big(\Delta_vf-v\cdot\nabla_xf+\nabla_v\cdot(vf)\Big)dt-\nabla_v\cdot(\sqrt{f}\circ dW_F(t)).
\end{equation}
Because the Stratonovich correction is singular at zero, equation
\eqref{SPDE-squareroot} is only a shorthand: the equation is defined by the
distribution-in-time kinetic identity below.}
{\color{black}The following definition records the weaker solution concept appropriate to the local compactness available in this section.}

{\color{black}
\begin{definition}\label{def-prob-weak-square-root}
A \emph{probabilistically weak renormalized kinetic solution} of \eqref{SPDE-squareroot} consists of a filtered probability space
$(\widetilde\Omega,\widetilde{\mathcal F},
(\widetilde{\mathcal F}_t)_{t\in[0,T]},\widetilde{\mathbb P})$,
a correlated Wiener process
$\widetilde W_F=\sum_{k\geq1}\widetilde B^k f_k$ adapted to this filtration,
an $\widetilde{\mathcal F}_0$-measurable initial random variable
$\widetilde f_0$, a nonnegative progressively measurable process $f$, and a
nonnegative random Radon measure $p$ on
$\mathbb R^{2d}\times(0,\infty)\times[0,T]$.  The stochastic basis satisfies
the usual conditions, the processes $(\widetilde B^k)_{k\geq1}$ are
independent $d$-dimensional Brownian motions, and
\begin{align*}
 \operatorname{Law}_{\widetilde{\mathbb P}}(\widetilde f_0)
 =\operatorname{Law}_{\mathbb P}(f_0),
\end{align*}
\begin{align*}
 f&\in L^1\bigl(\widetilde\Omega\times[0,T];
 L^1_{\rm loc}(\mathbb R^{2d})\bigr)
 \cap L^2\bigl(\widetilde\Omega\times[0,T];L^2(\mathbb R^{2d})\bigr),\\
 \nabla_vf&\in L^2\bigl(\widetilde\Omega\times[0,T];
 L^2(\mathbb R^{2d})\bigr).
\end{align*}
For $\widetilde{\mathbb P}$-almost every $\omega$ and almost every
$t\in[0,T]$,
\begin{align}\label{square-root-mass-inequality}
 \int_{\mathbb R^{2d}}f(t,z)\,dz
 \leq \int_{\mathbb R^{2d}}\widetilde f_0(z)\,dz
 <\infty.
\end{align}
Moreover, $p$ is locally finite almost surely and
\begin{align*}
 \delta_{f(t,z)}(d\zeta)|\nabla_vf(t,z)|^2\,dz\,dt
 \leq p(dz,d\zeta,dt)
\end{align*}
as measures.  For every
$\varphi\in C_c(\mathbb R^{2d}\times(0,\infty))$, the left-continuous
adapted process
\begin{align*}
 t\longmapsto
 \int_{\mathbb R^{2d}\times(0,\infty)\times[0,t)}\varphi(z,\zeta)
 \,p(dz,d\zeta,ds)
\end{align*}
is predictable.

	To state the equation, put
	$\chi_f(z,\zeta,t):=\boldsymbol1_{\{0<\zeta<f(t,z)\}}$.  For every
	$\psi\in C_c^\infty(\mathbb R^{2d}\times(0,\infty))$ and
	$\eta\in C_c^\infty([0,T))$, the following identity holds almost surely:
	\begin{align}\label{square-root-distributional-kinetic-identity}
	 -\int_0^T\eta'(t)\langle\chi_f(t),\psi\rangle\,dt
	={}&\eta(0)\langle\chi_{\widetilde f_0},\psi\rangle\notag\\
	 &+\int_0^T\eta(t)\Bigg[
	 -\int_{\mathbb R^{2d}}\nabla_vf\cdot(\nabla_v\psi)(z,f)\,dz
	 +\int_{\mathbb R^{2d}}\!\int_0^\infty
	 v\chi_f\cdot\nabla_x\psi\,d\zeta\,dz\notag\\
	 &\qquad+\int_{\mathbb R^{2d}}\nabla_v\cdot(vf)\psi(z,f)\,dz
	 -\frac18\int_{\mathbb R^{2d}}F_1f^{-1}\nabla_vf\cdot
	 (\nabla_v\psi)(z,f)\,dz\notag\\
	 &\qquad{\color{black}+\frac12\int_{\mathbb R^{2d}}
	 F_3f(\partial_\zeta\psi)(z,f)\,dz}\Bigg]dt\notag\\
	 &+\sum_{k\geq1}\int_0^T\eta(t)\mathcal G_\psi^k(f(t))\cdot
	 d\widetilde B^k(t)
	 -\int_{\mathbb R^{2d}\times(0,\infty)\times[0,T]}
	 \eta(t)\partial_\zeta\psi(z,\zeta)\,p(dz,d\zeta,dt),
	\end{align}
	where $\langle\chi,\psi\rangle$ denotes integration in $(z,\zeta)$.
	The factor $f^{-1}$ is evaluated only where $f$ lies in the compact
	$\zeta$-support of $\psi$, which is separated from zero.  The stochastic
	coefficient in \eqref{square-root-distributional-kinetic-identity} is
	defined by
\begin{align}\label{square-root-stochastic-integrand}
 \sum_{k\geq1}\int_0^T\eta(t)\,\mathcal G_\psi^k(f(t))
 \cdot d\widetilde B^k(t),\qquad
 \mathcal G_\psi^k(f)
 &:=-\int_{\mathbb R^{2d}}\psi(z,f)
       \nabla_v\bigl(\sqrt f\,f_k\bigr)\,dz .
\end{align}
The last expression is understood by spatial integration by parts.  If
\begin{align*}
 H_\psi(z,r):=\int_0^r\psi(z,a)\frac{da}{2\sqrt a},
\end{align*}
then, equivalently,
\begin{align}\label{square-root-stochastic-primitive}
 \mathcal G_\psi^k(f)
 ={}&\int_{\mathbb R^{2d}}f_k(z)
       \int_0^{f(z)}\nabla_v\psi(z,a)\frac{da}{2\sqrt a}\,dz\notag\\
 &+\int_{\mathbb R^{2d}}
       \bigl(H_\psi(z,f(z))-\psi(z,f(z))\sqrt{f(z)}\bigr)
       \nabla_v f_k(z)\,dz.
\end{align}
This formula contains no derivative of $f$ and makes the stochastic
integral well defined.  All singular coefficients in the remaining terms
are evaluated on the compact $\zeta$-support of $\psi$, which is separated
from zero.  No continuity in time, global $L^1$ convergence, mass equality,
or uniqueness is included in this definition.
\end{definition}}

{\color{black}Precisely, we have the following result.}
\begin{theorem}\label{thm-existence-squareroot}
	{\color{black}Assume that $f_0$ satisfies Assumption~\ref{Assump-initialdata-entropy}. Then there exists a probabilistically weak renormalized kinetic solution of \eqref{SPDE-squareroot} with initial data $f_0$ in the sense of Definition~\ref{def-prob-weak-square-root}.}
\end{theorem}
{\color{black}Let $(\sigma_n)_{n\geq1}$ be the smooth approximations from
Lemma~\ref{lem-sigma}.  We choose the construction in
\cite[Lemma~5.18]{FG24} so that the constants in the following bounds are
independent of $n$:
\begin{align}\label{square-root-uniform-coefficients}
 |\sigma_n(r)|&\leq C\sqrt r,&
 |\sigma_n(r)\sigma_n'(r)|&\leq C,&
 |\sigma_n'(r)|^2&\leq \frac{C}{r},
 \qquad r>0.
\end{align}
The bound on $(\sigma_n\sigma_n')'$ in
Assumption~\ref{Assump-sigma-n} is allowed to depend on $n$; it cannot be
uniform as the regularizations approach the square root at zero.  On every
compact interval $I\Subset(0,\infty)$, Lemma~\ref{lem-sigma} gives
$\sigma_n\to\sqrt{\cdot}$ in $C^1(I)$.  No uniform estimate or convergence
of $(\sigma_n\sigma_n')'$ will be used.  We use \eqref{SPDE-6} with
$V=0$ and coefficient $\sigma_n$ as an approximation of
\eqref{SPDE-squareroot}. No moment assumption on the initial datum is
imposed.}

The proof of Theorem~\ref{thm-existence-squareroot} is organized into five parts:
\begin{itemize}
  \item \textbf{Subsection~\ref{subsec-6-1}.} We establish $L^2$-estimates for \eqref{SPDE-6} with diffusion coefficient $(\sigma_n(\cdot))_{n\geq1}$, uniformly in $n\geq1$.
  
  \item \textbf{Subsection~\ref{subsec-6-2}.} We introduce a sequence of truncated equations associated with \eqref{SPDE-6}, which regularize the singularity near the zero-value of the solution.
  
  \item \textbf{Subsection~\ref{subsec-6-3}.} We derive uniform Besov regularity estimates for the truncated equations.
  
  \item {\color{black}\textbf{Subsection~\ref{subsec-6-4}.} We establish tightness in $L^1([0,T];L^1_{\rm loc}(\mathbb R^{2d}))$ by combining the Besov estimates with the truncation metric of \cite{FG24}.}
   
  \item {\color{black}\textbf{Subsection~\ref{subsec-6-5}.} We pass to the limit on a new stochastic basis and prove probabilistically weak existence.}
\end{itemize}

{\color{black}Throughout this section, $f_0$ satisfies Assumption~\ref{Assump-initialdata-entropy}. Let $\vartheta_n\in C_c^\infty(\mathbb R^{2d})$ satisfy $0\leq\vartheta_n\leq1$ and $\vartheta_n=1$ on $B_n$, and let $(\eta_{1/n})_{n\geq1}$ be standard nonnegative convolution kernels. Set
\begin{align}\label{square-root-initial-approximation}
 f_{0,n}:=(\vartheta_nf_0)*\eta_{1/n}.
\end{align}
Then $f_{0,n}\geq0$, every $f_{0,n}$ has finite spatial moments and hence is admissible for Theorem~\ref{thm-existence}, while
\begin{align*}
 f_{0,n}\longrightarrow f_0
 \quad\text{in }L^1(\mathbb R^{2d})\cap L^2(\mathbb R^{2d})
 \quad\text{almost surely}.
\end{align*}
Moreover, the $L^1$ and $L^2$ bounds used below are uniform in $n$; no estimate is allowed to depend on the moments of $f_{0,n}$.}

\subsection{\textcolor{black}{Uniform $L^2$ estimates}}\label{subsec-6-1}
{\color{black}We first record precisely the estimate that is uniform in the
square-root regularization.  Notice that it is a bound of
$\sup_t\mathbb E\|f_n(t)\|_2^2$, rather than the stronger quantity
$\mathbb E\sup_t\|f_n(t)\|_2^2$.}

\begin{proposition}\label{prop-square-root-uniform-L2}
{\color{black}Assume that $f_0$ satisfies
Assumption~\ref{Assump-initialdata-entropy}, and choose $(\sigma_n)_{n\geq1}$
as in \eqref{square-root-uniform-coefficients}.  For each $n$, let $f_n$ be
the weak solution supplied by Theorem~\ref{thm-existence} for $V=0$,
initial datum $f_{0,n}$, and coefficient $\sigma_n$.  There is a constant
$C_T$, independent of $n$ and of every spatial moment of $f_{0,n}$, such
that}
\begin{align}\label{square-root-uniform-L2}
\sup_{n \geq 1} \left[
\sup_{t \in [0,T]}\mathbb{E}\|f_n(t)\|_{L^2(\mathbb R^{2d})}^2
+ \mathbb{E}\int_0^T\|\nabla_v f_n(t)\|_{L^2(\mathbb R^{2d})}^2\,dt
\right]
\leq C_T\left(1+\mathbb E\|f_0\|_{L^2}^2
+\mathop{\rm ess\,sup}_{\omega\in\Omega}\|f_0(\omega)\|_{L^1}\right).
\end{align}
\end{proposition}

{\color{black}
\begin{proof}
Apply It\^o's formula to $\frac12\|f_n(t)\|_2^2$, first with the
spatial cutoffs used in the proof of Proposition~\ref{prp-uniforml2} and
then let the $x$-radius and the $v$-radius tend to infinity in this order.
This iterated limit treats the transport cutoff without using a spatial
moment.  The transport term then vanishes, the damping term is bounded by
$C\|f_n\|_2^2$, and the martingale has zero expectation.

We spell out the noise cancellation because a bound involving
$\|\sigma_n'\|_\infty$ would not be uniform in $n$.  Put
\begin{align*}
 B_n(r):=\int_0^r s|\sigma_n'(s)|^2\,ds.
\end{align*}
The bounds in \eqref{square-root-uniform-coefficients} give
$0\leq B_n(r)\leq Cr$, uniformly in $n$.  In the
whole-space identity the $F_1|\sigma_n'|^2|\nabla_vf_n|^2$ term from the
It\^o correction cancels its quadratic-variation counterpart. At the cutoff level, the potentially
nonuniform $F_1$ remainder is
\begin{align*}
 -\frac12\int F_1 f_n|\sigma_n'(f_n)|^2
             \nabla_vf_n\cdot\nabla_v\alpha_R
 =\frac12\int B_n(f_n)\nabla_v\cdot(F_1\nabla_v\alpha_R).
\end{align*}
Using ${\color{black}\nabla_vF_1=0}$, $F_1\in L^\infty$, and
$B_n(f_n)\leq Cf_n$, this converges to zero uniformly in $n$. Hence the only
remaining zeroth-order noise contribution is bounded uniformly in $n$ by
\begin{align*}
 C\int_{\mathbb R^{2d}}F_3(z)\sigma_n(f_n(z))^2\,dz
 \leq C\|F_3\|_{L^\infty}\|f_n\|_{L^1},
\end{align*}
where we used \eqref{square-root-uniform-coefficients}.  Mass preservation
for the regular approximating equation and Young's convolution inequality
give
$\|f_n(t)\|_1=\|f_{0,n}\|_1\leq\|f_0\|_1$.  Gronwall's inequality now
proves \eqref{square-root-uniform-L2}, with no dependence on
$\|\sigma_n'\|_\infty$ or on a spatial moment of $f_{0,n}$.
\end{proof}}

\subsection{The truncation equations}\label{subsec-6-2} 
For each fixed $\delta \in (0,1)$, let $\psi_{\delta} \in C^{\infty}([0,\infty))$ be a smooth, nondecreasing cutoff function satisfying $0 \leq \psi_{\delta} \leq 1$ and
\begin{align*}
 \psi_{\delta}(\zeta) =
  \begin{cases}
    1, & \text{if } \zeta \in [\delta, +\infty), \\
    0, & \text{if } \zeta \in [0, \delta/2], \\
    \text{smoothly interpolated}, & \text{otherwise}.
  \end{cases}
\end{align*}
By construction, there exists a constant $c > 0$, independent of $\delta$, such that
$$
\left|\psi_{\delta}'(\zeta)\right| \leq \frac{c}{\delta} \quad \text{for all } \zeta \in [0,\infty).
$$

\begin{definition}\label{def-hdelta}
For every $\delta \in (0,1)$, define the function $h_{\delta} \in C^{\infty}([0,\infty))$ by
$$
h_{\delta}(\zeta) := \psi_{\delta}(\zeta) \zeta, \quad \forall \zeta \in [0,\infty).
$$
{\color{black}When It\^o's formula is applied below, $h_\delta$ denotes any
smooth extension to $\mathbb R$ with the same values on $[0,\infty)$;
the choice of the extension is immaterial because $f_n\geq0$.}
\end{definition}

From the above definition, it immediately follows that the support of the derivative $h_{\delta}'$ is contained in the interval $[\frac{\delta}{2}, \infty)$. Moreover, we have the explicit expression
\begin{align}\label{eq-7.1}
{\color{black}|h_{\delta}'(\zeta)| = |\psi_{\delta}'(\zeta) \zeta + \psi_{\delta}(\zeta)| \leq c I_{\{\zeta \geq \frac{\delta}{2}\}},}
\end{align}
where $c > 0$ is a constant independent of $\delta$. Similarly, the second derivative satisfies
\begin{align}\label{eq-7.2}
{\color{black}|h_{\delta}''(\zeta)| = |\psi_{\delta}''(\zeta) \zeta + 2 \psi_{\delta}'(\zeta)| \leq c(\delta) I_{\{\frac{\delta}{2} \leq \zeta \leq \delta\}},}
\end{align}
for some constant $c(\delta) > 0$ that depends on $\delta$.

{\color{black}For each $n\geq1$, let $f_n$ be the weak solution of \eqref{SPDE-6} with $V=0$, initial datum $f_{0,n}$, and diffusion coefficient $\sigma_n$. Our goal is to derive the evolution equation satisfied by $h_\delta(f_n)$.}

Let $(\eta_{\gamma})_{\gamma \in (0,1)}$ be a sequence of standard convolution kernels on $\mathbb{R}^d_x$. For every test function $\varphi\in C^{\infty}_c(\mathbb{R}^{2d})$, applying It\^o's formula to $\int_{\mathbb{R}^{2d}} h_{\delta}(\eta_{\gamma} \ast f_n)\varphi$, and then using the integration by parts formula and passing to the limits $\gamma \to 0$, we derive the following distributional identity 
\begin{align*}
d h_{\delta}(f_n) &= \bigl( \Delta_v - v \cdot \nabla_x \bigr) h_{\delta}(f_n)\, dt - h_{\delta}''(f_n) \left| \nabla_v f_n \right|^2 \, dt\\
&\quad + h_{\delta}'(f_n) \nabla_v \cdot (v f_n) \, dt \\
&\quad - h_{\delta}'(f_n) \nabla_v \cdot \bigl( \sigma_n(f_n) \, dW_F \bigr)  + {\color{black}\frac{1}{2} \nabla_v \cdot \bigl( h_{\delta}'(f_n) \sigma_n'(f_n)^2 \nabla_v f_n F_1 \bigr) \, dt+\frac{1}{2}h_{\delta}''(f_n)\sigma_n(f_n)^2\, F_3\,dt}.
\end{align*}

For the damping term, we have the following identities:
{\color{black}
\begin{align*}
h_{\delta}'(f_n) \nabla_v \cdot (v f_n)
&= d\,h_{\delta}'(f_n) f_n + h_{\delta}'(f_n) v \cdot \nabla_v f_n \\
&= d\,h_{\delta}'(f_n) f_n + v \cdot \nabla_v h_{\delta}(f_n) \\
&= d\,h_{\delta}'(f_n) f_n - d\,h_{\delta}(f_n) + \nabla_v \cdot \bigl( v h_{\delta}(f_n) \bigr).
\end{align*}}

Collecting all the above computations, the evolution equation for $h_{\delta}(f_n)$ can be expressed as
\begin{align}\label{SPDE-hdelta}
d h_{\delta}(f_n) &= \bigl( \Delta_v - v \cdot \nabla_x \bigr) h_{\delta}(f_n) \, dt + \nabla_v \cdot \bigl( v h_{\delta}(f_n) \bigr) \, dt \notag \\
&\quad - h_{\delta}''(f_n) |\nabla_v f_n|^2 \, dt + {\color{black}d\,h_{\delta}'(f_n) f_n \, dt - d\,h_{\delta}(f_n) \, dt} - h_{\delta}'(f_n) \nabla_v \cdot \bigl( \sigma_n(f_n) dW_F \bigr) \notag \\
&\quad {\color{black}+ \frac{1}{2} \nabla_v \cdot \bigl( h_{\delta}'(f_n) \sigma_n'(f_n)^2 \nabla_v f_n F_1 \bigr) \, dt+\frac{1}{2}h_{\delta}''(f_n)\sigma_n(f_n)^2\, F_3\,dt}.
\end{align}

\subsection{Uniform Besov regularity estimates}\label{subsec-6-3}
In this section, we fix an arbitrary $\delta \in (0,1)$ and introduce the notation $u_n := h_{\delta}(f_n)$ for each $n \geq 1$. With this notation, the equation \eqref{SPDE-hdelta} can be equivalently rewritten in the form of \eqref{lkSPDE}. For the reader's convenience, we reformulate the equation as follows: 
\begin{align}\label{SPDE-un}
du_n = (\Delta_v - v \cdot \nabla_x) u_n \, dt + &\nabla_v \cdot (v u_n) \, dt \notag \\
&\quad + g_{1,n} \, dt + g_{2,n} \, dt +  h_{1,n} dW_F+h_{2,n}\nabla_v\cdot dW_F,
\end{align}
{\color{black}The terms $g_{1,n}$, $g_{2,n}$, $h_{1,n}$, and $h_{2,n}$ are defined by
\begin{align*}
g_{1,n} &:=\frac12\nabla_v\cdot\bigl(
 h_\delta'(f_n)\sigma_n'(f_n)^2\nabla_vf_nF_1\bigr),\\
g_{2,n} &:=-h_\delta''(f_n)|\nabla_vf_n|^2
 +d\,h_\delta'(f_n)f_n-d\,h_\delta(f_n)
 +\frac12h_\delta''(f_n)\sigma_n(f_n)^2F_3,\\
h_{1,n} &:=-h_\delta'(f_n)\nabla_v\sigma_n(f_n),\qquad
h_{2,n}:=-h_\delta'(f_n)\sigma_n(f_n).
\end{align*}}

Next, we establish several crucial uniform estimates. By H\"older's inequality, for any $\alpha \leq 0$ and $p \in [1,2]$, it holds that
\begin{align*}
\| h_{\delta}'(f_n) (f_n)^{\alpha} \nabla_v f_n \|_{L^p(\mathbb{R}^{2d})}
&\lesssim \delta^{\alpha} \left\| I_{\{f_n \ge \frac{\delta}{2}\}} |\nabla_v f_n| \right\|_{L^p(\mathbb{R}^{2d})} \\
&\lesssim \delta^{\alpha} \left\| \nabla_v f_n \right\|_{L^2(\mathbb{R}^{2d})} \| I_{\{f_n \ge \frac{\delta}{2}\}} \|_{L^{p^*}(\mathbb{R}^{2d})} \\
&\lesssim \delta^{\alpha  - \frac{1}{p^*}} \left\| \nabla_v f_n \right\|_{L^2(\mathbb{R}^{2d})} \| f_n \|_{L^1(\mathbb{R}^{2d})}^{\frac{1}{p^*}},
\quad \text{with } \frac{1}{p} = \frac{1}{2} + \frac{1}{p^*},
\end{align*}
and similarly,
\begin{align*}
\| h_{\delta}'(f_n) \sqrt{f_n} \|_{L^p(\mathbb{R}^{2d})} \lesssim \| I_{\{f_n \ge \frac{\delta}{2}\}} \sqrt{f_n} \|_{L^p(\mathbb{R}^{2d})} \lesssim \| \sqrt{f_n} \|_{L^2(\mathbb{R}^{2d})} \| I_{\{f_n \ge \frac{\delta}{2}\}} \|_{L^{p^*}(\mathbb{R}^{2d})} \lesssim \delta^{-\frac{1}{p^*}} \| f_n \|_{L^1(\mathbb{R}^{2d})}^{\frac{1}{p}}.
\end{align*}
Here, we have applied Markov's inequality in the last step of both estimates:
\begin{align*}
\| I_{\{f_n \ge \frac{\delta}{2}\}} \|_{L^{p^*}(\mathbb{R}^{2d})} \lesssim \delta^{-\frac{1}{p^*}} \| f_n \|_{L^1(\mathbb{R}^{2d})}^{\frac{1}{p^*}}.
\end{align*}

Furthermore, for any $\beta \in \mathbb{R}$, we note that
\begin{align*}
\| h_{\delta}''(f_n) (f_n)^{\beta} |\nabla_v f_n|^2 \|_{L^1(\mathbb{R}^{2d})} \lesssim_{\delta} \| I_{\{\delta \ge f_n \ge \frac{\delta}{2}\}} |\nabla_v f_n|^2 \|_{L^1(\mathbb{R}^{2d})} \lesssim_{\delta} \left\|\nabla_v f_n \right\|_{L^2(\mathbb{R}^{2d})}^2.
\end{align*}

{\color{black}Therefore, using the preceding estimates,
\eqref{square-root-uniform-coefficients}, and
\eqref{square-root-uniform-L2}, there exists $c_\delta>0$, independent of
$n$, such that for all $p\in[1,2]$,}
\begin{align}\label{0606:06}
\begin{split}
\color{black}\sup_{n \geq 1} \Bigg( & \| g_{1,n} \|_{L^2(\Omega\times[0,T]; \mathbf{B}^{-1}_{p;\theta})}
+ \| g_{2,n} \|_{L^1(\Omega\times[0,T]\times\mathbb{R}^{2d})}  \\
& + \| h_{1,n} \|_{L^2(\Omega\times[0,T]; L^p(\mathbb{R}^{2d}))}
+ \| h_{2,n} \|_{L^2(\Omega\times[0,T]; L^p(\mathbb{R}^{2d}))}  \Bigg) \leq c_{\delta}.
\end{split}
\end{align}

In the sequel, we will employ these uniform bounds to derive uniform Besov regularity estimates for the sequence $\{u_n\}_{n \geq 1}$ governed by \eqref{SPDE-un}.
\begin{theorem}\label{thm-besov-p-es-un}
{\color{black}Assume that $f_0$ satisfies
Assumption~\ref{Assump-initialdata-entropy}, and set
\begin{align*}
 \kappa_0:=1+\mathbb E\|f_0\|_{L^2(\mathbb R^{2d})}^2
 +\operatorname*{ess\,sup}_{\omega\in\Omega}
 \|f_0(\omega)\|_{L^1(\mathbb R^{2d})}<\infty.
\end{align*}
Let $(\sigma_n)_{n\geq1}$ be the approximations from
Lemma~\ref{lem-sigma}, and let $f_n$ be the weak solution of
\eqref{SPDE-6} with $V=0$, initial datum $f_{0,n}$, and coefficient
$\sigma_n$.}

For every fixed $\delta \in (0,1)$, let $h_{\delta}$ be defined as in Definition \ref{def-hdelta}, and set $u_n := h_{\delta}(f_n)$. Then for every
$$
{\color{black}p \in \left(1, \frac{2d+1}{2d}\right) \quad \text{and} \quad
\beta \in \left(0,1\wedge\left(\frac{4d+2}{p} - 4d\right)\right),}
$$
there holds
\begin{align}\label{0425:02}
{\color{black}\sup_{n \geq 1} \mathbb E
 \| u_n \|_{L^p([0,T]; \mathbf{B}^\beta_{p;\theta})}
 \lesssim_{\beta,T,\delta,p,d,F_1,F_3} \kappa_0.}
\end{align}
Furthermore, for every $\alpha \in (0,1/2)$ and any cutoff function $\chi \in C^\infty_c(\mathbb{R}^{2d})$, we have
\begin{align}\label{0521:050}
{\color{black}\sup_{n \geq 1} \mathbb E\| u_n \chi \|_{W^{\alpha,p}([0,T]; \mathbf{B}^{-6}_{p;\theta})}
 \lesssim_{\alpha,T,\delta,\chi,p,d,F_1,F_3} \kappa_0.}
\end{align}
\end{theorem}

{\color{black}
\begin{proof}
The first-moment estimate (S) of Theorem~\ref{thm:Besov} applies to
\eqref{SPDE-un} with $q=1$ and with the artificial $x$-viscosity set equal
to zero.  Indeed,
the only properties of the auxiliary heat semigroup $H_t^n$ used there are
its $L^p$-contraction and its commutation with $P_t$ and the dyadic blocks;
after removing $n^{-1}\Delta_x$, one simply replaces $H_t^n$ by the
identity.  Thus the proof, with this replacement and
\eqref{0606:06}, gives \eqref{0425:02}.  Here the $g_{2,n}$ input is used
in $L^1(\Omega\times[0,T]\times\mathbb R^{2d})$, exactly as prescribed
in (S), while the other source terms have the required outer
$L^2(\Omega)$ bounds; moreover,
$\mathbb E\|h_\delta(f_{0,n})\|_{L^1}\leq \kappa_0$.  The proof yields regularity below
both the exponent $(4d+2)/p-4d$ and the one derivative supplied by the
drift and stochastic convolutions, which explains the minimum with $1$.

For the time estimate, multiply \eqref{SPDE-un} by $\chi$.  The
deterministic terms are bounded in
$L^1(\Omega\times[0,T];\mathbf B^{-6}_{p;\theta})$ by
\eqref{0606:06}; the commutators generated by $\chi$ are supported on a
fixed compact set and obey the same bound.  The
Burkholder--Davis--Gundy inequality and the last two bounds in
\eqref{0606:06} give, for every $\alpha<1/2$, the corresponding
$W^{\alpha,p}$ estimate in first moment for the stochastic primitive.  The
deterministic primitive belongs to $W^{\alpha,p}$ because $\alpha<1/p$
for the present range of $p$, and its expected norm is controlled by the
$L^1(\Omega\times[0,T])$ bound in \eqref{0606:06}.  This proves
\eqref{0521:050}.  In particular, no estimate
for $n^{-1}\Delta_xu_n$ is used.
\end{proof}}

\begin{corollary}\label{cor-besov-1-es-un}
Under the same assumptions as in Theorem \ref{thm-besov-p-es-un}, for any cutoff function $\chi \in C^\infty_c(\mathbb{R}^{2d})$ and every $\beta \in (0,1)$, it holds that
\begin{align}\label{0425:02-cor}
{\color{black}\sup_{n \geq 1} \| u_n \chi \|_{L^1([0,T] \times \Omega; \mathbf{B}^\beta_{1;\theta})}
 \lesssim_{\beta,T,\delta,\chi,d,F_1,F_3} \kappa_0.}
\end{align}
\end{corollary}

{\color{black}
\begin{proof}
Given $\beta<1$, choose $p>1$ sufficiently close to $1$ that
$\beta<(4d+2)/p-4d$.  Multiplication by a fixed smooth compactly supported
function is bounded on $\mathbf B^\beta_{p;\theta}$.  Since
$0<\beta<1$, the difference characterization of anisotropic Besov spaces,
together with H\"older's inequality on a fixed compact enlargement of
$\operatorname{supp}\chi$, gives
\begin{align*}
 \|\chi u_n\|_{\mathbf B^\beta_{1;\theta}}
 \leq C_{\chi,p,\beta}\|u_n\|_{\mathbf B^\beta_{p;\theta}}.
\end{align*}
The assertion follows from \eqref{0425:02} and the pathwise estimate
$\|a\|_{L^1(0,T)}\leq T^{1-1/p}\|a\|_{L^p(0,T)}$, followed by expectation.
\end{proof}}

The following compactness result is a direct consequence of Aubin-Lions lemma (see \cite[Corollary 4]{Si86}):

\begin{proposition}
Under the assumptions of Theorem \ref{thm-besov-p-es-un}, for every $\delta \in (0,1)$, the sequence $\{h_{\delta}(f_n)\}_{n \geq 1}$ is tight in $L^1([0,T]; L^1(D))$, for every bounded domain $D\subset \mathbb{R}^{2d}$.
\end{proposition}
{\color{black}
\begin{proof}
Choose bounded domains $D\Subset D'$ and a cutoff equal to one on $D$.
For the $p>1$ selected above and any $0<\beta<1$, the localized compact
embedding
\begin{align*}
 \mathbf B^\beta_{p;\theta}(D')\Subset L^1(D)
 \hookrightarrow \mathbf B^{-6}_{p;\theta}(D')
\end{align*}
holds.  Estimates {\color{black}\eqref{0425:02}} and \eqref{0521:050}, followed by
the Aubin--Lions--Simon compactness criterion, give tightness in
$L^1([0,T];L^1(D))$.
\end{proof}}

\subsection{\textcolor{black}{Local $L^1$ topology and tightness}}\label{subsec-6-4}
\begin{definition}\label{def-metric}
		{\color{black}Let $D\subset\mathbb{R}^{2d}$ be a bounded subset, put
		$\delta_k=(k+1)^{-1}$, and set
		\begin{align*}
		 L^1_+([0,T]\times D):=
		 \{u\in L^1([0,T]\times D):u\geq0\ \text{a.e.}\}.
		\end{align*}
		For $h_{\delta_k}$ as in Definition~\ref{def-hdelta}, define
		$\mathcal E:L^1_+([0,T]\times D)^2\to[0,\infty)$ by}
		\begin{align}\label{eq-7.13}
		{\color{black}\mathcal{E}(f,g)=\sum_{k=1}^{\infty}2^{-k}\left(\frac{\|h_{\delta_k}(f)-h_{\delta_k}(g)\|_{L^{1}\left([0,T];L^{1}\left(D\right)\right)}}{1+\|h_{\delta_k}(f)-h_{\delta_k}(g)\|_{L^{1}\left([0,T];L^{1}\left(D\right)\right)}}\right).}
		\end{align}
\end{definition}
Further, referring to \cite[Lemma 5.20]{FG24}, the following results hold.
\begin{lemma}\label{lem-equivalenttopology}
		{\color{black}Let $D\subset\mathbb{R}^{2d}$ be bounded. Then the
		function $\mathcal{E}$ in Definition~\ref{def-metric} is a metric on
		$L^1_+([0,T]\times D)$. Its metric topology agrees with the topology
		induced by the strong $L^1([0,T]\times D)$ norm.}

\end{lemma}

\begin{lemma}\label{L1-tight}
Under the same assumptions as in Theorem \ref{thm-besov-p-es-un}, let $D \subset \mathbb{R}^{2d}$ be a bounded domain. Then the sequence $\{f_n\}_{n \geq 1}$ is tight in $L^1([0,T]; L^1(D))$. 
\end{lemma}
\begin{proof}

{\color{black}Fix $\varepsilon>0$. For each $k\geq1$, tightness of
$\{h_{\delta_k}(f_n)\}_{n\geq1}$ yields a compact set
$C_k^\varepsilon\subset L^1([0,T];L^1(D))$ such that
\begin{align*}
 \sup_{n\geq1}\mathbb P\bigl(h_{\delta_k}(f_n)\notin C_k^\varepsilon\bigr)
 \leq \varepsilon2^{-k}.
\end{align*}
The compactness criterion in \cite[Lemma~5.20]{FG24}, applied simultaneously
to the truncations $h_{\delta_k}$, provides a compact set
$K^\varepsilon\subset L^1([0,T];L^1(D))$ containing every function whose
$k$th truncation belongs to $C_k^\varepsilon$ for all $k$. Hence
\begin{align*}
 \sup_{n\geq1}\mathbb P(f_n\notin K^\varepsilon)
 \leq\sum_{k=1}^\infty
 \sup_{n\geq1}\mathbb P\bigl(h_{\delta_k}(f_n)\notin C_k^\varepsilon\bigr)
 \leq\varepsilon.
\end{align*}
This proves tightness on $D$. Applying the same argument to a countable
exhaustion of $\mathbb R^{2d}$ gives tightness in
$L^1([0,T];L^1_{\rm loc}(\mathbb R^{2d}))$.}

\end{proof}

{\color{black}We next record the martingale estimate in a form that will also
be used to identify the limiting noise.  For
$\psi\in C_c^\infty(\mathbb R^{2d}\times(0,\infty))$, set
\begin{align*}
 H_{n,\psi}(z,r):=\int_0^r\psi(z,a)\sigma_n'(a)\,da
\end{align*}
and
\begin{align}\label{approximating-stochastic-primitive}
 \mathcal G_{n,\psi}^k(g)
 :={}&\int_{\mathbb R^{2d}}f_k(z)
       \int_0^{g(z)}\nabla_v\psi(z,a)\sigma_n'(a)\,da\,dz\notag\\
 &+\int_{\mathbb R^{2d}}
 \bigl(H_{n,\psi}(z,g(z))-\psi(z,g(z))\sigma_n(g(z))\bigr)
 \nabla_vf_k(z)\,dz.
\end{align}
Spatial integration by parts shows that this is precisely
$-\int\psi(z,g)\nabla_v(\sigma_n(g)f_k)\,dz$.}

\begin{lemma}\label{martingale-tightness}
{\color{black}Under the same assumptions as in Theorem~\ref{thm-besov-p-es-un},
for every $\psi\in C^{\infty}_c(\mathbb{R}^{2d}\times(0,\infty))$ and
$\gamma\in(0,1/2)$, the laws of the martingales}
\begin{align}\label{martingale}
{\color{black}M_n^{\psi}(t):=
\sum_{k\geq1}\int_0^t\mathcal G_{n,\psi}^k(f_n(s))\cdot dB^k(s)}
\end{align}
{\color{black}are tight on $C^{\gamma}([0,T])$.}
\end{lemma}
{\color{black}
\begin{proof}
Because the $z$-support of $\psi$ is compact and its $\zeta$-support is
separated from zero, \eqref{square-root-uniform-coefficients} implies that
the two primitives in \eqref{approximating-stochastic-primitive} are
bounded uniformly in $n$ and in their last argument.  Minkowski's
inequality, together with
$\sum_k f_k^2=F_1$ and $\sum_k|\nabla_vf_k|^2=F_3$, therefore gives
\begin{align*}
 \sup_{n\geq1}\sup_{g\geq0}
 \sum_{k\geq1}|\mathcal G_{n,\psi}^k(g)|^2\leq C_\psi.
\end{align*}
The Burkholder--Davis--Gundy inequality consequently yields, for every
$q\geq2$,
\begin{align*}
 \mathbb E|M_n^\psi(t)-M_n^\psi(s)|^q
 \leq C_{q,\psi}|t-s|^{q/2}.
\end{align*}
Given $\gamma<1/2$, choose $\gamma'$ such that
$\gamma<\gamma'<1/2$.  Kolmogorov's criterion, with $q$ sufficiently
large, gives a uniform bound in $C^{\gamma'}([0,T])$.  The compact
embedding $C^{\gamma'}([0,T])\Subset C^\gamma([0,T])$ then proves the
asserted tightness, in agreement with \cite[Proposition~5.23]{FG24}.
\end{proof}}

\subsection{\textcolor{black}{Probabilistically weak existence}}\label{subsec-6-5}
{\color{black}We now construct a solution in the sense of Definition~\ref{def-prob-weak-square-root}.}

{\color{black}
\begin{lemma}[Weighted weak lower semicontinuity]
\label{lem-weighted-defect-lsc}
Let $Q$ be a finite-measure subset of a Euclidean space, let
$r_n,r:Q\to[0,\infty)$ be measurable, and let
$G_n,G\in L^2(Q;\mathbb R^d)$.  Assume that
\begin{align*}
 r_n\longrightarrow r\quad\text{in measure on }Q,
 \qquad G_n\rightharpoonup G\quad\text{weakly in }L^2(Q;\mathbb R^d).
\end{align*}
If $\Phi:Q\times[0,\infty)\to[0,\infty)$ is bounded and continuous in
its last variable, uniformly on $Q$, then
\begin{align}\label{weighted-defect-lsc}
 \int_Q\Phi(q,r(q))|G(q)|^2\,dq
 \leq\liminf_{n\to\infty}
 \int_Q\Phi(q,r_n(q))|G_n(q)|^2\,dq.
\end{align}
The same conclusion holds for a continuous compactly supported
$\Phi(q,a)$ after restricting $Q$ to the projection of its support.
\end{lemma}

\begin{proof}
Put $a_n(q)=\Phi(q,r_n(q))^{1/2}$ and
$a(q)=\Phi(q,r(q))^{1/2}$.  Then $(a_n)_n$ is uniformly bounded and
$a_n\to a$ in measure.  For every $H\in L^2(Q;\mathbb R^d)$,
bounded convergence in measure gives
\begin{align*}
 a_nH\longrightarrow aH\quad\text{strongly in }L^2(Q;\mathbb R^d).
\end{align*}
Indeed, this follows first along every almost-everywhere convergent
subsequence by dominated convergence, and then for the full sequence by the
subsequence principle.  Hence
\begin{align*}
 \int_Q a_nG_n\cdot H\,dq
 =\int_QG_n\cdot(a_nH)\,dq
 \longrightarrow\int_QG\cdot(aH)\,dq.
\end{align*}
Thus $a_nG_n\rightharpoonup aG$ weakly in $L^2(Q;\mathbb R^d)$, and
\eqref{weighted-defect-lsc} follows from weak lower semicontinuity of the
$L^2$ norm.
\end{proof}}

{\color{black}
\begin{proof}[Proof of Theorem~\ref{thm-existence-squareroot}]
For every $n\geq1$, let $f_n$ be the mass-preserving weak solution of
\eqref{SPDE-6} with $V=0$, initial datum $f_{0,n}$, and coefficient
$\sigma_n$.  Its kinetic function and parabolic defect measure are
\begin{align*}
 \chi_{f_n}(z,\zeta,t)&:=\boldsymbol1_{\{0<\zeta<f_n(t,z)\}},\\
 p_n(dz,d\zeta,dt)&:=
 \delta_{f_n(t,z)}(d\zeta)|\nabla_vf_n(t,z)|^2\,dz\,dt.
\end{align*}
By Proposition~\ref{equivalence}, each $f_n$ satisfies
\eqref{kinetic-formula-equivalence} with $V=0$, coefficient $\sigma_n$, and
defect measure $p_n$.  Multiplying that identity by a test function in time
and integrating by parts gives the distribution-in-time identity used below.
Thus $p_n$ is a measure on
\begin{align*}
 E:=\mathbb R^{2d}\times(0,\infty)\times[0,T]
\end{align*}
in the order $(z,\zeta,t)$.  For every compact $K\Subset E$,
\eqref{square-root-uniform-L2} gives
\begin{align}\label{local-kinetic-measure-first-moment}
 \sup_{n\geq1}\mathbb E p_n(K)
 \leq \sup_{n\geq1}\mathbb E\int_0^T
       \|\nabla_vf_n(t)\|_2^2\,dt<\infty.
\end{align}
Only this first-moment estimate will be used.

Choose an integer $l>d+\frac32$.  For each $R,m\in\mathbb N$, choose a
countable set dense in
$C_c^\infty(B_R\times(m^{-1},m))$ in the $H^l$ topology, and enumerate the
union of these sets as $(\psi_q)_{q\geq1}$.  In this way, approximation of
a test function never destroys the positive distance of its
$\zeta$-support from zero.  Let
$\mathcal M_{\rm loc}^+(E)$ denote the nonnegative Radon measures on $E$
with the local vague topology and put
\begin{align*}
 \mathcal X_f&:=
 L^1([0,T];L^1_{\rm loc}(\mathbb R^{2d}))
 \cap L^2([0,T]\times\mathbb R^{2d}),\\
 \mathbb X&:=\mathcal X_f\times
 {\color{black}L^2([0,T]\times\mathbb R^{2d};\mathbb R^d)_{\rm weak}}
 \times\mathcal M_{\rm loc}^+(E)
 \times C([0,T])^{\mathbb N}.
\end{align*}
On $\mathcal X_f$ we use the join of the strong local $L^1$ topology and
the weak global $L^2$ topology.  The four coordinates of the random
variable in $\mathbb X$ are
\begin{align*}
 \left(f_n,\nabla_vf_n,p_n,(M_n^{\psi_q})_{q\geq1}\right).
\end{align*}

{\color{black}The first coordinate is tight in its join topology by
Lemma~\ref{L1-tight} together with the weakly compact global $L^2$ balls
provided by \eqref{square-root-uniform-L2}.  The second and fourth
coordinates are tight by \eqref{square-root-uniform-L2} and
Lemma~\ref{martingale-tightness}, respectively.}  To treat the third
coordinate, take a compact exhaustion $(K_j)_{j\geq1}$ of $E$.
Estimate \eqref{local-kinetic-measure-first-moment} and Markov's inequality
give numbers $a_j<\infty$ such that, with probability at least
$1-\varepsilon$, $p_n(K_j)\leq a_j$ for every $j$.  The set of measures
satisfying these bounds is relatively compact for the local vague
topology.  Thus $(p_n)_{n\geq1}$ is tight as an
$\mathcal M_{\rm loc}^+(E)$-valued sequence.  This avoids any unsupported
$L^2(\Omega;\mathcal M)$ weak-compactness assertion.

Adding the Brownian paths $(B^k)_{k\geq1}$ and $f_{0,n}$ as coordinates,
the preceding bounds and \eqref{square-root-initial-approximation} give
tightness on the corresponding product Jakubowski space.  The
Jakubowski--Skorokhod representation theorem supplies a subsequence, a
probability space $(\bar\Omega,\bar{\mathcal F},\bar{\mathbb P})$, and
copies
\begin{align*}
 (\bar f_n,\nabla_v\bar f_n,\bar p_n,
  (\bar M_n^{\psi_q})_{q\geq1},(\bar B_n^k)_{k\geq1},\bar f_{0,n})
\end{align*}
having the original laws, which converge almost surely to
\begin{align*}
 (\bar f,{\color{black}\bar g},\bar p,
  (\bar M^{\psi_q})_{q\geq1},(\bar B^k)_{k\geq1},\bar f_0).
\end{align*}
In particular,
\begin{align}\label{square-root-representation-convergences}
 \bar f_n&\longrightarrow\bar f
 &&\text{strongly in }L^1([0,T];L^1_{\rm loc}),\notag\\
 \bar f_n&\rightharpoonup\bar f,\qquad
 \nabla_v\bar f_n\rightharpoonup{\color{black}\bar g}
 &&\text{weakly in the corresponding global }L^2\text{ spaces},\\
 \bar p_n&\longrightarrow\bar p
 &&\text{locally vaguely on }E,\notag\\
 \bar M_n^{\psi_q}&\longrightarrow\bar M^{\psi_q},\qquad
 \bar B_n^k\longrightarrow\bar B^k
 &&\text{uniformly on }[0,T].\notag
\end{align}
Also $\bar f_{0,n}\to\bar f_0$ strongly in $L^1(\mathbb R^{2d})$ and
$\operatorname{Law}(\bar f_0)=\operatorname{Law}(f_0)$.
{\color{black}The weak limit $\bar g$ is the distributional $v$-gradient of
$\bar f$.  Indeed, for every
$\Phi\in C_c^\infty((0,T)\times\mathbb R^{2d};\mathbb R^d)$,
\begin{align*}
 \int_0^T\!\int_{\mathbb R^{2d}}\bar g\cdot\Phi\,dz\,dt
 &=\lim_{n\to\infty}\int_0^T\!\int_{\mathbb R^{2d}}
       \nabla_v\bar f_n\cdot\Phi\,dz\,dt\\
 &=-\lim_{n\to\infty}\int_0^T\!\int_{\mathbb R^{2d}}
       \bar f_n\nabla_v\cdot\Phi\,dz\,dt
 =-\int_0^T\!\int_{\mathbb R^{2d}}
       \bar f\nabla_v\cdot\Phi\,dz\,dt,
\end{align*}
where the last limit follows from the local strong $L^1$ convergence.
Thus $\bar g=\nabla_v\bar f$; from now on we use the latter notation.}
Because the distributional kinetic identity, the equality defining
$p_n$, and the martingale identities are Borel properties of the joint
coordinates, every represented copy satisfies the same identities as its
original counterpart.

{\color{black}Let $(\mathcal G_t)_{t\in[0,T]}$ be the usual augmentation of the
canonical filtration generated up to time $t$ by all the limiting
coordinates, including $\bar f_0$ and the restriction of $\bar p$ to
$\mathbb R^{2d}\times(0,\infty)\times[0,t)$.  We spell out why the Brownian
property survives the representation.  Fix $0\leq s<t\leq T$, indices
$k,\ell\geq1$ and $i,j\in\{1,\ldots,d\}$, and let $\Xi$ be a bounded
continuous cylinder functional of the limiting coordinates restricted to
$[0,s]$; for the measure coordinate, its test functions are taken with time
support in $[0,s)$.  Let $\Xi_n$ be the same functional of the represented
approximating coordinates.  Equality of laws gives
\begin{align*}
 \bar{\mathbb E}\big[(\bar B_n^{k,i}(t)-\bar B_n^{k,i}(s))\Xi_n\big]&=0,\\
 \bar{\mathbb E}\Big[\Big((\bar B_n^{k,i}(t)-\bar B_n^{k,i}(s))
 (\bar B_n^{\ell,j}(t)-\bar B_n^{\ell,j}(s))
 -\delta_{k\ell}\delta_{ij}(t-s)\Big)\Xi_n\Big]&=0.
\end{align*}
Almost-sure convergence of the joint coordinates and uniform integrability
of the Brownian increments permit passage to the limit in both identities.
A monotone-class argument, followed by L\'evy's characterization for every
finite subfamily, shows that $(\bar B^k)_{k\geq1}$ are independent
$d$-dimensional Brownian motions relative to $(\mathcal G_t)$.  A jointly
measurable version of the canonical $L^1_{\rm loc}$ coordinate is
progressively measurable with respect to this filtration; we use this
version for $\bar f$.  Moreover, for every compactly supported}
$\varphi$, the process
\begin{align*}
 t\longmapsto\int_{\mathbb R^{2d}\times(0,\infty)\times[0,t)}
 \varphi(z,\zeta)\,\bar p(dz,d\zeta,ds)
\end{align*}
{\color{black}is left-continuous and adapted, hence predictable.}

We next identify the martingales rather than invoking continuity of the
stochastic-integral map.  For every $q,r,k$ the approximating martingales
satisfy
\begin{align}\label{approx-bracket-identities}
 \langle M_n^{\psi_q},M_n^{\psi_r}\rangle_t
 &=\int_0^t\sum_{k\geq1}
 \mathcal G_{n,\psi_q}^k(f_n(s))\cdot
 \mathcal G_{n,\psi_r}^k(f_n(s))\,ds,\notag\\
 \langle M_n^{\psi_q},B^{k,i}\rangle_t
 &=\int_0^t\mathcal G_{n,\psi_q}^{k,i}(f_n(s))\,ds.
\end{align}
The primitive formula \eqref{approximating-stochastic-primitive}, local
strong convergence in \eqref{square-root-representation-convergences}, and
the local $C^1$ convergence of $\sigma_n$ imply
\begin{align}\label{stochastic-integrand-convergence}
 (\mathcal G_{n,\psi_q}^k(\bar f_n))_{k\geq1}
 \longrightarrow(\mathcal G_{\psi_q}^k(\bar f))_{k\geq1}
 \quad\text{in }L^2(\bar\Omega\times[0,T];\ell^2).
\end{align}
{\color{black}For completeness, let $K_\psi\Subset\mathbb R^{2d}$ contain
the $z$-support of $\psi_q$ and put
\begin{align*}
 U_{n,\psi}(z,r)&:=\int_0^r\nabla_v\psi_q(z,a)\sigma_n'(a)\,da,\\
 V_{n,\psi}(z,r)&:=H_{n,\psi_q}(z,r)
                  -\psi_q(z,r)\sigma_n(r),
\end{align*}
with $U_\psi,V_\psi$ defined by replacing $\sigma_n$ with $\sqrt{\cdot}$.
Because the $\zeta$-support of $\psi_q$ is a compact subset of
$(0,\infty)$, these four functions are uniformly bounded, and
$U_{n,\psi}\to U_\psi$, $V_{n,\psi}\to V_\psi$ uniformly in $(z,r)$.
Minkowski's and Cauchy--Schwarz inequalities, together with the definitions
of $F_1$ and $F_3$, give
\begin{align*}
 &\sum_{k\geq1}\left|
 \mathcal G_{n,\psi_q}^k(\bar f_n)
 -\mathcal G_{\psi_q}^k(\bar f)\right|^2\\
 &\qquad\leq C_{\psi_q}\int_{K_\psi}
 \left(
 |U_{n,\psi}(z,\bar f_n)-U_\psi(z,\bar f)|^2
 +|V_{n,\psi}(z,\bar f_n)-V_\psi(z,\bar f)|^2
 \right)dz.
\end{align*}
The local strong $L^1$ convergence of $\bar f_n$ implies convergence in
measure on $[0,T]\times K_\psi$.  Boundedness of the displayed
compositions and dominated convergence, first pathwise and then in
$\bar\Omega$, prove \eqref{stochastic-integrand-convergence}; in particular,
no convergence of $\nabla_v\bar f_n$ is needed here.

It remains to justify the martingale passage relative to the limiting
filtration.  The approximating processes $M_n^{\psi_q}$,
\begin{align*}
 M_n^{\psi_q}M_n^{\psi_r}
 -\int_0^\cdot\sum_{k\geq1}
 \mathcal G_{n,\psi_q}^k(f_n)
 \cdot\mathcal G_{n,\psi_r}^k(f_n)\,ds,
\end{align*}
and
\begin{align*}
 M_n^{\psi_q}B^{k,i}
 -\int_0^\cdot\mathcal G_{n,\psi_q}^{k,i}(f_n)\,ds
\end{align*}
are martingales.  Equivalently, each of their increments has zero expectation
after multiplication by every bounded continuous past-cylinder functional
 $\Xi_n$ introduced above.  {\color{black}The moment bound} in
Lemma~\ref{martingale-tightness} gives uniform integrability.  Hence uniform
path convergence, \eqref{stochastic-integrand-convergence}, and the same
monotone-class argument pass these conditional martingale identities to the
limit.  Thus $\bar M^{\psi_q}$ is a continuous $(\mathcal G_t)$-martingale
with the limiting brackets and cross-brackets in
\eqref{approx-bracket-identities}.  It follows that}
\begin{align*}
 \bar M^{\psi_q}(t)=
 \sum_{k\geq1}\int_0^t
 \mathcal G_{\psi_q}^k(\bar f(s))\cdot d\bar B^k(s).
\end{align*}
Indeed, the difference between the two sides is a continuous martingale
with zero quadratic variation.  This proves the required martingale
identification, including its cross-variations with the driving Brownian
motions.  Set $\bar W_F:=\sum_{k\geq1}\bar B^kf_k$.

It remains to pass to the deterministic terms and the defect measure.
All tests have compact $z$-support and $\zeta$-support contained in
$[a,b]\Subset(0,\infty)$.  On this set the relevant functions of
$\sigma_n$ converge uniformly to
\begin{align*}
 \sigma(r)=\sqrt r,\qquad
 \sigma'(r)^2=\frac1{4r},\qquad
 \sigma(r)\sigma'(r)=\frac12.
\end{align*}
{\color{black}No derivative of $\sigma_n\sigma_n'$ occurs in the kinetic
formulation, so no convergence of such derivatives is required.}
Consequently, the strong local convergence of $\bar f_n$ and the weak
$L^2$ convergence of $\nabla_v\bar f_n$ pass all drift and It\^o
correction terms to the limit.  The transport term has the positive sign
$+\int v\chi_{\bar f}\cdot\nabla_x\psi$.  Local vague convergence passes
the defect term after multiplication by a compactly supported smooth
function of time.  Thus, for every $q$ and
$\eta\in C_c^\infty([0,T))$, almost surely,
\begin{align}\label{passing-to-the-limits}
 -\int_0^T\eta'(t)\langle\chi_{\bar f}(t),\psi_q\rangle\,dt
={}&\eta(0)\langle\chi_{\bar f_0},\psi_q\rangle\notag\\
&+\int_0^T\eta(t)\Bigg[
 -\int\nabla_v\bar f\cdot(\nabla_v\psi_q)(z,\bar f)\,dz
 +\iint v\chi_{\bar f}\cdot\nabla_x\psi_q\,dz\,d\zeta\notag\\
&\qquad+\int\nabla_v\cdot(v\bar f)\psi_q(z,\bar f)\,dz
 -\frac18\int F_1\bar f^{-1}\nabla_v\bar f\cdot
       (\nabla_v\psi_q)(z,\bar f)\,dz\notag\\
&\qquad{\color{black}+\frac12\int F_3\bar f
       (\partial_\zeta\psi_q)(z,\bar f)\,dz}
 \Bigg]dt\notag\\
&+\sum_{k\geq1}\int_0^T\eta(t)
 \mathcal G_{\psi_q}^k(\bar f(t))\cdot d\bar B^k(t)
 -\int_E\eta(t)\partial_\zeta\psi_q(z,\zeta)
 \,\bar p(dz,d\zeta,dt).
\end{align}
All unlabelled spatial integrals in this formula are over
$\mathbb R^{2d}$, and $\langle\chi,\psi\rangle$ denotes integration in
$(z,\zeta)$.  The factor $\bar f^{-1}$ is evaluated only on
$\{a\leq\bar f\leq b\}$ and is harmless.  The choice of the dense family
and Sobolev embedding extend \eqref{passing-to-the-limits} to every
$\psi\in C_c^\infty(\mathbb R^{2d}\times(0,\infty))$.

We finally verify the defect inequality and the global-in-space
statements that survive the local limit.  If
$\Phi\in C_c(E)$ is nonnegative, let $Q\subset
[0,T]\times\mathbb R^{2d}$ be a compact set containing the $(t,z)$-projection
of its support.  The local strong $L^1$ convergence implies
$\bar f_n\to\bar f$ in measure on $Q$, while
$\nabla_v\bar f_n\rightharpoonup\nabla_v\bar f$ weakly in $L^2(Q)$.
Lemma~\ref{lem-weighted-defect-lsc}, applied with
$r_n=\bar f_n$, $r=\bar f$, and $G_n=\nabla_v\bar f_n$, therefore gives
\begin{align*}
 \int\Phi(z,\bar f,t)|\nabla_v\bar f|^2\,dz\,dt
 &\leq\liminf_{n\to\infty}
 \int\Phi(z,\bar f_n,t)|\nabla_v\bar f_n|^2\,dz\,dt\\
 &=\int_E\Phi\,d\bar p.
\end{align*}
In the last equality we used the definition of $p_n$ and the local vague
convergence $\bar p_n\to\bar p$.  Thus the weighted lower-semicontinuity
step does not require strong $L^2$ convergence of the gradients.
Hence
\begin{align}\label{limiting-kinetic-measure-inequality}
 \delta_{\bar f(t,z)}(d\zeta)|\nabla_v\bar f(t,z)|^2\,dz\,dt
 \leq\bar p(dz,d\zeta,dt).
\end{align}
Weak lower semicontinuity also yields
\begin{align*}
 \bar{\mathbb E}\int_0^T
 \bigl(\|\bar f(t)\|_2^2+\|\nabla_v\bar f(t)\|_2^2\bigr)\,dt<\infty.
\end{align*}
{\color{black}Choose the further deterministic subsequence so that, on every
ball in a countable exhaustion, the probabilities that the local
$L^1_{t,z}$ error exceeds $2^{-j}$ are summable.  The Borel--Cantelli lemma
and then Fubini's theorem show that the local $L^1_z$ convergence holds for
$\bar{\mathbb P}$-almost every $\omega$ and almost every $t$.}  For such
$(\omega,t)$,
nonnegativity, Fatou's lemma on expanding balls, mass preservation of the
regular approximations, and the global $L^1$ convergence of the initial
data give
\begin{align*}
 \int_{\mathbb R^{2d}}\bar f(t,z)\,dz
 &\leq\liminf_{n\to\infty}\int_{\mathbb R^{2d}}\bar f_n(t,z)\,dz\\
 &=\lim_{n\to\infty}\int_{\mathbb R^{2d}}\bar f_{0,n}(z)\,dz
 =\int_{\mathbb R^{2d}}\bar f_0(z)\,dz.
\end{align*}
This is only a mass inequality; no global $L^1$ convergence of the
solutions has been used.

We have verified every item of
Definition~\ref{def-prob-weak-square-root}.  In particular,
\begin{align*}
 \bigl((\bar\Omega,\bar{\mathcal F},(\mathcal G_t)_{t\in[0,T]},
 \bar{\mathbb P}),\bar W_F,\bar f_0,\bar f,\bar p\bigr)
\end{align*}
is a probabilistically weak renormalized kinetic solution.  We neither
invoke uniqueness nor return to the original probability space.
\end{proof}}
\endgroup

\noindent{\bf  Acknowledgements}\quad This work is supported by the US Army Research Office, grant W911NF2310230, and the DFG through the CRC 1283 ``Taming uncertainty and profiting from randomness and low regularity in analysis, stochastics and their applications."

\bibliographystyle{alpha}
\bibliography{stochastic-kinetic.bib}

\end{document}